%% file: kiselyov.tex
\documentclass[twoside]{book}
\usepackage{book-ru}

\newcommand{\spell}[2]{#1} 

\hypersetup{pdftitle={Геометрия по Киселёву}}

\makeindex

\makeatletter
\newcommand{\rindex}[2][\imki@jobname]{%
  \index[#1]{\detokenize{#2}}%
}
\makeatother

\spell{}{\makeatletter\usepackage{environ}\RenewEnviron{figure}[1][]{ [РИС.] }\RenewEnviron{figure*}[1][]{ [РИС.] }\RenewEnviron{wrapfigure}[1][]{ [РИС.] }\usepackage{emptypage}\pagestyle{empty}\makeatother\usepackage[none]{hyphenat}\def\emph{\textit}\def\so{\textit}\renewcommand{\footnote}[1]{\textup{\ (#1)}}\renewcommand{\frac}[2]{(#1)/(#2)}\renewcommand{\tfrac}[2]{(#1)/(#2)}\everydisplay{\textstyle}}

\begin{document}

\cleardoublepage
\frontmatter
\title{Геометрия по Киселёву}
\author{Андрей Петрович Киселёв}
\date{}
\maketitle

\thispagestyle{empty}

Под редакцией Н. А. Ершова, А. М. Петрунина и С. Л. Табачникова.

На обложке: портрет  Андрея Петровича Киселёва, написанный его дочерью Еленой Андреевной Киселёвой, 1906 год.

\vfill

\noindent
\includegraphics[scale=.25]{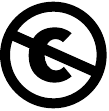}\ \ 
Это произведение находится в общественном достоянии.

\ 

\noindent\texttt{https://anton-petrunin.github.io/kiselyov/} 

\ 


\mainmatter

\input{predislovie.tex}
\input{vvedenie.tex}

\cleardoublepage
\part{Планиметрия}

\chapter{Прямая линия}
\input{2D/ugly.tex}
\input{2D/mat-predlozheniya.tex}
\input{2D/treugi.tex}
\input{2D/zadachi-na-postr.tex}

\chapter{Параллельные прямые}
\input{2D/parallelnye.tex}
\input{2D/parallelogrammy.tex}
\input{2D/postulat5.tex}
\input{2D/okruzhnost.tex}

\input{2D/dugi.tex}
\input{2D/raspolozheniya.tex}
\input{2D/vpisannye-ugly.tex}
\input{2D/vpis-opis-mnougi.tex}
\input{2D/zam-toch-trig.tex}

\chapter{Измерение отрезков}
\input{2D/izmereniya.tex}
\input{2D/nesoimerimost.tex}

\chapter{Подобные фигуры}
\input{2D/podobie-trig.tex}
\input{2D/podobie-fig.tex}
\input{2D/proportzii.tex}
\input{2D/teorema-pifagora.tex}
\input{2D/zadach-na-vych.tex}
\input{2D/proportzii-v-kruge.tex}
\input{2D/algebra-k-geomtrii.tex}

\chapter[Правильные многоугольники]{Правильные многоугольники и длина окружности}
\input{2D/pravilnye-mnogougi.tex}
\input{2D/dlina-okr.tex}

\chapter{Измерение площадей}
\input{2D/ponyatie-o-ploschadi.tex}
\input{2D/ploschadi-mnogugov.tex}

\input{2D/ploschad-kruga.tex}

\cleardoublepage
\part{Стереометрия}
\input{3D/pryamye_i_ploskosti.tex}

\input{3D/mnogogranniki.tex}
\input{3D/kruglye_tela.tex}

{
\bookmarksetupnext{startatroot,level=0}
\small
\printindex
}

\newpage
\cleardoublepage
\phantomsection
\pdfbookmark[0]{\contentsname}{bm:toc}
\tableofcontents


\end{document}

%% file: predislovie.tex
\section*{От редакции}

Перед вами один из лучших учебников по геометрии —
многие темы этого учебника просто невозможно объяснить более доходчиво.

Предмет «Геометрия» трактуется сегодня иначе;
во времена Киселёва геометрия включала в себя начала анализа и даже элементы теории чисел.
Этим учебник Киселёва существенно отличается от современных учебников, и поэтому он может помочь современным школьникам.
Например, для них может оказаться полезным геометрическое понимание предела, непрерывности, интеграла, алгоритма Евклида и иррациональных чисел.

Наша цель --- сделать удобным использование этого учебника сегодня;
мы заменили термины на современные и внесли незначительные уточнения, в основном исторические. 
Мы используем лицензию CC0 (то есть отказываемся от авторских прав); сам учебник Киселёва перешёл в общественное достояние несколько лет назад.
В частности, любой желающий может использовать любую часть этого учебника для написания своего под любой лицензией —
возможно, таким образом удастся продолжить традиции, и облегчить создание хороших учебников. 

Автор учебника, Андрей Петрович Киселёв, сочетал в себе хорошее понимание современной ему математики и талант школьного учителя — такое сочетание редко и ценно во все времена.
В учебнике раскрывается несколько тем, которые были новыми во время его написания. 
Этот учебник (очевидно) написан на основе учебника Августа Юльевича Давидова, структура взята практически без изменения, но с существенными улучшениями.
Практически весь учебник строится на строгих доказательствах, и при этом Киселёв смог обойтись без зауми;
очевидно, что в этом помог его преподавательский опыт.
Доказательства в учебнике всегда строятся по наиболее наглядному пути (пусть даже чуть более сложному) — однажды поняв такое доказательство, его уже невозможно забыть. 

На основе этого учебника (прямо или косвенно) были написаны все современные русские учебники по геометрии.
Среди авторов последующих учебников были математики и педагоги высокого уровня — 
Андрей Николаевич Колмогоров, 
Алексей Васильевич Погорелов, 
Левон Сергеевич Атанасян, 
Александр Данилович Александров, 
Игорь Фёдорович Шарыгин.
Но учебник Киселёва остаётся непревзойдённым по одному из основных критериев —
он предъявляет минимальные требования к преподавателю.

\begin{flushright}
Н. А. Ершов,
А. М. Петрунин, 
С. Л. Табачников.        
\end{flushright}
\clearpage

%% file: vvedenie.tex
\section*{Введение}

\paragraph{Геометрические фигуры.}\label{1938/1}
Часть пространства, ограниченная со всех сторон, называется \rindex{геометрическое тело}\textbf{геометрическим телом}.

Геометрическое тело отделяется от окружающего пространства \rindex{поверхность}\textbf{поверхностью}.

Часть поверхности отделяется от смежной части \rindex{линия}\textbf{линией}.

Часть линии отделяется от смежной части \rindex{точка}\textbf{точкой}.

Геометрическое тело, поверхность, линия и точка не существуют раздельно.
Однако при помощи отвлечения мы можем рассматривать поверхность независимо от геометрического тела, линию — независимо от поверхности и точку — независимо от линии.
При этом поверхность мы должны представить себе не имеющей толщины, линию — не имеющей ни толщины, ни ширины и точку — не имеющей ни длины, ни ширины, ни толщины.

Совокупность каких бы то ни было точек, линий, поверхностей или тел, расположенных известным образом в пространстве, называется вообще геометрической \rindex{фигура}\textbf{фигурой}.
Геометрические фигуры могут перемещаться в пространстве, не подвергаясь никаким изменениям.
Две геометрические фигуры называются \rindex{равные фигуры}\textbf{равными}, если перемещением одной из них в пространстве её можно совместить со второй фигурой так, что обе фигуры совместятся во всех своих частях.

\paragraph{Геометрия.}\label{1938/2}
Наука, рассматривающая свойства геометрических фигур, называется \textbf{геометрией}, что в переводе с греческого языка означает \textbf{землемерие}.
Возможно такое название этой науке было дано потому, что в древнее время главной целью геометрии было измерение расстояний и площадей на земной поверхности.

\subsection*{Плоскость}

\paragraph{Плоскость.}\label{1938/3}
Из различных поверхностей наиболее знакомая нам есть плоская поверхность, или просто \so{плоскость}, представление о которой даёт нам, например, поверхность хорошего оконного стекла или поверхность спокойной воды в пруде.

Укажем следующее свойство плоскости:
\textit{Всякую часть плоскости можно наложить всеми её точками на другое место этой или другой плоскости, причём накладываемую часть можно предварительно перевернуть другой стороной.}

\subsection*{Прямая линия}

\paragraph{Прямая линия.}\label{1938/4}
Самой простой линией является прямая.
Представление о прямой линии, или просто о прямой, всем хорошо знакомо.
Представление о ней даёт туго натянутая нить или луч света, выходящий из малого отверстия.
С этим представлением согласуется следующее основное свойство прямой.

\textit{Через всякие две точки пространства можно провести прямую и притом только одну.}

Из этого свойства следует:

\textit{Если две прямые наложены одна на другую так, что какие-нибудь две точки одной прямой совпадают с двумя точками другой прямой, то эти прямые сливаются и во всех остальных точках} (потому что в противном случае через две точки можно было бы провести две различные прямые, что невозможно).

По той же причине \textit{две прямые могут пересечься только в одной точке}.

Прямая линия может лежать на плоскости.
При этом плоскость обладает следующим свойством.

\textit{Если на плоскости взять какие-нибудь две точки и провести через них прямую линию, то все точки этой прямой будут находиться в этой плоскости.}

\paragraph{Луч и отрезок.}\label{1938/5}
Если прямую представляют продолженной в обе стороны бесконечно, то её называют \textbf{бесконечной} (или \textbf{неограниченной}) прямой.

\begin{wrapfigure}{o}{33 mm}
\vskip-2mm
\centering
\includegraphics{mppics/ris-1}
\caption{}\label{1938/ris-1}
\bigskip
\includegraphics{mppics/ris-2}
\caption{}\label{1938/ris-2}
\bigskip
\includegraphics{mppics/ris-3}
\caption{}\label{1938/ris-3}
\end{wrapfigure}

Прямую обозначают обыкновенно двумя большими буквами, поставленными у двух каких-либо её точек.
Так, говорят:
«прямая $AB$» или «$BA$» (рис.~\ref{1938/ris-1}).

Часть прямой, ограниченная с обеих сторон, называется \rindex{отрезок}\textbf{отрезком};
отрезок обыкновенно обозначается двумя буквами, поставленными у его концов (отрезок $CD$, рис.~\ref{1938/ris-2}).
Иногда прямую или отрезок обозначают и одной буквой (малой);
например, говорят: «прямая $a$, отрезок $b$».

Иногда рассматривают прямую, ограниченную только с одной стороны, например в точке $E$ (рис.~\ref{1938/ris-3}).
О~такой прямой говорят, что она исходит из точки $E$;
её называют \rindex{луч}\textbf{лучом} или \rindex{полупрямая}\textbf{полупрямой}. 

\paragraph{Равенство и неравенство отрезков.}\label{1938/6}
\emph{Два отрезка равны, если они могут быть наложены один на другой так, что их концы совпадут.}
Положим, например, что мы накладываем отрезок $AB$ на
отрезок $CD$ (рис.~\ref{1938/ris-4}) так, чтобы точка $A$ совпала с точкой $C$ и чтобы прямая $AB$ пошла по прямой $CD$, если при этом концы $B$ и $D$ совпадут, то отрезки $AB$ и $CD$ равны;
в противном случае отрезки будут не равны, причём меньшим считается тот, который составит часть другого.

\begin{figure}
\centering
\includegraphics{mppics/ris-4}
\caption{}\label{1938/ris-4}
\end{figure}

Чтобы на какой-нибудь прямой отложить отрезок, равный данному отрезку, употребляют \textbf{циркуль} — прибор, известный учащимся из опыта.

\paragraph{Сумма отрезков.}\label{1938/7}
\rindex{сумма!отрезков}
Суммой нескольких данных отрезков $AB$, $CD$, $EF,\dots$
(рис.~\ref{1938/ris-5}) называется такой отрезок, который получится следующим образом.
На какой-нибудь прямой берём произвольную точку $M$ и откладываем от неё отрезок $MN$, равный $AB$, затем от точки $N$ в том же направлении откладываем отрезок $NP$, равный $CD$, и отрезок $PQ$, равный $EF$.
Тогда отрезок $MQ$ и будет суммой отрезков $AB$, $CD$ и $EF$ (которые по отношению к этой сумме называются слагаемыми).
Подобным образом можно получить сумму какого угодно числа отрезков.

\begin{figure}[!ht]
\centering
\includegraphics{mppics/ris-5}
\caption{}\label{1938/ris-5}
\end{figure}

Сумма отрезков обладает всеми свойствами суммы чисел;
так, она не зависит от порядка слагаемых (\so{переместительный закон}) и не изменяется, если некоторые слагаемые будут заменены их суммой (\so{сочетательный закон}).
Так:
\[AB+CD+EF=AB+EF+CD=EF+CD+AB=\dots\]
и
\[AB+CD+EF=AB+(CD+EF)=CD+(AB+EF)=\dots\]

\paragraph{Действия над отрезками.}\label{1938/8}
Из понятия о сумме выводятся понятия о разности отрезков, умножении и делении отрезков на число.
Так, разность отрезков $AB$ и $CD$ (если $AB>CD$) есть такой третий отрезок, сумма которого с $CD$ равна $AB$;
произведение отрезка $AB$ на число $3$ есть сумма трёх отрезков, из которых каждый равен $AB$;
частное от деления отрезка $AB$ на число $3$ есть третья часть $AB$ и так далее.

Если данные отрезки измерены какой-нибудь линейной единицей (например, сантиметром), и длины их выражены соответствующими числами, то длина суммы отрезков выразится суммой чисел, измеряющих эти отрезки, разность выразится разностью чисел и~т.~д.

\subsection*{Понятие об окружности}

\paragraph{Окружность.}\label{1938/9}
Если дадим циркулю произвольный раствор и, поставив одну его ножку остриём в какую-нибудь точку $O$ плоскости (рис.~\ref{1938/ris-6}), станем вращать циркуль вокруг этой точки, то другая его ножка, снабжённая карандашом или пером, прикасающимся к плоскости, опишет на плоскости непрерывную линию, все точки которой одинаково удалены от точки $O$.
Эта линия называется \rindex{окружность}\textbf{окружностью}, и точка $O$ — её \rindex{центр!окружности}\textbf{центром}.
Отрезки $OA$, $OB$, $OC,\dots$, соединяющие центр с какими-нибудь точками окружности, называются \rindex{радиус}\textbf{радиусами}.
Все радиусы одной окружности равны между собой.

\begin{wrapfigure}{o}{39 mm}
\vskip-3mm
\centering
\includegraphics{mppics/ris-6}
\caption{}\label{1938/ris-6}
\end{wrapfigure}

Окружности, описанные одинаковыми радиусами, равны, так как они при совмещении их центров совмещаются всеми своими точками.
Прямая ($MN$, рис.~\ref{1938/ris-6}), проходящая через какие-нибудь две точки окружности, называется \rindex{секущая}\textbf{секущей}.

Отрезок ($EF$), соединяющий две какие-нибудь точки окружности, называется \rindex{хорда}\textbf{хордой}.

Всякая хорда ($AD$), проходящая через центр, называется \rindex{диаметр}\textbf{диаметром}.
Диаметр равен сумме двух радиусов, и потому все диаметры одной окружности равны между собой.
Какая-нибудь часть окружности (например, $EmF$) называется \rindex{дуга}\textbf{дугой}.

О хорде ($EF$), соединяющей концы какой-нибудь дуги, говорят, что она \textbf{стягивает} эту дугу.

Дуга обозначается иногда знаком $\smallsmile$;
например, вместо «дуга $EmF$» пишут «${\smallsmile} EmF$».
Часть плоскости, ограниченная окружностью, называется \rindex{круг}\textbf{кругом}%
\footnote{Иногда слово «круг» употребляют в том же смысле, как и окружность.
Но этого следует избегать, так как употребление одного и того же термина для разных понятий может приводить к ошибкам.}%
.

Часть круга, заключённая между двумя радиусами (часть $COB$, покрытая штрихами на рис.~\ref{1938/ris-6}), называется \rindex{сектор}\textbf{сектором}, а часть, отсекаемая от круга какой-нибудь секущей (часть $EmF$), называется \rindex{сегмент}\textbf{сегментом}.

\paragraph{Равенство и неравенство дуг.}\label{1938/10}
Две дуги одной и той же окружности (или равных окружностей) равны между собой, если они могут быть совмещены так, что их концы совпадут.
Положим, например, что мы накладываем дугу $AB$ (рис.~\ref{1938/ris-7}) на дугу $CD$ так, чтобы точка $A$ совпала с точкой $C$ и дуга $AB$ пошла по дуге $CD$;
если при этом концы $B$ и $D$ совпадут, то совпадут и все промежуточные точки этих дуг, так как они находятся на одинаковом расстоянии от центра, значит, ${\smallsmile} AB={\smallsmile} CD$;
если же $B$ и $D$ не совпадут, то дуги не равны, причём та считается меньше, которая составит часть другой.

\begin{wrapfigure}{r}{39 mm}
\vskip-4mm
\centering
\includegraphics{mppics/ris-7}
\caption{}\label{1938/ris-7}
\end{wrapfigure}

\paragraph{Сумма дуг.}\label{1938/11}
\rindex{сумма!дуг}
Суммой нескольких данных дуг одинакового радиуса называется такая дуга того же радиуса, которая составлена из частей, соответственно равных данным дугам.
Так, если от произвольной точки $M$ (рис.~\ref{1938/ris-7}) окружности отложим часть $MN$, равную $AB$, и затем от точки $N$ в том же направлении отложим часть $NP$, равную $CD$, то дуга $MP$ будет суммой дуг $AB$ и $CD$.
Подобно этому можно составить сумму трёх и более дуг.

При сложении дуг одинакового радиуса их сумма может не уместиться на одной окружности, одна из дуг может частично покрыть другую.
В таком случае суммой дуг будет являться дуга, б\'{о}льшая целой окружности.
Так, например, при сложении дуги $AmB$ с дугой $CnD$ (рис.~\ref{1938/ris-8}) получаем дугу, состоящую из целой окружности и дуги $AD$.

\begin{figure}[!ht]
\centering
\includegraphics{mppics/ris-8}
\caption{}\label{1938/ris-8}
\end{figure}

Сумма дуг, как и сумма отрезков, обладает свойствами \textbf{переместительным} и \textbf{сочетательным}.

Из понятия о сумме дуг выводятся понятия о разности дуг, умножении и делении дуги на число, так же как и для отрезков.

\paragraph{Разделение геометрии.}\label{1938/12}
Геометрия разделяется на две части:
\textbf{планиметрию} и \textbf{стереометрию}.
Первая рассматривает свойства таких фигур, все части которых помещаются на одной плоскости;
вторая — свойства таких фигур, у которых не все части помещаются на одной плоскости.

%% file: 2D/ugly.tex
\section{Углы} 

\subsection*{Предварительные понятия}

\paragraph{Угол.}\label{1938/13}
Фигура, образованная двумя лучами, исходящими из одной точки, называется \rindex{угол}\textbf{углом}.
Полупрямые, образующие угол, называются \rindex{сторона!угла}\textbf{сторонами}, а точка, из которой они исходят, — \rindex{вершина!угла}\textbf{вершиной угла}.
Стороны следует представлять себе неограниченно продолженными от вершины.

\begin{wrapfigure}{o}{31 mm}
\centering
\includegraphics{mppics/ris-9}
\caption{}\label{1938/ris-9}
\end{wrapfigure}

Угол обыкновенно обозначается тремя большими буквами, из которых средняя ставится у вершины, а крайние — у каких-нибудь точек сторон;
например, говорят:
«угол $AOB$» или «угол $BOA$» (рис.~\ref{1938/ris-9}).
Но можно обозначить угол и одной буквой, поставленной у вершины, если при этой вершине нет других углов.
Мы иногда будем обозначать угол цифрой, поставленной внутри угла около вершины.

Стороны угла разделяют всю плоскость, в которой лежит этот угол, на две области.
Одна из них называется \rindex{внутренняя область угла}\textbf{внутренней областью угла}, другая — \rindex{внешняя область угла}\textbf{внешней его областью}.
Обычно внутренней областью считается та, в которой целиком помещается отрезок, соединяющий две любые точки, взятые на сторонах угла, например точки $A$ и $B$ на сторонах угла $AOB$ (рис.~\ref{1938/ris-9}).
Но иногда приходится считать внутренней областью угла другую часть плоскости.
В этих случаях обычно делается специальное указание, какая область плоскости считается внутренней областью угла.

На рис.~\ref{1938/ris-10} представлены раздельно оба случая.
Внутренней областью угла в каждом случае служит заштрихованная часть плоскости.
Если из вершины угла (рис.~\ref{1938/ris-9}) проведём внутри него какие-нибудь прямые $OD, OE,\dots$, то образовавшиеся при этом углы $AOD$, $DOE$, $EOB,\dots$, рассматриваются как части угла $AOB$.

\begin{figure}[!ht]
\centering
\includegraphics{mppics/ris-10}
\caption{}\label{1938/ris-10}
\end{figure}

Слово «угол» при записи заменяется часто знаком $\angle$.
Например, вместо «угол $AOB$» обычно пишут:
$\angle AOB$.

\paragraph{Равенство и неравенство углов.}\label{1938/14}
В соответствии с общим определением равенства геометрических фигур (§~\ref{1938/1}) \emph{два угла считаются равными, если при наложении они могут совместиться.}

\begin{figure}[!ht]
\centering
\includegraphics{mppics/ris-11}
\caption{}\label{1938/ris-11}
\end{figure}

{\sloppy
Положим, например, что мы накладываем угол $AOB$ на угол $A_1O_1B_1$ (рис.~\ref{1938/ris-11}) так, чтобы вершина $O$ совпала с $O_1$, сторона $OB$ пошла по $O_1B_1$ и чтобы внутренние области обоих углов были расположены по одну сторону от прямой $O_1B_1$.
Если при этом сторона $OA$ совместится с $O_1A_1$, то углы равны;
если же $OA$ пойдёт внутри угла $A_1O_1B_1$ или же вне его, то углы не равны, причём тот из них будет меньше, который составит часть другого угла.

}

\begin{figure}[!ht]
\centering
\includegraphics{mppics/ris-12}
\caption{}\label{1938/ris-12}
\end{figure}

\paragraph{Сумма углов.}\label{1938/15}
\rindex{сумма!углов}
Суммой углов $AOB$ и $A_1O_1B_1$ (рис.~\ref{1938/ris-12}) называется такой угол, который получится следующим образом.
Строим угол $MNP$, равный первому данному углу $AOB$, и к нему пристраиваем угол $PNQ$, равный другому данному углу $A_1O_1B_1$ так, чтобы у обоих углов оказалась общая вершина $N$ и общая сторона $NP$ и чтобы внутренние области углов были расположены по разные стороны от общей стороны $NP$.
Тогда угол $MNQ$ называется суммой углов $AOB$ и $A_1O_1B_1$.
Внутренней областью этого угла служит та область плоскости, которая образована совокупностью внутренних областей складываемых углов.
Это та область, в которой лежит общая сторона ($NP$) складываемых углов.
Подобным образом может быть составлена сумма трёх и более углов.

\begin{wrapfigure}{r}{35 mm}
\vskip-4mm
\centering
\includegraphics{mppics/ris-ru-13}
\caption{}\label{1938/ris-13}
\vskip-4mm
\end{wrapfigure}

Сумма углов, как и сумма отрезков, обладает свойствами переместительным и сочетательным.
Часто приходится говорить о такой полупрямой, которая делит данный угол пополам;
такой полупрямой дали особое название:
\rindex{биссектриса}\textbf{биссектриса} (рис.~\ref{1938/ris-13}).

\paragraph{Расширение понятия об угле.}\label{1938/16}
При нахождении суммы углов могут представиться некоторые особые случаи, которые полезно рассмотреть отдельно.

1.
Может случиться, что после сложения нескольких углов, например трёх:
$AOB$, $BOC$ и $COD$ (рис.~\ref{1938/ris-14}), сторона $OD$ угла $COD$ составит продолжение стороны $OA$ угла $AOB$.
Мы получим тогда фигуру, образованную двумя полупрямыми ($OA$ и $OD$), исходящими из одной точки ($O$) и составляющими продолжение одна другой.
Такую фигуру принято тоже называть углом (\rindex{развёрнутый угол}\textbf{развёрнутым}).

\begin{figure}[!ht]
\begin{minipage}{.48\textwidth}
\centering
\includegraphics{mppics/ris-14}
\end{minipage}\hfill
\begin{minipage}{.48\textwidth}
\centering
\includegraphics{mppics/ris-15}
\end{minipage}

\medskip

\begin{minipage}{.48\textwidth}
\centering
\caption{}\label{1938/ris-14}
\end{minipage}\hfill
\begin{minipage}{.48\textwidth}
\centering
\caption{}\label{1938/ris-15}
\end{minipage}
\vskip-4mm
\end{figure}

2.
Может случиться, что после сложения нескольких углов, например пяти углов:
$AOB$, $BOC$, $COD$, $DOE$ и $EOA$ (рис.~\ref{1938/ris-15}), сторона $OA$ угла $EOA$ совместится со стороной $OA$ угла $AOB$.

Фигура, образованная такими совпавшими полупрямыми (рассматриваемая вместе со всей плоскостью, расположенной вокруг общей вершины $O$), также называется углом (\rindex{полный угол}\textbf{полным}).

3.
Наконец, может случиться, что, строя сумму углов, мы не только заполним всю плоскость вокруг их общей вершины, но даже будем вынуждены налагать углы один на другой, покрывая плоскость вокруг общей вершины во второй раз, в третий раз и~т.~д.
Такая сумма углов равна одному полному углу, сложенному с некоторым углом, или равна двум полным углам, сложенным с некоторым углом, и~т.~д.

\subsection*{Измерение углов}

{

\begin{wrapfigure}{r}{30 mm}
\vskip-2mm
\centering
\includegraphics{mppics/ris-16}
\caption{}\label{1938/ris-16}
\bigskip
\includegraphics{mppics/ris-17}
\caption{}\label{1938/ris-17}
\end{wrapfigure}

\paragraph{Центральный угол.}\label{1938/17}
Угол ($AOB$, рис.~\ref{1938/ris-16}), образованный двумя радиусами окружности, называется \rindex{центральный угол}центральным углом;
о таком угле и дуге, заключённой между его сторонами, говорят, что они \textbf{соответствуют} друг другу.

Центральные углы по отношению к соответствующим им дугам обладают следующими двумя свойствами.

\textbf{\emph{В одном круге или в равных кругах.}}

1) \textbf{\emph{Если центральные углы равны, то и соответствующие им дуги равны.}}

2) Обратно.
\textbf{\emph{Если дуги равны, то и соответствующие им центральные углы равны.}}

Пусть $\angle AOB=\angle COD$ (рис.~\ref{1938/ris-17}), покажем, что дуги $AB$ и $CD$ также равны.
Вообразим, что сектор $AOB$ мы повернули вокруг центра $O$ в направлении, указанном стрелкой, настолько, чтобы радиус $OA$ совпал с $OC$.
Тогда, вследствие равенства углов, радиус $OB$ совместится с $OD$;
значит, совместятся и дуги $AB$ и $CD$, то есть
они будут равны.

Второе свойство также легко обнаружить наложением.

}

\paragraph{Градусы дуговой и угловой.}\label{1938/18}
Вообразим, что какая-нибудь окружность разделена на 360 равных частей и ко всем точкам деления проведены радиусы.
Тогда вокруг центра образуются 360 центральных углов, которые, как соответствующие равным дугам, должны быть равны между собой.
Каждая из полученных таким образом на окружности дуг называется \rindex{градус}\textbf{дуговым градусом}, а каждый из образовавшихся при центре углов называется \textbf{угловым градусом}.
Значит, можно сказать, что дуговой градус есть $\tfrac1{360}$ часть окружности,
а угловой градус есть центральный угол, соответствующий дуговому градусу.

Градусы (дуговые и угловые) подразделяются ещё на 60 равных частей, называемых \rindex{минута}\textbf{минутами}, а минуты подразделяются ещё на 60 равных частей, называемых \rindex{секунда}\textbf{секундами}.

\begin{figure}[!ht]
\centering
\includegraphics{mppics/ris-18}
\caption{}\label{1938/ris-18}
\end{figure}

Пусть $AOB$ есть какой-нибудь угол (рис.~\ref{1938/ris-18}).
Опишем между его сторонами с центром в вершине $O$ дугу $CD$ произвольного радиуса;
тогда угол $AOB$ будет центральным углом, соответствующим дуге $CD$.
При этом угол измеряется соответствующей ему дугой; то есть в угле содержится столько угловых градусов, минут и секунд, сколько в соответствующей ему дуге содержится дуговых градусов, минут и секунд.
Если, например, в дуге $CD$ содержится $20{,}57$ градусов или, что то же самое, 20 градусов 34 минуты 12 секунд,%
\footnote{Поскольку \[20{,}57=20+\tfrac{57}{100}=20+\tfrac{34}{60}+\tfrac{12}{60\cdot 60}.\]}
 то и в угле $AOB$ заключается 20 градусов 34 минуты 12 секунд угловых, что принято выражать так:
$\angle AOB = 20{,}57\degree = 20\degree 34'12''$ обозначая знаками $\degree$, $'$ и $''$ соответственно градусы, минуты и секунды.

Величина углового градуса не зависит от радиуса окружности.
Действительно, если сложить 360 угловых градусов по правилу, указанному в §~\ref{1938/15}, то получим полный угол при центре окружности.
Каков бы ни был радиус окружности, этот полный угол будет один и тот же.
Значит, можно сказать, что угловой градус составляет $\tfrac1{360}$ часть полного угла;
это мера угла, определяющая его величину независимо от радиуса окружности.
Число угловых градусов в данном угле принимают за меру наклона одной стороны угла к другой.

Поскольку число 360 имеет много делителей,
градусы удобны для работы с дробными долями полного угла.
Тем не менее основной единицей измерения углов считается \textbf{радиан};
он будет определён в §~\ref{extra/radians}. 
Кроме градусов и радиан употребляются и другие единицы для измерения углов, например \textbf{прямой угол}, \textbf{оборот} равный полному углу и
\textbf{град} равный сотой доле прямого угла.

\begin{figure}[!ht]
\centering
\ \ \ \ \includegraphics{mppics/transportir-19}
\caption{}\label{1938/ris-19}
\end{figure}

\paragraph{Транспортир.}\label{1938/20}
Для измерения углов употребляется особый прибор — \rindex{транспортир}\textbf{транспортир}.
Этот прибор (рис.~\ref{1938/ris-19}) представляет собой полукруг, дуга которого разделена на $180\degree $.
Чтобы измерить угол $DCE$, накладывают на него прибор так, чтобы центр полукруга совпадал с вершиной угла, а радиус $CB$ был расположен по стороне $CE$.
Тогда число градусов, содержащееся в дуге, заключённой между сторонами угла $DCE$, покажет его величину.
При помощи транспортира можно также начертить угол, содержащий данное число градусов.

\paragraph{Прямой, острый и тупой углы.}\label{1938/21}
Угол в $90\degree$ (составляющий, следовательно, половину развёрнутого угла или четверть полного угла) называют \rindex{прямой угол}\textbf{прямым углом};
угол, меньший прямого, называют \rindex{острый угол}\textbf{острым}, а угол, больший прямого, но меньший развёрнутого, называют \rindex{тупой угол}\textbf{тупым} (рис.~\ref{1938/ris-20})

\begin{figure}[!ht]
\centering
\includegraphics{mppics/ris-20}
\caption{}\label{1938/ris-20}
\end{figure}

Конечно, \textbf{\emph{все прямые углы}}, как содержащие одинаковое число градусов, \textbf{\emph{равны между собой}}.
На чертежах прямой угол принято обозначать квадратиком, как показано на рис.~\ref{1938/ris-20}.

\subsection*{Смежные и вертикальные углы}

\paragraph{Смежные углы и их свойства.}\label{1938/22}
Два угла ($AOB$ и $BOC$, рис.~\ref{1938/ris-21}) называются \rindex{смежные углы}\textbf{смежными}, если одна сторона у них общая, а две другие стороны составляют продолжение одна другой.

Так как такие углы в сумме составляют развёрнутый угол, то \textbf{\emph{сумма двух смежных углов равна $\bm{180\degree}$}}.

\begin{figure}[!ht]
\begin{minipage}{.48\textwidth}
\centering
\includegraphics{mppics/ris-21}
\end{minipage}\hfill
\begin{minipage}{.48\textwidth}
\centering
\includegraphics{mppics/ris-22}
\end{minipage}

\medskip

\begin{minipage}{.48\textwidth}
\centering
\caption{}\label{1938/ris-21}
\end{minipage}\hfill
\begin{minipage}{.48\textwidth}
\centering
\caption{}\label{1938/ris-22}
\end{minipage}
\vskip-4mm
\end{figure}

Для каждого данного угла можно построить два смежных с ним угла.
Например, для угла $AOB$ (рис.~\ref{1938/ris-22}), продолжив сторону $AO$, мы получим один смежный угол $BOC$.
\emph{Два угла, $BOC$ и $AOD$, смежные с одним и тем же углом $AOB$, равны между собой,} так как каждый из них дополняет угол $AOB$ до $180\degree$.

\begin{wrapfigure}{r}{35 mm}
\vskip-0mm
\centering
\includegraphics{mppics/ris-23}
\caption{}\label{1938/ris-23}
\end{wrapfigure}

Если угол $AOB$ прямой (рис.~\ref{1938/ris-23}), то есть если он содержит $90\degree$, то и каждый из смежных с ним углов $COB$ и $AOD$ должен быть также прямой, так как он содержит в себе $180\degree-90\degree$, то есть $90\degree$;
четвёртый угол $COD$ тоже должен быть прямым, так как три угла $AOB$, $BOC$ и $AOD$ составляют в сумме $270\degree$ и, следовательно, от $360\degree$ на долю четвёртого угла $COB$ остаётся тоже $90\degree$.
Таким образом:
\emph{если при пересечении двух прямых \emph{($AC$ и $BD$, рис.~\ref{1938/ris-23})} один из четырёх углов окажется прямой, то и остальные три угла должны быть прямые.}

\paragraph{Перпендикуляр и наклонная.}\label{1938/23}
Общая сторона ($OB$) двух смежных углов называется \rindex{наклонная}\textbf{наклонной} к прямой ($AC$), на которой лежат две другие стороны, в том случае, когда смежные углы не равны между собой (рис.~\ref{1938/ris-24});
в том же случае, когда смежные углы равны (рис.~\ref{1938/ris-25}) и когда, следовательно, каждый из углов есть прямой, общая сторона называется \rindex{перпендикулярность}\textbf{перпендикуляром} к прямой, на которой лежат две другие стороны.
Общая вершина ($O$) в первом случае называется \rindex{основание!наклонной}\textbf{основанием наклонной}, во втором случае — \rindex{основание!перпендикуляра}\textbf{основанием перпендикуляра}.

\begin{figure}[!ht]
\begin{minipage}{.48\textwidth}
\centering
\includegraphics{mppics/ris-24}
\end{minipage}\hfill
\begin{minipage}{.48\textwidth}
\centering
\includegraphics{mppics/ris-25}
\end{minipage}

\medskip

\begin{minipage}{.48\textwidth}
\centering
\caption{}\label{1938/ris-24}
\end{minipage}\hfill
\begin{minipage}{.48\textwidth}
\centering
\caption{}\label{1938/ris-25}
\end{minipage}
\vskip-4mm
\end{figure}

Две прямые ($AC$ и $BD$, рис.~\ref{1938/ris-23}), пересекающиеся между собой под прямым углом, называются взаимно \rindex{перпендикулярность}\textbf{перпендикулярными}.
Что прямая $AC$ перпендикулярна к прямой $BD$, записывают так: $AC\perp BD$.

{\small

\smallskip

\mbox{\so{Замечания}.}
1) Если перпендикуляр к прямой $AC$ (рис.~\ref{1938/ris-25}) приходится проводить из точки $O$, лежащей на этой прямой, то говорят, что этот перпендикуляр надо  \so{восстановить} к прямой $AC$, а если требуется перпендикуляр провести из точки $B$, лежащей вне прямой, то говорят, что его надо \so{опустить} на прямую (всё равно: вниз или вверх, или вбок).

2) Очевидно, что из всякой точки данной прямой можно к этой прямой восстановить перпендикуляр и притом только один.

}

\begin{wrapfigure}{r}{38 mm}
\centering
\includegraphics{mppics/ris-26}
\caption{}\label{1938/ris-26}
\end{wrapfigure}

\paragraph{}\label{1938/24}
Докажем, что \textbf{\emph{из всякой точки, лежащей вне прямой, можно опустить на эту прямую перпендикуляр и только один}}.

Пусть дана какая-нибудь прямая $AB$ (рис.~\ref{1938/ris-26}) и вне её произвольная точка $M$;
требуется показать, что, во-первых, из этой точки можно опустить на прямую $AB$ перпендикуляр и, во-вторых, что этот перпендикуляр может быть только один.

Вообразим, что мы перегнули чертёж по прямой $AB$ так, чтобы верхняя его часть упала на нижнюю.
Тогда точка $M$ займёт некоторое положение $N$.
Отметив это положение, приведём чертёж в прежний вид и затем соединим точки $M$ и $N$ прямой.
Убедимся, что полученная прямая $MN$ перпендикулярна к $AB$, а всякая иная прямая, исходящая из $M$, например $MD$, не перпендикулярна к $AB$.
Для этого перегнём чертёж вторично.
Тогда точка $M$ снова совместится с $N$, а точки $C$ и $D$ останутся на своих местах;
следовательно, прямая $MC$ совпадёт с $NC$, а $MD$ с $ND$.
Из этого следует, что $\angle MCB = \angle BCN$ и $\angle MDC = \angle CDN$.

Но углы $MCB$ и $BCN$ смежные и, как теперь видим, равные;
следовательно, каждый из них есть \so{прямой}, а потому $MN\z\perp AB$.
Так как линия $MDN$ не прямая (потому что не может быть двух различных прямых, соединяющих точки $M$ и $N$), то сумма двух равных углов $MDC$ и $CDN$ \so{не равна} $180\degree$;
поэтому угол $MDC$ не есть прямой, и, значит, $MD$ не перпендикулярна к $AB$.
Таким образом, другого перпендикуляра из точки $M$ на прямую $AB$ опустить нельзя.

\paragraph{Чертёжный угольник.}\label{1938/25} 
Для построения перпендикуляра к данной прямой очень удобен угольник — чертёжный инструмент в виде треугольника у которого один из углов делается прямым.
Чтобы провести перпендикуляр к прямой $AB$ (рис.~\ref{1938/ris-27}) через точку $D$, взятую вне прямой, приставляют линейку к прямой $AB$, к линейке угольник, а затем, придерживая линейку рукой, двигают угольник вдоль линейки до тех пор, пока другая сторона прямого угла не пройдёт через точку $D$, затем проводят прямую $CE$.
Точно также поступают если точка $D$ лежит на прямой $AB$.

\begin{figure}[!ht]
\begin{minipage}{.48\textwidth}
\centering
\includegraphics{mppics/ris-wood-27}
\end{minipage}\hfill
\begin{minipage}{.48\textwidth}
\centering
\includegraphics{mppics/ris-28}
\end{minipage}

\medskip

\begin{minipage}{.48\textwidth}
\centering
\caption{}\label{1938/ris-27}
\end{minipage}\hfill
\begin{minipage}{.48\textwidth}
\centering
\caption{}\label{1938/ris-28}
\end{minipage}
\vskip-4mm
\end{figure}

\paragraph{Вертикальные углы и их свойство.}\label{1938/26}
Два угла называются \rindex{вертикальный!угол}\textbf{вертикальными}, если стороны одного составляют продолжение сторон другого.
Так, при пересечении двух прямых $AB$ и $CD$ (рис.~\ref{1938/ris-28}) образуются две пары вертикальных углов:
$AOD$ и $COB$, $AOC$ и $DOB$ (и четыре пары смежных углов).

{\sloppy

\textbf{\emph{Два вертикальных угла равны между собой}} (например, $\angle AOD \z= \angle BOC$), так как каждый из них есть смежный с одним и тем же углом (с~$\angle DOB$ или с~$\angle AOC$), а такие углы, как мы видели (§~\ref{1938/22}), равны друг другу.

}

\paragraph{Замечания об углах, имеющих общую вершину.}\label{1938/27}
Об углах, имеющих общую вершину, полезно помнить следующие простые истины.

\begin{figure}[!ht]
\begin{minipage}{.48\textwidth}
\centering
\includegraphics{mppics/ris-29}
\end{minipage}\hfill
\begin{minipage}{.48\textwidth}
\centering
\includegraphics{mppics/ris-30}
\end{minipage}

\medskip

\begin{minipage}{.48\textwidth}
\centering
\caption{}\label{1938/ris-29}
\end{minipage}\hfill
\begin{minipage}{.48\textwidth}
\centering
\caption{}\label{1938/ris-30}
\end{minipage}
\vskip-4mm
\end{figure}

1) \emph{Если сумма нескольких углов ($AOB$, $BOC$, $COD$, $DOE$, рис.~\ref{1938/ris-29}), имеющих общую вершину, составляет развёрнутый угол, то она равна $180\degree$.}

2) \emph{Если сумма нескольких углов ($AOB$, $BOC$, $COD$, $DOE$, $EOA$ рис.~\ref{1938/ris-30}), имеющих общую вершину, составляет полный угол, то она равна  $360\degree$.}

3) \emph{Если два угла ($AOB$ и $BOC$, рис.~\ref{1938/ris-24}) имеют общую вершину ($O$) и общую сторону ($OB$) и в сумме составляют $180\degree$, то их две другие стороны ($AO$ и $OC$) составляют продолжение одна другой} (то есть такие углы будут смежными).

{\small

\subsection*{Упражнения}

\begin{enumerate}[noitemsep]
\item
Некоторый угол равен $38\degree 20'$;
найти величину смежного с ним угла.

\item
Два угла $ABC$ и $CBD$, имея общую вершину $B$ и общую сторону $BC$, расположены так, что они не покрывают друг друга;
угол $ABC \z= 100\degree20'$ , а угол $CBD = 79\degree 40'$.
Составляют ли стороны $AB$ и $BD$ прямую или ломаную.

\item
Построить какой-нибудь угол и при помощи транспортира и линейки провести его биссектрису.

\end{enumerate}

\smallskip
\so{Доказать}, что:

\begin{enumerate}[resume,noitemsep]
\item
Биссектрисы двух смежных углов взаимно перпендикулярны.

\item
Биссектрисы двух вертикальных углов составляют продолжение одна другой.

\item
Если при точке $O$ прямой $AB$ (рис.~\ref{1938/ris-28}) построим по разные стороны от $AB$ равные углы $AOD$ и $BOC$, то стороны их $OD$ и $OC$ составляют одну прямую.

\item
Если полупрямые $OA$, $OD$, $OB$, $OC$ (рис.~\ref{1938/ris-28})  разделяют плоскость на четыре угла такие, что $\angle AOC = \angle DOB$ и $\angle AOD=\angle COB$, то $OD$ есть продолжение $OC$ и $OB$ есть продолжение $OA$.

\smallskip
\so{Указание}.
Надо применить §~\ref{1938/27}, 2 и 3.

\end{enumerate}

}

%% file: 2D/mat-predlozheniya.tex
\section{Математические предложения}

\paragraph{Теоремы, аксиомы, определения.}\label{1938/28}
Из того, что было изложено, можно заключить, что некоторые геометрические истины мы считаем вполне очевидными (например, свойства плоскости и прямой в §§~\ref{1938/3} и \ref{1938/4}), а другие устанавливаем путём рассуждений (например, свойства смежных углов в §~\ref{1938/22} и вертикальных в §~\ref{1938/26}).
Такие рассуждения являются в геометрии главным средством обнаружить свойства геометрических фигур.
Поэтому для дальнейшего полезно заранее познакомиться с теми видами рассуждений, которые применяются в геометрии.
Все истины, которые устанавливаются в геометрии, выражаются в виде предложений.

Эти предложения бывают следующих видов.

\textbf{Определения.}\rindex{определение}
Определениями называют предложения, в которых разъясняется, какой смысл придают тому или другому названию или выражению.
Например, мы уже встречали определения центрального угла, прямого угла и перпендикуляра.

\textbf{Аксиомы.}\rindex{аксиома}
Аксиомами называют истины, которые принимаются без доказательства.
Таковы, например, предложения, встречавшиеся нам ранее (§~\ref{1938/4}):
через всякие две точки можно провести прямую и притом только одну;
если две точки прямой лежат в данной плоскости, то и все точки этой прямой лежат в той же плоскости.

Укажем ещё следующие аксиомы, относящиеся ко всякого рода величинам:

если две величины равны порознь одной и той же третьей величине, то они равны и между собой.

если к равным величинам прибавим поровну или от равных величин отнимем поровну, то равенство не нарушится.

если к неравным величинам прибавим поровну или от неравных величин отнимем поровну, то смысл неравенства не изменится, то есть б\'{о}льшая величина останется б\'{о}льшей.

\textbf{Теоремы.}\rindex{теорема}
Теоремами называются предложения, истинность которых обнаруживается только после некоторого рассуждения (доказательства).
Примером могут служить следующие предложения.

если в одном круге или в равных кругах центральные углы равны, то и соответствующие им дуги равны.

если при пересечении двух прямых между собой один из четырёх углов окажется прямой, то и остальные три угла прямые.

\textbf{Следствия.}\rindex{следствие}
Следствиями называются предложения, которые составляют непосредственный вывод из аксиомы или из теоремы.
Например, из аксиомы:
«через две точки можно провести только одну прямую» следует, что «две прямые могут пересечься только в одной точке».

\paragraph{Состав теоремы.}\label{1938/29}
Во всякой теореме можно различить две части:
условие и заключение.
\textbf{Условие} выражает то, что предполагается данным;
\textbf{заключение} — то, что требуется доказать.
Например, в теореме:
«если центральные углы равны, то и соответствующие им дуги равны» условием служит первая часть теоремы:
«если центральные углы равны», а заключением — вторая часть:
«то и соответствующие им дуги равны»;
другими словами, нам дано (нам известно), что центральные углы равны, а требуется доказать, что при этом условии и соответствующие дуги также равны.

Условие и заключение теоремы могут иногда состоять из нескольких отдельных условий и заключений;
например, в теореме:
«если число делится на 2 и на 3, то оно разделится и на 6» условие состоит из двух частей:
«если число делится на 2» и «если число делится на 3».

Полезно заметить, что всякую теорему можно подробно выразить словами так, что её условие будет начинаться словом «если», а заключение — словом «то».
Например, теорему:
«вертикальные углы равны» можно подробнее высказать так:
«если два угла вертикальные, то они равны».

\paragraph{Обратная теорема.}\label{1938/30}\rindex{обратная теорема}
Теоремой, обратной данной теореме, называется такая, в которой условием поставлено заключение (или часть заключения), а заключением — условие (или часть условия) данной теоремы.
Например, следующие две теоремы обратны друг другу.

\medskip

{\sloppy

\columnratio{0.5}
\setlength{\columnseprule}{.2pt}
\begin{paracol}{2}
\textbf{\emph{Если центральные углы равны, то и соответствующие им дуги равны.}}
\switchcolumn
\textbf{\emph{Если дуги равны, то и соответствующие им центральные углы равны.}}
\end{paracol}

}

\medskip

Если одну из этих теорем назовём \textbf{прямой}, то другую следует назвать обратной.
В этом примере обе теоремы, и прямая, и обратная, оказываются верными.
Но так бывает не всегда.
Например, теорема:
«если два угла вертикальные, то они равны» верна, но обратное предложение:
«если два угла равны, то они вертикальные» неверно.

В самом деле, допустим, что в каком-либо углу проведена его биссектриса (рис.~\ref{1938/ris-13}).
Она разделит данный угол на два меньших угла.
Эти углы будут равны между собой, но они не будут вертикальными.

{\sloppy

\paragraph{Противоположная теорема.}\label{1938/31}\rindex{противоположная теорема}
Теоремой, противоположной данной теореме, называется такая, условие и заключение которой представляют отрицание условия и заключения данной теоремы.
Например, теореме:
«если сумма цифр делится на 9, то число делится на 9» соответствует такая противоположная:
«если сумма цифр не делится на 9, то число не делится на 9».

}

Заметим, что верность прямой теоремы ещё не служит доказательством верности противоположной:
например, противоположное предложение:
«если каждое слагаемое не делится на одно и то же число, то и сумма не разделится на это число» — неверно, тогда как прямое предложение верно.

\paragraph{Зависимость между теоремами: прямой, обратной и противоположной.}\label{1938/32}
Для лучшего уяснения этой зависимости выразим теоремы сокращённо так (буквой $A$ мы обозначим условие теоремы, а буквой $B$ — её заключение).

1) \textbf{Прямая:}
если есть $A$, то есть и $B$.

2) \textbf{Обратная:}
если есть $B$, то есть и $A$.

3) \textbf{Противоположная прямой:}
если нет $A$, то нет и $B$.

4) \textbf{Противоположная обратной:}
если нет $B$, то нет и $A$.

Рассматривая эти предложения, легко заметить, что первое из них находится в таком же отношении к четвёртому, как второе к третьему, а именно:
предложения первое и четвёртое обратимы одно в другое, равно как второе и третье.
Действительно, из предложения:
«если есть $A$, то есть и $B$» непосредственно следует:
«если нет $B$, то нет и $A$» (так как если бы $A$ было, то, согласно первому предложению, было бы и $B$);
обратно, из предложения:
«если нет $B$, то нет и $A$» выводим:
«если есть $A$, то есть и $B$» (так как если бы $B$ не было, то не было бы и $A$).
Совершенно так же убедимся, что из второго предложения следует третье, и наоборот.

Таким образом, чтобы иметь уверенность в справедливости всех четырёх теорем, нет надобности доказывать каждую из них отдельно, а достаточно ограничиться доказательством только двух:
прямой и обратной, или прямой и противоположной.

%% file: 2D/treugi.tex
\section{Треугольники}

\subsection*{Понятие о многоугольнике и треугольнике}

\paragraph{Ломаная линия.}\label{1938/33}
Линия, образуемая отрезками, не лежащими на одной прямой и расположенными так, что конец первого служит началом второго, конец второго — началом третьего и~т.~д., называется \rindex{ломаная}\textbf{ломаной линией} (рис.~\ref{1938/ris-31} и \ref{1938/ris-32}).
Эти отрезки называются \rindex{сторона!ломаной}\textbf{сторонами} ломаной или её \rindex{звено}\textbf{звеньями}, а вершины углов, образуемых соседними отрезками, — \rindex{вершина!ломаной}\textbf{вершинами} её.
Ломаная линия обозначается рядом букв, поставленных у её вершин и концов;
например говорят:
ломаная $ABCDE$.
Ломаная линия называется \rindex{выпуклая ломаная}\textbf{выпуклой}, если она вся расположена по одну сторону от каждого входящего в её состав отрезка, продолженного неограниченно в обе стороны.
Такова, например, линия, изображённая на рис.~\ref{1938/ris-31}, тогда как ломаная на рис.~\ref{1938/ris-32} не будет выпуклой (она расположена не по одну сторону от прямой $BC$), рис.~\ref{1938/ris-31}.

\begin{figure}[h]
\begin{minipage}{.48\textwidth}
\centering
\includegraphics{mppics/ris-31}
\end{minipage}\hfill
\begin{minipage}{.48\textwidth}
\centering
\includegraphics{mppics/ris-32}
\end{minipage}

\medskip

\begin{minipage}{.48\textwidth}
\centering
\caption{}\label{1938/ris-31}
\end{minipage}\hfill
\begin{minipage}{.48\textwidth}
\centering
\caption{}\label{1938/ris-32}
\end{minipage}
\vskip-4mm
\end{figure}

Когда концы ломаной сходятся в одну точку, то она называется \rindex{замкнутая ломаная}\textbf{замкнутой} (например, линия $ABCDE$ на рис.~\ref{1938/ris-33}).

\paragraph{Многоугольник.}\label{1938/34}\rindex{многоугольник}
Фигура, образованная замкнутой ломаной линией вместе с частью плоскости, ограниченной этой линией, называется многоугольником (рис.~\ref{1938/ris-33}).
Стороны этой ломаной называются \rindex{сторона!многоугольника}\textbf{сторонами} многоугольника, углы, составленные каждыми двумя соседними сторонами, — \rindex{угол!многоугольника}\textbf{углами} многоугольника, а их вершины — \rindex{вершина!многоугольника}\textbf{вершинами} его.

\begin{figure}[!ht]
\centering
\includegraphics{mppics/ris-33}
\caption{}\label{1938/ris-33}
\end{figure}

При этом внутренней областью угла многоугольника считается (рис. \ref{1938/ris-33}) та, к которой принадлежит непосредственно примыкающая к вершине внутренняя область самого многоугольника.
Так, для многоугольника $MNPQRS$ (рис.~\ref{1938/ris-33}) углом при вершине $P$ является угол, больший двух прямых (с заштрихованной внутренней областью).
Сама ломаная линия, ограничивающая многоугольник, называется \rindex{контур многоугольника}\textbf{контуром} его, а отрезок, равный сумме всех его сторон, — \rindex{периметр}\textbf{периметром}.

Многоугольник называется \rindex{выпуклый многоугольник}\textbf{выпуклым}, если он ограничен выпуклой ломаной линией;
таков, например, многоугольник $ABCDE$, изображённый на рис.~\ref{1938/ris-33} (многоугольник $MNPQRS$ нельзя назвать выпуклым);
мы будем рассматривать, главным образом, выпуклые многоугольники.

Всякая прямая (как $AD$, $BE$, $MR,\dots$, рис.~\ref{1938/ris-33}), которая соединяет вершины двух углов многоугольника, не прилежащих к одной стороне, называется \rindex{диагональ}\textbf{диагональю} многоугольника.

Наименьшее число сторон в многоугольнике — три.
По числу сторон многоугольник называется \textbf{треугольником}, \textbf{четырёхугольником}, \textbf{пятиугольником} и~т.~д. Многоугольник с $n$ сторонами называется \rindex{$n$-угольник}\textbf{$\bm{n}$-угольником}. 

Для краткости «треугольник» обозначается символом $\triangle$;
например вместо «треугольник $ABC$» пишут «$\triangle ABC$».

\begin{figure}[!ht]
\begin{minipage}{.32\textwidth}
\centering
\includegraphics{mppics/ris-34}
\end{minipage}\hfill
\begin{minipage}{.32\textwidth}
\centering
\includegraphics{mppics/ris-35}
\end{minipage}\hfill
\begin{minipage}{.32\textwidth}
\centering
\includegraphics{mppics/ris-36}
\end{minipage}

\medskip

\begin{minipage}{.32\textwidth}
\centering
\caption{}\label{1938/ris-34}
\end{minipage}\hfill
\begin{minipage}{.32\textwidth}
\centering
\caption{}\label{1938/ris-35}
\end{minipage}\hfill
\begin{minipage}{.32\textwidth}
\centering
\caption{}\label{1938/ris-36}
\end{minipage}
\vskip-4mm
\end{figure}

\begin{wrapfigure}[11]{r}{50mm}
\vskip-4mm

\begin{minipage}{24mm}
\centering
\includegraphics{mppics/ris-ru-37}
\end{minipage}
\hfill
\begin{minipage}{24mm}
\centering
\includegraphics{mppics/ris-38}
\end{minipage}

\medskip

\begin{minipage}{24mm}
\centering
\caption{}\label{1938/ris-37}
\end{minipage}
\hfill
\begin{minipage}{24mm}
\centering
\caption{}\label{1938/ris-38}
\end{minipage}
\end{wrapfigure}

\paragraph{Типы треугольников.}\label{1938/35}
Треугольники разделяются по сравнительной длине их сторон или по величине их углов.
Относительно длины сторон они бывают:
\rindex{разносторонний треугольник}\textbf{разносторонние} (рис.~\ref{1938/ris-34}), когда все стороны различной длины, и \rindex{равнобедренный треугольник}\textbf{равнобедренные} (рис. \ref{1938/ris-35}), когда две стороны одинаковы;
в частности, равнобедренный треугольник называется \rindex{равносторонний треугольник}\textbf{равносторонним} (рис.~\ref{1938/ris-36}), когда все три его стороны равны между собой.

Относительно величины углов треугольники бывают:
\rindex{остроугольный треугольник}\textbf{остроугольные} (рис.~\ref{1938/ris-34}), когда все углы острые, \rindex{прямоугольный треугольник}\textbf{прямоугольные} (рис. \ref{1938/ris-37}), когда в числе углов есть прямой, и \rindex{тупоугольный треугольник}\textbf{тупоугольные} (рис.~\ref{1938/ris-38}), когда в числе углов есть тупой.

В прямоугольном треугольнике стороны, образующие прямой угол, называются \rindex{катет}\textbf{катетами}, а сторона, лежащая против прямого угла, — \rindex{гипотенуза}\textbf{гипотенузой}.

\paragraph{Основные линии в треугольнике.}\label{1938/36}
Одну из сторон треугольника иногда называют \rindex{основание!треугольника}\textbf{основанием}, тогда вершину противоположного угла называют \rindex{вершина}\textbf{вершиной} треугольника, а перпендикуляр, опущенный из вершины на основание или на его продолжение, — \rindex{высота!треугольника}\textbf{высотой} его.

Так, если в $\triangle ABC$ (рис.~\ref{1938/ris-39}) за основание взята сторона $AC$, то $B$ будет вершина, $BD$ — высота треугольника.

В равнобедренном треугольнике основанием называют обыкновенно ту сторону, которая не принадлежит к равным;
тогда вершина равнобедренного треугольника будет вершиной того угла его, который образован равными сторонами.
Отрезок $BE$ (рис.~\ref{1938/ris-39}), соединяющий вершину какого-нибудь угла треугольника с серединой противоположной стороны, называется \rindex{медиана}\textbf{медианой}.
Отрезок $BF$ (рис.~\ref{1938/ris-39}а), делящий какой-нибудь угол треугольника пополам, называется его \rindex{биссектриса!треугольника}\textbf{биссектрисой} (биссектриса, вообще говоря, не совпадает ни с медианой, ни с высотой).

\begin{figure}[!ht]
\centering
\includegraphics{mppics/ris-ru-39}
\caption{}\label{1938/ris-39}
\end{figure}

Из вершины каждого угла треугольника можно опустить перпендикуляр на противоположную сторону или её продолжение;
следовательно, каждый треугольник имеет три высоты.
Вершину каждого угла треугольника можно соединить прямой с серединой противоположной стороны, следовательно, каждый треугольник имеет три медианы.
Точно так же ясно, что каждый треугольник имеет три биссектрисы.

\subsection*{Осевая симметрия}

\begin{wrapfigure}{R}{49mm}
\centering
\includegraphics{mppics/ris-40}
\caption{}\label{1938/ris-40}
\end{wrapfigure}

\paragraph{}\label{1938/37}
При изучении свойств треугольников, многоугольников и других геометрических фигур часто встречается случай особого расположения на плоскости двух равных фигур или двух равных отрезков, или двух точек по отношению к какой-либо прямой.
Если какие-нибудь две точки $A$ и $A'$ (рис.~\ref{1938/ris-40}) расположены по разные стороны от прямой $MN$ на одном и том же перпендикуляре к этой прямой и на одинаковом расстоянии от основания перпендикуляра ($Aa=A'a$), то такие точки называются \rindex{симметрия относительно прямой}\textbf{симметричными} относительно прямой $MN$.

Две фигуры (или две части одной и той же фигуры) называются симметричными относительно прямой $MN$, если каждой точке $A$, $B$, $C$, $D$, $E,\dots$
(рис.~\ref{1938/ris-40}) одной фигуры (или одной части фигуры) соответствуют симметричные точки $A'$, $B'$, $C'$, $D'$, $E',\dots$ другой фигуры (или другой части фигуры), и обратно.
Прямая $MN$ в таком случае называется \rindex{ось симметрии}\textbf{осью симметрии}. 
Здесь слово «ось» применяется потому, что если часть плоскости, лежащую по одну сторону от прямой $MN$ (например, левую часть), станем вращать вокруг $MN$, как около оси, до тех пор, пока эта часть плоскости не упадёт на ту часть, которая лежит по другую сторону от прямой $MN$ (на правую часть), то симметричные фигуры совместятся, так как точка $A$ совпадёт при этом с точкой $A'$, точка $B$ — с точкой $B'$ и~т.~д.

Обратно, если вращением вокруг некоторой прямой мы можем фигуру, лежащую по одну сторону от этой прямой, совместить с фигурой, лежащей по другую её сторону, то эти фигуры симметричны относительно оси вращения.
Из сказанного следует, что

\textbf{\emph{всякие две фигуры, симметричные относительно какой-либо оси, равны между собой.}}

Симметрия относительно оси называется \rindex{осевая симметрия}\textbf{осевой симметрией}. 

\begin{wrapfigure}{o}{36mm}
\centering
\includegraphics{mppics/ris-41}
\caption{}\label{1938/ris-41}
\end{wrapfigure}

\begin{figure}[!ht]
\begin{minipage}{.48\textwidth}
\centering
\includegraphics{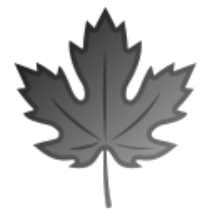}
\end{minipage}\hfill
\begin{minipage}{.48\textwidth}
\centering
\includegraphics{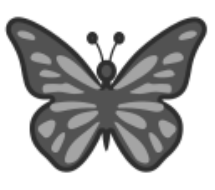}
\end{minipage}

\medskip

\begin{minipage}{.48\textwidth}
\centering
\caption{}\label{1938/ris-42}
\end{minipage}\hfill
\begin{minipage}{.48\textwidth}
\centering
\caption{}\label{1938/ris-43}
\end{minipage}
\vskip-4mm
\end{figure}

{\small
\smallskip
\mbox{\so{Замечание}.}
Хотя симметричные фигуры вращением вокруг оси симметрии могут быть приведены в совмещение, однако они, вообще говоря, не тождественны в своём расположении на плоскости.
Это нужно понимать в следующем смысле:
\emph{чтобы совместить две симметричные фигуры, необходимо одну из них перевернуть другой стороной и, следовательно, на время вывести её из плоскости.}
Если же не выводить фигуры из плоскости, то, вообще говоря, никаким перемещением в этой плоскости нельзя привести её к совпадению с фигурой, ей симметричной относительно оси.

На рис.~\ref{1938/ris-41} изображены два узора, симметричные относительно прямой $AB$.
Вращая правый узор около прямой $AB$, его можно совместить с левым узором.

При этом правый узор будет перевёрнут другой стороной.
Но если не отрывать правого узора от плоскости, а перемещать его так, чтобы он скользил по плоскости, то никаким передвижением его не удастся совместить с левым узором.

Осевая симметрия часто встречается в обыденной жизни.
Узоры на декоративных тканях и на комнатных обоях, архитектурные украшения на зданиях в виде плоских рисунков и самые фасады зданий имеют обычно форму, симметричную относительно некоторой оси.
В природе также часто встречаются симметричные формы.
Так, листья деревьев и лепестки цветов имеют форму, симметричную относительно среднего стебля.
Таков изображённый на рис.~\ref{1938/ris-42} лист клёна.
Крылья бабочки и их расцветка имеют форму, симметричную относительно оси её туловища (рис.~\ref{1938/ris-43}).

}

\subsection*{Свойства равнобедренного треугольника}

\paragraph{}\label{1938/38}
\so{Теоремы}.
1) \textbf{\emph{В равнобедренном треугольнике биссектриса угла при вершине есть одновременно и медиана и высота.}}

2) \textbf{\emph{В равнобедренном треугольнике углы при основании равны.}}

\begin{wrapfigure}{O}{36mm}
\centering
\includegraphics{mppics/ris-44}
\caption{}\label{1938/ris-44}
\end{wrapfigure}

Пусть $\triangle ABC$ (рис.~\ref{1938/ris-44}) равнобедренный и прямая $BD$ делит пополам угол $B$ при вершине его.
Требуется доказать, что эта биссектриса $BD$ есть также и медиана и высота.

Вообразим, что $\triangle ABD$ повёрнут вокруг стороны $BD$, как около оси, так, чтобы он упал на $\triangle BDC$.
Тогда, вследствие равенства углов 1 и 2, сторона $AB$ пойдёт по $BC$, а вследствие равенства этих сторон точка $A$ совпадёт с $C$.
Поэтому $BA$ совместится с $BC$, угол 4 совместится с углом 3 и угол 5 — с углом 6;
значит,
\[DA =DC,\quad \angle 4 = \angle 3\quad \text{и}\quad \angle 5 = \angle 6.\]

Из того, что $DA=DC$, следует, что $BD$ есть медиана;
из того, что углы 3 и 4 равны, вытекает, что эти углы прямые и, следовательно, $BD$ есть высота треугольника, и, наконец, углы 5 и 6 при основании треугольника равны.

\paragraph{}\label{1938/39}
\mbox{\so{Следствие}.}
Мы видим, что в равнобедренном треугольнике $ABC$ (рис.~\ref{1938/ris-44}) одна и та же прямая $BD$ обладает четырьмя свойствами:
она есть биссектриса угла при вершине, медиана, проведённая к основанию, высота, опущенная на основание, и, наконец, перпендикуляр к основанию, восстановленный из его середины. 
Перпендикуляр к отрезку, восстановленный из его середины называется \rindex{срединный перпендикуляр}\textbf{срединным перпендикуляром} к этому отрезку.

Так как каждое из этих четырёх свойств вполне определяет положение прямой $BD$, то существование одного из них влечёт все остальные.
Например, \emph{высота, опущенная на основание равнобедренного треугольника, служит одновременно биссектрисой угла при вершине, медианой, проведённой к основанию, и срединным перпендикуляром к основанию.} 

\paragraph{Симметрия равнобедренного треугольника.}\label{1938/40}
Мы видели, что равнобедренный $\triangle ABC$ (рис.~\ref{1938/ris-44}) делится биссектрисой $BD$ на такие два треугольника (левый и правый), которые вращением вокруг биссектрисы могут быть совмещены один с другим.
Из этого можно заключить, что какую бы точку на одной половине равнобедренного треугольника мы ни взяли, всегда можно на другой его половине найти точку, симметричную с первой относительно оси $BD$.
Возьмём, например, на стороне $AB$ точку $M$ (рис.~\ref{1938/ris-44}).
Опустим из неё на $BD$ перпендикуляр $MK$ и продолжим его до пересечения со стороной $BC$.
Мы получим тогда на этой стороне точку $M'$, симметричную с точкой $M$ относительно оси $BD$.
Действительно, если, вращая $\triangle ABD$ вокруг $BD$, мы его совместим с $\triangle BCD$, то при этом $KM$ пойдёт по $KM'$ (по равенству прямых углов), а сторона $BA$ пойдёт по стороне $BC$ (по равенству углов при вершине $B$);
значит, точка $M$, которая лежит на $KM$ и на $BA$, совпадёт с точкой $M'$, которая лежит и на $KM'$, и на $BC$.
Отсюда видно, что $KM=KM'$.
Таким образом, точки $M$ и $M'$ лежат по разные стороны от биссектрисы $BD$, на одном к ней перпендикуляре и на равных расстояниях от основания этого перпендикуляра;
значит, эти точки симметричны относительно оси $BD$.
Таким образом, \emph{в равнобедренном треугольнике биссектриса угла при вершине есть его ось симметрии.}

\subsection*{Признаки равенства треугольников}

\paragraph{Предварительные понятия.}\label{1938/41}
Две геометрические фигуры, например два треугольника, как мы знаем, называются равными, если они при наложении могут быть совмещены.
В совмещающихся треугольниках, конечно, должны быть соответственно равны все элементы их, то есть стороны, углы, высоты, медианы и биссектрисы.
Однако для того, чтобы утверждать, что два треугольника равны, нет необходимости устанавливать равенства всех их элементов, достаточно убедиться в равенстве только некоторых из них.

\paragraph{Признаки равенства треугольников.}\label{1938/42}\ 

1) \textbf{\emph{Если две стороны и угол, заключённый
между ними, одного треугольника соответственно равны
двум сторонам и углу, заключённому между ними, другого треугольника, то такие треугольники равны.}}

2) \textbf{\emph{Если два угла и прилежащая к ним сторона одного треугольника соответственно равны двум углам и прилежащей к ним стороне другого треугольника, то такие треугольники равны.}}

3) \textbf{\emph{Если три стороны одного треугольника равны трём сторонам другого треугольника, то такие треугольники равны.}}

1) Пусть $ABC$ и $A_1B_1C_1$ — два треугольника (рис.~\ref{1938/ris-45}), у которых
$AC\z=A_1C_1$, $AB \z= A_1B_1$, $\angle A \z= \angle A_1$.
Требуется доказать, что эти треугольники равны.
\begin{figure}[h]
\begin{minipage}{.48\textwidth}
\centering
\includegraphics{mppics/ris-45}
\end{minipage}\hfill
\begin{minipage}{.48\textwidth}
\centering
\includegraphics{mppics/ris-46}
\end{minipage}

\medskip

\begin{minipage}{.48\textwidth}
\centering
\caption{}\label{1938/ris-45}
\end{minipage}\hfill
\begin{minipage}{.48\textwidth}
\centering
\caption{}\label{1938/ris-46}
\end{minipage}
\vskip-4mm
\end{figure}

Наложим%
\footnote{Для выполнения указанных в этом параграфе наложений иногда приходится накладываемый треугольник перевернуть другой стороной.} 
$\triangle ABC$ на $\triangle A_1B_1C_1$ так, чтобы точка $A$ совпала с $A_1$ и сторона $AC$ пошла по $A_1C_1$.
Тогда, вследствие равенства этих сторон, точка $C$ совместится с $C_1$, вследствие равенства углов $A$ и $A_1$ сторона $AB$ пойдёт по $A_1B_1$, а вследствие равенства этих сторон точка $B$ совпадёт с $B_1$, поэтому сторона $CB$ совместится с $C_1B_1$ (так как две точки можно соединить только одной прямой), и треугольники совпадут;
значит, они равны.

2) Пусть $ABC$ и $A_1B_1C_1$ (рис.~\ref{1938/ris-46}) — два треугольника, у которых $\angle C= \angle C_1$,
$\angle B=\angle B_1$,
и
$CB = C_1B_1$.
Требуется доказать, что эти треугольники равны.

Наложим $\triangle ABC$ на $\triangle A_1B_1C_1$ так, чтобы точка $C$ совпала с $C_1$ и сторона $CB$ пошла по $C_1B_1$.
Тогда, вследствие равенства этих сторон, точка $B$ совпадёт с $B_1$, а вследствие равенства углов $B$ и $B_1$, $C$ и $C_1$, сторона $BA$ пойдёт по $B_1A_1$ и сторона $CA$ — по $C_1A_1$.

Так как две прямые могут пересечься только в одной точке, то вершина $A$ должна совпасть с $A_1$.
Таким образом, треугольники совместятся;
значит, они равны.

\begin{wrapfigure}{o}{35mm}
\vskip-0mm
\centering
\includegraphics{mppics/ris-47}
\caption{}\label{1938/ris-47}
\end{wrapfigure}

3) Пусть $ABC$ и $A_1B_1C_1$ (рис.~\ref{1938/ris-47}) — два треугольника, у которых
$AB \z= A_1B_1$,
$BC \z= B_1C_1$,
и
$CA = C_1A_1$.
Требуется доказать, что эти треугольники равны.

Доказывать этот признак равенства наложением, как мы это делали для первых признаков, было бы неудобно, так как, не зная ничего о величине углов, мы не можем утверждать, что при совпадении двух равных сторон совпадут и остальные стороны.
Вместо наложения применим здесь \textbf{приложение}.

Приложим $\triangle ABC$ к $\triangle A_1B_1C_1$ так, чтобы у них совместились равные стороны $AC$ и $A_1C_1$.
Тогда $\triangle ABC$ займёт положение $\triangle A_1C_1B_2$.

Соединив прямой точки $B_1$ и $B_2$, мы получим два равнобедренных треугольника $A_1B_1B_2$ и $B_1C_1B_2$ с общим основанием $B_1B_2$.
Но в равнобедренном треугольнике углы при основании равны (§~\ref{1938/38});
следовательно, $\angle 1 = \angle 2$ и $\angle 3 \z= \angle 4$, а потому $\angle A_1B_1C_1 = \angle A_1B_2C_1 = \angle B$.
Но в таком случае данные треугольники должны быть равны, так как две стороны и угол, заключённый между ними, одного треугольника соответственно равны двум сторонам и углу, заключённому между ними, другого треугольника.

\begin{wrapfigure}{r}{63mm}
\vskip-0mm
\begin{minipage}{31mm}
\centering
\includegraphics{mppics/ris-1914-40}
\end{minipage}\hfill
\begin{minipage}{31mm}
\centering
\includegraphics{mppics/ris-1914-41}
\end{minipage}
\medskip
\begin{minipage}{31mm}
\centering
\caption{}\label{1914/ris-40}
\end{minipage}\hfill
\begin{minipage}{31mm}
\centering
\caption{}\label{1914/ris-41}
\end{minipage}
\end{wrapfigure}

\medskip

Может случиться, что прямая $B_1B_2$ не пересечётся со стороной $A_1C_1$,
а пойдёт вне треугольников (если сумма углов $C$ и $C_1$, большее $180\degree$),
или сольётся с линией $B_1C_1B_2$ (если $\angle C +\angle  C_1 \z= 180\degree$). Доказательство остаётся то же самое, с той разницей, что углы $B_1$ и $B_2$
будут равны друг другу, не как \so{суммы} равных углов, а как их разности, (рис.~\ref{1914/ris-40}), или как углы при основании равнобедренного треугольника (рис.~\ref{1914/ris-41}).

{\small
\smallskip
\mbox{\so{Замечание}.}
\emph{В равных треугольниках против равных сторон лежат равные углы, и обратно, против равных углов лежат равные стороны.}

Доказанные теоремы о равенстве треугольников и умение распознавать равные треугольники по указанным признакам чрезвычайно облегчают решение многих геометрических задач и необходимы для доказательства многих теорем.
Теоремы о равенстве треугольников являются главным средством для обнаружения свойств сложных геометрических фигур.
Учащиеся убедятся в этом при дальнейшем прохождении предмета.
}

\subsection*{Внешний угол треугольника}

\paragraph{}\label{1938/43}
\mbox{\so{Определение}.}
Угол, смежный с каким-нибудь углом треугольника (или многоугольника), называется \rindex{внешний угол}\textbf{внешним} углом этого треугольника (или многоугольника). 

Таковы, например, углы $BCD$, $CBE$, $BAF$ (рис.~\ref{1938/ris-48}).
В отличие от внешних углы самого треугольника (или многоугольника) называются \rindex{внутренний угол}\textbf{внутренними}. 

При каждом угле треугольника (или многоугольника) можно построить по два внешних угла (продолжив одну или другую сторону угла).
Эти два угла равны, как углы вертикальные.

\begin{figure}[!ht]
\begin{minipage}{.48\textwidth}
\centering
\includegraphics{mppics/ris-48}
\end{minipage}\hfill
\begin{minipage}{.48\textwidth}
\centering
\includegraphics{mppics/ris-49}
\end{minipage}

\begin{minipage}{.48\textwidth}
\centering
\caption{}\label{1938/ris-48}
\end{minipage}\hfill
\begin{minipage}{.48\textwidth}
\centering
\caption{}\label{1938/ris-49}
\end{minipage}
\end{figure}

\paragraph{}\label{1938/44}
\mbox{\so{Теорема}.}
\textbf{\emph{Внешний угол треугольника больше каждого внутреннего угла его, не смежного с этим внешним.}}

Например, докажем, что внешний угол $BCD$ треугольника $ABC$
(рис.~\ref{1938/ris-49}) больше каждого из внутренних углов $A$ и $B$, не смежных с этим внешним.

Через середину $E$ стороны $BC$ проведём медиану $AE$ и на её продолжении отложим отрезок $EF=AE$.
Точка $F$, очевидно, будет лежать внутри угла $BCD$.
Соединим $F$ с $C$ прямой.
Треугольники $ABE$ и $EFC$ (покрытые штрихами) равны, так как при точке $E$ они имеют по равному углу, заключённому между двумя соответственно равными сторонами.
Из равенства их заключаем, что углы $B$ и $ECF$, лежащие против равных сторон $AE$ и $EF$, равны.
Но угол $ECF$ составляет часть внешнего угла $BCD$ и потому меньше его;
следовательно, и угол $B$ меньше угла $BCD$.

Продолжив сторону $BC$ за точку $C$, мы получим внешний угол $ACH$, равный углу $BCD$.
Если из вершины $B$ проведём к стороне $AC$ медиану и продолжим её на такую же длину за сторону $AC$, то совершенно так же докажем, что угол $A$ меньше угла $ACH$, то есть меньше угла $BCD$.

{

\begin{wrapfigure}[8]{r}{57mm}
\vskip-0mm
\begin{minipage}{28mm}
\centering
\includegraphics{mppics/ris-50}
\end{minipage}\hfill
\begin{minipage}{28mm}
\centering
\includegraphics{mppics/ris-51}
\end{minipage}
\medskip
\begin{minipage}{28mm}
\centering
\caption{}\label{1914/ris-50}
\end{minipage}\hfill
\begin{minipage}{28mm}
\centering
\caption{}\label{1914/ris-51}
\end{minipage}

\end{wrapfigure}

\paragraph{}\label{1938/45}
\mbox{\so{Следствие}.}
\emph{Если в треугольнике один угол прямой или тупой, то два других угла острые.}
Действительно, допустим, что угол $C$ треугольника $ABC$ 
(рис.~\ref{1914/ris-50} и \ref{1914/ris-51}) будет прямой или тупой, тогда смежный с ним внешний угол $BCD$ должен быть прямой или острый;
следовательно, углы $A$ и $B$, которые, по доказанному, меньше этого внешнего угла, должны быть оба острые.

}

\subsection*{Стороны и углы треугольника}

\paragraph{}\label{1938/46}
\so{Теоремы}.
\textbf{\emph{Во всяком треугольнике}}

\textbf{\emph{1) против равных сторон лежат равные углы.}}

\textbf{\emph{2) против б\'{о}льшей стороны лежит больший угол.}}

\begin{wrapfigure}{o}{33mm}
\vskip-4mm
\centering
\includegraphics{mppics/ris-52}
\caption{}\label{1938/ris-52}
\bigskip
\includegraphics{mppics/ris-53}
\caption{}\label{1938/ris-53}
\end{wrapfigure}

1) Если две стороны треугольника равны, то он равнобедренный, тогда углы, лежащие против этих сторон, должны быть равны, как углы при основании равнобедренного треугольника (§~\ref{1938/38}).

2) Пусть в $\triangle ABC$ (рис.~\ref{1938/ris-52}) сторона $AB$ больше $BC$;
требуется доказать, что угол $C$ больше угла $A$.

Отложим на б\'{о}льшей стороне $BA$ от вершины $B$ отрезок $BD$, равный меньшей стороне $BC$, и соединим $D$ с $C$ прямой.
Тогда получим равнобедренный $\triangle BDC$, у которого углы при основании равны, то есть $\angle BDC=\angle BCD$.
Но угол $BDC$ как внешний по отношению к $\angle ADC$, больше угла $A$, следовательно, и угол $BCD$ больше угла $A$, а потому и подавно угол $BCA$ больше угла $A$, что и требовалось доказать.

\paragraph{}\label{1938/47}
\mbox{\so{Обратные теоремы}.}
\textbf{\emph{Во всяком треугольнике.}}

1) \textbf{\emph{против равных углов лежат равные стороны;}}

2) \textbf{\emph{против б\'{о}льшего угла лежит б\'{о}льшая сторона.}}

1) Пусть в $\triangle ABC$ углы $A$ и $C$ равны (рис.~\ref{1938/ris-53});
требуется доказать, что $BA \z= BC$.

Предположим противное, то есть что стороны $AB$ и $BC$ не равны.
Тогда одна из этих сторон должна быть больше другой, и, следовательно, согласно прямой теореме, один из углов $A$ и $C$ должен быть больше другого.
Но это противоречит условию, что $\angle A = \angle C$;
значит, нельзя допустить, что стороны $AB$ и $BC$ не равны;
остаётся принять, что $AB=BC$.

2) Пусть в $\triangle ABC$ (рис.~\ref{1938/ris-54})
угол $C$ больше угла $A$;
требуется доказать, что $AB > BC$.

{

\begin{wrapfigure}[10]{r}{30mm}
\centering
\includegraphics{mppics/ris-54}
\caption{}\label{1938/ris-54}
\end{wrapfigure}

Предположим противное, то есть что $AB$ не больше $BC$.
Тогда могут представиться два случая:
или $AB=BC$, или $AB<BC$.

В первом случае, согласно прямой теореме, угол $C$ был бы равен углу $A$, во втором случае угол $C$ был бы меньше угла $A$;
и то и другое противоречит условию;
значит, оба эти случая исключаются.
Остаётся один возможный случай, что $AB>BC$.

\smallskip
\mbox{\so{Следствия}.}
1) \emph{В равностороннем треугольнике все углы равны.}

2) \emph{В равноугольном треугольнике все стороны равны.}

}

{\sloppy
\paragraph{Доказательство от противного.}\label{1938/48}
Способ, которым мы только что доказали обратные теоремы, называется \rindex{доказательство от противного}\textbf{доказательством от противного}, или приведением к \rindex{противоречие}\textbf{противоречию}.
Первое название этот способ получил потому, что в начале рассуждения делается предположение, противное (противоположное) тому, что требуется доказать.
Приведением к противоречию он называется вследствие того, что, рассуждая на основании сделанного предположения, мы приходим к противоречию (к абсурду).
Получение такого вывода заставляет нас отвергнуть сделанное вначале допущение и принять то, которое требовалось доказать.

}

Этот приём очень часто употребляется для доказательства теорем.

\paragraph{Замечание об обратных теоремах.}\label{1938/49}
Начинающие изучать геометрию часто делают одну характерную ошибку.
Она заключается в том, что правильность обратной теоремы считают само собой разумеющейся, если доказана прямая теорема.
Отсюда представление, что доказательство обратных теорем вообще излишне.
Ошибочность такого заключения может быть показана в ряде примеров.
В частности, такой пример был приведён в §~\ref{1938/30}.
Поэтому обратные теоремы, когда они верны, всегда доказываются особо.

\subsection*{Длина ломаной}

\paragraph{}\label{1938/50}
\so{Теорема}.
\textbf{\emph{В треугольнике каждая сторона меньше суммы двух других сторон.}}

Если в треугольнике возьмём сторону не самую большую, то, конечно, она окажется менее суммы двух других сторон.
Значит, нам надо доказать, что даже \so{наибольшая} сторона треугольника меньше суммы двух других сторон.

Пусть в $\triangle ABC$ (рис.~\ref{1938/ris-55}) наибольшая сторона есть $AC$.
Продолжив сторону $AB$, отложим $BD=BC$ и проведём $DC$.
Так как $\triangle BDC$ равнобедренный, то $\angle D = \angle DCB$;
поэтому угол $D$ меньше угла $DCA$, и, следовательно, в $\triangle ADC$ сторона $AC$ меньше $AD$ (§~\ref{1938/47}), то есть
$AC \z< AB + BD$.
Заменив $BD$ на $BC$, получим:
\[AC < AB + BC.\]

\begin{wrapfigure}{o}{33mm}
\vskip-8mm
\centering
\includegraphics{mppics/ris-55}
\caption{}\label{1938/ris-55}
\end{wrapfigure}

\smallskip
\mbox{\so{Следствие}.}
Отнимем от обеих частей выведенного неравенства по $AB$ или по $BC$:
\begin{align*}
AC-AB&<BC;
\\
AC-BC&<AB.
\end{align*}
Читая эти неравенства справа налево, видим, что каждая из сторон $BC$ и $AB$ больше разности двух других сторон;
так как это же можно, очевидно, сказать и о третьей, наибольшей стороне $AC$, то, значит, \emph{в треугольнике каждая сторона больше разности двух других сторон.}

\paragraph{Неравенство треугольника.}\label{extra/3inq}
\emph{Для любых трёх точек $A$, $B$ и $C$ выполняется неравенство
\[AC \le AB + BC;\]
то есть отрезок $AC$ может быть равен сумме $AB + BC$ или меньше неё.
При этом равенство достигается только в случае если $B$ лежит на отрезке $AC$.}
Это неравенство называется \rindex{неравенство треугольника}\textbf{неравенством треугольника}.

Случай когда $B$ не лежит на прямой $AC$ (то есть если $A$, $B$ и $C$ являются вершинами некоторого треугольника) уже доказан (§~\ref{1938/50}). 
Остаётся рассмотреть случай когда $B$ лежит на прямой $AC$.

Равенство 
\[AC = AB + BC;\]
очевидно выполняется в случае если $B$ лежит на отрезке $AC$.
В случае если $B$ лежит на продолжении отрезка $AC$, то очевидно
$AC=AB\z-BC<AB$ или $AC=CB-BA<CB$;
в обоих случаях имеем
\[AC < AB + BC.\] 

\paragraph{}\label{1938/51}
\mbox{\so{Теорема}.}
\textbf{\emph{Отрезок, соединяющий две какие-нибудь точки, меньше всякой ломаной, соединяющей эти же точки.}}

Если ломаная состоит только из двух сторон, то теорема уже была доказана в предыдущем параграфе.
Рассмотрим случай, когда ломаная состоит более чем из двух сторон.
Пусть $AE$ (рис.~\ref{1938/ris-56}) есть отрезок, соединяющий точки $A$ и $E$, а $ABCDE$ — какая-нибудь ломаная, соединяющая те же точки.
Требуется доказать, что 
\[AE<AB\z+BC\z+CD\z+DE.\eqno(1)\]

{

\begin{wrapfigure}{r}{33mm}
\vskip0mm
\centering
\includegraphics{mppics/ris-56}
\caption{}\label{1938/ris-56}
\end{wrapfigure}

Соединив $A$ с $C$ и $D$, находим, согласно неравенству треугольника (§~\ref{extra/3inq})
\begin{align*}
AE&\le AD+DE;
\\
AD&\le AC +CD;
\\
AC&< AB+BC.
\end{align*}
Последнее неравенство строгое, поскольку $B$ не может лежать на отрезке $AC$.

}

Сложим почленно эти неравенства и затем от обеих частей полученного неравенства отнимем по $AD$ и $AC$; тогда получим неравенство~(1).

\paragraph{Треугольники с двумя соответственно равными сторонами.}\label{1938/52}\ 

\smallskip
\so{Теоремы}.
\textbf{\emph{Если две стороны одного треугольника соответственно равны двум сторонам другого треугольника, то:}}

1) \textbf{\emph{против большего из углов, заключённых между ними, лежит б\'{о}льшая сторона.}}

2) \so{обратно}:
\textbf{\emph{против б\'{о}льшей из неравных сторон лежит больший угол.}}

\begin{figure}[!ht]
\centering
\includegraphics{mppics/ris-57}
\caption{}\label{1938/ris-57}
\end{figure}

1) Пусть в треугольниках $ABC$ и $A_1B_1C_1$ дано (рис.~\ref{1938/ris-57}):
$AC\z=A_1C_1$, $AB=A_1B_1$ и $\angle A > \angle A_1$.
Требуется доказать, что $BC\z>B_1C_1$.
Наложим $\triangle A_1B_1C_1$ на $\triangle ABC$ так, чтобы сторона $A_1C_1$ совпадала с $AC$.
Так как $\angle A_1 < \angle BAC$, то сторона $A_1B_1$ пойдёт внутри угла $BAC$;
пусть $\triangle A_1B_1C_1$ займёт положение $AB_2C$ (вершина $B_2$ может оказаться или вне $\triangle ABC$, или внутри него, или же на стороне $BC$;
доказательство может быть применено ко всем этим случаям).
Проведём биссектрису $AD$ угла $BAB_2$ и соединим $D$ с $B_2$;
тогда получим два треугольника:
$ABD$ и $DAB_2$, которые равны, потому что у них $AB$ — общая сторона, $AB=AB_2$ по условию и $\angle BAD=\angle DAB_2$ по построению.
Из равенства треугольников следует:
$BD=DB_2$.
Из $\triangle DCB_2$ выводим:
$B_2C < B_2D \z+ DC$ (§~\ref{1938/50}), или (заменив $B_2D$ на $BD$):
\[B_2C <BD +DC,\quad\text{значит,}\quad B_1C_1 < BC.\]

2) Пусть в тех же треугольниках $ABC$ и $A_1B_1C_1$ дано:
$AB\z=A_1B_1$;
$AC=A_1C_1$ и $BC>B_1C_1$;
докажем, что $\angle A > \angle A_1$.

Допустим противное, то есть что угол $A$ не больше угла $A_1$, тогда могут представиться два случая:
или $\angle BAC = \angle A_1$, или $\angle BAC \z< \angle A_1$.
В первом случае треугольники были бы равны и, следовательно, сторона $BC$ равнялась бы $B_1C_1$, что противоречит условию;
во втором случае сторона $BC$ (согласно теореме 1) была бы меньше $B_1C_1$, что также противоречит условию.
Значит, оба эти случая исключаются;
остаётся один возможный случай, что $\angle A > \angle A_1$.

\subsection*{Перпендикуляр и наклонная}

\paragraph{}\label{1938/53}
\mbox{\so{Теорема}.}
\textbf{\emph{Перпендикуляр, опущенный из какой-нибудь точки на прямую, меньше всякой наклонной, проведённой из той же точки на эту прямую%
\footnote{В §§~\ref{1938/53}, \ref{1938/54} и \ref{1938/55} ради краткости термины «перпендикуляр» и «наклонная» употребляются вместо «отрезок перпендикуляра, ограниченный данной точкой и основанием перпендикуляра» и «отрезок наклонной, ограниченный данной точкой и основанием наклонной».}.%
}}

\begin{wrapfigure}{r}{32mm}
\vskip -7mm
\centering
\includegraphics{mppics/ris-58}
\caption{}\label{1938/ris-58}
\end{wrapfigure}

Пусть $AB$ (рис.~\ref{1938/ris-58}) есть перпендикуляр, опущенный из точки $A$ на прямую $MN$, и $AC$ — какая-нибудь наклонная, проведённая из той же точки $A$ к прямой $MN$;
требуется доказать, что $AB<AC$.

В $\triangle ABC$ угол $B$ прямой, а угол $C$ острый (§~\ref{1938/45});
значит, $\angle C<\angle B$ и потому $AB<AC$, что и требовалось доказать.

{\small
\smallskip
\mbox{\so{Замечание}.}
Когда говорят:
«расстояние от точки до прямой», имеется в виду \so{кратчайшее} расстояние, измеряемое по перпендикуляру, опущенному из этой точки на прямую.
}

\paragraph{}\label{1938/54}
\mbox{\so{Теорема}.}
\textbf{\emph{Если из одной и той же точки, взятой вне прямой, проведены к этой прямой перпендикуляр и какие-нибудь наклонные, то:}}

1) \textbf{\emph{если основания двух наклонных одинаково удалены от основания перпендикуляра, то такие наклонные равны.}}

2) \textbf{\emph{если основания двух наклонных неодинаково удалены от основания перпендикуляра, то та из наклонных больше, основание которой дальше отстоит от основания перпендикуляра.}}

1) Пусть $AC$ и $AD$ (рис.~\ref{1938/ris-59}) будут две наклонные, проведённые из точки $A$ к прямой $MN$, основания которых $C$ и $D$ одинаково удалены от основания перпендикуляра $AB$, то есть $CB\z=BD$;
требуется доказать, что $AC = AD$.

\begin{wrapfigure}{o}{42mm}
\vskip2mm
\centering
\includegraphics{mppics/ris-59}
\caption{}\label{1938/ris-59}
\end{wrapfigure}

В треугольниках $ABC$ и $ABD$ есть общая сторона $AB$ и сверх того $BC\z=BD$ (по условию) и $\angle ABC \z= \angle ABD$ (как углы прямые);
значит, эти треугольники равны, и потому $AC = AD$.

2) Пусть $AC$ и $AE$ (рис.~\ref{1938/ris-59}) будут две такие наклонные, проведённые из точки $A$ к прямой $MN$, основания которых неодинаково удалены от основания перпендикуляра;
например, пусть $BE>BC$.
Требуется доказать, что $AE>AC$.

Отложим $BD=BC$ и проведём $AD$.
По доказанному выше $AD \z= AC$.
Сравним $AE$ с $AD$.
Угол $ADE$ есть внешний по отношению $\triangle ABD$, и потому он больше прямого угла $ABD$;
следовательно, угол $ADE$ тупой, и потому угол $AED$ должен быть острый (§~\ref{1938/45}), значит, $\angle ADE>\angle AED$ и, следовательно, $AE>AD$, и потому $AE>AC$.

\paragraph{}\label{1938/55}
\so{Обратные теоремы}.
\textbf{\emph{Если из одной и той же точки, взятой вне прямой}} (рис.~\ref{1938/ris-59}), \textbf{\emph{проведены к этой прямой перпендикуляр и какие-нибудь наклонные, то:}}

1) \textbf{\emph{если две наклонные равны, то их основания одинаково удалены от основания перпендикуляра.}}

2) \textbf{\emph{если две наклонные не равны, то основание б\'{о}льшей из них дальше отстоит от основания перпендикуляра.}}

Предоставляем учащимся самим доказать эти теоремы (способом от противного).

\subsection*{Признаки равенства прямоугольных треугольников}

\paragraph{Признаки, не требующие особого доказательства.}\label{1938/56}
Так как в прямоугольных треугольниках углы, содержащиеся между катетами, всегда равны, как углы прямые, то \textbf{\emph{прямоугольные треугольники равны.}}

1) \textbf{\emph{если катеты одного треугольника соответственно равны катетам другого;}}

2) \textbf{\emph{если катет и прилежащий к нему острый угол одного треугольника соответственно равны катету и прилежащему к нему острому углу другого треугольника.}}

Эти два признака не требуют особого доказательства, так как они представляют частные случаи общих признаков.
Докажем ещё два следующих признака, относящихся только к прямоугольным треугольникам.

\paragraph{Признаки, требующие особого доказательства.}\label{1938/57}\ 

\so{Теоремы}.
\textbf{\emph{Прямоугольные треугольники равны:}}

1) \textbf{\emph{если гипотенуза и острый угол одного треугольника соответственно равны гипотенузе и острому углу другого}} или

2) \textbf{\emph{если гипотенуза и катет одного треугольника соответственно равны гипотенузе и катету другого.}}

1) Пусть $ABC$ и $\triangle A_1B_1C_1$ (рис.~\ref{1938/ris-60}) — два прямоугольных треугольника, у которых $AB=A_1B_1$ и $\angle A = \angle A_1$;
требуется доказать, что эти треугольники равны.

Наложим $\triangle ABC$ на $\triangle A_1B_1C_1$ так, чтобы у них совместились равные гипотенузы.
Тогда по равенству углов $A$ и $A_1$ катет $AC$ пойдёт по $A_1C_1$.
При этом точка $C$ должна совпадать с точкой $C_1$, потому что если предположим, что она не совпадёт с точкой $C_1$, то тогда катет $BC$ занял бы положение $B_1C_2$ или $B_1C_3$, что невозможно, так как из одной точки $B_1$ нельзя на прямую $A_1C_1$ опустить два перпендикуляра ($B_1C_1$ и $B_1C_2$ или $B_1C_1$ и $B_1C_3$).

\begin{figure}[!ht]
\begin{minipage}{.48\textwidth}
\centering
\includegraphics{mppics/ris-60}
\end{minipage}\hfill
\begin{minipage}{.24\textwidth}
\centering
\includegraphics{mppics/ris-61}
\end{minipage}\hfill
\begin{minipage}{.24\textwidth}
\centering
\includegraphics{mppics/ris-62}
\end{minipage}

\medskip

\begin{minipage}{.48\textwidth}
\centering
\caption{}\label{1938/ris-60}
\end{minipage}\hfill
\begin{minipage}{.24\textwidth}
\centering
\caption{}\label{1938/ris-61}
\end{minipage}\hfill
\begin{minipage}{.24\textwidth}
\centering
\caption{}\label{1938/ris-62}
\end{minipage}
\vskip-4mm
\end{figure}

2) Пусть (рис.~\ref{1938/ris-61} и \ref{1938/ris-62}) в прямоугольных треугольниках дано:
$AB=A_1B_1$ и $BC=B_1C_1$;
требуется доказать, что треугольники равны.
Наложим $\triangle ABC$ на $\triangle A_1B_1C_1$ так, чтобы у них совместились равные катеты $BC$ и $B_1C_1$.
Тогда по равенству прямых углов $CA$ пойдёт по $C_1A_1$.
При этом гипотенуза $AB$ не может не совместиться с гипотенузой $A_1B_1$, потому что, если бы она заняла положение $A_2B_1$ или $A_3B_1$, то тогда мы имели бы две равные наклонные ($A_1B_1$ и $A_2B_1$ или $A_1B_1$ и $A_3B_1$), которые неодинаково удалены от основания перпендикуляра, что невозможно (§~\ref{1938/54}).

\subsection*{Срединный перпендикуляр и биссектриса} 

\paragraph{}\label{1938/58}
Свойство срединного перпендикуляра к отрезку, и свойство биссектрисы угла очень сходны между собой.
Чтобы лучше видеть сходство этих свойств, мы изложим их параллельно.

\columnratio{0.5}
\setlength{\columnseprule}{.2pt}
\begin{paracol}{2}

{\sloppy

\textbf{\emph{1) Если какая-нибудь точка}} ($K$, рис.~\ref{1938/ris-63}) \textbf{\emph{лежит на срединном перпендикуляре}} ($MN$), \textbf{\emph{к отрезку}} ($AB$), \textbf{\emph{то она одинаково удалена от концов этого отрезка}} (то есть $KA=KB$).

Так как \mbox{$MN\perp AB$} и $AO\z=OB$, то $AK$ и $KB$ — наклонные к $AB$, основания которых одинаково удалены от основания перпендикуляра, значит $KA=KB$.

\medskip

{\centering
\includegraphics{mppics/ris-63}
\captionof{figure}{}
\label{1938/ris-63}
\addtocounter{figure}{1}
}

\medskip

2) \so{Обратная теорема}.

\textbf{\emph{Если какая-нибудь точка}} ($K$, рис.~\ref{1938/ris-63}) \textbf{\emph{одинаково удалена от концов отрезка $\bm{AB}$}} (то есть если $KA=KB$), \textbf{\emph{то она лежит на перпендикуляре, проведённом к отрезку $\bm{AB}$ через его середину.}}

Проведём через $K$ прямую $MN\perp AB$;
тогда мы получим два прямоугольных треугольника $KAO$ и $KBO$, которые, имея общий катет $KO$ и равные гипотенузы, равны, а потому $AO=OB$.
Значит, прямая $MN$, проведённая нами через $K$ перпендикулярно к $AB$, делит отрезок $AB$ пополам.

}

\switchcolumn

{\sloppy

\textbf{\emph{1) Если какая-нибудь точка}} ($K$, рис.~\ref{1938/ris-64}) \textbf{\emph{лежит на биссектрисе}} ($OM$ угла $AOB$), \textbf{\emph{то она одинаково удалена от сторон этого угла}} (то есть перпендикуляры $KD$ и $KC$ равны).

Так как $OM$ делит угол пополам, то прямоугольные треугольники $OCK$ и $ODK$, имея общую гипотенузу и равные острые углы при вершине $O$, равны и потому $KC=KD$.

\medskip

{\centering
\addtocounter{figure}{1}
\includegraphics{mppics/ris-64}
\captionof{figure}{}
\label{1938/ris-64}
\addtocounter{figure}{1}
}

\medskip

2) \so{Обратная теорема}.

\textbf{\emph{Если какая-нибудь точка}} ($K$, рис.~\ref{1938/ris-64}) \textbf{\emph{внутри угла одинаково удалена от его сторон}} (то есть если перпендикуляры $KC$ и $KD$ равны), \textbf{\emph{то она лежит на биссектрисе этого угла.}}

Через $O$ и $K$ проведём прямую $OM$.
Тогда получим два прямоугольных треугольника $OCK$ и $ODK$, которые, имея общую гипотенузу и равные катеты $CK$ и $DK$, равны, а потому равны и углы при вершине $O$.
Значит, прямая $OM$, проведённая через точку $K$, будет биссектрисой угла $AOB$.

}
\end{paracol}

\paragraph{}\label{1938/59}
\so{Следствие}.
Из двух доказанных теорем (прямой и обратной) можно вывести следующие \so{противоположные теоремы:}

\begin{paracol}{2}

{\sloppy

\textbf{\emph{Если какая-нибудь точка не лежит на срединном перпендикуляре к отрезку, то она неодинаково удалена от концов этого отрезка.}}

\medskip

{\centering
\includegraphics{mppics/ris-1931-62}
\captionof{figure}{}
\label{1931/ris-62}
\addtocounter{figure}{1}
}

}

\switchcolumn

{\sloppy

\textbf{\emph{Если какая-нибудь точка внутри угла не лежит на его биссектрисе, то она неодинаково удалена от сторон этого угла.}}

\medskip

{\centering
\includegraphics{mppics/ris-1931-63}
\captionof{figure}{}
\label{1931/ris-63}
\addtocounter{figure}{1}
}

}

\end{paracol}
\setlength{\columnseprule}{0pt}

\medskip

Предоставляем самим учащимся доказать эти теоремы (способом от противного).

\paragraph{Геометрическое место точек.}\label{1938/60}\rindex{геометрическое место точек}
Геометрическим местом точек, обладающих некоторым свойством, называется такая линия (или поверхность в пространстве) или вообще такая совокупность точек, которая содержит в себе все точки, обладающие этим свойством, и не содержит ни одной точки, не обладающей им.

Например, геометрическое место точек, находящихся на данном расстоянии $r$ от данной точки $C$, есть окружность с центром в точке $C$ и радиусом $r$.
Из теорем предыдущих параграфов следует:

\emph{Геометрическое место точек, одинаково удалённых от двух данных точек, есть срединный перпендикуляр к отрезку соединяющему эти точки.}

\emph{Геометрическое место точек внутри угла, одинаково удалённых от его сторон, есть биссектриса этого угла.}

\begin{wrapfigure}{r}{40mm}
\vskip3mm
\centering
\includegraphics{mppics/ris-extra-1}
\caption{}\label{extra/1}
\end{wrapfigure}

\smallskip
\mbox{\so{Упражнение}.} Докажите, что геометрическое место точек, одинаково удалённых от двух данных пресекающихся прямых, состоит из пары прямых, делящих пополам углы, образованные данными прямыми.

Такие пары прямых всегда перпендикулярны поскольку угол между такими прямыми составлен из половин смежных углов.

%% file: 2D/zadachi-na-postr.tex
\section{Задачи на построение}

\paragraph{Предварительное замечание.}\label{1938/61}
Теоремы, доказанные нами ранее, позволяют решать некоторые задачи на \so{построение}.
Заметим, что в элементарной геометрии рассматриваются только такие построения, которые могут быть выполнены с помощью \so{линейки и циркуля}.
Употребление чертёжного треугольника и некоторых других приборов хотя и допускается ради сокращения времени, но не является необходимым.

\paragraph{}\label{1938/62}
\so{Задача 1}.
\emph{Построить треугольник по трём его сторонам $a$, $b$ и $c$} (рис.~\ref{1938/ris-65}).

\begin{figure}[!ht]
\centering
\includegraphics{mppics/ris-65}
\caption{}\label{1938/ris-65}
\end{figure}

На какой-нибудь прямой $MN$ откладываем отрезок $CB$, равный одной из данных сторон, например $a$.
Описываем две небольшие дуги с центрами в точках $C$ и $B$, одну радиусом, равным $b$, другую радиусом, равным $c$.
Точку $A$, в которой эти дуги пересекаются, соединяем с $B$ и С;
$\triangle ABC$ будет искомый.

{\small

\smallskip
\so{Замечание}.
Чтобы три отрезка могли служить сторонами треугольника, необходимо и достаточно, чтобы больший из них был меньше суммы двух остальных (необходимость доказана в §~\ref{1938/50}, а условие равносильное достаточности будет принято за очевидное в §~\ref{1914/118}).

}

\paragraph{}\label{1938/63}
\so{Задача 2}.
\emph{Построить угол, равный данному углу $ABC$, одной из сторон которого является данная прямая и вершина которого находится в данной точке $O$} (точка $O$ расположена на прямой $MN$, рис.~\ref{1938/ris-66}).

\begin{figure}[!ht]
\centering
\includegraphics{mppics/ris-66}
\caption{}\label{1938/ris-66}
\end{figure}

Описываем произвольным радиусом с центром в вершине $B$ между сторонами данного угла дугу $EF$;
затем, не изменяя раствора циркуля, переносим его остриё в точку $O$ и описываем дугу $PQ$.
Далее описываем дугу $ab$ с центром в точке $P$ радиусом, равным расстоянию между точками $E$ и $F$.
Наконец, через точки $O$ и $R$ (пересечение двух дуг) проводим прямую.
Угол $ROP$ равен углу $ABC$, потому что треугольники $ROP$ и $FBE$, имеющие соответственно равные стороны, равны.

\begin{wrapfigure}{o}{40mm}
\vskip-4mm
\centering
\includegraphics{mppics/ris-67}
\caption{}\label{1938/ris-67}
\bigskip
\includegraphics{mppics/ris-68}
\caption{}\label{1938/ris-68}
\end{wrapfigure}

\paragraph{}\label{1938/64}
\mbox{\so{Задача 3}.}
\emph{Разделить данный угол $ABC$ пополам (рис.~\ref{1938/ris-67}), другими словами, построить биссектрису данного угла или провести его ось симметрии.}
С центром в вершине $B$ произвольным радиусом опишем между сторонами угла дугу $DE$.
Затем, взяв произвольный раствор циркуля, больший, однако, половины расстояния между точками $E$ и $D$ (смотри замечание к задаче 1), описываем этим раствором небольшие дуги с центрами в точках $D$ и $E$, которые пересекутся в некоторой точке $F$.
Проведя прямую $BF$, мы получим биссектрису угла $ABC$.

Для доказательства соединим прямыми точку $F$ с $D$ и $E$;
тогда получим два треугольника $BEF$ и $BDF$, которые равны, так как у них $BF$ — общая сторона\footnote{На чертежах общая сторона треугольников обычно обозначается значком похожим на $S$.}, $BD=BE$ и $DF=EF$ по построению.
Из равенства треугольников следует:
$\angle ABF = \angle CBF$.

\paragraph{}\label{1938/65}
\so{Задача 4}.
Из данной точки $C$ прямой $AB$ восстановить к этой прямой перпендикуляр (рис.~\ref{1938/ris-68}).

Отложим на $AB$ по обе стороны от данной точки $C$ равные отрезки (произвольной длины) $CD$ и $CE$.
С центрами в точках $E$ и $D$ одним и тем же раствором циркуля (б\'{о}льшим, однако, $CD$) опишем две небольшие дуги, которые пересекутся в некоторой точке $F$.
Прямая, проведённая через точки $C$ и $F$, будет искомым перпендикуляром.

Действительно, как видно из построения, точка $F$ одинаково удалена от точек $D$ и $E$;
следовательно, она должна лежать на срединном перпендикуляре к отрезку $DE$ (§~\ref{1938/58});
но середина этого отрезка есть $C$, а через точки $C$ и $F$ можно провести только одну прямую;
значит, $FC \perp DE$.

\paragraph{}\label{1938/66}
\so{Задача 5}.
\emph{Из данной точки $A$ опустить перпендикуляр на данную прямую $BC$} (рис.~\ref{1938/ris-69}).

\begin{wrapfigure}{o}{41mm}
\centering
\includegraphics{mppics/ris-69}
\caption{}\label{1938/ris-69}
\bigskip
\includegraphics{mppics/ris-70}
\caption{}\label{1938/ris-70}
\end{wrapfigure}

С центром в точке $A$ произвольным раствором циркуля (б\'{о}льшим, однако, расстояния от $A$ до $BC$) опишем дугу, которая пересечётся с $BC$ в каких-нибудь точках $D$ и $E$.
С центрами в этих точках произвольным, но одним и тем же раствором циркуля (б\'{о}льшим, однако, $\tfrac12\cdot DE$), проводим две небольшие дуги, которые пересекутся между собой в некоторой точке $F$.
Прямая $AF$ будет искомым перпендикуляром.

Действительно, как видно из построения, каждая из точек $A$ и $F$ одинаково удалена от $D$ и $E$, а такие точки лежат на срединном перпендикуляре к отрезку $DE$ (§~\ref{1938/58}).

\paragraph{}\label{1938/67}
\mbox{\so{Задача 6}.}
\emph{Провести срединный перпендикуляр к данному отрезку \emph{($AB$)}} (рис.~\ref{1938/ris-70}).

Произвольным, но одинаковым раствором циркуля (б\'{о}льшим половины $AB$) описываем две дуги с центрами в точках $A$ и $B$, которые пересекутся между собой в некоторых точках $C$ и $D$.
Прямая $CD$ будет искомым перпендикуляром.
Действительно, как видно из построения, каждая из точек $C$ и $D$ одинаково удалена от $A$ и $B$;
следовательно, эти точки должны лежать на оси симметрии отрезка $AB$.

\smallskip
\so{Задача 7}.
\emph{Разделить пополам данный отрезок} (рис.~\ref{1938/ris-70}).
Решается так же, как предыдущая задача.

\paragraph{Пример более сложной задачи.}\label{1938/68}
При помощи этих основных задач можно решать задачи более сложные.
Для примера решим следующую задачу.

\smallskip
\so{Задача}.
\emph{Построить треугольник, зная его основание $b$, угол $\alpha$, прилежащий к основанию, и сумму $s$ двух боковых сторон} (рис.~\ref{1938/ris-71}).

\begin{figure}[!ht]
\centering
\includegraphics{mppics/ris-71}
\caption{}\label{1938/ris-71}
\end{figure}

Чтобы составить план решения, предположим, что задача решена, то есть что найден такой $\triangle ABC$, у которого основание $AC = b$, $\angle A=\alpha$ и $AB+BC=s$.
Рассмотрим полученный чертёж.
Сторону $AC$, равную $b$, и угол $A$, равный $\alpha$, мы построить умеем.
Значит, остаётся найти на другой стороне угла $A$ такую точку $B$, чтобы сумма $AB+BC$ равнялась $s$.
Продолжив $AB$, отложим отрезок $AD$, равный $s$.

Вопрос сводится к тому, чтобы на прямой $AD$ отыскать такую точку $B$, которая была бы одинаково удалена от $C$ и $D$.
Такая точка, как мы знаем (§~\ref{1938/58}), должна лежать на срединном перпендикуляре к отрезку $CD$. 
Точка $B$ найдётся в пересечении этого срединного перпендикуляра с $AD$. 

Итак, вот решение задачи:
строим (рис.~\ref{1938/ris-71}) угол $A$, равный данному углу $\alpha$;
на сторонах его откладываем $AC=b$ и $AD=s$ и соединяем точку $D$ с $C$.
Проведём срединный перпендикуляр к отрезку $CD$;
пересечение его с $AD$, то есть точку $B$, соединяем с $C$.
Треугольник $ABC$ будет  искомый, так как он удовлетворяет всем требованиям задачи:
у него $AC=b$, $\angle A = \alpha$ и $AB+BC=s$ (потому что $BD=BC$).

Рассматривая построение, мы замечаем, что решения задачи может не быть.
Действительно, если сумма задана слишком малой сравнительно с $b$, то перпендикуляр $BE$ может не пересечь отрезка $AD$ (или пересечёт его продолжение за точку $A$ или за точку $D$);
в этом случае задача не имеет решения.
И независимо от построения можно видеть, что решения нет, если $s<b$ или $s=b$, потому что не может быть такого треугольника, у которого сумма двух сторон была бы меньше или равна третьей стороне.

Если задача имеет решение, то оно единственно;
то есть существует только один треугольник, удовлетворяющий требованиям задачи, так как перпендикуляр $BE$ может пересечься с прямой $AD$ только в одной точке.

{\small

\paragraph{}\label{1938/69}
\so{Замечание}.
Из приведённого примера видно, что решение сложной задачи на построение состоит из следующих четырёх частей.

1) Предположив, что задача решена, делают от руки приблизительный чертёж искомой фигуры и затем, внимательно рассматривая начерченную фигуру, стремятся найти такие зависимости между данными задачи и искомыми, которые позволили бы свести задачу к другим, известным ранее.
Эта самая важная часть решения задачи, имеющая целью составить план решения, носит название \textbf{анализа}.

2) Когда таким образом план решения найден, выполняют сообразно ему \textbf{построение}.

3) Для проверки правильности плана показывают затем на основании известных теорем, что полученная фигура удовлетворяет всем требованиям задачи.
Эта часть называется \textbf{синтезом}.

4) Затем задаются вопросом, при всяких ли данных задача имеет решение, допускает ли она одно решение или несколько, и нет ли в задаче каких-либо особенных случаев, когда построение упрощается или, наоборот, усложняется.
Эта часть решения называется \textbf{исследованием} задачи.

Когда задача очень проста и не может быть сомнения относительно существования её решения, то обыкновенно анализ и исследование опускают, а указывают прямо построение и приводят доказательство.
Так мы делали, излагая решение первых семи задач этой главы;
так же будем делать и впоследствии, когда нам придётся излагать решения несложных задач.

}

{\small

\subsection*{Упражнения}

\begin{center}\so{Доказать теоремы}
\end{center}

\begin{enumerate}[noitemsep]

\item
В равнобедренном треугольнике две медианы равны, две биссектрисы равны, две высоты равны.

\item
Если к каждой из равных сторон равнобедренного треугольника восстановим срединные перпендикуляры до пересечения с другой из равных сторон, то эти перпендикуляры будут равны. 

\item
Прямая, перпендикулярная к биссектрисе угла, отсекает от его сторон равные отрезки.

\item
Медиана треугольника меньше его полупериметра.

\item
Медиана треугольника меньше полусуммы сторон, между которыми она заключается.

\smallskip
\so{Указание}.
Продолжить медиану на расстояние, равное ей, полученную точку соединить с одним концом стороны, к которой проведена медиана, и рассмотреть образовавшуюся фигуру.

\item
Сумма медиан треугольника меньше периметра, но больше полупериметра.

\smallskip
\so{Указание}.
См. предыдущее упражнение, а также следствие в §~\ref{1938/50}.

\item
Сумма диагоналей четырёхугольника меньше его периметра, но больше полупериметра.

\item
Доказать как прямую теорему, что всякая точка, не лежащая на срединном перпендикуляре к отрезку, неодинаково удалена от концов этого отрезка, а именно: 
она ближе к тому концу, с которым она расположена по одну сторону от перпендикуляра.

\item
Доказать как прямую теорему, что всякая точка, не лежащая на биссектрисе угла, неодинаково отстоит от сторон его.

\item 
Медиана, исходящая из какой-нибудь вершины треугольника, равно отстоит от двух других его вершин.

\item
На одной стороне угла $A$ отложены отрезки $AB$ и $AC$ и на другой стороне отложены отрезки $AB'=AB$ и $AC' = AC$.
Доказать, что прямые $BC'$ и $B'C$ пересекаются на биссектрисе угла $A$.

\item
Вывести отсюда способ построения биссектрисы угла.

\item
Если $A'$ и $A$, $B'$ и $B$ — две пары точек, симметричных относительно какой-нибудь прямой $XY$, то четыре точки $A'$, $A$, $B'$, $B$ лежат на одной окружности.

\item
Дан острый угол $XOY$ и точка $A$ внутри этого угла.
Найти на стороне $OX$ точку $B$ и на стороне $OY$ точку $C$ так, чтобы периметр $\triangle ABC$ был наименьший.

\smallskip
\so{Указание}.
Надо взять точки, симметричные с $A$ относительно сторон $OX$ и $OY$.

\end{enumerate}

\begin{center}
\so{Задачи на построение}
\end{center}

\begin{enumerate}[resume,noitemsep]

\item
Построить сумму двух, трёх и более углов.

\item
Построить разность двух углов.

\item
По данной сумме и разности двух углов найти эти углы.

\item
Разделить угол на 4, 8 и 16 равных частей.

\item
Через вершину данного угла провести вне его такую прямую, которая со сторонами угла образовала бы равные углы.

\item
Построить треугольник:
а) по двум сторонам и углу между ними;
б) по стороне и двум прилежащим углам;
в) по двум сторонам и углу, лежащему против б\'{о}льшей из них;
г) по двум сторонам и углу, лежащему против меньшей из них (в этом случае получаются два решения, или одно, или ни одного).

\item
Построить равнобедренный треугольник:
а) по основанию и боковой стороне;
б) по основанию и прилежащему углу;
в) по боковой стороне и углу при вершине;
г) по боковой стороне и углу при основании.

\item
Построить прямоугольный треугольник:
а) по двум катетам;
б) по катету и гипотенузе;
в) по катету и прилежащему острому углу.

\item
Построить равнобедренный треугольник:
а) по высоте и боковой стороне;
б) по высоте и углу при вершине;
в) по основанию и перпендикуляру, опущенному из конца основания на боковую сторону.

\item
Построить прямоугольный треугольник по гипотенузе и острому углу.

\item
Через точку, данную внутри угла, провести такую прямую, которая отсекла бы от сторон угла равные части.

\item
По данной сумме и разности двух отрезков найти эти отрезки.

\item
Разделить данный отрезок на 4, 8, 16 равных частей.

\item
На данной прямой найти точку, одинаково удалённую от двух данных точек (вне прямой).

\item
Найти точку, равно отстоящую от трёх вершин треугольника.

\item
На прямой, пересекающей стороны угла, найти точку, одинаково удалённую от сторон этого угла.

\item
Найти точку, одинаково удалённую от трёх сторон треугольника.

\item
На бесконечной прямой $AB$ найти такую точку $C$, чтобы полупрямые $CM$ и $CN$, проведённые из $C$ через данные точки $M$ и $N$, расположенные по одну сторону от $AB$, составляли с полупрямыми $CA$ и $CD$ равные углы.

\smallskip
\so{Указание}.
Построить точку $M'$, симметричную с $M$ относительно оси $AB$, и соединить $M'$ с $N$.

\item
Построить прямоугольный треугольник по катету и сумме гипотенузы с другим катетом.

\item
Построить треугольник по основанию, углу, прилежащему к основанию, и разности двух других сторон.
(Рассмотреть два случая:
1) когда дан меньший из двух углов, прилежащих к основанию;
2) когда дан больший из них.)

\smallskip
\so{Указание}.
См. задачу §~\ref{1938/68}.

\item
Построить прямоугольный треугольник по катету и разности двух других сторон.

\item
Дан угол $A$ и точки $B$ и $C$, расположенные одна на одной стороне угла, другая — на другой.
Найти:
1) точку $M$, равно отстоящую от сторон угла, и такую, чтобы $MC=MB$;
2) точку $N$, равно отстоящую от сторон угла так, чтобы $NC=CB$.

\item
По соседству с железной дорогой расположены две деревни $A$ и $B$.
Найти на линии железной дороги (имеющей прямолинейную форму) место для станции, которая была бы одинаково удалена от $A$ и $B$.

\item
Дан угол $A$ и точка $B$ на одной из его сторон.
Найти на другой стороне такую точку $C$, чтобы сумма $CA+CB$ была равна данному отрезку~$\ell$.

\end{enumerate}

}

%% file: 2D/parallelnye.tex
\subsection*{Основные теоремы}

{

\begin{wrapfigure}{o}{27mm}
\vskip-5mm
\centering
\includegraphics{mppics/ris-extra-3}
\caption{}\label{extra/ris-3}
\end{wrapfigure}

\paragraph{}\label{1938/70}
\mbox{\so{Определение}.}
Две прямые называются \rindex{параллельность}\textbf{параллельными}, если они лежат в одной плоскости и \textbf{не пересекаются}, сколько бы их ни продолжали.
При этом прямая также считается параллельной самой себе.

Параллельность прямых обозначается письменно знаком $\parallel$.
Так, если прямые $AB$ и $CD$ параллельны, то пишут:
$AB \parallel CD$. 
На чертежах параллельные прямые принято отмечать одинаковыми стрелочками как на рис. \ref{extra/ris-3}.

Существование параллельных прямых обнаруживается следующей теоремой.

}

{

\begin{wrapfigure}{o}{34mm}
\vskip-6mm
\centering
\includegraphics{mppics/ris-72}
\caption{}\label{1938/ris-72}
\end{wrapfigure}

\paragraph{}\label{1938/71}
\mbox{\so{Теорема}.}
\textbf{\emph{Два перпендикуляра}} ($AB$ и $CD$, рис.~\ref{1938/ris-72}) \textbf{\emph{к одной и той же прямой}} ($MN$) \textbf{\emph{не могут пересечься, сколько бы мы их ни продолжали.}}

Действительно, если бы эти перпендикуляры пересеклись в какой-нибудь точке $P$, то из этой точки на прямую $MN$ были бы опущены два перпендикуляра, что невозможно (§~\ref{1938/24}).
Таким образом, \emph{два перпендикуляра к одной прямой параллельны между собой.}

}

\paragraph{Названия углов, получаемых при пересечении двух прямых третьей.}\label{1938/72}
Пусть две прямые $AB$ и $CD$ (рис.~\ref{1938/ris-73}) пересечены третьей прямой $MN$.
Тогда получаются 8 углов (мы их обозначили цифрами), которые попарно носят следующие названия.

\rindex{соответственные углы}\textbf{соответственные углы:}
1 и 5, 4 и 8, 2 и 6, 3 и 7.

\rindex{накрест лежащие углы}\textbf{накрест лежащие углы:}
3 и 5, 4 и 6 (внутренние);
1 и 7, 2 и 8 (внешние).

\rindex{односторонние углы}\textbf{односторонние углы:}
4 и 5, 3 и 6 (внутренние);
1 и 8, 2 и 7 (внешние).

\begin{wrapfigure}[12]{r}{42mm}
\vskip-2mm
\centering
\includegraphics{mppics/ris-73}
\caption{}\label{1938/ris-73}
\end{wrapfigure}

\paragraph{Признаки параллельности двух прямых.}\label{1938/73}
\textbf{\emph{Если при пересечении двух прямых}} ($AB$ и $CD$, рис.~\ref{1938/ris-74}) \textbf{\emph{третьей прямой}} ($MN$) \textbf{\emph{окажется, что:}}

1) \textbf{\emph{какие-нибудь соответственные углы равны, или}}

2) \textbf{\emph{какие-нибудь накрест лежащие углы равны, или}}

3) \textbf{\emph{сумма каких-нибудь двух внутренних или двух внешних односторонних углов равна $\bm{180\degree}$,}}

\textbf{\emph{то эти две прямые параллельны.}}

Пусть, например, дано, что соответственные углы 2 и 6 равны;
требуется доказать, что в таком случае $AB \parallel CD$.
Предположим противное, то есть что прямые $AB$ и $CD$ не параллельны;
тогда эти прямые пересекутся в какой-нибудь точке $P$, лежащей направо от $MN$, или в какой-нибудь точке $P'$, лежащей налево от $MN$.

\begin{figure}[!ht]
\centering
\includegraphics{mppics/ris-74}
\caption{}\label{1938/ris-74}
\end{figure}

Если пересечение будет в $P$, то образуется треугольник, в котором угол 2 будет внешним, а угол 6 — внутренним, не смежным с внешним углом 2, и, значит, тогда угол 2 должен быть больше угла 6 (§~\ref{1938/44}), что противоречит условию;
значит, пересечься в какой-нибудь точке $P$, лежащей направо от $MN$, прямые $AB$ и $CD$ не могут.

Если предположим, что пересечение будет в точке $P'$, то тогда образуется треугольник, у которого угол 4, равный углу 2, будет внутренним, а угол 6 — внешним, не смежным с внутренним углом 4;
тогда угол 6 должен быть больше угла 4 и, следовательно, больше угла 2, что противоречит условию.
Значит, прямые $AB$ и $CD$ не могут пересечься и в точке, лежащей налево от $MN$;
следовательно, эти прямые нигде не пересекаются, то есть они параллельны.

Подобным же образом доказывается, что $AB \parallel CD$, если $\angle 1 = \angle 5$ или $\angle 3 = \angle 7$ и~т.~д.

Пусть ещё дано, что $\angle 4 + \angle 5 = 180\degree$.
Тогда мы должны заключить, что $\angle 4 = \angle 6$, так как угол $6$ в сумме с углом $5$ тоже составляет $180\degree$.
Но если $\angle 4 = \angle 6$, то прямые не могут пересечься, так как в противном случае углы 4 и 6 не могли бы быть равными (один был бы внешний, а другой внутренний, не смежный с ним).

\paragraph{}\label{1938/74}
\so{Задача}.
\emph{Через данную точку $M$} (рис.~\ref{1938/ris-75}) \emph{провести прямую, параллельную данной прямой $AB$.}

Наиболее простое решение этой задачи состоит в следующем:
описываем произвольным радиусом с центром в точке $M$ дугу $CD$, далее описываем с центром в точке $C$ тем же радиусом дугу $ME$.
Затем, дав циркулю раствор, равный расстоянию от $E$ до $M$, описываем небольшую дугу с центром в точке $C$, которая пересечётся с $CD$ в некоторой точке $F$.
Прямая $MF$ будет параллельна $AB$.

Для доказательства проведём вспомогательную прямую $MC$;
образовавшиеся при этом углы 1 и 2 равны по построению (ибо треугольники $EMC$ и $MCD$ равны по трём сторонам), а если накрест лежащие углы равны, то линии параллельны.

\begin{figure}[!ht]
\begin{minipage}{.40\textwidth}
\centering
\includegraphics{mppics/ris-75}
\end{minipage}\hfill
\begin{minipage}{.56\textwidth}
\centering
\includegraphics{mppics/ris-wood-76}
\end{minipage}

\medskip

\begin{minipage}{.32\textwidth}
\centering
\captionof{figure}{}\label{1938/ris-75}
\end{minipage}\hfill
\begin{minipage}{.54\textwidth}
\centering
\caption{}\label{1938/ris-76}
\end{minipage}
\vskip-4mm
\end{figure}

Для построения параллельных прямых удобно пользоваться угольником и линейкой, как это видно из рис.~\ref{1938/ris-76}.

\paragraph{Аксиома параллельных линий.}\label{1938/75}
\textbf{\emph{Через одну и ту же точку нельзя провести двух различных прямых, параллельных одной и той же прямой.}}

Так, если (рис.~\ref{1938/ris-77}) $CE\parallel AB$, то никакая другая прямая $CE'$, проведённая через точку $C$, не может быть параллельной $AB$, то есть $CE'$ при продолжении пересечётся с $AB$.

Доказать это предложение, то есть вывести его как следствие из ранее принятых аксиом, оказывается невозможным.
Поэтому приходится принимать его как некоторое новое допущение (постулат или аксиому).

\begin{figure}[!ht]
\begin{minipage}{.48\textwidth}
\centering
\includegraphics{mppics/ris-77}
\end{minipage}\hfill
\begin{minipage}{.48\textwidth}
\centering
\includegraphics{mppics/ris-78}
\end{minipage}

\medskip

\begin{minipage}{.48\textwidth}
\centering
\caption{}\label{1938/ris-77}
\end{minipage}\hfill
\begin{minipage}{.48\textwidth}
\centering
\caption{}\label{1938/ris-78}
\end{minipage}
\vskip-4mm
\end{figure}

\paragraph{}\label{1938/76}
\mbox{\so{Следствия}.}
1) \emph{Если $CE\z\parallel AB$ \emph{(рис.~\ref{1938/ris-77})} и какая-нибудь третья прямая $CE'$ пересекается с одной из этих двух параллельных, то она пересекается и с другой.}
В противном случае через одну и ту же точку $C$ проходили бы две различные прямые $CE'$ и $CE$, параллельные $AB$, что невозможно.

2) \emph{Если каждая из двух прямых $A$ и $B$ \emph{(рис.~\ref{1938/ris-78})} параллельна одной и той же третьей прямой $C$, то они параллельны между собой.}

Действительно, если бы мы предположили, что прямые $A$ и $B$ пересекаются в некоторой точке $M$, то тогда через эту точку проходили бы две различные прямые, параллельные $C$, что невозможно.

\paragraph{Пара параллельных и их секущая.}\label{1938/77}\ 

\smallskip
\mbox{\so{Теорема} (обратная теорема, §~\ref{1938/73}).}\\
\textbf{\emph{Если две параллельные прямые}} ($AB$ и $CD$, рис.~\ref{1938/ris-79}) \textbf{\emph{пересечены какой-нибудь прямой}} ($MN$), \textbf{\emph{то:}}

1) \textbf{\emph{соответственные углы равны;}}

2) \textbf{\emph{накрест лежащие углы равны;}}

3) \textbf{\emph{сумма внутренних односторонних углов равна $\bm{180\degree}$;}}

4) \textbf{\emph{сумма внешних односторонних углов равна $\bm{180\degree}$.}}

\begin{wrapfigure}[12]{r}{41mm}
\vskip-3mm
\centering
\includegraphics{mppics/ris-79}
\caption{}\label{1938/ris-79}
\end{wrapfigure}

Докажем, например, что если $AB\z\parallel CD$, то соответственные углы $\alpha$ и $\beta$ равны.

Предположим противное, то есть что эти углы не равны (например, пусть $\alpha \z> \beta$).
Построив $\angle MEB_1 = \beta$, мы получим тогда прямую $A_1B_1$, не сливающуюся с $AB$, и, следовательно, будем иметь две различные прямые, проходящие через точку $E$ и параллельные одной и той же прямой $CD$, именно:
$AB\parallel CD$, согласно условию теоремы, и $A_1B_1\z\parallel CD$ вследствие равенства соответственных углов $\angle MEB_1$ и $ \beta$.
Так как это противоречит аксиоме параллельных линий, то наше предположение, что углы $\alpha$ и $\beta$ не равны, должно быть отброшено;
остаётся принять, что $ \alpha =  \beta$.

Таким же путём можно доказать и остальные заключения теоремы.
Из доказанных выше предложений непосредственно вытекает следующая теорема.

\textbf{\emph{Перпендикуляр к одной из двух параллельных прямых есть также перпендикуляр и к другой.}}

\begin{wrapfigure}[11]{o}{41mm}
\centering
\includegraphics{mppics/ris-80}
\caption{}\label{1938/ris-80}
\end{wrapfigure}

Действительно, если $AB\parallel CD$ (рис.~\ref{1938/ris-80}) и $ME\perp AB$, то, во-первых, $ME$, пересекаясь с $AB$, пересекается и с $CD$ в некоторой точке $F$, во-вторых, соответственные углы $\alpha$ и $\beta$ равны.
Но угол $\alpha$ прямой, значит, и угол $\beta$ прямой, то есть
$ME\perp CD$.

\paragraph{Признаки непараллельности прямых.}
\label{1938/78}
Из двух теорем:
прямой (§~\ref{1938/73}) и ей обратной (§~\ref{1938/77}) можно вывести заключение, что \so{противоположные теоремы} также верны, то есть:

\emph{Если при пересечении двух прямых третьей окажется, что 1) соответственные углы не равны или 2) внутренние накрест лежащие углы \so{не равны} и~т.~д., то прямые \so{не параллельны};
если две прямые \so{не параллельны}, то при пересечении их третьей прямой 1) соответственные углы \so{не равны}, 2) внутренние накрест лежащие углы \so{не равны}} и~т.~д.
Из этих признаков непараллельности (легко доказываемых способом от противного) полезно обратить особое внимание на следующий.

\begin{wrapfigure}[8]{o}{55mm}
\vskip-0mm
\centering
\includegraphics{mppics/ris-81}
\caption{}\label{1938/ris-81}
\end{wrapfigure}

\emph{если сумма внутренних односторонних углов \emph{($\alpha$ и $\beta$, рис.~\ref{1938/ris-81})} не равна $180\degree$, то прямые \emph{($AB$ и $CD$)} при достаточном продолжении пересекаются,} так как если бы эти прямые не пересекались, то они были бы параллельны, и тогда сумма внутренних односторонних углов равнялась бы $180\degree$, что противоречит условию.

Это предложение (дополненное утверждением, что прямые пересекутся по ту сторону от секущей линии, по которой сумма внутренних односторонних углов \so{меньше} $180\degree$) было принято знаменитым греческим геометром Евклидом (жившим в III веке до нашей эры) в его «Началах» геометрии без доказательства как аксиома параллельных линий, и потому оно известно под именем постулата Евклида.
В настоящее время предпочитают принимать за такую аксиому более простое предложение (§~\ref{1938/75}).

Укажем ещё два следующих признака непараллельности, которые понадобятся нам впоследствии:

{

1) \emph{Перпендикуляр \emph{($AB$, рис.~\ref{1938/ris-82})} и наклонная \emph{($CD$)} к одной и той же прямой \emph{($EF$)} пересекаются,} потому что сумма внутренних односторонних углов 1 и 2 не равна $180\degree$.

\begin{figure}[!ht]
\begin{minipage}{.48\textwidth}
\centering
\includegraphics{mppics/ris-82}
\end{minipage}\hfill
\begin{minipage}{.48\textwidth}
\centering
\includegraphics{mppics/ris-83}
\end{minipage}

\medskip

\begin{minipage}{.48\textwidth}
\centering
\caption{}\label{1938/ris-82}
\end{minipage}\hfill
\begin{minipage}{.48\textwidth}
\centering
\caption{}\label{1938/ris-83}
\end{minipage}
\vskip-4mm
\end{figure}

2) \emph{Две прямые \emph{($AB$ и $CD$, рис.~\ref{1938/ris-83}),} перпендикулярные к двум пересекающимся прямым \emph{($FE$ и $FG$),} пересекаются.}

Действительно, если предположим противное, то есть
что $AB \parallel CD$, то прямая $FG$, будучи перпендикулярна к одной из параллельных (к $CD$), была бы перпендикулярна и к другой параллельной (к $AB$), и тогда из одной точки $F$ к прямой $AB$ были бы проведены два перпендикуляра:
$FB$ и $FD$, что невозможно.

}

\subsection*{Углы с соответственно параллельными или перпендикулярными сторонами}

\paragraph{}\label{1938/79}
\so{Теорема}.
\textbf{\emph{Если стороны одного угла соответственно параллельны сторонам другого угла, то такие углы или равны, или в сумме составляют два прямых.}}

\begin{wrapfigure}{o}{43mm}
\vskip-0mm
\centering
\includegraphics{mppics/ris-84}
\caption{}\label{1938/ris-84}
\end{wrapfigure}

Рассмотрим особо следующие три случая (рис.~\ref{1938/ris-84}).

1) Пусть стороны угла 1 соответственно параллельны сторонам угла 2 и, сверх того, имеют \so{одинаковое направление от вершины} (на чертеже направления указаны стрелками).

Продолжив одну из сторон угла 2 до пересечения с непараллельной ей стороной угла 1, мы получим угол 3, равный и углу 1, и углу 2 (как соответственные при параллельных прямых);
следовательно, $\angle 1 = \angle 2$.

2) Пусть стороны угла 1 соответственно параллельны сторонам угла 4, но имеют противоположное направление от вершины.

Продолжив обе стороны угла 4, мы получим угол 2, который равен углу 1 (по доказанному выше) и углу 4 (как вертикальный ему);
следовательно, $\angle 4 = \angle 1$.

3) Пусть, наконец, стороны угла 1 соответственно параллельны сторонам углов 5 и 6, причём две из этих сторон имеют одинаковое направление, а две другие — противоположное.

Продолжив одну сторону угла 5 или угла 6, мы получим угол 2, который равен (по доказанному) углу 1;
но $\angle 5$ (или $\angle 6) +\angle 2 \z= 180\degree$ (по свойству смежных углов);
следовательно, и $\angle 5$ (или $\angle 6) \z+ \angle 1 \z= 180\degree$.

Таким образом, углы с параллельными сторонами оказываются равными, когда их стороны имеют или одинаковое, или противоположное направление от вершины;
если же это условие не выполнено, то углы составляют в сумме $180\degree$.

{\small
\smallskip
\so{Замечание}.
Можно было бы сказать, что углы с параллельными сторонами равны тогда, когда они оба острые или оба тупые;
но бывают случаи, когда по виду углов трудно определить, острые ли они или тупые;
поэтому приходится сравнить направления сторон углов.
}

\paragraph{}\label{1938/80}
\so{Теорема}.
\textbf{\emph{Если стороны одного угла соответственно перпендикулярны к сторонам другого угла, то такие углы или равны, или в сумме составляют два прямых.}}

Пусть угол $ABC$, обозначенный цифрой 1 (рис.~\ref{1938/ris-85}), есть один из данных углов;
за другой данный угол возьмём какой-нибудь из четырёх углов:
2, 3, 4 или 5, образованных двумя пересекающимися прямыми, из которых одна перпендикулярна к стороне $AB$, а другая — к стороне $BC$ (общая вершина их может находиться в любой точке плоскости).

\begin{wrapfigure}[11]{o}{35mm}
\centering
\includegraphics{mppics/ris-85}
\caption{}\label{1938/ris-85}
\end{wrapfigure}

Проведём из вершины угла 1 две вспомогательные прямые:
$BD\perp BC$ и $BE\z\perp BA$.
Образованный ими угол 6 равен углу 1 по следующей причине:
углы $DBC$ и $EBA$ равны, так как оба они прямые;
отняв от каждого из них по одному и тому же углу $EBC$, получим: $\angle 6 = \angle 1$.

Заметим, что стороны вспомогательного угла 6 параллельны пересекающимся прямым, образующим углы 2, 3, 4 и 5 (потому что два перпендикуляра к одной прямой параллельны, §~\ref{1938/71}). 
Следовательно, эти углы или равны углу 6, или составляют с ним в сумме $180\degree$.
Заменив угол 6 равным ему углом 1, получим то, что требовалось доказать.

\subsection*{Сумма углов треугольника и многоугольника}

\begin{wrapfigure}[8]{r}{41mm}
\vskip-2mm
\centering
\includegraphics{mppics/ris-86}
\caption{}\label{1938/ris-86}
\end{wrapfigure}

\paragraph{}\label{1938/81}
\mbox{\so{Теорема}.}
\textbf{\emph{Сумма углов треугольника равна $\bm{180\degree}$.}}

Пусть $ABC$ (рис.~\ref{1938/ris-86}) — какой-нибудь треугольник;
требуется доказать, что сумма углов $A$, $B$ и $C$ равна $180\degree$.
Продолжив сторону $AC$ и проведя $CE\parallel AB$, найдём:
$\angle A \z= \angle ECD$ (как углы, соответственные при параллельных), $\angle B = \angle BCE$ (как углы, накрест лежащие при параллельных). 
Отсюда
\[\angle A + \angle B+\angle C = \angle ECD + \angle BCE + \angle C =180\degree.\]

\smallskip
\so{Следствия}.
1) \emph{Всякий внешний угол треугольника равен сумме двух внутренних углов, не смежных с ним} (так, $\angle BCD \z= \angle A + \angle B$).

2) \emph{Если два угла одного треугольника соответственно равны двум углам другого, то и третьи углы равны.}

{

\begin{wrapfigure}{o}{38mm}
\centering
\includegraphics{mppics/ris-87}
\caption{}\label{1938/ris-87}
\end{wrapfigure}

3) \emph{Сумма двух острых углов прямоугольного треугольника равна одному прямому углу, то есть $90\degree$.}

4) \emph{В равнобедренном прямоугольном треугольнике каждый острый угол равен $45\degree$.}

5) \emph{В равностороннем треугольнике каждый угол равен $60\degree$.}

6) \emph{Если в прямоугольном треугольнике $ABC$ \emph{(рис.~\ref{1938/ris-87})} один из острых углов \emph{(например $\angle B$)} равен $30\degree$, то лежащий против него катет составляет половину гипотенузы.}

}

Заметив, что в таком треугольнике другой острый угол равен $60\degree$, пристроим к треугольнику $ABC$ другой треугольник $ABD$, равный данному.
Тогда мы получим треугольник $DBC$, у которого каждый угол равен $60\degree$.
Такой равноугольный треугольник должен быть равносторонним (§~\ref{1938/47}), и потому $BC\z=DC$.
Но $AC=\tfrac12DC$, значит, $AC=\tfrac12BC$.

{\sloppy 

Предоставляем самим учащимся доказать обратное предложение:
\emph{если катет равен половине гипотенузы, то противолежащий ему острый угол равен $30\degree$.}

}

\begin{wrapfigure}{r}{30mm}
\vskip-5mm
\centering
\includegraphics{mppics/ris-88}
\caption{}\label{1938/ris-88}
\bigskip
\includegraphics{mppics/ris-89}
\caption{}\label{1938/ris-89}
\end{wrapfigure}

\paragraph{}\label{1938/82}
\mbox{\so{Теорема}.}
\textbf{\emph{Сумма углов выпуклого $\bm{n}$-угольника равна $\bm{180\degree\cdot(n-2)}$.}}

Взяв внутри выпуклого $n$-угольника произвольную точку $O$ (рис. \ref{1938/ris-88}), соединим её со всеми вершинами.
Тогда $n$-угольник разобьётся на $n$ треугольников — к каждой его стороне примыкает один треугольник.

Сумма углов каждого треугольника равна $180\degree$, следовательно, сумма углов всех треугольников равна $180\degree\cdot n$.
Эта величина, очевидно, превышает сумму углов $n$-угольника на сумму всех тех углов, которые расположены вокруг точки $O$;
но эта сумма равна $360\degree$;
следовательно, сумма углов $n$-угольника равна:
\[180\degree\cdot n -360\degree = 180\degree\cdot(n - 2).\]

{\small 

\smallskip
\mbox{\so{Замечание}.}
Эту теорему можно доказать ещё и так:
Из вершины какого-нибудь выпуклого $n$-угольника проведём его диагонали (рис.~\ref{1938/ris-89}).
Тогда $n$-угольник разобьётся на $n-2$ треугольника.
Действительно, если не будем считать двух сторон, образующих угол, из вершины которого проведены диагонали, то на каждую из остальных сторон придётся по одному треугольнику.
Но в каждом треугольнике сумма углов равна $180\degree$.
Значит, сумма углов всех треугольников будет $180\degree\cdot(n-2)$;
но эта сумма и есть сумма всех углов $n$-угольника.

{\small

\medskip
\smallskip
\so{Примечание}.
Доказанная теорема верна и для вогнутых многоугольников.
При этом, если внутри $n$-угольника можно найти такую точку, что отрезки, соединяющие её с вершинами $n$-угольника, лежат внутри него, то теорему можно доказать, повторяя дословно те рассуждения, которые мы приводили выше при первом способе доказательства.
Если же такой точки найти нельзя, то следует весь $n$-угольник разбить на треугольники, проведя некоторые его диагонали — можно доказать, что такое разбиение всегда существует и при этом всех треугольников будет $n-2$.
Подсчитав затем сумму углов в каждом из треугольников и сложив эти суммы, получим ту же формулу $180\degree\cdot(n - 2)$.

}

}

{

\begin{wrapfigure}{o}{33mm}
\centering
\includegraphics{mppics/ris-90}
\caption{}\label{1938/ris-90}
\end{wrapfigure}

\paragraph{}\label{1938/83}
\mbox{\so{Теорема}.}
\textbf{\emph{Если из вершины каждого угла выпуклого многоугольника проведём продолжение одной из сторон этого угла, то сумма всех образовавшихся при этом внешних углов многоугольника равна $\bm{360\degree}$}} (независимо от числа сторон многоугольника).

Каждый из таких внешних углов (рис.~\ref{1938/ris-90}) составляет дополнение до $180\degree$ к смежному с ним внутреннему углу многоугольника.

Следовательно, если к сумме всех внутренних углов $n$-угольника прибавим сумму всех внешних углов, то получим $180\degree\cdot n$;
но сумма внутренних углов, как мы видели, равна $180\degree\cdot (n - 2)$;
следовательно, сумма внешних углов равна разности:
\begin{align*}
&180\degree\cdot n-180\degree\cdot (n - 2)=
\\
=~&180\degree\cdot n-180\degree\cdot n+360\degree=
\\
=~&360\degree.
\end{align*}

}

\subsection*{Центральная симметрия}

\paragraph{}\label{1938/84}
В §~\ref{1938/37} был рассмотрен случай симметричного расположения двух равных фигур относительно прямой.
Выведенные выше свойства параллельных прямых позволяют изучить ещё один замечательный вид расположения двух равных фигур, или двух равных отрезков, или двух точек по отношению к некоторой точке на плоскости.

\begin{wrapfigure}{o}{25mm}
\vskip-0mm
\centering
\includegraphics{mppics/ris-91}
\caption{}\label{1938/ris-91}
\end{wrapfigure}

\textbf{Если две какие-либо точки $\bm{A}$ и $\bm{A'}$} (рис.~\ref{1938/ris-91}) \textbf{расположены на одной прямой с точкой $\bm{O}$ по разные стороны от неё и на одинаковом от неё расстоянии} ($OA=OA'$), \textbf{то такие точки называются симметричными относительно точки} ($O$).

Чтобы построить точку, симметричную с данной точкой $A$ относительно другой данной точки $O$, следует соединить точки $A$ и $O$ прямой, продолжить эту прямую за точку $O$ и отложить на ней от точки $O$ отрезок $OA'$, равный $OA$, таким образом, чтобы точки $A$ и $A'$ были расположены по разные стороны относительно точки $O$.
Точка $A'$ будет искомой.

\paragraph{}\label{1938/85}
\so{Теорема}.
\textbf{\emph{Если для двух точек $\bm{A}$ и $\bm{B}$ какой-либо прямой $\bm{AB}$ построить симметричные им точки $\bm{A'}$ и $\bm{B'}$ относительно некоторой точки $\bm{O}$, то:}}

\begin{wrapfigure}{r}{35mm}
\vskip-4mm
\centering
\includegraphics{mppics/ris-92}
\caption{}\label{1938/ris-92}
\end{wrapfigure}

1) \textbf{\emph{Прямая, соединяющая точки $\bm{A'}$ и $\bm{B'}$, будет параллельна данной прямой $\bm{AB}$, причём отрезок $\bm{AB}$ равен отрезку $\bm{A'B'}$.}}

2) \textbf{\emph{Каждой точке данной прямой $\bm{AB}$ соответствует симметричная ей точка на построенной прямой $\bm{A'B'}$.}}

\smallskip

1) Треугольники $AOB$ и $A'OB'$ равны (рис.~\ref{1938/ris-92}), потому что у них $AO=A'O$ и $BO\z=B'O$ (по построению), $\angle AOB=\angle A'OB'$ (как вертикальные углы).
Из равенства этих треугольников следует:
$AB=A'B'$ и $\angle OAB \z= \angle OA'B'$;
значит, $AB\parallel A'B'$ (§~\ref{1938/73}, 2-й случай).

2) Возьмём на прямой $AB$ какую-либо точку $D$ (рис.~\ref{1938/ris-92}).
Рассмотрим прямую, соединяющую точку $D$ с точкой $O$.
Эта прямая пересечёт прямую $A'B'$ в некоторой точке $D'$.
Треугольники $AOD$ и $A'OD'$ равны, потому что у них $AO\z=A'O$, $\angle 1 = \angle 2$ (как накрест лежащие при параллельных прямых) и $\angle 3 = \angle 4$ (как вертикальные).
Из равенства этих треугольников следует:
$OD = OD'$.
Значит, точки $D$ и $D'$ симметричны относительно точки $O$.

\paragraph{Симметричные фигуры.}\label{1938/86}
\textbf{\emph{Две фигуры называются симметричными относительно данной точки $\bm{O}$, если каждой точке одной фигуры соответствует симметричная ей точка другой фигуры.}} 

Точка $O$ называется \rindex{центр!симметрии}\textbf{центром симметрии} данных фигур.
Сама симметрия называется \rindex{центральная симметрия}\textbf{центральной} в отличие от осевой, с которой мы уже встречались раньше (§~\ref{1938/37}).
Если каждой точке данной фигуры соответствует симметричная ей точка той же самой фигуры (относительно некоторого центра), то говорят, что данная фигура имеет центр симметрии.
Примером такой фигуры служит окружность.
Центром её симметрии является её центр.

\textbf{Каждую фигуру можно совместить с фигурой, ей симметричной, путём вращения её вокруг центра симметрии.}

В самом деле, возьмём, например, два треугольника $ABC$ и $A'B'C'$ (рис.~\ref{1938/ris-93}), симметричные относительно некоторого центра $O$.

Всю фигуру $OABC$, не отрывая от плоскости, будем вращать вокруг точки $O$ как вокруг центра до тех пор, пока прямая $OA$ не пойдёт по $OA'$.

Так как $\angle 1 = \angle 2$ и $\angle 3 = \angle 4$, то прямая $OB$ пойдёт по $OB'$, а прямая $OC$ по $OC'$.

\begin{wrapfigure}{o}{45mm}
\vskip-2mm
\centering
\includegraphics{mppics/ris-93}
\caption{}\label{1938/ris-93}
\end{wrapfigure}

Так как $OA = OA'$, $OB=OB'$, $OC\z=OC'$, то точка $A$ совпадёт с $A'$, точка $B$ с $B'$ и точка $C$ с $C'$.
Таким образом, треугольник $ABC$ совместится с треугольником $A'B'C'$.

Очевидно, что при таком повороте каждая прямая $OA$, $OB$, $OC$, а также каждая сторона треугольника $ABC$ повернётся на $180\degree$.
Если фигура имеет центр симметрии, то после поворота её вокруг центра симметрии на $180\degree$ эта фигура совместится сама с собой.

\begin{figure}[!ht]
\begin{minipage}{.68\textwidth}
\centering
\begin{lpic}[t(1 mm),b(1 mm),r(0 mm),l(0 mm)]{jpg/Manning_propeller(.8)}
\lbl[b]{54,19;$O$}
\end{lpic}
\end{minipage}
\hfill
\begin{minipage}{.28\textwidth}
\centering
\includegraphics[scale=.19]{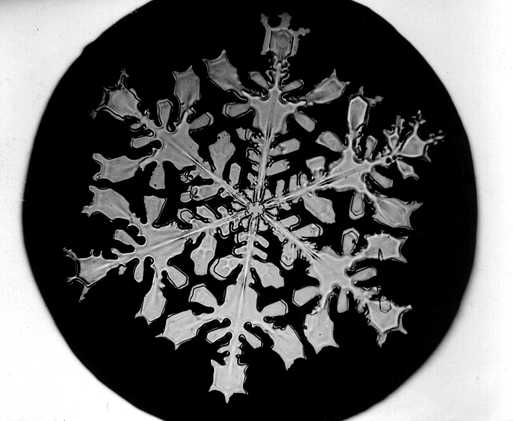}
\end{minipage}

\begin{minipage}{.68\textwidth}
\centering
\caption{}\label{1938/ris-94}
\end{minipage}
\hfill
\begin{minipage}{.28\textwidth}
\centering
\caption{}\label{1938/ris-95}
\end{minipage}
\vskip-4mm
\end{figure}

{\small

\medskip

\so{Замечание}.
При вращении, которое мы произвели для совмещения треугольников $ABC$ и $A'B'C'$, треугольник $ABC$ скользил по плоскости.
Таким образом, фигуры, симметричные относительно центра, можно совместить, не выводя их из плоскости.
Этим центральная симметрия существенно отличается от осевой (§~\ref{1938/37}), где для совмещения симметричных фигур необходимо было одну из них перевернуть другой стороной.

Центральная симметрия фигур, так же как и осевая, весьма часто встречается в природе и в обыденной жизни.
На рис.~\ref{1938/ris-94} приведено изображение пропеллера самолёта.
Оно имеет центром симметрии точку $O$.
На рис.~\ref{1938/ris-95} дано изображение снежинки, оно также обладает центром симметрии.

}

%% file: 2D/parallelogrammy.tex
\section{Параллелограммы и трапеции}

\subsection*{Параллелограммы}

{

\begin{wrapfigure}{r}{40mm}
\vskip-7mm
\centering
\includegraphics{mppics/ris-96}
\caption{}\label{1938/ris-96}
\end{wrapfigure}

\paragraph{Параллелограмм.}\label{1938/87}
Четырёхугольник, у которого противоположные стороны попарно параллельны, называется \rindex{параллелограмм}параллелограммом.
Такой четырёхугольник ($ABCD$, рис. \ref{1938/ris-96}) получится, например, если какие-нибудь две параллельные прямые $KL$ и $MN$ пересечём двумя другими параллельными прямыми $RS$ и $PQ$.

}

\paragraph{}\label{1938/88}
\mbox{\so{Теорема}} (выражающая свойство сторон и углов параллелограмма).
\textbf{\emph{Во всяком параллелограмме противоположные стороны равны, противоположные углы равны и сумма углов, прилежащих к одной стороне, равна $\bm{180\degree}$}} (рис.~\ref{1938/ris-97}).

\begin{wrapfigure}{o}{35mm}
\centering
\includegraphics{mppics/ris-97}
\caption{}\label{1938/ris-97}
\end{wrapfigure}

Проведя диагональ $BD$, мы получим два треугольника:
$ABD$ и $BCD$, которые равны, потому что у них $BD$ — общая сторона, $\angle 1 \z= \angle 4$ и $\angle 2 = \angle 3$ (как накрест лежащие при параллельных прямых).
Из равенства треугольников следует:
$AB=CD$, $AD=BC$ и $\angle A \z= \angle C$.
Противоположные углы $B$ и $D$ также равны, так как они представляют собой суммы равных углов.

Наконец, углы, прилежащие к одной стороне, например углы $A$ и $B$, дают в сумме $180\degree$, так как это углы внутренние односторонние при параллельных прямых.

{\small
\smallskip
\so{Замечание}.
Равенство противоположных сторон параллелограмма иногда кратко выражают другими словами, так:
\emph{отрезки параллельных, отсекаемые параллельными, равны.}

}

{

\begin{wrapfigure}{r}{41mm}
\centering
\vskip-6mm
\includegraphics{mppics/ris-98}
\caption{}\label{1938/ris-98}
\end{wrapfigure}

\smallskip
\mbox{\so{Следствие}.}
\emph{Если две прямые параллельны, то все точки каждой из них одинаково удалены от другой параллельной;
короче:
параллельные прямые \emph{($AB$ и $CD$ рис.~\ref{1938/ris-98})} везде одинаково удалены одна от другой.}

Действительно, если из каких-нибудь двух точек $M$ и $N$ прямой $CD$ опустим на $AB$ перпендикуляры $MP$ и $NQ$, то эти перпендикуляры параллельны (§~\ref{1938/71}), и потому фигура $MNQP$ — параллелограмм;
отсюда следует, что $MP=NQ$, то есть точки $M$ и $N$ одинаково удалены от прямой~$AB$.

}

\paragraph{Два признака параллелограммов.}\label{1938/89}\ 

\smallskip
\mbox{\so{Теорема}.}
\textbf{\emph{Если в выпуклом четырёхугольнике:}}

1) \textbf{\emph{противоположные стороны равны между собой или}}

2) \textbf{\emph{две противоположные стороны равны и параллельны, то такой четырёхугольник есть параллелограмм.}}

Рассмотрим четырёхугольник $ABCD$ (рис.~\ref{1938/ris-99}), у которого.
\[AB=CD\quad \text{и}\quad BC=AD.\]
Требуется доказать, что эта фигура — параллелограмм, то есть что $AB\parallel CD$ и $BC \parallel AD$.
Проведя диагональ $BD$, мы получим два треугольника, которые равны, так как у них $BD$ — общая сторона, $AB\z=CD$ и $BC = AD$ (по условию).
Из равенства этих треугольников следует:
$\angle 1 = \angle 4 $ и $\angle 2 = \angle 3$ (в равных треугольниках против равных сторон лежат равные углы);
вследствие этого $AB \parallel CD$ и $BC\parallel AD$ (если накрест лежащие углы равны, то прямые параллельны).

\begin{wrapfigure}{o}{36mm}
\centering
\includegraphics{mppics/ris-99}
\caption{}\label{1938/ris-99}
\end{wrapfigure}

2) Пусть в четырёхугольнике ($ABCD$, рис.~\ref{1938/ris-99}) $BC\parallel AD$ и $BC = AD$.
Требуется доказать, что $ABCD$ есть параллелограмм, то есть что $AB \parallel CD$.

Треугольники $ABD$ и $BCD$ равны, потому что у них $BD$ — общая сторона, $BC = AD$ (по условию) и $\angle 2 = \angle 3$ (как накрест лежащие углы при параллельных прямых).
Из равенства треугольников следует:
$\angle 1 = \angle 4$;
поэтому $AB\parallel CD$.

\paragraph{}\label{1938/90}
\so{Теорема} (выражающая свойство диагоналей параллелограмма).
\textbf{\emph{Если четырёхугольник}} ($ABCD$, рис.~\ref{1938/ris-100}) \textbf{\emph{— параллелограмм, то его диагонали, пересекаясь, делятся пополам.}}

{

\begin{wrapfigure}[9]{r}{44mm}
\vskip-4mm
\centering
\includegraphics{mppics/ris-100}
\caption{}\label{1938/ris-100}
\end{wrapfigure}

Обратно:
\textbf{\emph{если в четырёхугольнике диагонали точкой их пересечения делятся пополам, то данный четырёхугольник — параллелограмм.}}

1) Треугольники $BOC$ и $AOD$ равны, потому что у них:
$BC\z=AD$ (как противоположные стороны параллелограмма), $\angle 1 = \angle 2$ и $\angle 3 = \angle 4$ (как накрест лежащие при параллельных прямых).

Из равенства треугольников следует:
$OC=OA$ и $OB=OD$.

2) Если $AO=OC$ и $BO=OD$, то треугольники $AOD$ и $BOC$ равны (по двум сторонам и углу между ними).
Из равенства треугольников следует:
$\angle 1 = \angle 2$ и $\angle 3 = \angle 4$.
Следовательно, $BC \parallel AD$ (углы накрест лежащие равны) и $BC=AD$;
поэтому фигура $ABCD$ есть параллелограмм.

}

\paragraph{Центр симметрии.}\label{1938/91}
\textbf{\emph{Параллелограмм имеет центр симметрии}}, причём центром симметрии служит точка пересечения диагоналей (рис.~\ref{1938/ris-100}).

Действительно, так как $BO=OD$ и $OC=OA$, то отрезки $BC$ и $AD$ симметричны относительно точки $O$ и каждой точке $P$ отрезка $BC$ соответствует симметричная ей точка $Q$ отрезка $AD$ (§~\ref{1938/85}).
Таким же образом убеждаемся, что отрезки $AB$ и $CD$ симметричны относительно той же точки $O$.
Если параллелограмм повернуть вокруг точки пересечения его диагоналей на $180\degree$, то новое положение параллелограмма совпадёт с первоначальным.
При этом каждая из его вершин поменяется местом с противоположной вершиной 
(на рис.~\ref{1938/ris-100} вершина $A$ с $C$ и $B$ с $D$).

\subsection*{Прямоугольник, ромб и квадрат}

\paragraph{Прямоугольник и его свойства.}\label{1938/92}
Если один из углов параллелограмма прямой, то три остальных его угла также прямые (§~\ref{1938/88}).
Параллелограмм, у которого все углы прямые, называется \rindex{прямоугольник}\textbf{прямоугольником}.

Так как прямоугольник есть параллелограмм, то он обладает всеми свойствами параллелограмма;
например, диагонали его делятся пополам и точка их пересечения есть центр симметрии.
Но у прямоугольника есть ещё свои особые свойства.

1) \emph{В прямоугольнике \emph{($ABCD$, рис.~\ref{1938/ris-101})} диагонали равны.}

Прямоугольные треугольники $ACD$ и $ABD$ равны, потому что у них:
$AD$ — общий катет и $AB=CD$ (как противоположные стороны параллелограмма).
Из равенства треугольников следует:
$AC = BD$.

2) \emph{Прямоугольник имеет две оси симметрии.}
Именно, каждая прямая, проходящая через центр симметрии прямоугольника и параллельная двум его противоположным сторонам, есть ось его симметрии.
Эти оси симметрии перпендикулярны между собой (смотри рис. \ref{1938/ris-102}).

\begin{figure}[h]
\begin{minipage}{.32\textwidth}
\centering
\includegraphics{mppics/ris-101}
\end{minipage}\hfill
\begin{minipage}{.32\textwidth}
\centering
\includegraphics{mppics/ris-102}
\end{minipage}\hfill
\begin{minipage}{.32\textwidth}
\centering
\includegraphics{mppics/ris-103}
\end{minipage}

\medskip

\begin{minipage}{.32\textwidth}
\centering
\caption{}\label{1938/ris-101}
\end{minipage}\hfill
\begin{minipage}{.32\textwidth}
\centering
\caption{}\label{1938/ris-102}
\end{minipage}\hfill
\begin{minipage}{.32\textwidth}
\centering
\caption{}\label{1938/ris-103}
\end{minipage}
\vskip-4mm
\end{figure} 

\paragraph{Ромб и его свойства.}\label{1938/93}
Параллелограмм, у которого все стороны равны, называется \rindex{ромб}ромбом.
Конечно, ему принадлежат все свойства параллелограмма, но у него есть следующие два особых свойства.

1) \emph{Диагонали ромба \emph{($ABCD$, рис.~\ref{1938/ris-103})} взаимно перпендикулярны и делят углы ромба пополам.}
Треугольники $AOB$ и $COB$ равны, потому что у них:
$BO$ — общая сторона, $AB=BC$ (так как у ромба все стороны равны) и $AO=OC$ (так как диагонали всякого параллелограмма делятся пополам).
Из равенства треугольников следует:
\[\angle 1 =\angle 2,
\quad
\text{то есть}
\quad
BD\perp AC
\quad
\text{и}
\quad
\angle 3 =\angle 4,
\]

\begin{wrapfigure}[20]{o}{34mm}
\vskip-6mm
\centering
\includegraphics{mppics/ris-104}
\captionof{figure}{}\label{1938/ris-104}
\bigskip
\includegraphics{mppics/ris-105}
\caption{}\label{1938/ris-105}
\end{wrapfigure}

\noindent
то есть угол $B$ делится диагональю пополам.
Из равенства треугольников $BOC$ и $COD$ заключаем, что угол $C$ делится диагональю пополам.

2) \emph{Каждая диагональ ромба есть его ось симметрии.}

Диагональ $BD$ (рис.~\ref{1938/ris-104}) является осью симметрии ромба $ABCD$,
то есть вращая $\triangle ABD$ вокруг $BD$, мы можем совместить его с $\triangle BCD$.
В самом деле, диагональ $BD$ делит углы $B$ и $D$ пополам и, кроме того, $AB=BC$ и $AD=CD$.

То же самое можно сказать о диагонали $AC$.

\paragraph{Квадрат и его свойства.}\label{1938/94}\rindex{квадрат}
Квадратом называется параллелограмм, у которого все стороны равны и все углы прямые;
можно также сказать, что квадрат — это прямоугольник, у которого стороны равны, или ромб, у которого углы прямые.
Поэтому квадрату принадлежат все свойства параллелограмма, прямоугольника и ромба.
Например, у квадрата имеется четыре оси симметрии (рис.~\ref{1938/ris-105}):
две, проходящие через середины противоположных сторон (как у прямоугольника), и две, проходящие через вершины противоположных углов (как у ромба).

\subsection*{Теоремы, основанные на свойствах параллелограмма}

\paragraph{}\label{1938/95}
\so{Теорема}.
\textbf{\emph{Если на одной стороне угла}} (например, на стороне $BC$ угла $ABC$, рис.~\ref{1938/ris-106}) \textbf{\emph{отложим равные между собой отрезки}} ($DE=EF=\dots$) \textbf{\emph{и через их концы проведём параллельные прямые}} ($DM$, $EN$, $FP,\dots$) \textbf{\emph{до пересечения с другой стороной угла, то и на этой стороне отложатся равные между собой отрезки}} ($MN=NP=\dots$).
Проведём вспомогательные прямые $DK$ и $EL$, параллельные $AB$.
Полученные при этом треугольники $DKE$ и $ELF$ равны, так как у них:
$DE=EF$ (по условию), $\angle KDE=\angle LEF$ и $\angle KED = \angle LFE$ (как углы, соответственные при параллельных прямых).
Из равенства этих треугольников следует:
$DK=EL$.
Но $DK\z=MN$ и $EL\z=NP$ (как противоположные стороны параллелограммов);
значит, $MN=NP$.

\begin{figure}[!ht]
\centering
\includegraphics{mppics/ris-106}
\caption{}\label{1938/ris-106}
\end{figure}

{\small
\smallskip
\mbox{\so{Замечание}.}
Равные отрезки могут быть откладываемы и от вершины угла $B$, то есть
так:
$BD=DE= EF=\dots$
Тогда и на другой стороне равные отрезки надо считать от вершины угла, то есть так:
$BM=MN=NP=\dots$

}

\paragraph{}\label{1938/96}
\mbox{\so{Следствие}.}
\emph{Прямая \emph{($DE$, рис.~\ref{1938/ris-107}),} проведённая через середину стороны ($AB$) треугольника параллельно другой его стороне ($AC$), делит третью сторону ($BC$) пополам.}

\begin{wrapfigure}[8]{r}{30mm}
\vskip-4mm
\centering
\includegraphics{mppics/ris-107}
\caption{}\label{1938/ris-107}
\end{wrapfigure}

Действительно, мы видим, что на стороне угла $B$ отложены равные отрезки $BD=DA$ и через точки деления $D$ и $A$ проведены параллельные прямые $DE$ и $AC$ до пересечения со стороной $BC$;
значит, по доказанному, на этой стороне тоже отложатся равные отрезки $BE=EC$, и потому $BC$ разделится в точке $E$ пополам.

{\small
\smallskip
\mbox{\so{Замечание}.}
Отрезок, соединяющий середины двух сторон треугольника, называется его \rindex{средняя линия!треугольника}\textbf{средней линией}.

}

\paragraph{}\label{1938/97}
\so{Теорема} (о средней линии треугольника).
\textbf{\emph{Прямая}} ($DE$, рис. \ref{1938/ris-107}), \textbf{\emph{проведённая через середины двух сторон треугольника, параллельна третьей его стороне;
отрезок этой прямой, лежащий внутри треугольника, равен половине третьей стороны.}}

Для доказательства вообразим, что через середину $D$ стороны $AB$ мы провели прямую, параллельную стороне $AC$.
Тогда, по доказанному в предыдущем параграфе, эта прямая разделит сторону $BC$ пополам и, следовательно, сольётся с прямой $DE$, соединяющей середины сторон $AB$ и $BC$.

Проведя ещё $EF \parallel AD$, найдём, что сторона $AC$ также разделится пополам в точке $F$;
значит, $AF=FC$ и, кроме того, $AF=DE$
(как противоположные стороны параллелограмма $ADEF$), откуда следует:
\[DE=\tfrac12AC.\]

\renewcommand{\bottomtitlespace}{.11\textheight}

\subsection*{Трапеции}

\begin{wrapfigure}{r}{40mm}
\vskip-6mm
\centering
\includegraphics{mppics/ris-ru-108}
\caption{}\label{1938/ris-108}
\end{wrapfigure}

\paragraph{}\label{1938/98}
Четырёхугольник, у которого две противоположные стороны параллельны, а две другие не параллельны, называется \rindex{трапеция}\textbf{трапецией}.
Параллельные стороны трапеции ($AD$ и $BC$) называются \rindex{основание!трапеции}\textbf{её основаниями}, непараллельные ($AB$ и $CD$) — \rindex{боковая сторона трапеции}\textbf{боковыми сторонами} (рис.~\ref{1938/ris-108}).
Если боковые стороны равны, трапеция называется \rindex{равнобочная трапеция}\textbf{равнобочной}.

\paragraph{Средняя линия трапеции.}\label{1938/99}
Прямая, соединяющая середины боковых сторон трапеции, называется \rindex{средняя линия!трапеции}\textbf{её средней линией}.
Линия эта обладает следующим свойством.

\smallskip
\mbox{\so{Теорема}.}
\textbf{\emph{Средняя линия}} ($EF$, рис.~\ref{1938/ris-109}) \textbf{\emph{трапеции параллельна основаниям, и равна их полусумме.}}

Через точки $B$ и $F$ проведём прямую до пересечения с продолжением стороны $AD$ в некоторой точке $G$.
Тогда получим два треугольника:
$BCF$ и $DFG$, которые равны, так как у них:
$CF=FD$ (по условию), $\angle BFC\z=\angle DFG$ (как углы вертикальные) и $\angle BCF \z= \angle FDG$ (как углы накрест лежащие при параллельных прямых).
Из равенства треугольников следует:
$BF=FG$ и $BC=DG$.
То есть в треугольнике $ABG$ прямая $EF$ соединяет середины двух сторон, значит (§~\ref{1938/97}), $EF \parallel AG$ и $EF = \tfrac12(AD+DG)$, другими словами, $EF\parallel AD$ и $EF \z= \tfrac12(AD + BC)$.

\begin{figure}[h]
\begin{minipage}{.54\textwidth}
\centering
\includegraphics{mppics/ris-109}
\end{minipage}\hfill
\begin{minipage}{.44\textwidth}
\centering
\includegraphics{mppics/ris-110}
\end{minipage}

\medskip

\begin{minipage}{.54\textwidth}
\centering
\caption{}\label{1938/ris-109}
\end{minipage}\hfill
\begin{minipage}{.44\textwidth}
\centering
\caption{}\label{1938/ris-110}
\end{minipage}
\vskip-4mm
\end{figure} 

\paragraph{}\label{1938/100}
\mbox{\so{Задача}.}
\emph{Данный отрезок \emph{($AB$, рис.~\ref{1938/ris-110})} разделить на данное число равных частей} (например, на~3).

Из конца $A$ проводим прямую $AC$, образующую с $AB$ какой-нибудь угол;
откладываем на $AC$ от точки $A$ три произвольной величины и равных между собой отрезка:
$AD$, $DE$ и $EF$;
точку $F$ соединяем с $B$;
наконец, из $E$ и $D$ проводим прямые $EN$ и $DM$, параллельные $FB$.
Тогда отрезок $AB$, по доказанному, разделится в точках $M$ и $N$ на три равные части.

{\small

\subsection*{Задачи на построение}

\paragraph{Метод параллельного перенесения.}\label{1938/101}
На свойствах параллелограмма основан особый приём решения задач на построение, известный под названием метода параллельного перенесения.
Его сущность лучше всего выяснить на примере.

\begin{wrapfigure}{r}{45mm}
\vskip-7mm
\centering
\includegraphics{mppics/ris-111}
\caption{}\label{1938/ris-111}
\end{wrapfigure}

\smallskip
\mbox{\so{Задача}.}
Построить четырёхугольник $ABCD$ (рис.~\ref{1938/ris-111}), зная все его стороны и отрезок $EF$, соединяющий середины противоположных сторон.

Чтобы сблизить между собой данные линии, перенесём параллельно самим себе стороны $AD$ и $BC$ в положения $ED'$ и $EC'$.
Тогда сторона $DD'$ будет равна и параллельна $AE$, а сторона $CC'$ равна и параллельна $BE$, но так как $AE=BE$, то $DD'=CC'$ и $DD'\parallel CC'$.
Вследствие этого треугольники $DD'F$ и $CC'F$ будут равны (так как у них:
$DD' = CC'$, $DF=FC$ и $\angle D'DF=\angle FCC'$), значит, $\angle D'FD=\angle C'FC$ и потому линия $D'FC'$ должна быть прямая, то есть
фигура $ED'FC'$ окажется треугольником.
В этом треугольнике известны две стороны ($ED'=AD$ и $EC'=BC$) и медиана $EF$, проведённая к третьей стороне.
По этим данным легко построить треугольник $ED'C'$ (если на продолжении медианы $EF$ за точку $E$ отложим отрезок, равный $EF$, и полученную точку соединим с $D'$ и $C'$, то получим параллелограмм, у которого известны стороны и одна диагональ).

Найдя $\triangle ED'C'$, строим затем треугольники $D'DF$ и $C'CF$, а затем и весь четырёхугольник $ABCD$.

Предоставляем самим учащимся с помощью этого метода решить следующие задачи:

\medskip

1.
Построить трапецию по одному её углу, двум диагоналям и средней линии.

2.
Построить четырёхугольник по трём сторонам $a$, $b$, $c$ и двум углам $\alpha$ и $\beta$, прилежащим к неизвестной стороне.

3.
Построить трапецию по четырём данным её сторонам.

\paragraph{Метод симметрии.}\label{1938/102}
Свойства осевой симметрии также могут быть использованы при решении задач на построение.
Иногда искомый приём построения легко обнаруживается, если перегнём часть чертежа вокруг некоторой прямой так, чтобы эта часть заняла симметричное положение по другую сторону от этой прямой.
Приведём пример:

\smallskip
\so{Задача}.
\emph{На прямой $AB$ \emph{(рис.~\ref{1938/ris-112})} найти точку $x$, чтобы сумма расстояний от $x$ до данных точек $M$ и $N$ была наименьшая.}

\begin{wrapfigure}{o}{40mm}
\vskip-0mm
\centering
\includegraphics{mppics/ris-112}
\caption{}\label{1938/ris-112}
\end{wrapfigure}

Если, перегнув чертёж вокруг $AB$, приведём точку $M$ в симметричное относительно $AB$ положение $M'$, то расстояние от точки $M$ до какой угодно точки прямой $AB$ равно расстоянию от точки $M'$ до той же точки прямой $AB$.
Поэтому суммы $Mx\z+xN,  Mx_1\z+x_1N,\dots $ равны соответственно суммам $M'x+xN, M'x_1\z+x_1N, \dots$;
но из последних сумм наименьшая будет та, при которой линия $M'xN$ — прямая.
Отсюда становится ясным приём построения.

То же самое построение решает и другую задачу:
\emph{на прямой $AB$ найти такую точку $x$, чтобы прямые $xM$ и $xN$, проведённые от неё к данным точкам $M$ и $N$, составляли с $AB$ равные углы.}

Предоставляем учащимся решить методом симметрии следующие задачи:

\medskip

1.
Построить по четырём сторонам четырёхугольник $ABCD$, зная, что его диагональ $AC$ делит угол о пополам.

2.
На прямоугольном бильярде дано положение двух шаров $A$ и $B$.
В каком направлении надо толкнуть шар $A$, чтобы он, отразившись последовательно от всех четырёх бортов, ударил затем шар $B$.

3.
Дан угол и внутри него точка.
Построить треугольник наименьшего периметра, такой, чтобы одна его вершина лежала в данной точке, а две другие на сторонах угла.

}

{\small

\subsection*{Упражнения}

\begin{center}
\so{Доказать теоремы}
\end{center}

\begin{enumerate}[noitemsep]

\item
Соединив последовательно середины сторон какого-нибудь четырёхугольника, получим параллелограмм.

\item
В прямоугольном треугольнике медиана, проведённая к гипотенузе, равна её половине.

\smallskip
\so{Указание}.
Продлите медиану на расстояние, равное её длине.

\item
Обратно:
если медиана равна половине стороны, к которой она проведена, то треугольник прямоугольный.

\item
В прямоугольном треугольнике медиана и высота, приведённые к гипотенузе, образуют угол, равный разности острых углов треугольника.

\smallskip
\so{Указание}.
См.
задачу 2.

\item
В $\triangle ABC$ биссектриса угла $A$ встречает сторону $BC$ в точке $D$;
прямая, проведённая из $D$ параллельно $CA$, встречает $AB$ в точке $E$;
прямая, проведённая из $E$ параллельно $BC$, встречает $AC$ в $F$.
Доказать, что $EA=FC$.

\item
Внутри данного угла построен другой угол, стороны которого параллельны сторонам данного и равно отстоят от них.
Доказать, что биссектриса построенного угла лежит на биссектрисе данного угла.

\item
Всякая прямая, соединяющая какую-нибудь точку нижнего основания трапеции с какой-нибудь точкой верхнего основания, делится средней линией пополам.

\item
В треугольнике через точку пересечения биссектрис углов, прилежащих к основанию, проведена прямая параллельно основанию.
Доказать, что отрезок, заключённый между боковыми сторонами треугольника, равен сумме отрезков боковых сторон, считая их от основания.

\item
Через вершины углов треугольника проведены прямые, параллельные противоположным сторонам.
Доказать, что образованный ими треугольник составлен из четырёх треугольников, равных данному, и что каждая сторона его в два раза более соответствующей стороны данного треугольника.

\item
В равнобедренном треугольнике сумма расстояний от каждой точки основания до боковых сторон есть величина постоянная, а именно:
она равна высоте, опущенной на боковую сторону.

\item
Как изменится эта теорема, если взять точку на продолжении основания.

\item
В равностороннем треугольнике сумма расстояний всякой точки, взятой внутри этого треугольника, до сторон его есть величина постоянная, равная высоте треугольника.

\item
Всякий параллелограмм, у которого диагонали равны, есть прямоугольник.

\item
Всякий параллелограмм, у которого диагонали взаимно перпендикулярны, есть ромб.

\item
Всякий параллелограмм, у которого диагональ делит угол пополам, есть ромб.

\item
Из точки пересечения диагоналей ромба опущены перпендикуляры на стороны ромба.
Доказать, что основания этих перпендикуляров образуют вершины прямоугольника.

\smallskip
\so{Указание}.
См.
задачу 13.

\item
Биссектрисы углов прямоугольника своим пересечением образуют квадрат.

\item
Пусть $A'$, $B'$, $C'$ и $B'$ будут середины сторон $CD$, $DA$, $AB$ и $BC$ квадрата.
Доказать, что отрезки $AA'$, $CC'$, $BB'$ и $DD'$ образуют своим пересечением квадрат, сторона которого равна $\tfrac25$ каждого из этих отрезков.

\item
Дан квадрат $ABCD$.
На сторонах его отложены равные части $AA_1$, $BB_1$, $CC_1$ и $DD_1$.
Точки $A_1$, $B_1$, $C_1$, $D_1$ соединены последовательно прямыми.
Доказать, что $A_1B_1C_1D_1$ есть квадрат.

\item
Если середины сторон какого угодно четырёхугольника взять за вершины нового четырёхугольника, то последний есть параллелограмм.
Определить, при каких условиях этот параллелограмм будет:
1) прямоугольником, 2) ромбом, 3) квадратом.

\end{enumerate}

\begin{center}
\so{Найти геометрические места}
\end{center}

\begin{enumerate}[noitemsep]

\item
Середин всех отрезков, проведённых из данной точки к различным точкам данной прямой.

\item
Точек, равно отстоящих от двух параллельных прямых.

\item
Вершины треугольников, имеющих общее основание и равные высоты.
\end{enumerate}
\begin{center}
\so{Задачи на построение}
\end{center}

\begin{enumerate}[resume,noitemsep]

\item
Даны два угла треугольника;
построить третий.

\item
Дан острый угол прямоугольного треугольника;
построить другой острый угол.

\item
Провести прямую, параллельную данной прямой и находящуюся от неё на данном расстоянии.

\item
Разделить пополам угол, вершина которого не помещается на чертеже.

\item
Через данную точку провести прямую под данным углом к данной прямой.

\item
Даны две прямые $XY$ и $X'Y'$ и точка $F$;
провести через эту точку такую секущую, чтобы часть её, заключённая между данными прямыми, делилась точкой $F$ пополам.

\item
Через данную точку провести прямую так, чтобы отрезок её, заключённый между двумя данными параллельными прямыми, равнялся данному отрезку.

\item
Между сторонами данного острого угла поместить отрезок данной длины так, чтобы он был перпендикулярен к одной стороне угла.

\item
Между сторонами данного угла поместить отрезок данной длины параллельно заданной прямой, пересекающей обе стороны данного угла.

\item
Между сторонами данного угла поместить отрезок данной длины так, чтобы он отсекал от сторон угла равные отрезки.

\item
Построить прямоугольный треугольник по данным:
острому углу и противолежащему катету.

\item
Построить треугольник по двум углам и стороне, лежащей против одного из них.

\item
Построить \so{равнобедренный} треугольник по углу при вершине и основанию.

\item
То же — по углу при основании и высоте, опущенной на боковую сторону.

\item
То же — по боковой стороне и высоте, опущенной на неё.

\item
Построить равносторонний треугольник по его высоте.

\item
Разделить прямой угол на три равные части (другими словами, построить угол, равный  $30\degree$).

\item
Построить треугольник по основанию, высоте и боковой стороне.

\item
То же — по основанию, высоте и углу при основании.

\item
То же — по углу и двум высотам, опущенным из стороны этого угла.

\item
То же — по стороне, сумме двух других сторон и высоте, опущенной на одну из этих сторон.

\item
То же — по высоте, периметру и углу при основании.

\item
Провести в треугольнике прямую, параллельную основанию, так, чтобы отрезок, заключённый между боковыми сторонами, был равен сумме отрезков боковых сторон, считая от основания.

\item
Построить многоугольник, равный данному.

\smallskip
\so{Указание}.
Диагоналями разбивают данный многоугольник на треугольники.

\item
Построить \so{четырёхугольник} по трём его углам и двум сторонам, образующим четвёртый угол.

\smallskip
\so{Указание}.
Надо найти четвёртый угол.

\item
То же — по трём сторонам и двум диагоналям.

\item
Построить \so{параллелограмм} по двум неравным сторонам и одной диагонали.

\item
То же — по стороне и двум диагоналям.

\item
То же — по двум диагоналям и углу между ними.

\item
То же — по основанию, высоте и диагонали.

\item
Построить \so{прямоугольник} по диагонали и углу между диагоналями.

\item
Построить \so{ромб} по стороне и диагонали.

\item
То же — по двум диагоналям.

\item
То же — по высоте и диагонали.

\item
То же — по углу и диагонали, проходящей через этот угол.

\item
То же — по диагонали и противолежащему углу.

\item
То же — по сумме диагоналей и углу, образованному диагональю со стороной.

\item
Построить \so{квадрат} по данной диагонали.

\item
Построить \so{трапецию} по основанию, прилежащему к нему углу и двум непараллельным сторонам (могут быть два решения, одно и ни одного).

\item
То же — по разности оснований, двум боковым сторонам и одной диагонали.

\item
То же — по четырём сторонам (всегда ли задача имеет решение?).

\item
То же — по основанию, высоте и двум диагоналям (условие существования решения).

\item
То же — по двум основаниям и двум диагоналям (условие существования решения).

\item
Построить \so{квадрат} по сумме стороны с диагональю.

\item
То же — по разности диагонали и стороны.

\item
Построить \so{параллелограмм} по двум диагоналям и высоте.

\item
То же — по стороне, сумме диагоналей и углу между ними.

\item
Построить треугольник по двум сторонам и медиане, проведённой к третьей стороне.

\item
То же — по основанию, высоте и медиане, проведённой к боковой стороне.

\item
Построить \so{прямоугольный треугольник} по гипотенузе и сумме катетов (исследовать).

\item
То же — по гипотенузе и разности катетов.

\item
Даны две точки $A$ и $B$, расположенные по одну сторону от данной прямой $XY$.
Расположить на этой прямой отрезок $MN$ данной длины $\ell$ так, чтобы ломаная $AM+MN+NB$ была наименьшей длины.

\smallskip
\so{Указание}.
Приблизим точку $B$ к точке $A$, двигая её по прямой, параллельной $XY$, на расстояние, равное $MN$.

\end{enumerate}

}

%% file: 2D/postulat5.tex
{\small

\section{Об аксиоме параллельных}

\paragraph{Две аксиомы параллельных.}\label{1914/91} 
Принятая нами в §~\ref{1938/75} аксиома параллельных была предложена древнегреческим математиком Проклом (IV век нашей эры) но часто называется \rindex{аксиома!Плейфэра}\textbf{аксиомой Плейфэра} в честь шотландского математика Джона Плейфэра благодаря которому многие узнали эту формулировку.

Легко показать, что пятый постулат Евклида (§~\ref{1938/78}) и аксиома Плейфэра обратимы одна в другую. 
То есть из аксиомы Плейфэра можно вывести, как логическое следствие, аксиому Евклида (что и сделано в §~\ref{1938/78}) и, обратно, из этого постулата можно логически получить постулат Плейфэра.
Последнее можно выполнить, например, так:

\begin{wrapfigure}[10]{r}{34mm}
\vskip-7mm
\centering
\includegraphics{mppics/ris-1914-88}
\caption{}\label{1914/ris-88}
\end{wrapfigure}

Пусть через точку $E$ (рис. \ref{1914/ris-88}), взятую вне прямой $CD$, проведены какие-нибудь две прямые $AB$ и $A_1B_1$.
Докажем, исходя из постулата Евклида, что эти прямые не могут быть обе параллельны прямой $CD$.
Для этого проведём через $E$ какую-нибудь секущую прямую $MN$;
обозначим внутренние односторонние углы, образуемые этой секущей
с прямыми $CD$ и $AB$, буквами $\alpha$ и $\beta$.
Тогда одно из двух: или сумма $\alpha+\beta$ не равна $180\degree$,
или она равна $180\degree$.
В первом случае согласно постулату Евклида,
прямая $AB$ должна пересечься с $CD$ и, следовательно, она не может быть параллельной $CD$. 
Во втором случае $\alpha+\angle B_1EN\ne 180\degree$, (так как $\angle B_1EN\ne \angle BEN$).
Значит, тогда, согласно тому же постулату, прямая $A_1B_1$ должна пересечься
с $CD$ и, следовательно, эта прямая не может быть параллельной
$CD$.
Таким образом, одна из прямых $AB$ или $A_1B_1$ непременно окажется непараллельной прямой $CD$; следовательно, через одну точку нельзя провести двух различных прямых, параллельных одной и той же прямой.

\paragraph{Другие предложения равносильные аксиоме параллельных.}\label{1914/92}
Есть много других предложений, также логически обратимых с постулатом Евклида (и, следовательно, ему логически равносильных).
Приведём несколько знаменитых примеров:

\emph{Существует по крайней мере один треугольник, у которого сумма углов равна $180\degree$} (французский математик Адриен Мари Лежандр, начало XIX столетия).

\emph{Существует выпуклый четырёхугольник} (прямоугольник), \emph{у которого все четыре угла прямые} (французский математик Клод  Клеро, XVIII столетие).

\emph{Существует треугольник, подобный, но не равный, другому треугольнику} (итальянский математик Джироламо Саккери, начало XVIII столетия).

\emph{Через всякую точку, взятую внутри угла, меньшего $180\degree$, можно
провести прямую, пересекающую обе стороны этого угла} (немецкий
математик Иоганн Фридрих Лоренц, конец XVIII столетия).

Многие делали попытки доказать постулат Евклида
(или какой-нибудь другой, ему равносильный), то есть вывести его,
как логическое следствие, из других аксиом геометрий.
Все эти попытки оказались неудачными: в каждом из таких «доказательств», после подробного разбора его, находили логическую ошибку.

\paragraph{Открытие новой геометрии.}\label{1914/93} Неудачи в поисках доказательств
постулата Евклида привели к мысли, что этот постулат (как и любой ему равносильный) и не может быть выведен из других аксиом геометрии, а представляет собою независимое от них самостоятельное допущение о свойствах пространства.
Сначала эта мысль высказывалась только в частной переписке;
до нас дошло такое письмо 1816 года, написанное Карлом Гауссом;
ещё более уверенное письмо 1818 года написано Фердинандом Швейкартом.
Однако эти математики воздерживались от публикации своих результатов, сейчас трудно понять тому истинную причину.

Независимо, те же идеи были развиты Николаем Ивановичем Лобачевским в 
сочинении изданным Казанским университетом в 1836—1838 годах.
Чуть позже, также независимо, те же результаты были опубликованы венгерским математиком Яношем Бояи.

\begin{wrapfigure}[8]{r}{37mm}
\vskip-5mm
\centering
\includegraphics{mppics/ris-1914-89}
\caption{}\label{1914/ris-89}
\end{wrapfigure}

В своём сочинении Лобачевский обнародовал особую геометрию, названную им «воображаемой», а теперь обычно называемой «геометрией Лобачевского».
В её основание положены те же геометрические аксиомы, на которых основана
геометрия Евклида, за исключением только постулата параллельных линий, вместо которого Лобачевский взял следующее допущение:
\emph{через точку, лежащую вне прямой, можно провести бесчисленное множество параллельных этой прямой}.

То есть, он допустил, что если $AB$ (рис. \ref{1914/ris-89}) есть прямая и $C$ какая-нибудь точка вне её, то при этой точке существует некоторый угол $DEC$, обладающий следующим свойствам:
1) всякая прямая, проведённая через $C$ внутри этого угла (например, прямая $CF$) не пересекается с $AB$,
2) то же верно и для продолжений сторон $DE$ и $EC$ угла,
а при этом 3) всякая прямая, проведённая через $C$ вне этого угла, пересекается с $AB$.

Понятно, что такое допущение отрицает аксиому параллельности Евклида.
Несмотря однако на это отрицание, геометрия Лобачевского представляет собою такую же стройную
систему геометрических теорем как и геометрия Евклида.
Конечно, теоремы геометрии Лобачевского существенно отличаются от теорем геометрии Евклида, но
в ней, как и в геометрии Евклида, не встречается никаких логических противоречий ни теорем с аксиомами, положенными в основание этой геометрии, ни одних теорем с другими теоремами.

Между тем, если бы постулат Евклида мог быть доказан, то есть если бы он представлял собою
некоторое, хотя бы и очень отдалённое, логическое следствие из других геометрических аксиом, то тогда отрицание этого постулата, положенное в основу геометрии вместе с принятием всех других аксиом, непременно привело бы к логически противоречивым следствиям.

Отсутствие таких противоречий в геометрии Лобачевского служит указанием на независимость пятого
постулата Евклида от прочих геометрических аксиом и, следовательно, на невозможность доказать его.
Заметим, однако, что одно только отсутствие противоречий в геометрии Лобачевского ещё не служит доказательством независимости пятого постулата от других аксиом геометрии.
Ведь всегда можно возразить, что это отсутствие противоречий есть только
случайное явление, происходящее, быть может, от того, что в геометрии Лобачевского не сделано ещё достаточного количества выводов, что со временем, быть может, и удастся кому-нибудь получить
такой логический вывод в этой геометрии, который окажется в противоречии с каким-нибудь другим выводом той же геометрии.

Однако доказано следующее: \emph{если бы в геометрии Лобачевского нашлось противоречие, то нашлось бы соответствующее противоречие и в Евклидовой геометрии};
также верно и обратное — \emph{если есть противоречие в Евклидовой геометрии то  было бы  противоречие и в геометрии Лобачевского}.
То есть удаётся доказать, что в логическом смысле геометрия Лобачевского «не хуже и не лучше» геометрии Евклида.

\paragraph{Неевклидовы геометрии.}\label{1914/94} 
Позже немецкий математик Бернхард Риман (1826—1866) построил ещё особую, также лишённую противоречий, геометрию (названную потом \textbf{геометрией Римана}), 
в которой вместо постулата Евклида принимается допущение, что
через точку, взятую вне прямой, нельзя провести ни одной параллельной этой прямой (другими словами, все прямые плоскости пересекаются).
Такие геометрии (как геометрии Лобачевского и Римана), в которых в основание положено какое-нибудь допущение о параллельных линиях, не согласное с постулатом Евклида, носят общее название
\textbf{неевклидовых геометрий}.

Другой пример неевклидовой геометрии — это так называемая \rindex{абсолютная геометрия}\textbf{абсолютная геометрия}, независимая от пятого постулата.
Другими словами, в этой геометрии пятый постулат может выполняться, а может и не выполняться.
Эта геометрия как раз и рассматривалась в упомянутом сочинении Яноша Бояи,
она включает геометрии Евклида и Лобачевского как частные случаи.

\paragraph{Теоремы геометрии Лобачевского.}\label{1914/95} Приведём некоторые теоремы геометрии Лобачевского, резко различающиеся от соответствующих теорем геометрии Евклида:

\emph{Два перпендикуляра к одной и той же прямой, по мере удаления
от этой прямой, расходятся неограниченно.}

\emph{Сумма углов треугольника меньше $180\degree$} (в геометрии Римана она
больше $180\degree$), причём эта сумма не есть величина постоянная для разных треугольников.
\emph{Чем больше площадь треугольника, тем больше сумма его углов разнится от $180\degree$.}

\begin{wrapfigure}[6]{r}{25mm}
\vskip-5mm
\centering
\includegraphics{mppics/ris-extra-5}
\caption{}\label{extra/ris-5}
\end{wrapfigure}

\emph{Если в выпуклом четырёхугольнике три угла прямые, то четвёртый угол острый.}
В частности в этой геометрии нет прямоугольников.
Четырёхугольник с тремя прямыми углами называется \rindex{четырёхугольник Ламберта}\textbf{четырёхугольником Ламберта} в честь шведского математика Иоганна Генриха Ламберта (1728—1777).
Попытка изобразить четырёхугольник Ламберта (рис. \ref{extra/ris-5}) не очень удачная, поскольку сам лист бумаги не похож на плоскость Лобачевского. 

\emph{Если углы одного треугольника соответственно равны углам другого треугольника, то такие треугольники равны} (следовательно, в геометрии Лобачевского не существует подобия).

\emph{Геометрическое место точек плоскости, равноотстоящих от какой-нибудь прямой этой плоскости, не является прямой.}

}

%% file: 2D/okruzhnost.tex
\chapter{Окружность}

\section{Форма и положение окружности}

\paragraph{}\label{1938/103}
\so{Предварительное замечание}.
Очевидно, что через одну точку ($A$, рис.~\ref{1938/ris-113}) можно провести сколько угодно окружностей:
центры их можно брать произвольно.
Через две точки ($A$ и $B$, рис.~\ref{1938/ris-114}) тоже можно провести сколько угодно окружностей, но центры их нельзя брать произвольно, так как точки, одинаково удалённые от двух точек $A$ и $B$, должны лежать на срединном перпендикуляре к отрезку $AB$ (§~\ref{1938/58}). 

\begin{figure}[!ht]
\begin{minipage}{.48\textwidth}
\centering
\includegraphics{mppics/ris-113}
\end{minipage}
\hfill
\begin{minipage}{.48\textwidth}
\centering
\includegraphics{mppics/ris-114}
\end{minipage}

\medskip

\begin{minipage}{.48\textwidth}
\centering
\caption{}\label{1938/ris-113}
\end{minipage}
\hfill
\begin{minipage}{.48\textwidth}
\centering
\caption{}\label{1938/ris-114}
\end{minipage}
\vskip-4mm
\end{figure}

Посмотрим, можно ли провести окружность через три точки.

\paragraph{}\label{1938/104}
\so{Теорема}.
\textbf{\emph{Через три точки, не лежащие на одной прямой, можно провести окружность и притом только одну.}}

Через три точки $A$, $B$ и $C$ (рис.~\ref{1938/ris-115}), только тогда можно провести окружность, если существует такая четвёртая точка $O$, которая одинаково удалена от точек $A$, $B$ и $C$.

\begin{wrapfigure}{o}{45mm}
\centering
\includegraphics{mppics/ris-115}
\caption{}\label{1938/ris-115}
\end{wrapfigure}

Докажем, что если $A$, $B$ и $C$ не лежат на одной прямой 
(другими словами, если точки $A$, $B$ и $C$ являются вершинами треугольника),
то такая точка $O$ существует и притом только одна.
Для этого примем во внимание, что всякая точка, одинаково удалённая от точек $A$ и $B$, должна лежать на срединном перпендикуляре $MN$, проведённом к стороне $AB$ (§~\ref{1938/58}); 
точно так же всякая точка, одинаково удалённая от точек $B$ и $C$, должна лежать на срединном перпендикуляре $PQ$, проведённом к стороне $BC$.
Значит, если существует точка, одинаково удалённая от трёх точек $A$, $B$ и $C$, то она должна лежать одновременно и на $MN$, и на $PQ$, что возможно только тогда, когда она совпадает с точкой пересечения этих двух прямых.
Прямые $MN$ и $PQ$ всегда пересекаются, так как они перпендикулярны к пересекающимся прямым $AB$ и $BC$ (§~\ref{1938/78}).
Точка $O$ их пересечения и будет точкой, одинаково удалённой от $A$, от $B$ и от $C$;
значит, если примем эту точку за центр, а за радиус возьмём отрезок $OA$ (или $OB$, или $OC$), то окружность пройдёт через точки $A$, $B$ и $C$.
Так как прямые $MN$ и $PQ$ могут пересечься только в одной точке, то центр такой окружности может быть только один, и длина её радиуса может быть только одна;
значит, искомая окружность — единственная.

{\small

\smallskip
\so{Замечание}.
Если бы три точки $A$, $B$ и $C$ (рис.~\ref{1938/ris-115}) лежали на одной прямой то перпендикуляры $MN$ и $PQ$, были бы параллельны и, значит, не могли бы пересечься.
Следовательно, через три точки, лежащие на одной прямой, нельзя провести окружности.

}

\smallskip
\so{Следствие}.
Точка $O$ (рис.~\ref{1938/ris-115}), находясь на одинаковом расстоянии от $A$ и от $C$, должна также лежать на срединном перпендикуляре $RS$, проведённом к стороне $AC$. 
Таким образом:
\emph{три срединных перпендикуляра к сторонам треугольника пересекаются в одной точке.}

\paragraph{}\label{1938/105}
\mbox{\so{Теорема}.}
\textbf{\emph{Диаметр}} ($AB$, рис.~\ref{1938/ris-116}), \textbf{\emph{перпендикулярный к хорде}} ($CD$), \textbf{\emph{делит эту хорду и обе стягиваемые ею дуги пополам.}}
Перегнём чертёж по диаметру $AB$ так, чтобы его левая часть упала на правую.
Тогда левая полуокружность совместится с правой полуокружностью, и перпендикуляр $KC$ пойдёт по $KD$.
Из этого следует, что точка $C$, представляющая собой пересечение полуокружности с $KC$, совпадёт с $D$;
поэтому $CK=KD$,
${\smallsmile} BC={\smallsmile} BD$ и
${\smallsmile} AC={\smallsmile} AD$.

\paragraph{}\label{1938/106}
\mbox{\so{Обратные теоремы}.}

1.
\textbf{\emph{Диаметр}} ($AB$), \textbf{\emph{проведённый через середину хорды}} ($CD$)\textbf{\emph{, не проходящей через центр, перпендикулярен к этой хорде и делит дугу, стягиваемую ею, пополам}} (рис.~\ref{1938/ris-116}).

\begin{wrapfigure}{r}{33mm}
\centering
\includegraphics{mppics/ris-116}
\caption{}\label{1938/ris-116}
\end{wrapfigure}

2.
\textbf{\emph{Диаметр}} ($AB$), \textbf{\emph{проведённый через середину дуги}} ($CBD$), \textbf{\emph{перпендикулярен к хорде, стягивающей эту дугу, и делит её пополам.}}

Оба эти предложения легко доказываются от противного.

\paragraph{}\label{1938/107}
\mbox{\so{Теорема}.}
\textbf{\emph{Дуги}} ($AC$ и $BD$, рис.~\ref{1938/ris-117}), \textbf{\emph{заключённые между параллельными хордами}} ($AB$ и $CD$), \textbf{\emph{равны.}}

Перегнём чертёж по диаметру $EF\z\perp AB$.
Тогда на основании предыдущей теоремы можно утверждать, что точка $A$ совпадёт с $B$, точка $C$ совпадёт с $D$ и, следовательно, дуга $AC$ совместится с дугой $BD$, то есть эти дуги равны.

\begin{figure}[h]
\begin{minipage}{.32\textwidth}
\centering
\includegraphics{mppics/ris-117}
\end{minipage}\hfill
\begin{minipage}{.32\textwidth}
\centering
\includegraphics{mppics/ris-118}
\end{minipage}\hfill
\begin{minipage}{.32\textwidth}
\centering
\includegraphics{mppics/ris-119}
\end{minipage}

\medskip

\begin{minipage}{.32\textwidth}
\centering
\caption{}\label{1938/ris-117}
\end{minipage}\hfill
\begin{minipage}{.32\textwidth}
\centering
\caption{}\label{1938/ris-118}
\end{minipage}\hfill
\begin{minipage}{.32\textwidth}
\centering
\caption{}\label{1938/ris-119}
\end{minipage}
\vskip-4mm
\end{figure} 

\paragraph{}\label{1938/108}
\mbox{\so{Задачи}.}
1) \emph{Разделить данную дугу \emph{($AB$, рис.~\ref{1938/ris-118})} пополам.}

Соединив концы дуги хордой $AB$, опускаем на неё перпендикуляр из центра и продолжаем его до пересечения с дугой.
По доказанному в предыдущей теореме, дуга $AB$ разделится этим перпендикуляром пополам.

Если же центр не известен, тогда к хорде $AB$ следует провести срединный перпендикуляр. 

2) \emph{Найти центр данной окружности} (рис.~\ref{1938/ris-119}).

Взяв на данной окружности какие-нибудь три точки $A$, $B$ и $C$, проводят через них две хорды, например $AB$ и $CB$, и проводят к ним срединные перпендикуляры $MN$ и $PQ$.

Искомый центр, будучи одинаково удалён от $A$, $B$ и $C$, должен лежать и на $MN$ и на $PQ$, следовательно, он находится в их пересечении, то есть в точке $O$.

%% file: 2D/dugi.tex
\section{Дуги, хорды и расстояния до центра}

\begin{wrapfigure}{r}{40mm}
\centering
\includegraphics{mppics/ris-120}
\caption{}\label{1938/ris-120}
\end{wrapfigure}

\paragraph{}\label{1938/109}
\mbox{\so{Теоремы}.}
\textbf{\emph{В одном круге или в равных кругах:}}

1) \textbf{\emph{если дуги равны, то стягивающие их хорды равны и одинаково удалены от центра.}}

2) \textbf{\emph{если две дуги, меньшие полуокружности, не равны, то б\'{о}льшая из них стягивается б\'{о}льшей хордой и из обеих хорд б\'{о}льшая расположена ближе к центру.}}

1) Пусть дуга $AB$ равна дуге $CD$ (рис.~\ref{1938/ris-120}), требуется доказать, что хорды $AB$ и $CD$ равны, а также равны перпендикуляры $OE$ и $OF$, опущенные из центра на хорды.

Повернём сектор $AOB$ вокруг центра $O$ в направлении, указанном стрелкой, на столько, чтобы радиус $OB$ совпал с $OC$.
Тогда дуга $BA$ пойдёт по дуге $CD$ и вследствие их равенства эти дуги совместятся.
Значит, хорда $AB$ совместится с хордой $CD$ и перпендикуляр $OE$ совпадёт с $OF$ (из одной точки можно опустить на прямую только один перпендикуляр), то есть
$AB=CD$ и $OE=OF$.

\begin{wrapfigure}{o}{40mm}
\centering
\includegraphics{mppics/ris-121}
\caption{}\label{1938/ris-121}
\end{wrapfigure}

2) Пусть дуга $AB$ (рис.~\ref{1938/ris-121}) меньше дуги $CD$, и притом обе дуги меньше полуокружности;
требуется доказать, что хорда $AB$ меньше хорды $CD$, а перпендикуляр $OE$ больше перпендикуляра $OF$.
Отложим на дуге $CD$ дугу $CK$, равную $AB$, и проведём вспомогательную хорду $CK$, которая, по доказанному, равна хорде $AB$ и одинаково с ней удалена от центра.
У треугольников $COD$ и $COK$ две стороны одного равны двум сторонам другого (как радиусы), а углы, заключённые между этими сторонами, не равны;
в этом случае, как мы знаем (§~\ref{1938/52}), против большего из углов, то есть
$\angle COD$, должна лежать б\'{о}льшая сторона;
значит, $CD>CK$, и потому $CD>AB$.

Для доказательства того, что $OE>OF$, проведём $OL\perp CK$ и примем во внимание, что, по доказанному, $OE=OL$;
следовательно, нам достаточно сравнить $OF$ с $OL$.
В прямоугольном треугольнике $OFM$ (покрытом на чертеже штрихами) гипотенуза $OM$ больше катета $OF$;
но $OL>OM$;
значит, и подавно $OL>OF$, и потому $OE>OF$.

Теорема, доказанная нами для одного круга, остаётся верной и для равных кругов, потому что такие круги один от другого отличаются только положением.

\paragraph{}\label{1938/110}
\so{Обратные теоремы}.
Так как в предыдущем параграфе рассмотрены всевозможные взаимно исключающие случаи относительно сравнительной величины двух дуг одного радиуса, причём получились взаимно исключающие выводы относительно сравнительной величины хорд и расстояний их от центра, то обратные предложения должны быть верны, а именно.

\textbf{\emph{В одном круге или в равных кругах:}}

1) \textbf{\emph{равные хорды одинаково удалены от центра и стягивают равные дуги;}}

2) \textbf{\emph{хорды, одинаково удалённые от центра, равны и стягивают равные дуги;}}

3) \textbf{\emph{из двух неравных хорд б\'{о}льшая ближе к центру и стягивает б\'{о}льшую дугу;}}

4) \textbf{\emph{из двух хорд, неодинаково удалённых от центра, та, которая ближе к центру, больше и стягивает б\'{о}льшую дугу.}}

Эти предложения легко доказываются от противного.
Например, для доказательства первого из них рассуждаем так:
если бы данные хорды стягивали неравные дуги, то, согласно прямой теореме, они были бы не равны, что противоречит условию;
значит, равные хорды должны стягивать равные дуги;
а если дуги равны, то, согласно прямой теореме, стягивающие их хорды одинаково удалены от центра.

\paragraph{}\label{1938/111}
\mbox{\so{Теорема}.}
\textbf{\emph{Диаметр есть наибольшая из хорд.}}

\begin{wrapfigure}{o}{40mm}
\centering
\includegraphics{mppics/ris-122}
\caption{}\label{1938/ris-122}
\end{wrapfigure}

Если соединим с центром $O$ концы какой-нибудь хорды, не проходящей через центр, например хорды $AB$ (рис.~\ref{1938/ris-122}), то получим треугольник $AOB$, в котором одна сторона есть эта хорда, а две другие — радиусы.
Но в треугольнике каждая сторона менее суммы двух других сторон;
следовательно, хорда $AB$ менее суммы двух радиусов, тогда как всякий диаметр $CD$ равен сумме двух радиусов.
Значит, диаметр больше всякой хорды, не проходящей через центр.
Но так как диаметр есть тоже хорда, то можно сказать, что диаметр есть наибольшая из хорд.

%% file: 2D/raspolozheniya.tex
\section{Расположения прямой и окружности}

\paragraph{Точки внутри круга и точки вне его.}\label{1914/118}
Окружность разделяет все точки плоскости на три области:

\begin{wrapfigure}{r}{33mm}
\vskip-6mm
\centering
\includegraphics{mppics/ris-1914-113}
\caption{}\label{1914/ris-113}
\end{wrapfigure}

1) точки \rindex{внешние точки круга}\textbf{вне круга}, расстояние от которых до центра больше радиуса; 

2) точки на окружности, расстояние от которых до центра равно радиусу;

3) точки \rindex{внутренние точки круга}\textbf{внутри круга}, расстояние от которых до центра меньше радиуса. 

Следующее утверждение мы принимаем за очевидное:

\emph{Если отрезок (или дуга) соединяет точку \emph{($A$, рис.~\ref{1914/ris-113})} внутри круга с какою-нибудь точкой \emph{($B$)} вне его, то он пересекается где-нибудь с окружностью \emph{(в точке $X$ для отрезка $AB$ и в точках $Y$ и $Z$ для двух дуг с концами в $A$ и $B$)}.}


\paragraph{}\label{1938/112}
Прямая и окружность могут, очевидно, находиться только в следующих трёх относительных положениях.

1) \emph{Расстояние \emph{($OC$)} от центра до прямой \emph{($AB$) (то есть длина перпендикуляра $OC$, опущенного из центра на прямую)} больше радиуса окружности} (рис.~\ref{1938/ris-123}).
Тогда точка $C$ прямой удалена от центра больше, чем на радиус, и потому лежит вне круга.
Так как все остальные точки прямой удалены от $O$ ещё более, чем точка $C$ (наклонные длиннее перпендикуляра), то они все лежат вне круга;
значит, тогда прямая не имеет никаких точек, общих с окружностью.

2) \emph{Расстояние \emph{($OC$)} от центра до прямой меньше радиуса.}
В этом случае (рис.~\ref{1938/ris-124}) точка $C$ лежит внутри круга. 
Но на прямой $AB$, по обе стороны от точки $C$, можно найти такие точки $D$ и $E$, которые удалены от $O$ более, чем на радиус%
\footnote{Если, например, на прямой $AB$ отложим от точки $C$ по обе стороны от неё, отрезки, большие радиуса, то расстояния от их концов до центра будут больше радиуса, так как перпендикуляр короче наклонной.}
и которые, следовательно, лежат вне круга.
Но тогда каждый из двух отрезков: $CD$ и $CE$, соединяя внутреннюю точку с внешней, должен пересечься с окружностью (§~\ref{1914/118}).
Следовательно, в этом случае прямая имеет с окружностью 2 общие точки.

\begin{figure}
\begin{minipage}{.32\textwidth}
\centering
\includegraphics{mppics/ris-123}
\end{minipage}\hfill
\begin{minipage}{.32\textwidth}
\vskip6mm
\centering
\includegraphics{mppics/ris-124}
\end{minipage}\hfill
\begin{minipage}{.32\textwidth}
\vskip4mm
\centering
\includegraphics{mppics/ris-125}
\end{minipage}

\medskip

\begin{minipage}{.32\textwidth}
\centering
\caption{}\label{1938/ris-123}
\end{minipage}
\begin{minipage}{.32\textwidth}
\vfill
\centering
\caption{}\label{1938/ris-124}
\end{minipage}
\begin{minipage}{.32\textwidth}
\vfill
\centering
\caption{}\label{1938/ris-125}
\end{minipage}
\vskip-4mm
\end{figure} 

3) \emph{Расстояние \emph{($OC$)} от центра до прямой равно радиусу.}
Тогда точка $C$ (рис.~\ref{1938/ris-125}) принадлежит и прямой и окружности, все же остальные точки прямой, будучи удалены от $O$ более, чем точка $C$, лежат вне круга.
Значит, в этом случае прямая и окружность имеют только одну общую точку, именно ту, которая служит основанием перпендикуляра, опущенного из центра на прямую.

Такая прямая, которая с окружностью имеет только одну общую точку, называется \rindex{касательная}\textbf{касательной} к окружности;
общая точка называется \rindex{точка!касания}\textbf{точкой касания}.

\begin{wrapfigure}{o}{45mm}
\vskip-4mm
\centering
\includegraphics{mppics/ris-1914-134}
\caption{}\label{1914/ris-134}
\end{wrapfigure}

\medskip

{\small

\smallskip
\mbox{\so{Замечание}.}
Когда речь идёт о произвольной гладкой кривой, то за определение касательной принимают \emph{предельное положение} $MT$, к которому стремится секущая $MP$, когда точка $P$ скользит по кривой неограниченно приближаясь к $M$.
Для окружности оба определения равносильны, но в общем случае, определяемая таким образом касательная может иметь с кривой более одной общей точки (рис.~\ref{1914/ris-134}).

} 

\paragraph{}\label{1938/113}
Относительно касательной мы докажем следующие \so{две теоремы} (прямую и обратную):

1) \textbf{\emph{если прямая}} ($MN$ рис.~\ref{1938/ris-126}) \textbf{\emph{перпендикулярна к радиусу}} ($OA$) \textbf{\emph{в конце его}} ($A$), \textbf{\emph{лежащем на окружности, то она касается окружности,}} и обратно,

2) \textbf{\emph{если прямая касается окружности, то радиус, проведённый в точку касания, перпендикулярен к ней.}}

\begin{wrapfigure}{o}{33mm}
\centering
\includegraphics{mppics/ris-126}
\caption{}\label{1938/ris-126}
\end{wrapfigure}

1) Точка $A$, как конец радиуса, лежащий на окружности, принадлежит этой окружности;
в то же время она принадлежит и прямой $MN$.
Значит, эта точка есть общая у окружности и прямой.
Все же остальные точки прямой $MN$, как $B$, $C$ и другие, отстоят от центра $O$ дальше, чем на радиус (так как отрезки $OB$, $OC,\dots$, как наклонные, больше перпендикуляра $OA$), и потому они лежат вне окружности.
Таким образом, у прямой $MN$ есть только одна точка ($A$), общая с окружностью, и, значит, прямая $MN$ есть касательная.

2) Пусть $MN$ касается окружности в точке $A$; требуется доказать, что $MN\perp OA$.
Предположим противное, то есть что радиус $OA$ не перпендикулярен к $MN$, а представляет собою наклонную к этой прямой.
В таком случае какая-нибудь другая прямая, например $OC$, будет перпендикуляром, опущенным из центра $O$ на касательную $MN$ (§~\ref{1938/24}).
Значит на $MN$ найдётся другая точка $A'$, лежащая от $C$ на том же расстоянии, что и $A$.
В этом случае, $OA=OA'$ (§~\ref{1938/54});
то есть $MN$ имеет две общие точки с окружностью и значит не может её касаться.

\begin{figure}
\begin{minipage}{.32\textwidth}
\centering
\includegraphics{mppics/ris-1931-149}
\end{minipage}\hfill
\begin{minipage}{.32\textwidth}
\hfill
\centering
\includegraphics{mppics/ris-127}
\end{minipage}\hfill
\begin{minipage}{.32\textwidth}
\hfill
\centering
\includegraphics{mppics/ris-128}
\end{minipage}

\medskip

\begin{minipage}{.32\textwidth}
\centering
\caption{}\label{1931/ris-149}
\end{minipage}
\begin{minipage}{.32\textwidth}
\vfill
\centering
\caption{}\label{1938/ris-127}
\end{minipage}
\begin{minipage}{.32\textwidth}
\vfill
\centering
\caption{}\label{1938/ris-128}
\end{minipage}
\vskip-4mm
\end{figure} 

\smallskip
\mbox{\so{Следствие}.}
\emph{Две касательные, проведённые к окружности из точки вне её, равны и образуют равные углы с прямой, соединяющей эту точку с центром,} что следует из равенства прямоугольных треугольников $AOB$ и $AOB_1$ (рис.~\ref{1931/ris-149}).

\paragraph{}\label{1938/114}
\mbox{\so{Теорема}.}
\textbf{\emph{Если касательная параллельна хорде, то точка касания делит дугу, стягиваемую хордой, пополам.}}

Пусть прямая $AB$ касается окружности в точке $M$ (рис.~\ref{1938/ris-127}) и параллельна хорде $CD$;
требуется доказать, что ${\smallsmile}CM={\smallsmile}MD$.

Проведя через точку касания диаметр $ME$, будем иметь:
$EM\perp AB$ и, следовательно, $EM\z\perp CD$;
поэтому ${\smallsmile}CM\z={\smallsmile}MD$.

\paragraph{}\label{1938/115}
\mbox{\so{Задача}.}
\emph{Провести касательную к данной окружности с центром $O$ параллельно данной прямой $AB$} (рис. \ref{1938/ris-128}).

Опускаем на $AB$ из центра $O$ перпендикуляр $OC$ и через точку $D$, в которой этот перпендикуляр пересекается с окружностью, проводим $EF\parallel AB$.
Искомая касательная будет $EF$.
Действительно, так как $OC\perp AB$ и $EF\parallel AB$, то $EF\perp OD$, а прямая, перпендикулярная к радиусу в конце его, лежащем на окружности, есть касательная.

\paragraph{}\label{1931/115}
\mbox{\so{Задача}.}
\emph{Через данную точку провести касательную к данной окружности.}
Если данная точка (например, точка $C$, рис. \ref{1938/ris-125}) находится на окружности, то проводят через неё радиус и через конец радиуса перпендикулярную прямую.
Эта прямая и будет искомой касательной.
Другой касательной через ту же точку окружности провести нельзя, так как касательная должна быть перпендикулярна к радиусу в конце его, лежащем на окружности, а двух различных перпендикуляров к одному и тому же радиусу через одну и ту же точку провести нельзя.

Рассмотрим случай, когда точка дана вне круга.
Пусть требуется (рис. \ref{1931/ris-124}) провести к окружности с центром $O$ касательную через точку $A$.
Для этого из точки $A$ как из центра описываем дугу радиуса $AO$, а из точки $O$ как из центра пересекаем эту дугу в точках $B$ и $C$ раствором циркуля, равным удвоенному радиусу данной окружности.
Проведя затем хорды $OB$ и $OC$, соединим точку $A$ с точками $D$ и $E$, в которых эти хорды пересекаются с
данной окружностью.
Прямые $AD$ и $AE$ и будут касательными к данной окружности.

{

\begin{wrapfigure}{o}{32mm}
\vskip-4mm
\centering
\includegraphics{mppics/ris-1931-124}
\caption{}\label{1931/ris-124}
\end{wrapfigure}

Действительно, из построения видно, что $\triangle AOB$ и $\triangle AOC$ равнобедренные ($AO=AB\z=AC$) с основаниями $OB$ и $OC$, равными удвоенному радиусу окружности.
Так как $OD$ и $OE$ — радиусы, то $D$ есть середина $OB$, а $E$ — середина $OC$.
Значит, прямые $AD$ и $AE$ являются медианами, проведёнными к основаниям равнобедренных треугольников и потому перпендикулярны к этим основаниям (§~\ref{1938/38}).
Если же прямые $AD$ и $AE$ перпендикулярны к радиусам $OD$ и $OE$ в их концах, лежащих
на окружности, то они касательные (§~\ref{1938/113}).

}

{\small

\smallskip
\mbox{\so{Замечание}.} Ниже (§~\ref{1938/128}) будет указан другой приём проведения касательной.

}

{

\begin{wrapfigure}{r}{48mm}
\vskip-3mm
\centering
\includegraphics{mppics/ris-146}
\caption{}\label{1938/ris-146}
\end{wrapfigure}

\paragraph{}\label{1938/129}
\mbox{\so{Задача}.}
\emph{К двум окружностям с центрами $O$ и $O_1$ провести общую касательную} (рис.~\ref{1938/ris-146}).

1) \so{Анализ}.
Предположим, что задача решена.
Пусть $AB$ будет общая касательная, $A$ и $B$ — точки касания.
Очевидно, что если мы найдём одну из этих точек, например $A$, то затем легко найдём и другую.
Проведём радиусы $OA$ и $O_1B$.
Эти радиусы, будучи перпендикулярны к общей касательной, параллельны между собой;
поэтому, если из $O_1$ проведём $O_1C\parallel BA$, то треугольник $OCO_1$ будет прямоугольный с прямым углом при вершине $C$;
вследствие этого, если опишем с центром в точке $O$ радиусом $OC$ окружность, то она будет касаться прямой $O_1C$ в точке $C$.
Радиус этой вспомогательной окружности известен:
он равен $OA-CA=OA-O_1B$, то есть
он равен разности радиусов данных окружностей.

}

\smallskip
\so{Построение}.
Таким образом, построение можно выполнить так:
описываем окружность с центром в точке $O$ радиусом, равным разности данных радиусов;
из $O_1$ проводим к этой окружности касательную $O_1C$ (способом, указанным в предыдущей задаче);
через точку касания $C$ проводим радиус $OC$ и продолжаем его до встречи с данной окружностью в точке $A$.
Наконец, из $A$ проводим $AB$ параллельно $CO_1$.

Совершенно таким же способом мы можем построить другую общую касательную $A_1B_1$.
Прямые $AB$ и $A_1B_1$ называются \rindex{внешние касательные}\textbf{внешними касательными} двух окружностей.

Аналогично можно построить и \rindex{внутренние касательные}\textbf{внутренние касательные} (рис. \ref{1938/ris-147}).

\begin{wrapfigure}{o}{48mm}
\centering
\includegraphics{mppics/ris-147}
\caption{}\label{1938/ris-147}
\end{wrapfigure}

2) \so{Анализ}.
Предположим, что задача решена.
Пусть $AB$ будет искомая касательная.
Проведём радиусы $OA$ и $O_1B$ в точки касания $A$ и $B$.
Эти радиусы, будучи оба перпендикулярны к общей касательной, параллельны между собой.

Поэтому, если из $O_1$ проведём $O_1C\parallel BA$ и продолжим $OA$ до пересечения с $O_1C$ в точке $C$, то $OC$ будет перпендикулярен к $O_1C$, вследствие этого окружность, описанная радиусом $OC$ с центром в точке $O$, будет касаться прямой $O_1C$ в точке $C$.
Радиус этой вспомогательной окружности известен, он равен:
$OA+AC=OA+O_1B$, то есть
он равен сумме радиусов данных окружностей.

\smallskip
\so{Построение}.
Таким образом, построение может быть выполнено так:
описываем окружность с центром в точке $O$ радиусом, равным сумме данных радиусов;
из $O_1$ проводим к этой окружности касательную $O_1C$;
точку касания $C$ соединяем с $O$;
наконец, через точку $A$, в которой $OC$ пересекается с данной окружностью, проводим $AB \parallel CO_1$.
Подобным же способом можно построить и другую общую внутреннюю касательную $A_1B_1$.

\section{Расположения двух окружностей}

\paragraph{}\label{1938/117}
\so{Определение}.
Если две окружности имеют только одну общую точку, то говорят, что они \rindex{касающиеся окружности}\textbf{касаются};
если же две окружности имеют две общие точки, то говорят, что они \rindex{пересекающиеся окружности}\textbf{пересекаются}.

Трёх общих точек две несливающиеся окружности иметь не могут, потому что в противном случае через три точки можно было бы провести две различные окружности, что невозможно (§~\ref{1938/104}).

Будем называть \rindex{линия!центров}\textbf{линией центров} прямую, проходящую через центры двух окружностей. 

\paragraph{}\label{1938/118}
\so{Теорема}.
\textbf{\emph{Если две окружности}} (рис.~\ref{1938/ris-131}) \textbf{\emph{имеют общую точку}} ($A$), \textbf{\emph{расположенную вне линии центров, то они имеют ещё и другую общую точку}} ($A_1$), \textbf{\emph{симметричную первой относительно линии центров}} (и, следовательно, такие окружности пересекаются).

\begin{figure}[!ht]
\centering
\includegraphics{mppics/ris-131}
\caption{}\label{1938/ris-131}
\end{figure}

Линия центров содержит в себе диаметры обеих окружностей и поэтому должна быть осью симметрии всей фигуры.
Поэтому общей точке $A$, лежащей вне линии центров, должна соответствовать симметричная общая точка $A_1$, расположенная по другую сторону от оси симметрии 
(на одном перпендикуляре к линии центров и на равном расстоянии от неё). 

\medskip

\smallskip
\so{Следствие}.
\emph{Общая хорда \emph{($AA_1$, рис.~\ref{1938/ris-131})} двух пересекающихся окружностей перпендикулярна к их линии центров и делится ею пополам.}

\paragraph{}\label{1938/119}
\so{Теорема}.
\textbf{\emph{Если две окружности имеют общую точку}} ($A$) \textbf{\emph{на линии их центров, то они касаются}} (рис.~\ref{1938/ris-132} и \ref{1938/ris-133}).

\begin{figure}[!ht]
\begin{minipage}{.55\textwidth}
\centering
\includegraphics{mppics/ris-132}
\end{minipage}
\hfill
\begin{minipage}{.41\textwidth}
\centering
\includegraphics{mppics/ris-133}
\end{minipage}

\medskip

\begin{minipage}{.55\textwidth}
\centering
\caption{}\label{1938/ris-132}
\end{minipage}
\hfill
\begin{minipage}{.41\textwidth}
\centering
\caption{}\label{1938/ris-133}
\end{minipage}
\vskip-4mm
\end{figure}

Окружности не могут иметь другой общей точки вне линии центров, потому что в противном случае они имели бы ещё третью общую точку по другую сторону от линии центров и, следовательно, должны были бы слиться.
Они не могут иметь другой общей точки и на линии центров, так как, имея на этой линии две общие точки, они должны были бы иметь и общую хорду, соединяющую эти точки.
Но хорда, проходящая через центры, должна быть диаметром;
если же окружности имеют общий диаметр, то они сливаются в одну окружность.

{\small
\sloppy

\smallskip
\so{Замечание}.
Касание двух окружностей называется \rindex{внешнее касание}\textbf{внешним}, если окружности расположены одна вне другой (рис.~\ref{1938/ris-132}), и \rindex{внутреннее касание}\textbf{внутренним}, если одна из окружностей лежит внутри другой (рис.~\ref{1938/ris-133}).

}

\paragraph{}\label{1938/120}
\so{Теорема} (обратная предыдущей).
\textbf{\emph{Если две окружности касаются}} (в точке $A$, рис.~\ref{1938/ris-132} и \ref{1938/ris-133}), \textbf{\emph{то точка касания лежит на линии центров.}}

Точка $A$ не может лежать вне линии центров, потому что в противном случае окружности имели бы ещё другую общую точку, что противоречит условию теоремы (§~\ref{1938/118}).

\paragraph{}\label{1938/121}
\so{Следствие}.
\emph{Две касающиеся окружности имеют общую касательную в точке касания}, потому что, если проведём через точку касания прямую $MN$ (рис.~\ref{1938/ris-132} и \ref{1938/ris-133}), перпендикулярную к радиусу $OA$, то эта прямая будет также перпендикулярна и к радиусу $O_1A$.

\begin{figure}[!ht]
\begin{minipage}{.54\textwidth}
\centering
\includegraphics{mppics/ris-134}
\end{minipage}
\hfill
\begin{minipage}{.42\textwidth}
\centering
\includegraphics{mppics/ris-135}
\end{minipage}

\medskip

\begin{minipage}{.54\textwidth}
\centering
\caption{}\label{1938/ris-134}
\end{minipage}
\hfill
\begin{minipage}{.42\textwidth}
\centering
\caption{}\label{1938/ris-135}
\end{minipage}
\vskip-4mm
\end{figure}

{\sloppy 
\paragraph{Различные случаи взаимного расположения двух окружностей.}\label{1938/122}
Обозначим радиусы двух окружностей буквами $R$ и $R_1$ и расстояние между их центрами буквой $d$.
Можно предположить, что $R_1\z\le R$.
Рассмотрим, какова зависимость между этими величинами в различных случаях взаимного расположения двух окружностей.
Этих случаев можно указать пять, а именно.

}

1) \so{Окружности лежат одна вне другой, не касаясь} (рис. \ref{1938/ris-134});
в этом случае, очевидно, $d>R+R_1$.

2) \so{Окружности имеют внешнее касание} (рис. \ref{1938/ris-135});
тогда $d\z=R+R_1$, так как точка касания лежит на линии центров.

3) \so{Окружности пересекаются} (рис.~\ref{1938/ris-131});
тогда $d\z<R\z+R_1$ и в то же время $d>R-R_1$, потому что в $\triangle OAO_1$%
\footnote{На рис.~\ref{1938/ris-131} провести прямые $OA$ и $O_1A$.}
сторона $OO_1$, равная $d$, меньше суммы, но больше разности двух других сторон, равных радиусам $R$ и $R_1$.

{

\begin{wrapfigure}{o}{47mm}
\vskip-6mm
\centering
\includegraphics{mppics/ris-136}
\caption{}\label{1938/ris-136}
\end{wrapfigure}

4) \so{Окружности имеют внутреннее касание} (рис. \ref{1938/ris-133});
в этом случае $d=R-R_1$, потому что точка касания лежит на линии центров.

{\sloppy

5) \so{Одна окружность лежит внутри другой, не касаясь} (рис.~\ref{1938/ris-136});
тогда, очевидно, $d\z<R-R_1$, и в частном случае $d= 0$, когда центры обеих окружностей сливаются (такие окружности называются \rindex{концентрические окружности}\textbf{концентрическими}).

}

}

{\small

\smallskip
\mbox{\so{Замечание}.}
Учащимся предлагается проверить правильность обратных предложений, а именно:

1) \emph{Если $d>R+R_1$, то окружности расположены одна вне другой, не касаясь.}

2) \emph{Если $d=R+R_1$, то окружности касаются извне.}

3) \emph{Если $d<R+R_1$, и в то же время $d>R-R_1$, то окружности пересекаются.} 

4) \emph{Если $d=R-R_1$, то окружности касаются изнутри.}

5) \emph{Если $d<R-R_1$, то одна окружность лежит внутри другой, не касаясь.}

\smallskip
\so{Указание:} Все эти предложения легко доказываются от противного;
в доказательстве 3), следует воспользоваться допущением в §~\ref{1914/118}. 

}


%% file: 2D/vpisannye-ugly.tex
\section{Вписанные углы}

\paragraph{Вписанный угол.}\label{1938/123}
Угол, образованный двумя хордами, исходящими из одной точки окружности, называется \rindex{вписанный угол}\textbf{вписанным} углом.
Таков, например, угол $ABC$ (рис.~\ref{1938/ris-137}).

О вписанном угле принято говорить, что он опирается на дугу, заключённую между его сторонами.
Так, угол $ABC$ {}\textbf{опирается} на дугу $AC$.

\paragraph{}\label{1938/124} 
\so{Теорема}.
\textbf{\emph{Вписанный угол измеряется половиной дуги, на которую он опирается.}}
Эту теорему надо понимать так:
вписанный угол содержит в себе столько угловых градусов, минут и секунд, сколько дуговых градусов, минут и секунд заключается в половине дуги, на которую он опирается.

\begin{figure}[h]
\begin{minipage}{.32\textwidth}
\centering
\includegraphics{mppics/ris-137}
\end{minipage}
\hfill
\begin{minipage}{.32\textwidth}
\centering
\includegraphics{mppics/ris-138}
\end{minipage}
\hfill
\begin{minipage}{.32\textwidth}
\centering
\includegraphics{mppics/ris-139}
\end{minipage}

\medskip

\begin{minipage}{.32\textwidth}
\centering
\caption{}\label{1938/ris-137}
\end{minipage}
\hfill
\begin{minipage}{.32\textwidth}
\centering
\caption{}\label{1938/ris-138}
\end{minipage}
\hfill
\begin{minipage}{.32\textwidth}
\centering
\caption{}\label{1938/ris-139}
\end{minipage}
\vskip-4mm
\end{figure}

При доказательстве теоремы рассмотрим особо три случая:

1) Центр $O$ (рис.~\ref{1938/ris-137}) лежит на стороне вписанного угла $ABC$.
Проведя радиус $AO$, мы получим $\angle AOB$, в котором $OA \z= OB$ (как радиусы), и, следовательно, $\angle ABO\z=\angle BAO$.
По отношению к этому треугольнику угол $AOC$ есть внешний, поэтому он равен сумме углов $ABO$ и $BAO$ и, значит, равен двум углам $ABO$;
поэтому угол $ABO$ равен половине центрального угла $AOC$.
Но угол $AOC$ измеряется дугой $AC$, то есть
он содержит в себе столько угловых градусов, минут и секунд, сколько дуговых градусов, минут и секунд содержится в дуге $AC$;
следовательно, вписанный угол $ABC$ измеряется половиной дуги $AC$.

2) Центр $O$ лежит между сторонами вписанного угла $ABC$ (рис.~\ref{1938/ris-138}).

Проведя диаметр $BD$, мы разделим угол $ABC$ на два угла, из которых, по доказанному, один измеряется половиной дуги $AD$, а другой — половиной дуги $DC$;
следовательно, угол $ABC$ измеряется суммой
\[\tfrac12{\smallsmile}AD\z+\tfrac12{\smallsmile}DC=\tfrac12({\smallsmile}AD\z+{\smallsmile}DC)=\tfrac12{\smallsmile}AC.\]

3) Центр $O$ лежит вне вписанного угла $ABC$ (рис.~\ref{1938/ris-139}).
Проведя диаметр $BD$, мы будем иметь:
\[\angle ABC=\angle ABD-\angle CBD.\]

Но углы $ABD$ и $CBD$ измеряются, по доказанному, половинами дуг $AD$ и $CD$;
следовательно, угол $ABC$ измеряется разностью
\[\tfrac12{\smallsmile}AD-\tfrac12{\smallsmile}CD=\tfrac12({\smallsmile}AD-{\smallsmile}CD)=\tfrac12{\smallsmile}AC.\]

\begin{figure}[h]
\begin{minipage}{.48\textwidth}
\centering
\includegraphics{mppics/ris-140}
\end{minipage}
\hfill
\begin{minipage}{.48\textwidth}
\centering
\includegraphics{mppics/ris-141}
\end{minipage}

\medskip

\begin{minipage}{.48\textwidth}
\centering
\caption{}\label{1938/ris-140}
\end{minipage}
\hfill
\begin{minipage}{.48\textwidth}
\centering
\caption{}\label{1938/ris-141}
\end{minipage}
\vskip-4mm
\end{figure}

\paragraph{}\label{1938/125}
\so{Следствия}. 
1) \emph{Все вписанные углы, опирающиеся на одну и ту же дугу, равны между собой} (рис.~\ref{1938/ris-140}), потому что каждый из них измеряется половиной одной и той же дуги.
Если величину одного из таких углов обозначим $\alpha$, то можно сказать, что сегмент $AmB$, покрытый на чертеже штрихами, \so{вмещает в себя угол}, равный~$\alpha$. 

2) \emph{Всякий вписанный угол, опирающийся на диаметр, есть прямой} (рис.~\ref{1938/ris-141}), потому что каждый такой угол измеряется половиной полуокружности и, следовательно, содержит $90\degree$.

\paragraph{}\label{1938/126}
\so{Задача}.
\emph{Построить прямоугольный треугольник по гипотенузе $a$ и катету $b$} (рис.~\ref{1938/ris-142}).

На какой-нибудь прямой $MN$ отложим $AB=a$, на $AB$ опишем полуокружность.
Затем проводим дугу радиусом, равным $b$, с центром в точке $A$ (или $B$).

Точку пересечения $C$ полуокружности и дуги соединим с концами диаметра $AB$.
Треугольник $ABC$ будет искомый, так как угол $C$ — прямой, $a$ является гипотенузой, а $b$ — катетом.

\begin{figure}[!ht]
\begin{minipage}{.48\textwidth}
\centering
\includegraphics{mppics/ris-142}
\end{minipage}
\hfill
\begin{minipage}{.48\textwidth}
\centering
\includegraphics{mppics/ris-143}
\end{minipage}

\medskip

\begin{minipage}{.48\textwidth}
\centering
\caption{}\label{1938/ris-142}
\end{minipage}
\hfill
\begin{minipage}{.48\textwidth}
\centering
\caption{}\label{1938/ris-143}
\end{minipage}
\vskip-4mm
\end{figure}

\paragraph{}\label{1938/127}
\so{Задача}.
\emph{Из конца $A$ \emph{(рис.~\ref{1938/ris-143})} данного луча $AB$, не продолжая его, восстановить к нему перпендикуляр.}

Возьмём вне прямой какую-либо точку $O$ так, чтобы окружность с центром в этой точке и радиусом, равным отрезку $OA$, пересекла луч $AB$ в какой-либо точке $C$.
Через эту точку $C$ проведём диаметр $CD$ и конец его $D$ соединим с $A$.
Прямая $AD$ есть искомый перпендикуляр, потому что угол $A$ прямой, как вписанный и опирающийся на диаметр.

\begin{wrapfigure}{r}{42mm}
\centering
\includegraphics{mppics/ris-145}
\caption{}\label{1938/ris-145}
\end{wrapfigure}

\paragraph{}\label{1938/128}
\mbox{\so{Задача}.}
\emph{Через данную точку провести касательную к данной окружности.}

Эта задача уже  была решена в §~\ref{1931/115};
мы приведём другое решение, случай, когда \so{данная точка} ($A$, рис.~\ref{1938/ris-145}) \so{лежит вне окружности} (с центром $O$).

Соединив $A$ с $O$, делим $AO$ пополам в точке $O_1$ и с центром в этой точке радиусом $OO_1$ описываем окружность.
Через точки $B$ и $B_1$, в которых эта окружность пересекается с данной, проводим прямые $AB$ и $AB_1$.
Эти прямые и будут касательными, так как углы $OBA$ и $OB_1A$, как опирающиеся на диаметр, — прямые.

\paragraph{}\label{1938/130}
\so{Теорема}.
1) \textbf{\emph{Угол}} ($ABC$, рис.~\ref{1938/ris-148}), \textbf{\emph{вершина которого лежит внутри круга, измеряется полусуммой двух дуг}} ($AC$ и $DE$), \textbf{\emph{из которых одна заключена между его сторонами, а другая — между продолжениями сторон.}}

2) \textbf{\emph{Угол}} ($ABC$, рис.~\ref{1938/ris-149}), \textbf{\emph{вершина которого лежит вне круга и стороны пересекаются с окружностью, измеряется полуразностью двух дуг}} ($AC$ и $ED$), \textbf{\emph{заключённых между его сторонами.}}

\begin{figure}[h]
\begin{minipage}{.48\textwidth}
\centering
\includegraphics{mppics/ris-148}
\end{minipage}
\hfill
\begin{minipage}{.48\textwidth}
\centering
\includegraphics{mppics/ris-149}
\end{minipage}

\medskip

\begin{minipage}{.48\textwidth}
\centering
\caption{}\label{1938/ris-148}
\end{minipage}
\hfill
\begin{minipage}{.48\textwidth}
\centering
\caption{}\label{1938/ris-149}
\end{minipage}
\vskip-4mm
\end{figure}

Проведя хорду $AD$ (на том и на другом чертеже), мы получим $\triangle ABD$, относительно которого рассматриваемый угол $ABC$ служит внешним, когда его вершина лежит внутри круга, и внутренним, когда его вершина лежит вне круга.
Поэтому в первом случае:
\[\angle ABC = \angle ADC+\angle DAE;\]
во втором случае:
\[\angle ABC = \angle ADC-\angle DAE.\]

Но углы $ADC$ и $DAE$, как вписанные, измеряются половинами дуг $AC$ и $DE$, поэтому угол $ABC$ измеряется:
в первом случае суммой
$\tfrac12{\smallsmile}AC+\tfrac12{\smallsmile}DE$, которая равна $\tfrac12({\smallsmile}AC+{\smallsmile}DE)$, а во втором случае разностью $\tfrac12{\smallsmile}AC-\tfrac12{\smallsmile}DE$, которая равна $\tfrac12({\smallsmile}AC-{\smallsmile}DE)$.

\paragraph{}\label{1938/131}
\so{Теорема}.
\textbf{\emph{Угол}} ($ACD$, рис.~\ref{1938/ris-150} и \ref{1938/ris-151}), \textbf{\emph{составленный касательной и хордой, измеряется половиной дуги, заключённой внутри него.}}

\begin{wrapfigure}{o}{30mm}
\vskip-3mm
\centering
\includegraphics{mppics/ris-150}
\caption{}\label{1938/ris-150}
\bigskip
\includegraphics{mppics/ris-151}
\caption{}\label{1938/ris-151}
\end{wrapfigure}

Предположим сначала, что хорда $CD$ проходит через центр $O$, то есть что эта хорда есть диаметр (рис.~\ref{1938/ris-150}).
Тогда угол $ACD$ — прямой и,
следовательно, равен $90\degree$.
Но и половина дуги $CmD$ также равна $90\degree$, так как целая дуга $CmD$, составляя полуокружность, содержит $180\degree$.
Значит, теорема справедлива в этом частном случае.

Рассмотрим общий случай (рис.~\ref{1938/ris-151}), когда хорда $CD$ не проходит через центр.
Проведя диаметр $CE$, мы будем иметь:
\[\angle ACD = \angle ACE - \angle DCE.\]

Угол $ACE$, как составленный касательной и диаметром, измеряется, по доказанному, половиной дуги $CDE$;
угол $DCE$, как вписанный, измеряется половиной дуги $DE$;
следовательно, угол $ACD$ измеряется разностью $\tfrac12{\smallsmile}CDE-\tfrac12{\smallsmile}DE$, то есть половиной дуги $CD$.

Подобным же образом можно доказать, что тупой угол $BCD$ (рис. \ref{1938/ris-151}), также составленный касательной и хордой, измеряется половиной дуги $CnED$;
разница в доказательстве только та, что этот угол надо рассматривать не как разность, а как сумму прямого угла $BCE$ и острого $ECD$.

\paragraph{}\label{1938/132}
\mbox{\so{Задача}.}
\emph{На данном отрезке $AB$ построить сегмент, вмещающий данный угол $\alpha$} (рис.~\ref{1938/ris-152}).

\smallskip
\mbox{\so{Анализ}.}
Предположим, что задача решена;
пусть сегмент $AmB$ \textbf{вмещает} в себя угол $\alpha$, то есть всякий вписанный в него угол $ACB$ равен $\alpha$.
Проведём вспомогательную прямую $AE$, касательную к окружности в точке $A$.
Тогда угол $BAE$, составленный касательной и хордой, должен равняться вписанному углу $ACB$, так как и тот, и другой угол измеряется половиной дуги $AnB$.
Примем во внимание, что центр $O$ окружности должен лежать на срединном перпендикуляре $DO$, проведённом к отрезку $AB$, и в то же время он должен лежать и на перпендикуляре $AO$, восстановленном к касательной $AE$ из точки касания.
Отсюда выводим следующее построение.

\begin{wrapfigure}{o}{33mm}
\centering
\includegraphics{mppics/ris-152}
\caption{}\label{1938/ris-152}
\end{wrapfigure}

\smallskip
\mbox{\so{Построение}.}
При конце отрезка $AB$ строим угол $DAE$, равный углу $\alpha$;
проводим срединный перпендикуляр $DO$ к $AB$ и из точки $A$ восстанавливаем перпендикуляр к $AE$; 
точку пересечения $O$ этих двух перпендикуляров принимаем за центр и радиусом $AO$ описываем окружность.

\smallskip
\mbox{\so{Доказательство}.}
Сегмент $AmB$ будет искомый, потому что всякий вписанный в него угол измеряется половиной дуги $AnB$, а половина этой дуги измеряет также $\angle DAE=\alpha$.

{\small
\smallskip
\mbox{\so{Замечание}.}
На рис.~\ref{1938/ris-152} построен сегмент, расположенный выше отрезка $AB$.
Такой же сегмент может быть построен и по другую сторону отрезка $AB$.
Таким образом, можно сказать, что геометрическое место точек, из которых данный отрезок $AB$ виден под данным углом $\alpha$, состоит из дуг двух сегментов, из которых каждый вмещает в себя угол $\alpha$ и один расположен по одну сторону отрезка $AB$, а другой — по другую сторону.

}

\subsection*{Задачи на построение}

\paragraph{Метод геометрических мест.}\label{1938/133}
Для решения многих задач на построение можно с успехом применять понятие о геометрическом месте и основанный на нём метод геометрических мест.
Этот метод, известный ещё со времён Платона (IV век до нашей эры), состоит в следующем.

Положим, что решение предложенной задачи сводится к нахождению некоторой точки, которая должна удовлетворять известным условиям.
Отбросим из этих условий какое-нибудь одно;
тогда задаче может удовлетворять бесчисленное множество точек.
Эти точки составят некоторое геометрическое место.
Построим его, если это окажется возможным.
Затем примем во внимание отброшенное нами условие и откинем какое-нибудь другое;
тогда задача будет снова удовлетворяться бесчисленным множеством точек, которые составят новое геометрическое место.
Построим его, если это возможно.
Искомая точка, удовлетворяя всем условиям, должна лежать на обоих геометрических местах, то есть она должна находиться в их пересечении.
Задача не имеет решения, если найденные геометрические места не пересекаются;
и задача будет иметь столько решений, сколько окажется точек пересечения.

Приведём на этот метод один пример, который вместе с тем покажет нам, как иногда приходится вводить в чертёж вспомогательные линии с целью принять во внимание все данные условия задачи.

\paragraph{}\label{1938/134}
\so{Задача}.
\emph{Построить треугольник по основанию $a$, углу при вершине $A$ и сумме $s$ боковых сторон.}

\begin{wrapfigure}{o}{53mm}
\centering
\includegraphics{mppics/ris-153}
\caption{}\label{1938/ris-153}
\end{wrapfigure}

Пусть $ABC$ будет искомый треугольник (рис.~\ref{1938/ris-153}).
Чтобы принять во внимание данную сумму боковых сторон, продолжим $BA$ и отложим $BM=s$.
Проведя $MC$, получим вспомогательный треугольник $BMC$.
Если мы построим этот треугольник, то затем легко построить искомый треугольник $ABC$.

Построение треугольника $BMC$ сводится к нахождению точки~$M$.

Заметив, что треугольник $AMC$ равнобедренный ($AM=AC$) и, следовательно, $\angle M \z=  \tfrac12\angle BAC$ (так как $\angle M+\angle MCA = \angle BAC$), мы видим, что точка $M$ должна удовлетворять двум условиям:
1) она удалена от $B$ на расстояние $s$, 
2) из неё данный отрезок $BC$ виден под углом, равным $\tfrac12\angle A$.
Отбросив второе условие, мы получим бесчисленное множество точек $M$, лежащих на окружности, описанной радиусом, равным $s$, с центром в точке $B$.
Отбросив первое условие, мы получим также бесчисленное множество точек $M$, лежащих на дуге сегмента, построенного на $BC$ и вмещающего угол, равный $\tfrac12\angle A$.
Таким образом, нахождение точки $M$ сводится к построению двух геометрических мест, из которых каждое мы построить умеем.
Задача не будет иметь решения, если эти геометрические места не будут иметь общих точек;
задача будет иметь одно или два решения, смотря по тому, касаются или же пересекаются эти геометрические места (на нашем чертеже получаются два треугольника $ABC$ и $A_1BC$, удовлетворяющие условиям задачи).

\medskip

Иногда задача сводится не к определению точки, а к нахождению прямой, удовлетворяющей нескольким условиям.
Если отбросим одно из условий, то получим бесчисленное множество прямых;
при этом может случиться, что эти прямые определяют некоторую линию (например, все они будут касательными к некоторой окружности).
Отбросив другое условие и приняв во внимание то, которое было откинуто ранее, мы получим снова бесчисленное множество прямых, которые, быть может, определят некоторую другую линию.
Построив, если возможно, эти две линии, мы затем легко найдём и искомую прямую.
Приведём пример.

\paragraph{}\label{1938/135}
\so{Задача}.
\emph{Провести секущую к двум данным окружностям с центрами $O$ и $O_1$ так, чтобы части секущей, заключённые внутри окружностей, равнялись соответственно данным отрезкам $a$ и~$a_1$.}

\begin{wrapfigure}{o}{53mm}
\vskip-2mm
\centering
\includegraphics{mppics/ris-extra-6}
\caption{}\label{extra/ris-6}
\end{wrapfigure}

Если возьмём только одно условие, например, чтобы часть секущей, лежащая внутри круга  с центром $O$, равнялась $a$, то получим бесчисленное множество секущих, которые все должны быть одинаково удалены от $O$ (так как равные хорды одинаково удалены от центра).

Поэтому, если в этом круге где-нибудь построим хорду, равную $a$, а затем радиусом, равным расстоянию от этой хорды до центра, опишем окружность, с центром в $O$, то все секущие, о которых идёт речь, должны касаться этой вспомогательной окружности.

Подобным образом, приняв во внимание только второе условие, мы увидим, что искомая секущая должна касаться второй вспомогательной окружности с центром $O_1$.
Значит, вопрос сводится к построению общей касательной к двум окружностям.
{\small

\subsection*{Упражнения}

\begin{center}
\so{Найти геометрические места}
\end{center}

\begin{enumerate}[noitemsep]

\item
Точек, из которых касательные, проведённые к данной окружности, равны данному отрезку.
 
\item

Точек, из которых данная окружность видна под данным углом (то есть две касательные, проведённые из данной точки к окружности, составляют между собой данный угол).

\item
Центров окружностей, описанных данным радиусом и касающихся данной прямой.

\item
Центров окружностей, описанных данным радиусом и касающихся данной окружности (два случая:
касание внешнее и касание внутреннее).

\item
Отрезок данной длины движется параллельно самому себе так, что один его конец скользит по окружности.
Найти геометрическое место, описанное другим концом.

\smallskip
\so{Указание}.
Возьмём две прямые, изображающие два положения, движущейся прямой, и через концы их, лежащие на окружности, проведём радиусы, а через другие концы проведём прямые, параллельные этим радиусам, до пересечения с прямой, проходящей через центр и параллельной движущейся линии.
Рассмотрим образовавшиеся параллелограммы.

\item
Отрезок данной длины движется так, что концы его скользят по сторонам прямого угла.
Найти геометрическое место, описываемое серединой этого отрезка.

\end{enumerate}

\begin{center}
\so{Доказать теоремы}
\end{center}

\begin{enumerate}[resume,noitemsep]

\item
Из всех хорд, проходящих через точку $A$, взятую внутри данного круга, наименьшей будет та, которая перпендикулярна к диаметру, проходящему через точку $A$.

\item
На хорде $AB$ взяты две точки $D$ и $E$ на равном расстоянии от середины $C$ этой хорды, и через эти точки восстановлены к $AB$ перпендикуляры $DF$ и $EG$ до пересечения с окружностью.
Доказать, что эти перпендикуляры равны.

\smallskip
\so{Указание}.
Перегнуть чертёж по диаметру.

\item
В круге проведены две хорды $CC'$ и $DD'$ перпендикулярно к диаметру $AB$.
Доказать, что прямая $MM'$, соединяющая середины хорд $CD$ и $C'D'$, перпендикулярна к $AB$.

\item
В круге с центром $O$ проведена хорда $AB$ и продолжена на расстояние $BC$, равное радиусу.
Через точку $C$ и центр $O$ проведена секущая $CD$ ($D$ — вторая точка пересечения с окружностью).
Доказать, что угол $AOD$ равен утроенному углу $ACD$.

\item
Если через центр окружности и данную точку вне её проведём секущую, то часть её, заключённая между данной точкой и ближайшей точкой пересечения, есть наименьшее расстояние, а часть, заключённая между данной точкой и другой точкой пересечения, есть наибольшее расстояние от данной точки до точек окружности.

\item
Кратчайшее расстояние между двумя окружностями, лежащими одна вне другой, есть отрезок линии центров, заключённый между окружностями.

\item
Если через точку пересечения двух окружностей будем проводить секущие, не продолжая их за окружности, то наибольшей из них окажется та, которая параллельна линии центров.

\item
Если к двум окружностям, касающимся извне, провести три общие касательные, то внутренняя из них делит пополам тот отрезок каждой внешней, который ограничен точками касания.

\item
Через точку $A$ окружности проведена хорда $AB$ и затем касательная в точке $B$.
Диаметр, перпендикулярный радиусу $OA$, встречает касательную и хорду (или её продолжение) соответственно в точках $C$ и $O$.
Доказать, что $BC=CD$.

\item
К двум окружностям с центрами $O$ и $O_1$, касающимся извне в точке $A$, проведена общая внешняя касательная $BC$ ($B$ и $C$ — точки касания);
доказать, что угол $BAC$ есть прямой.

\smallskip
\so{Указание}.
Провести в точке $A$ общую касательную и рассмотреть равнобедренные треугольники $ABD$ и $ADC$.

\item
Две прямые исходят из одной и той же точки $M$ и касаются окружности в точках $A$ и $B$.
Проведя радиус $OB$, продолжают его за точку $B$ на расстояние $BC=OB$.
Доказать, что $\angle AMC\z=3\angle BMC$.

\item
Две прямые, исходящие из точки $M$, касаются окружности в точках $A$ и $B$.
На меньшей из двух дуг, ограниченных точками $A$ и $B$, берут произвольную точку $C$ и через неё проводят третью касательную до пересечения с $MA$ и $MB$ в точках $D$ и $E$.
Доказать, что:
1) периметр $\triangle MDE$ и 2) угол $DOE$ не изменяются при изменении положения точки $C$.

\smallskip
\so{Указание}.
Периметр $DME=MA+MB$, $\angle DOE\z=\tfrac12 \angle AOB$.
\end{enumerate}


\begin{center}
\so{Задачи на построение}
\end{center}

\begin{enumerate}[resume,noitemsep]
\item
Разделить данную дугу на $4, 8, 16, \dots$
равных частей.

\item
По сумме и разности дуг одного и того же радиуса найти эти дуги.

\item
Описать такую окружность с центром в данной точке, которая разделила бы данную окружность пополам.

\item
На данной прямой найти точку, наименее удалённую от данной окружности.

\item
В круге дана хорда.
Провести другую хорду, которая делилась бы первой пополам и составляла бы с ней данный угол (при всяком ли данном угле задача имеет решение?).

\item
Через данную в круге точку провести хорду, которая делилась бы этой точкой пополам.

\item
С центром в точке на стороне данного угла, описать окружность, которая от другой стороны угла отсекла бы хорду данной длины.

\item
Данным радиусом описать окружность, центр которой лежал бы на стороне данного угла и которая от другой стороны его отсекла бы хорду данной длины.

\item
Данным радиусом описать окружность, которая касалась бы данной прямой в данной точке.

\item
Провести касательную к данной окружности параллельно данной прямой.

\item
Описать окружность, которая проходила бы через данную точку $A$ и касалась бы данной прямой в данной на ней точке $B$.

\item
Описать окружность, касательную к сторонам данного угла, причём к одной из них в данной точке.

\item
Между двумя параллельными прямыми дана точка;
провести окружность, проходящую через эту точку и касающуюся данных прямых.

\item
Провести к данной окружности касательную под данным углом к данной прямой (сколько решений?).

\item
Из точки, данной вне круга, провести к нему секущую так, чтобы её внутренняя часть равнялась данному отрезку (исследовать задачу).

\item
Данным радиусом описать окружность, проходящую через данную точку и касательную к данной прямой.

\item
На данной прямой найти такую точку, чтобы касательные, проведённые из неё к данной окружности, были данной длины.

\item
Построить треугольник, зная один угол и две высоты, из которых одна проведена из вершины данного угла.

\item
Даны две точки;
провести прямую так, чтобы перпендикуляры, опущенные на неё из этих точек, имели данные длины.

\item
Описать окружность, которая проходила бы через данную точку и касалась бы данной окружности в данной точке.

\item
Описать окружность, которая касалась бы двух данных параллельных прямых и круга, находящегося между ними.

\item
Данным радиусом описать окружность, которая касалась бы данного круга и проходила бы через данную точку.
(Рассмотреть три случая:
данная точка лежит 
1) вне круга, 
2) на окружности 
и 3) внутри круга.)

\item
Данным радиусом описать окружность, которая касалась бы данной прямой и данного круга.

\item
Данным радиусом описать окружность, которая от сторон данного угла отсекала бы хорды данной длины.

\item
Описать окружность, касающуюся данного круга в данной точке и данной прямой (два решения).

\item
Описать окружность, касающуюся данной прямой в данной точке и данного круга (два решения).

\item
Описать окружность, касающуюся двух данных кругов, причём одного из них в данной точке.
(Рассмотреть три случая:
1) искомый круг лежит вне данных кругов;
2) один из данных кругов лежит вне искомого, другой внутри;
3) оба данных круга лежат внутри искомого.)

\item
Описать окружность, касающуюся трёх равных кругов извне или внутри.

\item
В данный сектор вписать окружность, касающуюся радиусов, ограничивающих сектор, и дуги сектора.

\item
Вписать в данный круг три равных круга, которые касались бы попарно между собой и данного круга.

\item
Через точку, данную внутри круга, провести хорду так, чтобы разность её отрезков равнялась данному отрезку.

\smallskip
\so{Указание}.
Провести окружность, концентрическую данной и проходящую через данную точку.
В этой окружности от данной точки построить хорду данной длины.

\item
Через точку пересечения двух окружностей провести секущую так, чтобы часть её, заключённая внутри окружностей, равнялась данной длине.

\smallskip
\so{Указание}.
Построить прямоугольный треугольник, имеющий гипотенузой отрезок прямой, соединяющий центры данных окружностей с катетом, равным половине данного отрезка, и~т.~д.

\item
Из точки, данной вне круга, провести секущую так, чтобы внешняя часть её равнялась внутренней.

\smallskip
\so{Указание}.
Пусть $O$ — центр окружности, $R$ — её радиус, $A$~— данная точка.
Строим $\triangle AOB$, где $AB=R$, $OB=2R$.
Если $C$ — середина отрезка $OB$, то прямая $AC$ — искомая.

\item
Построить окружность, касающуюся данной прямой в данной точке, и проходящую через другую данную точку.

\item
Построить окружность, касающуюся двух прямых, если дана одна из точек касания.

\end{enumerate}

}

%% file: 2D/vpis-opis-mnougi.tex
\section{Вписанные и описанные многоугольники}

\paragraph{}\label{1938/136}
\mbox{\so{Определения}.}
Если все вершины многоугольника $ABCDE$ лежат на окружности (рис.~\ref{1938/ris-155}), то говорят, что этот многоугольник \rindex{вписанный многоугольник}\textbf{вписан в окружность}, или что окружность \rindex{описанная окружность}\textbf{описана вокруг него}.

\begin{wrapfigure}[10]{r}{40mm}
\centering
\includegraphics{mppics/ris-155}
\caption{}\label{1938/ris-155}
\end{wrapfigure}

Если все стороны какого-нибудь многоугольника ($MNPQ$, рис.~\ref{1938/ris-155}) касаются окружности, то говорят, что этот многоугольник \rindex{описанный многоугольник}\textbf{описан около окружности}, или что окружность \rindex{вписанная окружность}\textbf{вписана в него}.

{
\sloppy

\paragraph{}\label{1938/137}
\mbox{\so{Теоремы}.}
1) \textbf{\emph{Около всякого треугольника можно описать окружность и притом только одну.}}

}

2) \textbf{\emph{Во всякий треугольник можно вписать окружность и притом только одну.}}

1) Вершины $A$, $B$ и $C$ всякого треугольника не лежат на одной прямой, а через такие точки, как мы видели (§~\ref{1938/104}), всегда можно провести окружность и притом только одну.

\begin{wrapfigure}{o}{40mm}
\centering
\includegraphics{mppics/ris-156}
\caption{}\label{1938/ris-156}
\end{wrapfigure}

2) Если существует такая окружность, которая касалась бы всех сторон треугольника $ABC$ (рис.~\ref{1938/ris-156}),
то её центр должен быть точкой, одинаково удалённой от этих сторон.
Докажем, что такая точка существует.
Геометрическое место точек, равно отстоящих от сторон $AB$ и $AC$, есть биссектриса $AM$ угла $A$ (§~\ref{1938/60});
геометрическое место точек, равно отстоящих от сторон $BA$ и $BC$, есть биссектриса $BN$ угла $B$.
Эти две биссектрисы должны, очевидно, пересечься внутри треугольника в некоторой точке $O$.
Эта точка и будет равноудалённой от всех сторон треугольника, так как она находится на обоих геометрических местах.

Итак, чтобы вписать круг в треугольник, делим какие-нибудь два угла его, например $A$ и $B$, пополам и точку пересечения биссектрис берём за центр.
За радиус берём один из перпендикуляров $OP$, $OQ$ или $OR$, опущенных из центра на стороны треугольника.
Окружность коснётся сторон в точках $P$, $Q$ и $R$, так как стороны в этих точках перпендикулярны к радиусам в их концах, лежащих на окружности (§~\ref{1938/113}).
Другой вписанной окружности не может быть, так как две биссектрисы пересекаются только в одной точке, а из одной точки на прямую можно опустить только один перпендикуляр.

{\small

\smallskip
\so{Замечание}.
Оставляем самим учащимся убедиться, что центр описанной окружности лежит внутри треугольника только тогда, когда треугольник остроугольный;
в тупоугольном же треугольнике он лежит вне его, а в прямоугольном — на середине гипотенузы.
Центр вписанной окружности лежит всегда внутри треугольника.

}

\smallskip
\so{Следствие}.
Точка $O$ (рис.~\ref{1938/ris-156}), находясь на одинаковом расстоянии от сторон $CA$ и $CB$, должна лежать на биссектрисе угла $C$;
следовательно, \emph{биссектрисы трёх углов треугольника пересекаются в одной точке.}

\begin{wrapfigure}{r}{45mm}
\vskip-4mm
\centering
\includegraphics{mppics/ris-157}
\caption{}\label{1938/ris-157}
\end{wrapfigure}

{\small

\paragraph{Вневписанные окружности.}\label{1938/138}
Вневписанными называются окружности (рис.~\ref{1938/ris-157}), которые касаются одной стороны треугольника и \so{продолжений} двух других сторон (они лежат вне треугольника, вследствие чего и получили название \rindex{вневписанная окружность}\textbf{вневписанных}).

Таких окружностей для всякого треугольника может быть три.
Чтобы построить их, проводят биссектрисы внешних углов треугольника $ABC$ и точки их пересечений берут за центры.
Так, центром окружности, вписанной в угол $A$, служит точка $O_a$, то есть
пересечение биссектрис $BO_a$ и $CO_a$ внешних углов, не смежных с $A$;
радиус этой окружности есть перпендикуляр, опущенный из $O_a$ на какую-либо сторону треугольника.

\begin{wrapfigure}{r}{27mm}
\centering
\includegraphics{mppics/ris-158}
\caption{}\label{1938/ris-158}
\end{wrapfigure}

}

\paragraph{Свойства вписанного выпуклого четырехугольника.}\label{1938/139}
1) \textbf{\emph{В выпуклом вписанном четырёхугольнике сумма противоположных углов равна $\bm{180\degree}$.}}

2) \textbf{\emph{Обратно, если в выпуклом четырёхугольнике сумма противоположных углов равна $\bm{180\degree}$, то около него можно описать окружность.}}

1) Пусть $ABCD$ (рис.~\ref{1938/ris-158}) есть вписанный выпуклый четырёхугольник;
требуется доказать, что
\[\angle B+\angle D = 180\degree
\quad\text{и}\quad 
\angle A + \angle C = 180\degree.\]

Так как сумма всех четырёх углов всякого выпуклого четырёхугольника равна $360\degree$ (§~\ref{1938/82}), то достаточно доказать только одно из требуемых равенств.

Докажем, например, что $\angle B+\angle D = 180\degree$.

Углы $B$ и $D$, как вписанные, измеряются:
первый — половиной дуги $ADC$, второй — половиной дуги $ABC$;
следовательно, сумма $\angle B+\angle D$ измеряется суммой $\tfrac12{\smallsmile}ADC + \tfrac12{\smallsmile}ABC$, а эта сумма равна $\tfrac12({\smallsmile}ADC\z+{\smallsmile}ABC)$, то есть
равна половине окружности;
значит:
\[\angle B+\angle D=180\degree.\]

2) Пусть $ABCD$ (рис.~\ref{1938/ris-158}) есть такой выпуклый четырёхугольник, у которого $\angle B+\angle D = 180\degree$, и, следовательно, $\angle A \z+ \angle C =180\degree$.
Требуется доказать, что около такого четырёхугольника можно описать окружность.

Через какие-нибудь три его вершины, например через $A$, $B$ и $C$, проведём окружность (что всегда можно сделать).
Четвёртая вершина $D$ должна находиться на этой окружности, потому что в противном случае вершина угла $B$ лежала бы или внутри круга, или вне его, и тогда этот угол не измерялся бы половиной дуги $ABC$;
поэтому сумма $\angle B+\angle D$ не измерялась бы полусуммой дуг $ADC$ и $ABC$ (§~\ref{1938/130}) и, значит, сумма $\angle B+\angle D$ не равнялась бы $180\degree$, что противоречит условию.

\smallskip
\so{Следствия}.
1) \emph{Из всех параллелограммов только около прямоугольника можно описать окружность.}

2) \emph{Около трапеции можно описать окружность только тогда, когда она равнобочная.}

\paragraph{Свойство описанного четырехугольника.}\label{1938/140}
\textbf{\emph{В описанном четырёхугольнике суммы противоположных сторон равны.}}

\begin{wrapfigure}{o}{44mm}
\centering
\includegraphics{mppics/ris-159}
\caption{}\label{1938/ris-159}
\end{wrapfigure}

Пусть $ABCD$ (рис.~\ref{1938/ris-159}) будет описанный четырёхугольник, то есть стороны его касаются окружности;
требуется доказать, что 
\[AB+CD=BC+AD.\]

Обозначим точки касания буквами $M$, $N$, $P$ и $Q$.
Так как две касательные, проведённые из одной точки к окружности, равны, то 
\begin{align*}
AM&=AQ,& BM&=BN,
\\
CP&=CN, & DP &= DQ.
\end{align*}

Следовательно,
\begin{align*}
AM&+MB+CP+PD = 
\\
&=
AQ + QD+BN+NC,
\end{align*}
то есть 
\[AB+CD=AD+BC.\]

%% file: 2D/zam-toch-trig.tex
\section{Замечательные точки треугольника}

\paragraph{}\label{1938/141}
Мы видели, что:

1) \emph{три срединных перпендикуляра к сторонам треугольника пересекаются в одной точке} (которая есть центр описанного круга). 

2) \emph{три биссектрисы углов треугольника пересекаются в одной точке} (которая есть центр вписанного круга).

Следующие две теоремы указывают ещё две замечательные точки треугольника:

3) точку пересечения трёх высот и 

4) точку пересечения трёх медиан.

\paragraph{}\label{1938/142}
\so{Теорема}.
\textbf{\emph{Три высоты треугольника пересекаются в одной точке.}}

\begin{wrapfigure}{o}{50mm}
\centering
\includegraphics{mppics/ris-160}
\caption{}\label{1938/ris-160}
\end{wrapfigure}

Через каждую вершину $\triangle ABC$ (рис.~\ref{1938/ris-160}) проведём прямую, параллельную противоположной стороне его.
Тогда получим вспомогательный $\triangle A_1B_1C_1$, к сторонам которого высоты данного треугольника перпендикулярны.
Так как $C_1B=AC=BA_1$ (как противоположные стороны параллелограммов), то точка $B$ есть середина стороны $A_1C_1$.

Подобно этому убедимся, что $C$ есть середина $A_1B_1$ и $A$ — середина $B_1C_1$.
Таким образом, высоты $AD$, $BE$ и $CF$ являются
срединными перпендикулярами к сторонам $\triangle A_1B_1C_1$ и как мы знаем (§~\ref{1938/104}), пересекаются в одной точке.

{\small
\smallskip
\mbox{\so{Замечание}.}
Точка, в которой пересекаются высоты треугольника, называется его \rindex{ортоцентр}\textbf{ортоцентром}.
}

\begin{wrapfigure}{r}{37mm}
\vskip-4mm
\centering
\includegraphics{mppics/ris-161}
\caption{}\label{1938/ris-161}
\end{wrapfigure}

\paragraph{}\label{1938/143} 
\mbox{\so{Теорема}.}
\textbf{\emph{Три медианы треугольника пересекаются в одной точке;
эта точка отсекает от каждой медианы третью часть, считая от соответствующей стороны.}}

Возьмём в $\triangle ABC$ (рис.~\ref{1938/ris-161}) какие-нибудь две медианы, например $AE$ и $BD$, пересекающиеся в точке $Z$, и докажем, что
\[ZD=\tfrac13 BD\quad\text{и}\quad ZE = \tfrac13 AE.\]

Для этого, разделив $ZA$ и $ZB$ пополам в точках $F$ и $G$, построим четырёхугольник $DEGF$.
Так как отрезок $FG$ соединяет середины двух сторон $\triangle ABZ$, то $FG\parallel AB$ и $FG=\tfrac12 AB$.
Отрезок $DE$ также соединяет середины двух сторон $\triangle ABC$;
поэтому $DE\parallel AB$ и $DE=\tfrac12 AB$.
Отсюда выводим, что $DE \parallel FG$ и $DE=FG$;
следовательно, четырёхугольник $DEGF$ есть параллелограмм (§~\ref{1938/89}), и потому $ZF=ZE$ и $ZD=ZG$.
Отсюда следует, что
\[ ZE=\tfrac13AE\quad\text{и}\quad ZD=\tfrac13 BD.\]

Если теперь возьмём третью медиану с одной из медиан $AE$ или $BD$, то также убедимся, что точка их пересечения отсекает от каждой из них треть, считая от основания;
значит, и третья медиана должна пересечься с медианами $AE$ и $BD$ в той же точке $Z$.

Из физики известно, что пересечение медиан треугольника есть его \so{центр тяжести};
он всегда лежит внутри треугольника.

{\small

\subsection*{Упражнения}

\begin{center}
\so{Найти геометрические места}
\end{center}

\begin{enumerate}[noitemsep]

\item
Оснований перпендикуляров, опущенных из данной точки $A$ на прямые, проходящие через другую данную точку $B$.

\item
Середин хорд, проведённых в окружности через данную внутри неё точку.

\end{enumerate}

\begin{center}
\so{Доказать теоремы}
\end{center}

\begin{enumerate}[resume,noitemsep]

\item
Если две окружности касаются, то всякая секущая, проведённая через точку касания, отсекает от окружностей две противолежащие дуги одинакового числа градусов.

\item
Отрезки двух равных хорд, пересекающихся в одной окружности, соответственно равны.

\item
Две окружности пересекаются в точках $A$ и $B$, через $A$ проведена секущая, пересекающая окружности в точках $C$ и $D$;
доказать, что угол $CBD$ есть величина постоянная для всякой секущей, проведённой через точку $A$.

\smallskip
\so{Указание}.
Углы $ACB$ и $ADB$ имеют постоянную величину.

\item
Если через точку касания двух окружностей проведём две секущие и концы их соединим хордами, то эти хорды параллельны.

{\sloppy

\item
Если через точку касания двух окружностей проведём внутри них какую-либо секущую, то касательные, проведённые через концы этой секущей, параллельны.

}

\item
Если основания высот остроугольного треугольника соединим прямыми, то получим новый треугольник, для которого высоты первого треугольника служат биссектрисами. 

\item
На окружности, описанной около равностороннего $\triangle ABC$, взята произвольная точка $M$;
доказать, что наибольший из отрезков $MA$, $MB$, $MC$ равен сумме двух остальных.

\item
Из точки $F$ проведены к окружности две касательные $PA$ и $FB$ и через точку $B$ — диаметр $BC$.
Доказать, что прямые $CA$ и $OP$ параллельны ($O$ — центр окружности).

\item
Через одну из точек пересечения двух окружностей проводят диаметр в каждой из них.
Доказать, что прямая, соединяющая концы этих диаметров, проходит через вторую точку пересечения окружностей.

\item
Диаметр $AB$ и хорда $AC$ образуют угол в $30\degree$.
Через $C$ проведена касательная, пересекающая продолжение $AB$ в точке $D$.
Доказать, что $\triangle ACD$ равнобедренный.

\item
Если около треугольника опишем окружность и из произвольной точки её опустим перпендикуляры на стороны треугольника, то их основания лежат на одной прямой (прямая Симпсона).

\smallskip
\so{Указание}.
Доказательство основывается на свойствах вписанных углов (§~\ref{1938/124}) и углов вписанного четырёхугольника (§~\ref{1938/139}).

\end{enumerate}

\begin{center}
\so{Задачи на построение}
\end{center}

\begin{enumerate}[resume,noitemsep]

\item
На данной бесконечной прямой найти точку, из которой данный отрезок был бы виден под данным углом.

\item
Построить треугольник по основанию, углу при вершине и высоте.

{\sloppy

\item
К дуге данного сектора провести такую касательную, чтобы часть её, заключённая между продолженными радиусами (ограничивающими сектор), равнялась данному отрезку (свести эту задачу к предыдущей).

}

\item
Построить треугольник по основанию, углу при вершине и медиане, проведённой к основанию.

\item
Даны по величине и положению два отрезка $a$ и $b$.
Найти такую точку, из которой отрезок $a$ был бы виден под данным углом $\alpha$ и отрезок $b$ под данным углом $\beta$.

\item
В треугольнике найти точку, из которой его стороны были бы видны под равными углами.
(Не всякий треугольник допускает решение.)

\smallskip
\so{Указание}.
Обратить внимание на то, что каждый из этих углов должен равняться $120\degree$.

\item
Построить треугольник по углу при вершине, высоте и медиане, проведённой к основанию.

\smallskip
\so{Указание}.
Продолжив медиану на равное расстояние и соединив полученную точку с концами основания, рассмотреть образовавшийся параллелограмм.

\item
Построить треугольник, в котором даны:
основание, прилежащий к нему угол и угол, составленный медианой, проведённой из вершины данного угла, и стороной, к которой эта медиана проведена.

\item
Построить параллелограмм по двум его диагоналям и одному углу.

\item
Построить треугольник по основанию, углу при вершине и сумме или разности двух других сторон.

\item
Построить четырёхугольник по двум диагоналям, двум соседним сторонам и углу, образованному остальными двумя сторонами.

\item
Даны три точки $A$, $B$ и $C$.
Провести через $A$ такую прямую, чтобы расстояние между перпендикулярами, опущенными на эту прямую из точек $B$ и $C$, равнялось данному отрезку.

\item
В данный круг вписать треугольник, у которого два угла даны.

\item
Около данного круга описать треугольник, у которого два угла даны.

\item
Построить треугольник по радиусу описанного круга, углу при вершине и высоте.

\item
Вписать в данный круг треугольник, у которого известны:
сумма двух сторон и угол, противолежащий одной из этих сторон.

\item
Вписать в данный круг четырёхугольник, у которого даны сторона и два угла, не прилежащие к этой стороне.

\item
В данный ромб вписать круг.

\item
В равносторонний треугольник вписать три круга, которые попарно касаются друг друга и из которых каждый касается двух сторон треугольника.

\item
Построить четырёхугольник, который можно было бы вписать в окружность, по трём его сторонам и одной диагонали.

\item
Построить ромб по данной стороне и радиусу вписанного круга.

\item
Около данного круга описать равнобедренный прямоугольный треугольник.

\item
Построить равнобедренный треугольник по основанию и радиусу вписанного круга.

\item
Построить треугольник по основанию и двум медианам, исходящим из концов основания.
Указание: смотри §~\ref{1938/143}.

\item
То же по трём медианам.
Указание: смотри §~\ref{1938/143}.

\item
Дана окружность и на ней точки $A$, $B$ и $C$.
Вписать в эту окружность такой треугольник, чтобы его биссектрисы при продолжении встречали окружность в точках $A$, $B$ и $C$.

\item
Та же задача, с заменой биссектрис треугольника его высотами.

\item
Дана окружность и на ней три точки $M$, $N$ и $P$, в которых пересекаются с окружностью (при продолжении) высота, биссектриса и медиана, исходящие из одной вершины вписанного треугольника.
Построить этот треугольник.

\item
На окружности даны две точки $A$ и $B$.
Из этих точек провести две параллельные хорды, сумма которых дана.

\end{enumerate}

\begin{center}
\so{Задачи на вычисление}
\end{center}

\begin{enumerate}[resume,noitemsep]

\item
Вычислить вписанный угол, опирающийся на дугу, равную $\tfrac1{12}$ части окружности.

\item
Круг разделён на два сегмента хордой, делящей окружность на части в отношении 5:7.
Вычислить углы, которые вмещаются этими сегментами.

\item
Две хорды пересекаются под углом в $36\degree\, 15'\, 32''$.
Вычислить в градусах, минутах и секундах две дуги, заключённые между сторонами этого угла и их продолжениями, если одна из этих дуг относится к другой как 2:3.

\item
Угол, составленный двумя касательными, проведёнными из одной точки к окружности, равен $25\degree15'$.
Вычислить дуги, заключённые между точками касания.

\item
Вычислить угол, составленный касательной и хордой, если хорда делит окружность на две части, относящиеся как $3:7$.

\item
Две окружности одинакового радиуса пересекаются под углом $60\degree$;
определить в градусах меньшую из дуг, заключающихся между точками пересечения. (Два возможных решения.)

\smallskip
\so{Примечание}.
Углом двух пересекающихся дуг называется любой угол, составленный двумя касательными, проведёнными к этим дугам из точки пересечения.

\item
Из одного конца диаметра проведена касательная, а из другого — секущая, которая с касательной составляет угол в $20\degree 30'$.
Найти меньшую из дуг, заключённых между касательной и секущей.

\end{enumerate}

}

%% file: 2D/izmereniya.tex
\section{Измерение отрезков}

\paragraph{Задача измерения отрезка.}\label{1938/144}
До сих пор, сравнивая между собой два отрезка, мы могли определить, равны ли они между собой, и если не равны, то какой из них больше (§~\ref{1938/6}).
Нам приходилось это делать при изучении соотношений между сторонами и углами треугольника (§§~\ref{1938/46}, \ref{1938/47}), при сравнении отрезка прямой с ломаной (§~\ref{1938/50}, \ref{1938/51}) и в некоторых других случаях (§§~\ref{1938/54}, \ref{1938/55}, \ref{1938/63}).
Но такое сравнение отрезков между собой ещё не даёт точного представления о величине каждого из них.

Наша задача установить точное понятие о длине отрезка и найти способы выражать эту длину при помощи числа.

\paragraph{Понятие об измерении отрезков.}\label{1938/150} 
Чтобы составить ясное представление о величине данного отрезка, его сравнивают с другим, уже известным нам отрезком, например с метром.
Этот известный отрезок, с которым сравнивают другие отрезки, называется \rindex{единица длины}\textbf{единицей длины}.

\begin{figure}[!ht]
\centering
\includegraphics{mppics/ris-168}
\caption{}\label{1938/ris-168}
\end{figure}

Пусть, например, надо измерить отрезок $a$ (рис.~\ref{1938/ris-168}) при помощи единицы $b$.
Тогда поступают так:
положим, что мы желаем найти отрезки, которые отличались бы от $a$ меньше, чем на
$\tfrac1{10}$ единицы длины $b$.
Тогда делим единицу $b$ на 10 равных частей (рис.~\ref{1938/ris-168}) и одну такую долю откладываем на отрезке $a$ столько раз, сколько возможно.
Пусть она уложится 13 раз с некоторым остатком меньшим $\tfrac1{10}b$.
Тогда получим отрезок $a_1=\tfrac{13}{10}b$ и меньший, чем~$a$.
Отложив $\tfrac1{10}b$ ещё один раз, получим другой отрезок, $a_2=\tfrac{14}{10}b$,  больший, чем $a$, который разнится от $a$ менее чем на $\tfrac1{10}$ единицы.
Длины отрезков $a_1$ и $a_2$ выражаются числами $\tfrac{13}{10}$ и $\tfrac{14}{10}$.
Эти числа рассматриваются как \so{приближённые меры} длины отрезка $a$:
первое с недостатком, второе — с избытком.
При этом, так как отрезок $a$ разнится от $a_1$ и от $a_2$ менее чем на $\tfrac1{10}$ единицы, то принято говорить, что каждое из этих чисел выражает длину отрезка $a$ с точностью до $\tfrac1{10}$.

Вообще, чтобы найти приближённые меры длины отрезка $a$ с точностью до $\tfrac1n$ единицы, делят единицу $b$ на $n$ равных частей и узнают, сколько раз $\tfrac1n$-я доля единицы содержится в $a$;
если она содержится $m$ раз с некоторым остатком, меньшим $\tfrac1n b$, то числа $\tfrac mn$ и $\tfrac {m+1}n$ считаются приближёнными мерами длины отрезка $a$ с точностью до $\tfrac1n$-й, первое с недостатком, второе — с избытком.

Может случиться, что этим путём мы найдём точный результат, то есть если отрезок $\tfrac1n b$ уложится целое число раз в отрезке $a$.

Для получения того числа, которое можно было бы принять за точную меру длины отрезка $a$, поступают следующим образом.

Вычисляют последовательно приближённую меру длины отрезка $a$ с недостатком с точностью до $0{,}1$, затем ту же меру с недостатком с точностью до $0{,}01$, затем её же с точностью до $0{,}001$ и продолжают беспредельно этот процесс последовательного вычисления приближённой меры длины $a$, каждый раз повышая точность в 10 раз.
При таком процессе будут получаться последовательно десятичные дроби сначала с одним десятичным знаком, затем с двумя, тремя и дальше всё с б\'{о}льшим и б\'{о}льшим числом десятичных знаков.
Неограниченное продолжение описанного процесса построения десятичных дробей определяет бесконечную десятичную дробь.

Бесконечную десятичную дробь нельзя, конечно, полностью записать на листе бумаги, так как число её десятичных знаков бесконечно.
Тем не менее её считают известной, если известен способ, при помощи которого можно определить любое число её десятичных знаков.

\paragraph{Бесконечные десятичные дроби.}\label{1938/151}
Введение бесконечных десятичных дробей производится в алгебре на основе следующих определений.

1) Бесконечная десятичная дробь называется вещественным числом.

2) Две бесконечные десятичные дроби считаются равными, если их десятичные знаки одинакового порядка равны.

3) Из двух неравных бесконечных десятичных дробей считается б\'{о}льшим вещественным числом та дробь, в которой первый из неравных десятичных знаков одинакового порядка со второй дробью больше.

4) Если в бесконечной десятичной дроби все десятичные знаки, начиная с некоторого порядка, равны нулю, то дробь считается равной той конечной десятичной дроби, которая получится из данной зачёркиванием всех нулей, стоящих справа от последней значащей цифры.
Так, бесконечная десятичная дробь $7{,}8530078000\dots$
равна конечной дроби $7{,}8530078$.

5) Бесконечная периодическая дробь с периодом 9 всегда заменяется конечной десятичной дробью, получаемой из данной увеличением на единицу её последнего десятичного знака, отличного от $9$, и отбрасыванием всех последующих девяток.
Так, дробь $3{,}72999\dots$ заменяют конечной дробью $3{,}73$.

\paragraph{Приближённые значения бесконечной десятичной дроби.}\label{1938/152}
Если оборвать данную бесконечную десятичную дробь на её $n$-м знаке, то полученная конечная дробь называется приближённым значением бесконечной десятичной дроби с точностью до $\tfrac1{10^n}$ с недостатком.
Если же в этой дроби увеличить на единицу её последний десятичный знак, то есть
прибавить к ней $\tfrac1{10^n}$, то получится новая конечная дробь, которая называется приближённым значением бесконечной дроби с той же точностью с избытком.
Если приближённое значение вещественного числа $\alpha$ с $n$ десятичными знаками с недостатком обозначим через $\alpha_n$, а с избытком через  $\alpha_n'$, то  $\alpha_n'=\alpha_n+\tfrac1{10^n}$.
Из определения неравенства вещественных чисел следует, что 
\[\alpha_n\le \alpha<\alpha_n';\]
то есть приближённое значение взятое с недостатком не превосходит само число,
а само число меньше своего приближённого значения с избытком.

Пусть, например, дано вещественное число, определяющее  $\sqrt{2}  \z= 1{,}414\dots$;
его приближённое значение с точностью до $0{,}01$ с недостатком:
$1{,}41$, с избытком:
$1{,}42$;
так как
$1{,}41 = 1{,}41000$
и
$1{,}42 = 1{,}42000$,
то в силу определения неравенства вещественных чисел имеем:
\[
1{,}41000\ldots
< 1{,}414\ldots
< 1{,}42000\ldots,
\qquad\text{или}\qquad
1{,}41 <  \sqrt{2}  < 1{,}42.
\]

\paragraph{Сложение вещественных чисел.}\label{1938/153}

Пусть даны два вещественных числа $\alpha$ и $\beta$.
Возьмём их приближённые значения с произвольным числом $n$ десятичных знаков, сначала с недостатком, а затем с избытком.
Приближённые значения чисел $\alpha$ и $\beta$ с недостатком обозначим соответственно через $\alpha_n$ и $\beta_n$, а приближённые значения с избытком — через $\alpha_n'$ и $\beta_n'$.
При этом:
\[\alpha_n'=\alpha_n +\tfrac1{10^n},
\quad
 \beta_n'=\beta_n +\tfrac1{10^n}.\eqno(1)
\]
Составим теперь суммы $\alpha_n+\beta_n$ и $\alpha_n'+ \beta_n'$.
Каждая из них есть десятичная дробь, содержащая $n$ десятичных знаков.

Назовём первую $\gamma_n$, а вторую $\gamma_n'$:
\[\alpha_n+\beta_n=\gamma_n,\quad\alpha_n'+\beta_n'=\gamma_n'.\]
Складывая почленно равенства (1), получим:
\[\alpha_n'+\beta_n'= \alpha_n + \beta_n + \tfrac2{10^n},\]
или $\gamma_n'=\gamma_n+ \tfrac2{10^n}$.
Это равенство показывает, что дробь $\gamma_n$ получается из
дроби $\gamma_n$ прибавлением двух единиц к её последнему десятичному знаку.

Будем теперь увеличивать $n$;
в таком случае дробь $\gamma_n'$ приведёт к образованию бесконечной десятичной дроби, которую обозначим~$\gamma$.
Эта дробь может оказаться или периодической, или непериодической.

Допустим, что дробь $\gamma$ непериодическая.
В таком случае она должна содержать бесчисленное множество десятичных знаков, отличных от 9.
В этом случае в дроби $\gamma$ число десятичных знаков, отличных от 9, должно возрастать с возрастанием~$n$.
Так как прибавка в дроби $\gamma$ числа $\tfrac2{10^n}$ не может оказать влияния на её десятичные знаки, стоящие левее двух последних знаков, отличных от 9, то число общих первых десятичных знаков в дробях $\gamma_n$ и $\gamma_n'$ будет неограниченно возрастать с возрастанием~$n$.
Следовательно, дробь $\gamma_n'$ будет приводить к той же бесконечной десятичной дроби, что и дробь $\gamma_n$.
При этом из предыдущего следует, что при любом~$n$
\[\gamma_n\le \gamma<\gamma_n'.\eqno(2)\]

Допустим теперь, что дробь $\gamma$ периодическая.
В таком случае она представляет собой некоторое рациональное число.
Это число и только оно удовлетворяет неравенству (2) при всех~$n$. 

\smallskip
\so{Определение}.
\emph{Вещественное число $\gamma$, удовлетворяющее неравенствам (2) при всех $n$, называется суммой вещественных чисел $\alpha$ и $\beta$.}
\[\gamma=\alpha+\beta.\]

\paragraph{Другие действия с вещественными числами.}\label{1938/154}
Совершенно аналогичным образом можно определить разность двух вещественных чисел, их произведение и частное от деления одного вещественного числа на другое.
Более подробное изучение результатов этих действий показывает, что определённые таким образом сумма и произведение вещественных чисел подчиняются основным законам действий, имеющим место для чисел рациональных:
сложение подчиняется переместительному и сочетательному законам.

\[\alpha+\beta=\beta+\alpha,
\quad
(\alpha+\beta)+\gamma=\alpha+(\beta+\gamma),
\]
а умножение — переместительному, сочетательному и распределительному законам.
\[\alpha\beta=\beta\alpha,
\quad
(\alpha\beta)\gamma=\alpha(\beta\gamma),
\quad
(\alpha+\beta)\gamma=\alpha\gamma+\beta\gamma.
\]

В тех случаях, когда бесконечные десятичные дроби будут периодическими, определённые выше действия над ними будут приводить, как легко показать, к тем же результатам, что и действия над обыкновенными дробями, получаемыми после обращения периодических дробей в простые.

Таким образом, рациональные числа являются лишь частным видом вещественных чисел.
 
\paragraph{Длины отрезков и их отношения.}\label{1938/155}
Число, получаемое в результате измерения отрезка $a$, называется \rindex{длина}\textbf{длиной} этого отрезка.

Заметим, что равенство длин отрезков и равенство отрезков, определённое нами с помощью наложения (§~\ref{1938/6}), эквивалентны, конечно, если для измерения отрезков мы пользовались одной единицей длины.
То же верно и для сравнения, сложения и других действий над отрезками и их длинами (§§~\ref{1938/6}—\ref{1938/8}).

Например неравенство $AB>2\cdot CD$ может пониматься двояко — то, что отрезок $CD$ укладывается два раза в отрезке $AB$ с некоторым остатком и то, что длина отрезка $AB$ больше чем удвоенная длина отрезка $CD$, измеренная той же единицей длины.
При этом выражение «измеренная той же единицей длины» мы будем опускать, предполагая, что в каждой конкретной задаче измерения производятся только одной единицей.

Под отношением двух отрезков мы понимаем отношение их длин. 
Это же отношение равно длине первого, если второй взять за единицу длины.

Заметим, что отношение двух отрезков не зависит от того, как выбрана единица измерения.
В самом деле, если, например, вместо одной уже выбранной единицы измерения взять другую, в 3 раза меньшую, то в каждом отрезке эта новая единица уложится втрое большее число раз, чем прежняя.
В той дроби, которая представляет отношение отрезков, числитель и знаменатель оба увеличатся в 3 раза.
Величина же самой дроби от этого не изменится.

\paragraph{Пропорции.}\label{extra/proportions}
В геометрических задачах часто появляется уравнение типа  
\[\frac{a}{b}=\frac{c}{d}\eqno(1)\]
на длины $a$, $b$, $c$ и $d$ некоторых отрезков.
Такое уравнение называется \rindex{пропорция}\textbf{пропорцией}. 

Следующие наблюдения часто оказываются полезными:

1) Пропорцию (1) можно переписать как произведение:
\begin{align*}
a\cdot d&=b\cdot c.
\end{align*}

2) Пропорцию (1) можно продолжить, складывая или вычитая соответствующие члены:
\[\frac{a}{b}=\frac{c}{d}=\frac{a+c}{b+d}=\frac{a-c}{b-d}=\frac{2a+3c}{2b+3d}=\dots\]
конечно, если знаменатели в новых дробях не равны нулю.

3) По пропорции (1) можно написать другие пропорции, например:
\begin{align*}
\frac{a}{c}&=\frac{b}{d},
&
\frac{a+b}{b}&=\frac{d+c}{d},
&
\frac{a}{b-a}&=\frac{c}{d-c}\quad \text{и так далее.}
\end{align*}

Приведём пример использования этих наблюдений.

\smallskip
\so{Задача 1}. Предположим для точек $C$ и $C'$, лежащих на отрезке $AB$, выполняется пропорция
\[\frac{AC}{CB}=\frac{AC'}{C'B}.\eqno(2)\]
Доказать, что $C=C'$.

Из пропорции (2) можно написать другую
\[\frac{AC}{AC+CB}=\frac{AC'}{AC'+C'B}.\]
Поскольку точки $C$ и $C'$ лежат на отрезке $AB$, 
\[AC+CB=AC'+C'B=AB.\]
Значит
\[\frac{AC}{AB}=\frac{AC'}{AB};\]
следовательно $AC=AC'$ и $C=C'$.

Аналогично решается следующая задача:

\smallskip
\so{Задача 2}. Предположим для точек $C$ и $C'$, лежащих на продолжении отрезка $AB$, выполняется пропорция
\[\frac{AC}{CB}=\frac{AC'}{C'B}.\eqno(3)\]
Доказать, что $C=C'$.

Не умаляя общности можно предположить, что $AC>CB$;
тогда из пропорции (3) следует, что $AC'>C'B$,
то есть обе точки $C$ и $C'$ лежат на продолжении $AB$ за точку $B$ и, значит,
\[AC-CB=AC'-C'B=AB.\]
Из пропорции (3) можно написать другую
\[\frac{AC}{AC-CB}=\frac{AC'}{AC'-C'B}\]
или
\[\frac{AC}{AB}=\frac{AC'}{AB};\]
следовательно $AC=AC'$ и $C=C'$.

%% file: 2D/nesoimerimost.tex
\section{Соизмеримые и несоизмеримые отрезки}

\begin{wrapfigure}{r}{37mm}
\centering
\includegraphics{mppics/ris-162}
\caption{}\label{1938/ris-162}
\end{wrapfigure}

\paragraph{Общая мера.}\label{1938/145}\rindex{общая мера}
Общей мерой двух отрезков называется такой третий отрезок, который в каждом из первых двух содержится целое число раз без остатка.
Так, если отрезок $AM$ (рис.~\ref{1938/ris-162}) содержится 5 раз в $AB$ и 3 раза в $CD$, то $AM$ есть общая мера $AB$ и $CD$.
Подобно этому можно говорить об общей мере двух дуг одинакового радиуса, двух углов и вообще двух однородных величин.

{\small
\smallskip
\so{Замечание}.
Очевидно, что если отрезок $AM$ есть общая мера
отрезков $AB$ и $CD$, то, разделив $AM$ на 2, 3, 4 и так далее равные
части, мы получим меньшие общие меры для отрезков $AB$ и $CD$.
Таким образом, если два отрезка имеют какую-нибудь общую меру, то можно сказать, что они имеют бесчисленное множество общих мер.
Одна из них будет наибольшая.
}

\paragraph{Теоремы о наибольшей общей мере.}\label{1938/146}
Чтобы найти наибольшую общую меру двух отрезков, употребляют так называемый \rindex{алгоритм Евклида}\textbf{алгоритм Евклида} — способ последовательного отложения, подобный тому последовательному делению, каким в арифметике находят наибольший общий делитель двух целых чисел. 

Этот способ основывается на следующих теоремах.

\begin{wrapfigure}{o}{24mm}
\centering
\includegraphics{mppics/ris-163}
\caption{}\label{1938/ris-163}
\end{wrapfigure}

1.
\textbf{\emph{Если меньший из двух отрезков}} ($a$ и $b$, рис.~\ref{1938/ris-163}) \textbf{\emph{содержится в большем целое число раз без остатка, то наибольшая общая мера этих отрезков есть меньший из них.}}

Пусть, например, отрезок $b$ содержится в отрезке $a$ ровно $3$ раза;
так как при этом, конечно, отрезок $b$ содержится в самом себе ровно $1$ раз, то $b$ есть общая мера отрезков $a$ и $b$;
с другой стороны, эта мера есть и наибольшая, так как никакой отрезок, больший $b$, не может содержаться в $b$ целое число раз.

2.
\textbf{\emph{Если меньший из двух отрезков}} ($b$, рис.~\ref{1938/ris-164}) \textbf{\emph{содержится в большем}} ($a$) \textbf{\emph{целое число раз с некоторым остатком}} ($r$), \textbf{\emph{то наибольшая общая мера этих отрезков}} (если она существует) \textbf{\emph{должна быть и наибольшей общей мерой меньшего отрезка}} ($b$) \textbf{\emph{и остатка}} ($r$).

\begin{wrapfigure}{O}{30mm}
\centering
\includegraphics{mppics/ris-164}
\caption{}\label{1938/ris-164}
\end{wrapfigure}

Пусть, например.
\[a=b+b+b+r.\]
Из этого равенства мы можем вывести следующие два заключения.

1) Если существует отрезок, содержащийся без остатка в $b$ и $r$, то он содержится также без остатка и в $a$;
если, например, какой-нибудь отрезок содержится в $b$ ровно 5 раз и в $r$ содержится ровно 2 раза, то в $a$ он содержится 5 + 5 + 5 + 2, то есть 17 раз без остатка.

2) Обратно:
если существует отрезок, содержащийся без остатка в $a$ и $b$, то он содержится также без остатка и в $r$;
если, например, какой-нибудь отрезок содержится в $a$ ровно 17 раз и в $b$ ровно 5 раз, то в той части отрезка $a$, которая равна $3b$, он содержится 15 раз;
следовательно, в остающейся части отрезка $a$, то есть в $r$, он содержится $17-15$, то есть 2 раза.

Таким образом, у двух пар отрезков
\[\overbrace{a\ \text{и}\ b},\quad \overbrace{b\ \text{и}\ r}\]
должны быть одни и те же общие меры (если они существуют);
поэтому и \so{наибольшая} общая мера у них должна быть одна и та же.
К этим двум теоремам надо ещё добавить следующую \rindex{аксиома!измерения}\textbf{аксиому измерения} (так называемую \rindex{аксиома!Архимеда}аксиому Архимеда). 

\textbf{\emph{Как бы велик ни был больший отрезок}} ($a$) \textbf{\emph{и как бы мал ни был меньший отрезок}} ($b$), \textbf{\emph{всегда, откладывая меньший на большем последовательно 1, 2, 3 и так далее раз, мы получим, что после некоторого $\bm{m}$-го отложения или не получится никакого остатка, или получится остаток, меньший меньшего отрезка}} ($b$);
другими словами, всегда можно найти достаточно большое целое положительное число $m$, что $b \cdot  m \le a$, но $b \cdot  (m + 1) > a$. 

\paragraph{Алгоритм Евклида.}\label{1938/147}
Пусть требуется найти наибольшую общую меру двух данных отрезков $AB$ и $CD$ (рис.~\ref{1938/ris-165}).

\begin{wrapfigure}{o}{48mm}
\centering
\includegraphics{mppics/ris-165}
\caption{}\label{1938/ris-165}
\end{wrapfigure}

Для этого на большем отрезке откладываем (с помощью циркуля) меньший отрезок столько раз, сколько это возможно.
При этом, согласно аксиоме измерения, случится одно из двух:
или 
1) $CD$ уложится в $AB$ без остатка, тогда искомая мера, согласно теореме 1, будет $CD$, 
или 2) получится некоторый остаток $EB$, меньший $CD$ (как у нас на чертеже);
тогда, согласно теореме 2, вопрос приведётся к нахождению наибольшей общей меры двух меньших отрезков, именно $CD$ и первого остатка $EB$.
Чтобы найти её, поступаем точно так же, то есть откладываем $EB$ на $CD$ столько раз, сколько можно.
И опять произойдёт одно из двух:
или 1) $EB$ уложится в $CD$ без остатка, тогда искомая мера и будет $EB$, 
или 2) получится остаток $FD$, меньший $EB$ (как у нас на чертеже);
тогда вопрос приведётся к нахождению наибольшей общей меры двух меньших отрезков, именно $EB$ и второго остатка $FD$.

Продолжая этот приём далее, мы можем встретиться с такими двумя возможными случаями:

1) после некоторого отложения не получится никакого остатка или

2) процесс последовательного отложения не будет иметь конца (в предположении, что мы имеем возможность откладывать какие угодно малые отрезки, что, конечно, возможно только теоретически).

В первом случае последний остаток и будет наибольшей общей мерой данных отрезков.
Чтобы удобней вычислить, сколько раз эта наибольшая общая мера содержится в данных отрезках, выписываем ряд равенств, получаемых после каждого отложения.
Так, по нашему чертежу мы будем иметь:
\begin{align*}
&\text{после}
&&\text{первого}
&\text{отложения}&:
&AB &= 3\cdot CD + EB;
\\
&\text{\ —\textquotedbl—}
&&\text{второго}
&\text{—\textquotedbl—\ \ \ \ }&:
&CD &= 2\cdot EB + FD;
\\
&\text{\ —\textquotedbl—}
&&\text{третьего}
&\text{—\textquotedbl—\ \ \ \ }&:
&EB &= 4\cdot FD.
\end{align*}
Переходя в этих равенствах от нижнего к верхнему, последовательно находим:
\begin{align*}
EB&=4\cdot FD;
\\
CD&=(4\cdot FD)\cdot 2+FD=9\cdot FD;
\\
AB&=(9\cdot FD)\cdot 3+4\cdot FD=31\cdot FD.
\end{align*}
Подобно этому можно находить наибольшую общую меру двух дуг одинакового радиуса, а также двух углов.

\smallskip
\so{Во втором} случае данные отрезки не могут иметь общей меры.
Чтобы обнаружить это, предположим, что данные отрезки $AB$ и $CD$ имеют какую-нибудь общую меру.
Мера эта, как мы видели, должна содержаться целое число раз не только в $AB$ и в $CD$, но и в остатке $EB$, следовательно, и во втором остатке $FD$, и в третьем, и в четвёртом и~т.~д.
Так как остатки эти идут, последовательно уменьшаясь, то в каждом из них общая мера должна содержаться меньшее число раз, чем в предыдущем остатке.
Если, например, в $EB$ общая мера содержится $100$ раз (вообще $m$ раз), то в $FD$ она содержится менее $100$ раз (значит, не более $99$ раз);
в следующем остатке она должна содержаться менее $99$ раз (значит, не более $98$ раз) и~т.~д.
Так как ряд целых положительных уменьшающихся чисел:
$100, 99, 98, \dots$
(и вообще $m, m-1, m-2,\dots$) имеет конец (как бы велико ни было число $m$), то и процесс последовательного отложения, при достаточном его продолжении, должен дойти до конца, то есть мы дойдём до того, что уже не получится никакого остатка.
Значит, если последовательное отложение не имеет конца, то данные отрезки никакой общей меры иметь не могут.

\paragraph{Соизмеримые и несоизмеримые отрезки.}\label{1938/148}\rindex{соизмеримые отрезки}\rindex{несоизмеримые отрезки}
Два отрезка называются соизмеримыми, если они имеют общую меру, и несоизмеримыми, когда такой общей меры не существует.

На практике нет возможности убедиться в существовании несоизмеримых отрезков, потому что, продолжая последовательное отложение, мы всегда дойдём до столь малого остатка, который в предыдущем остатке, по-видимому, укладывается целое число раз.
Быть может при этом и должен был бы получиться некоторый остаток, но по причине неточности инструментов (циркуля) и несовершенства наших органов чувств (зрения) мы не в состоянии его заметить.
Однако, как мы сейчас докажем, несоизмеримые отрезки существуют.

{\small

\paragraph{}\label{1914/156}
\so{Теорема}.
\textbf{\emph{Если в равнобедренном треугольнике угол при основании равен $\bm{36\degree}$, то боковая сторона его несоизмерима с основанием.}}

Пусть $ABC$ равнобедренный треугольник (рис.~\ref{1914/ris-147}), у которого каждый из углов $A$ и $C$ равен $36\degree$; 
требуется доказать, что боковая сторона $AB$ несоизмерима с основанием $AC$.

\begin{wrapfigure}{o}{48mm}
\centering
\includegraphics{mppics/ris-1914-147}
\caption{}\label{1914/ris-147}
\end{wrapfigure}

Прежде всего определим, которая из этих сторон больше.
Для этого достаточно сравнить углы, против которых лежат эти стороны.
Так как, по условию, $\angle A=\angle C\z=36\degree$, то 
\[\angle B=180\degree-36\degree-36\degree=108\degree;\]
следовательно, $\angle B>\angle C$;
поэтому $AC\z>AB$.

Теперь найдём сколько раз в $AC$ уложится $AB$.
Так как $AC\z<AB+BC$ и $AB=BC$ то $AC<2\cdot AB$;
значит $AB$ может уложиться только один раз с некоторым остатком.

Таким образом, мы замечаем следующее свойство:
\emph{если в равнобедренном треугольнике угол при основании равен $36\degree$, то боковая его сторона содержится в основании только один раз и притом с некоторым остатком.}

Заметив это, приступим к последовательному отложению.
Отложим на $AC$ часть $AD$, равную $AB$ тогда получим остаток $DC$, который надо накладывать на $AB$, или, что все равно, на $BC$.

Чтобы узнать, сколько раз $DC$ уложится на $BC$, соединим $B$ с $D$ и рассмотрим $\triangle DBC$.
Найдём его углы.
Так как $\triangle ABD$ равнобедренный, то $\angle ABD = \angle ADB$;
следовательно, каждый из них равен 
\[\tfrac12(180\degree -\angle A)=\tfrac12(180\degree -36\degree)=72\degree.\]
Но угол $ABC$, как мы прежде нашли, равен $108\degree$; следовательно, 
\[\angle DBC=108\degree-72\degree=36\degree.\]
Таким образом, в треугольнике $DBC$ есть два равных угла при $BC$; следовательно, он равнобедренный, причём каждый угол при его основании $BC$ равен $36\degree$.

Вследствие этого, по доказанному выше, боковая сторона его $DC$
(или $BD$) уложится в основании $BC$ \emph{один раз с некоторым остатком}.
Пусть этот остаток будет $EC$.
Соединив $E$ с $D$, мы снова получим равнобедренный треугольник $CDE$, в котором каждый угол при основании $CD$ равен $36\degree$.
Отложив $EC$ (или $DE$) на $DC$ (от точки $D$), мы снова получим равнобедренный треугольник $CEF$, у которого каждый угол при основании $CE$ равен $36\degree$.

Таким образом, мы постоянно будем приходить к равнобедренному треугольнику (всё меньшему и меньшему) с углами при основании, равными $36\degree$;
следовательно, мы никогда в этом процессе \so{не дойдём до конца}.
Значит, стороны $AC$ и $AB$ не могут иметь общей меры.

{\small

\smallskip
{\so{Примечание}.}
Первый пример несоизмеримых отрезков приписывается пифагорейцу Гиппасу из Метапонта (приблизительно 470 год до нашей эры).
Есть основания полагать, что Гиппас доказал несоизмеримость стороны и диагонали правильного пятиугольника, хотя сколько-нибудь уверенных подтверждений тому нет.
Нетрудно видеть, что треугольник, отсекаемый диагональю от правильного пятиугольника, равнобедренный и имеет углы при основании равные $36\degree$.
То есть возможно приведённое построение повторяет то, что было сделано Гиппасом.

}

\paragraph{}\label{1938/149}
\so{Теорема}.
\textbf{\emph{Диагональ квадрата несоизмерима с его стороной.}}

\begin{wrapfigure}{o}{42mm}
\centering
\includegraphics{mppics/ris-166}
\caption{}\label{1938/ris-166}
\end{wrapfigure}

Так как диагональ делит квадрат на два равнобедренных прямоугольных треугольника, то теорему эту можно высказать иными словами так:
\textbf{\emph{гипотенуза равнобедренного прямоугольного треугольника несоизмерима с его катетом.}}

Предварительно докажем следующее свойство такого треугольника;
если на гипотенузе (рис.~\ref{1938/ris-166}) отложим отрезок $AD$, равный катету, и проведём $DE\perp AC$, то образовавшийся при этом прямоугольный треугольник $DEC$ будет равнобедренный, а отрезок $BE$ катета $BC$ окажется равным отрезку $DC$ гипотенузы.
Чтобы убедиться в этом, проведём прямую $BD$ и рассмотрим углы треугольников $DEC$ и $BED$.
Так как треугольник $ABC$ равнобедренный и прямоугольный, то $\angle 1 \z= \angle 4$, и, следовательно, $\angle 1 \z= 45\degree$, а потому в прямоугольном треугольнике $DEC$ и $\angle 2 = 45\degree$ и, значит, треугольник $DEC$ имеет два равных угла и потому его стороны $DC$ и $DE$ равны.

В треугольнике $BDE$ угол 3 равен прямому углу $B$ без угла $ABD$, а угол 5 равен прямому углу $ADE$ без угла $ADB$.
Но $\angle ADB \z= \angle ABD$ (так как $AB=AD$);
значит, и $\angle 3 = \angle 5$.
Но тогда треугольник $BDE$ должен быть равнобедренный, и потому $BE=ED=DC$.

Заметив это, станем находить общую меру отрезков $AB$ и $AC$.

Так как $AC>AB$ и $AC<AB+BC$, то есть
$AC<2AB$, то катет $AB$ отложится на гипотенузе $AC$ только один раз с некоторым остатком $DC$.
Теперь надо этот остаток откладывать на $AB$ или, что всё равно, на $BC$.
Но отрезок $BE$, по доказанному, равен $DC$.
Значит, надо $DC$ отложить ещё на $EC$.
Но $EC$ есть гипотенуза равнобедренного треугольника $DEC$.
Следовательно, процесс отложения для нахождения общей меры сводится теперь к откладыванию катета $DC$ прямоугольного равнобедренного треугольника $DEC$ на его гипотенузе $EC$.
В свою очередь это отложение сведётся к откладыванию катета на гипотенузе нового меньшего прямоугольного равнобедренного треугольника и~т.~д., очевидно, без конца.
А если процесс этот не может окончиться, то общей меры отрезков $AC$ и $AB$ не существует.

\paragraph{Цепные дроби.}\label{extra/tzepnye-drobi}
В дополнение к десятичным дробям мы рассмотрим другой способ, который позволяет записать число.
Он основан на последовательности частных, полученных в алгоритме Евклида (§~\ref{1938/147}).
Этот способ записи неудобен для складывания и умножения чисел, но очень удобен для отыскания приближений числа в виде дробей с небольшими знаменателями.

{

\begin{wrapfigure}{r}{48mm}
\vskip0mm
\centering
\includegraphics{mppics/ris-165}
\caption{}\label{1938/ris-165-1}
\bigskip
\includegraphics{mppics/ris-extra-7}
\caption{}\label{extra/ris-7}
\bigskip
\includegraphics{mppics/ris-extra-8}
\caption{}\label{extra/ris-8} 
\end{wrapfigure}

Вычисление отношения отрезков $\frac{AB}{CD}\z=\frac{31}{9}$ (рис.~\ref{1938/ris-165-1}) данное в §~\ref{1938/147},
можно компактно записать следующим образом:
\begin{align*}
\frac{AB}{CD}&=3+\frac{EB}{CD}=3+\frac{1}{\frac{CD}{EB}}=
\\
&=3+\frac{1}{2+\frac{FD}{EB}}=3+\frac{1}{2+\frac{1}{\frac{EB}{FD}}}
\\
&=3+\frac{1}{2+\frac{1}{4}}.
\end{align*}

Последнее выражение называется \rindex{цепная дробь}\textbf{цепной дробью}.
В этом случае процесс построения цепной дроби оборвался на третьем шаге, так как отрезок $FD$ вмещается 4 раза в $EB$ без остатка.
Вообще, любое рациональное число представляется конечной цепной дробью.

Для несоизмеримых отрезков этот процесс можно было бы продолжать сколь угодно долго, получая всё более длинные дроби.

Например для пары отрезков $AC$ и $AB$ в §~\ref{1914/156} мы получим бесконечную цепную дробь, составленную из одних единиц:
\begin{align*}
\frac{AC}{AB}&=1+\frac{DC}{AB}=1+\frac1{\frac{AB}{DC}}=
\\
&=1+\frac1{1+\frac{EC}{DC}}=1+\frac1{1+\frac1{\frac{DC}{EC}}}=
\\
&=\cdots=1+\frac1{1+\frac1{1+\cdots}}
\end{align*}

}

Из доказательства в §~\ref{1938/149} следует, что для отношения диагонали квадрата к его стороне получается цепная дробь из двоек:
\[\alpha=1+\frac1{2+\frac1{2+\cdots}}.\eqno(1)\]
Из теоремы Пифагора доказанной ниже (§~\ref{1938/191}) следует, что $\alpha\z=\sqrt2$;
но это же можно увидеть непосредственно:
в знаменателе выражения (1) стоит его копия плюс один; то есть $\alpha=1+\frac1{1+\alpha}$ и, значит, $\alpha^2=2$.

Для нахождения приближений числа с небольшим знаменателем,
достаточно оборвать цепную дробь и привести её к обычному виду.
Такие дроби называются \rindex{подходящая дробь}\textbf{подходящими}.
Например для $\sqrt2$ мы получаем следующие приближения:
\[1,\quad \tfrac 32=1+\tfrac12,\quad \tfrac75=1+\frac1{2+\frac12},\quad \tfrac{17}{12}=1+\frac1{2+\frac1{2+\frac12}},\quad \dots\]
Несложно проверить, что нечётные подходящие дроби дают приближение с недостатком, а чётные с избытком.
При этом оказывается, что подходящие дроби являются \so{лучшими} приближениями среди всех дробей с меньшими знаменателями.

\paragraph{О теории пропорций.}\label{extra/evdox}
В современной математике уравнение $\frac{a}{b}\z=\frac{c}{d}$
для длин $a$, $b$, $c$ и $d$ четырёх отрезков  понимается как равенство двух вещественных чисел; каждое из них — это отношение пары вещественных чисел.

Во времена Евклида, под числом понимались только положительные рациональные числа.
Длины же отрезков и промежутки времени назывались \rindex{величина}\textbf{величинами} и числом не считались.
Как мы знаем, 
если выбрана единица измерения (например секунда для времени или метр для длины), то величину можно выразить вещественным числом; но во времена Евклида это не было известно.

Величины одного рода можно было сравнивать, складывать и умножать на натуральное число. 
Чтобы делить или умножать одну величину на другую пользовались специальными приёмами;
например про произведение двух длин думали как про площадь прямоугольника с данными сторонами.
Была также развита теория пропорций по сути позволяющая делить одну величину на другую.

Для того чтобы выразить уравнение $\frac{a}{b}=\frac{c}{d}$, говорили «$a$ относится к $b$ также, как $c$ относится к $d$».
Изначально это выражение означало, что цепная дробь, получаемая при делении $a$ на $b$ в точности совпадает с цепной дробью при делении $c$ на $d$.
Иначе говоря, если применить алгоритм Евклида (§~\ref{1938/147}) к паре отрезков $a$ и $b$, то полученные частные будут в точности те же, что и для пары $c$ и $d$.

Позже Евдокс (живший в IV веке до нашей эры) дал следующее, более удобное определение, равносильное следующему: \emph{$a$ относится к $b$ также как $c$ относится к $d$, если для любых натуральных чисел $m$ и $n$ верно одно из трёх утверждений:
\begin{align*}
\text{либо}\quad m\cdot a&>n\cdot b\quad\text{и}\quad m\cdot c>n\cdot d,
\\
\text{либо}\quad m\cdot a&=n\cdot b\quad\text{и}\quad m\cdot c=n\cdot d,
\\
\text{либо}\quad m\cdot a&<n\cdot b\quad\text{и}\quad m\cdot c<n\cdot d.
\end{align*}
}
Это определение было включено в «Начала» Евклида.

Позже определение Евдокса было использовано немецким математиком Рихардом Дедекиндом для того чтобы дать точное определение вещественных чисел;
сегодня его определение считается основным.
Поскольку учащиеся знакомы с десятичными дробями, мы использовали их как определение вещественных чисел, но строгий вывод всех необходимых свойств десятичных дробей оказывается несоизмеримо сложней, чем при подходе  Евдокса и Дедекинда.

}

%% file: 2D/podobie-trig.tex
\section{Подобие треугольников}

\paragraph{Предварительные понятия.}\label{1938/156}
В окружающей нас жизни часто встречаются фигуры, имеющие различные размеры, но одинаковую форму.
Таковы, например, одинаковые фотографии одного и того же лица, изготовленные в различных размерах, или планы здания или целого города, вычерченные в различных масштабах. 
Такие фигуры принято называть подобными.
Умение измерять длины отрезков позволяет точно определить понятие о геометрическом подобии фигур и дать способы изменения размера фигуры без изменений её формы.

Изучение подобия фигур мы начнём с простейшего случая, именно с подобия треугольников.

\paragraph{}\label{1938/158}
\so{Определение}.
\emph{Два треугольника называются \rindex{подобные!треугольники}подобными, если:
1) углы одного соответственно равны углам другого и 
2) стороны одного пропорциональны соответственным сторонам другого.}

\begin{wrapfigure}{r}{51mm}
\centering
\includegraphics{mppics/ris-extra-4}
\caption{}\label{extra/ris-4}
\end{wrapfigure}

Подобие треугольников обозначается знаком $\sim$;
например,
\[\triangle ABC\z\sim \triangle DEF\]
означает, что треугольники $ABC$ и $DEF$ подобны.
При этом принято выписывать соответственные вершины треугольников в том же порядке;
то есть  $\triangle ABC\z\sim \triangle DEF$ обычно означает, что 
\[\angle A=\angle D,\quad
 \angle B=\angle E,\quad
 \angle C=\angle F
\]
и
\[\frac{AB}{DE}=\frac{BC}{EF}=\frac{AC}{DF}.\]

То, что подобные треугольники существуют, показывает следующая \rindex{лемма}лемма\footnote{Леммой называется вспомогательная теорема, которая излагается для того, чтобы при её помощи доказать следующую за ней теорему.}.

{\sloppy

\paragraph{}\label{1938/159}
\so{Лемма}.
\textbf{\emph{Прямая}} ($DE$, рис.~\ref{1938/ris-169}), \textbf{\emph{параллельная какой-нибудь стороне}} ($AC$) \textbf{\emph{треугольника}} ($ABC$), \textbf{\emph{отсекает от него треугольник}} ($DBE$), \textbf{\emph{подобный данному.}}

}

Пусть в треугольнике $ABC$ прямая $DE$ параллельна стороне $AC$.
Требуется доказать, что $\triangle DBE\sim \triangle ABC$.

Предстоит доказать, во-первых, равенство соответственных углов и, во-вторых, пропорциональность соответственных сторон треугольников $ABC$ и $DBE$.

1.
Углы треугольников соответственно равны, так как угол $B$ у них общий, а $\angle D = \angle A$ и $\angle E= \angle C$, как соответственные углы при параллельных $DE$ и $AC$ и секущих $AB$ и $CB$.

\begin{wrapfigure}[13]{r}{51mm}
\centering
\includegraphics{mppics/ris-169}
\caption{}\label{1938/ris-169}
\end{wrapfigure}

2.
Докажем, что стороны $\triangle DBE$ пропорциональны соответственным сторонам $\triangle ABC$, то есть что
\[\frac{BD}{BA}=\frac{BE}{BC}=\frac{DE}{AC}.\]

Для этого рассмотрим отдельно следующие два случая:

1.
\mbox{\so{Стороны}} $AB$ \so{и $BD$ имеют общую меру}.

Разделим $AB$ на части, равные этой общей мере.
Тогда $BD$ разделится на целое число таких частей.
Пусть этих частей содержится $m$ в $BD$ и $n$ в $AB$.
Проведём из точек деления ряд прямых, параллельных $AC$, и другой ряд прямых, параллельных $BC$.
Тогда $BE$ и $BC$ разделятся на равные части (§~\ref{1938/95}), которых будет $m$ в $BE$ и $n$ в $BC$.
Точно так же $DE$ разделится на $m$ равных частей, а $AC$ на $n$ равных частей, причём части $BE$ равны частям $AC$ (как противоположные стороны параллелограммов).
Теперь очевидно, что
\begin{align*}
\frac{BD}{BA}&=\frac mn,
&
\frac{BE}{BC}&=\frac mn,
&
\frac{DE}{AC}&=\frac mn.
\end{align*}

Следовательно
\[\frac{BD}{BA}=\frac{BE}{BC}=\frac{DE}{AC}.\]

2. \so{Стороны} $AB$ и $BD$ \so{не имеют общей меры} (рис.~\ref{1938/ris-170}).

Найдём приближённые значения каждого из отношений $\frac{BD}{BA}$ и $\frac{BE}{BC}$, сначала с точностью до $\tfrac1{10}$;
затем до $\tfrac1{100}$ и далее будем последовательно повышать степень точности в 10 раз.

\begin{wrapfigure}{o}{45mm}
\centering
\includegraphics{mppics/ris-170}
\caption{}\label{1938/ris-170}
\end{wrapfigure}

Для этого разделим сторону $AB$ сначала на 10 равных частей и через точки деления проведём прямые, параллельные $AC$.
Тогда сторона $BC$ разделится также на 10 равных частей.
Предположим, что $\tfrac1{10}$ доля $AB$ укладывается в $BD$ более $m$
раз, причём получается остаток, меньший $\tfrac1{10}AB$.

Тогда, как видно из рис.~\ref{1938/ris-170}, $\tfrac1{10}$ доля $BC$ укладывается в $BE$ также $m$ раз с остатком, меньшим $\tfrac1{10}BC$.
Следовательно, с точностью до $\tfrac1{10}$ имеем:
\[\frac{BD}{AB}=\frac{m}{10}; 
\qquad
\frac{BE}{BC}=\frac{m}{10}.\]
Далее, разделим $AB$ на 100 равных частей и предположим, что $\tfrac1{100}AB$ укладывается $m_1$ раз в $BD$.
Проводя опять через точки деления прямые, параллельные $AC$, убеждаемся, что $\tfrac1{100}BC$ укладывается в $BE$ также $m_1$ раз.
Поэтому с точностью до $\tfrac1{100}$ имеем:
\[\frac{BD}{AB}=\frac{m_1}{100}; 
\qquad
\frac{BE}{BC}=\frac{m_1}{100}.\]

Повышая далее степень точности в $10,100,\dots$ раз, убеждаемся, что приближённые значения соотношений $\frac{BD}{BA}$ и $\frac{BE}{BC}$, вычисленные с произвольной, но одинаковой десятичной точностью, равны.
Следовательно, значения этих отношений выражаются одной и той же бесконечной десятичной дробью;
значит
\[\frac{BD}{BA}=\frac{BE}{BC}.\]

Точно так же, проводя через точки деления стороны $AB$ прямые, параллельные стороне $BC$, найдём, что
$\frac{BD}{BA}=\frac{DE}{AC}$.
Таким образом, и в этом случае имеем:
\[\frac{BD}{BA}=\frac{BE}{BC}=\frac{DE}{AC}.\]

{\small

\paragraph{}\label{1938/160}
\so{Замечание}:
Доказанные соотношения представляют собой три следующие пропорции:
\[\frac{BD}{BA}=\frac{BE}{BC};
\quad
\frac{BD}{BA}=\frac{DE}{AC};
\quad
\frac{BE}{BC}=\frac{DE}{AC}.\]
Переставив в них средние члены, получим:
\[\frac{BD}{BE}=\frac{BA}{BC};
\quad
\frac{BD}{DE}=\frac{BA}{AC};
\quad
\frac{BE}{DE}=\frac{BC}{AC}.\]

Таким образом, если в треугольниках стороны пропорциональны, то отношение любых двух сторон одного треугольника равно отношению соответственных сторон другого треугольника.
}

\subsection*{Признаки подобия треугольников}

\paragraph{}\label{1938/161}
\so{Теоремы}.
\textbf{\emph{Если в двух треугольниках:}}

1) \textbf{\emph{два угла одного треугольника соответственно равны двум углам другого}} или

2) \textbf{\emph{две стороны одного треугольника пропорциональны двум сторонам другого и углы, лежащие между этими сторонами, равны}} или

3) \textbf{\emph{если три стороны одного треугольника пропорциональны трём сторонам другого, то такие треугольники подобны.}}

1) Пусть $ABC$ и $A_1B_1C_1$ (рис.~\ref{1938/ris-171}) будут два треугольника, у которых $\angle A = \angle A_1$, $\angle B=\angle B_1$ и, следовательно, $\angle C=\angle C_1$.
Требуется доказать, что такие треугольники подобны.

\begin{figure}[!ht]
\centering
\includegraphics{mppics/ris-171}
\caption{}\label{1938/ris-171}
\end{figure}

Отложим на $AB$ отрезок $BD$, равный $A_1B_1$, и проведём $DE\z\parallel AC$.
Согласно доказанной выше лемме, $\triangle DBE\sim\triangle ABC$.

С другой стороны, $\triangle DBE= \triangle A_1B_1C_1$, потому что у них:
$BD\z=A_1B_1$ (по построению), $\angle B=\angle B_1$ (по условию) и $\angle D \z= \angle A_1$ (потому что $\angle D = \angle A$ и $\angle A = \angle A_1$).
Но очевидно, что если из двух равных треугольников один подобен третьему, то и другой ему подобен;
следовательно, 
\[\triangle A_1B_1C_1\sim\triangle ABC.\]

2) Пусть в треугольниках $ABC$ и $A_1B_1C_1$ (рис.~\ref{1938/ris-172}) дано:

\[\angle B=\angle B_1
\quad
\text{и}
\quad
\frac{AB}{A_1B_1}=\frac{BC}{B_1C_1}.\eqno(1)\]

\begin{figure}[!ht]
\centering
\includegraphics{mppics/ris-172}
\caption{}\label{1938/ris-172}
\end{figure}

Требуется доказать, что такие треугольники подобны.
Отложим снова на $AB$ отрезок $BD$, равный $A_1B_1$, и проведём $DE\parallel AC$.
Тогда получим вспомогательный $\triangle DBE$, подобный $\triangle ABC$.
Докажем, что он равен $\triangle A_1B_1C_1$.
Из подобия треугольников $ABC$ и $DBE$ следует:
\[\frac{AB}{DB}=\frac{BC}{BE}\eqno(2)\]

Сравнивая эту пропорцию с данной пропорцией (1), замечаем, что первые отношения обеих пропорций одинаковы ($DB\z=A_1B_1$ по построению);
следовательно, остальные отношения этих пропорций также равны, то есть 
\[\frac{BC}{B_1C_1}=\frac{BC}{BE}\]
Но если в пропорции предыдущие члены равны, то должны быть равны и последующие члены, значит
\[B_1C_1=BE.\]

Теперь видим, что треугольники $DBE$ и $A_1B_1C_1$ имеют по равному углу ($\angle B=\angle B_1$), заключённому между соответственно равными сторонами;
значит, эти треугольники равны.
Но $\triangle DBE$ подобен $\triangle ABC$, поэтому и $\triangle A_1B_1C_1$ подобен $\triangle ABC$.

3) Пусть в треугольниках $ABC$ и $A_1B_1C_1$ (рис.~\ref{1938/ris-173}) дано:
\[
\frac{AB}{A_1B_1}=\frac{BC}{B_1C_1}=\frac{AC}{A_1C_1}.\eqno(1)\]
Требуется доказать, что такие треугольники подобны.

\begin{figure}[!ht]
\centering
\includegraphics{mppics/ris-173}
\caption{}\label{1938/ris-173}
\end{figure}

Сделав построение такое же, как и прежде, покажем, что $\triangle DBE\z=\triangle A_1B_1C_1$.
Из подобия треугольников $ABC$ и $DBE$ следует:
\[\frac{AB}{DB}=\frac{BC}{BE}=\frac{AC}{DE}\eqno(2)\]

Сравнивая этот ряд отношений с данным рядом (1), замечаем, что первые отношения у них равны, следовательно, и остальные отношения равны, и потому
\[\frac{BC}{B_1C_1}=\frac{BC}{BE},
\qquad\text{откуда}\qquad
B_1C_1=BE,\]
и
\[\frac{AC}{A_1C_1}=\frac{AC}{DE},
\qquad\text{откуда}\qquad
A_1C_1=DE.\]

То есть треугольники $DBE$ и $A_1B_1C_1$ имеют по три соответственно равные стороны;
значит, они равны.
Но один из них, именно $\triangle DBE$, подобен $\triangle ABC$;
следовательно, и другой $\triangle A_1B_1C_1$ подобен $\triangle ABC$.

{\small
\paragraph{Замечания о приёме доказательства.}\label{1938/162}
Полезно обратить внимание на то, что приём доказательства, употреблённый нами в трёх предыдущих теоремах, один и тот же, а именно:
отложив на стороне большего треугольника отрезок, равный соответственной стороне меньшего, и проведя прямую, параллельную другой стороне, мы образуем вспомогательный треугольник, подобный большему данному.
После этого, в силу условия доказываемой теоремы и свойства подобных треугольников, мы обнаруживаем равенство вспомогательного треугольника меньшему данному и, наконец, делаем заключение о подобии данных треугольников.
}

\renewcommand{\bottomtitlespace}{.15\textheight}

\subsection*{Признаки подобия прямоугольных треугольников}

\renewcommand{\bottomtitlespace}{.1\textheight}

\paragraph{Признаки, не требующие особого доказательства.}\label{1938/163}
Так как прямые углы всегда равны друг другу, то на основании доказанных признаков подобия треугольников мы можем утверждать, что если в двух прямоугольных треугольниках:

1) \textbf{\emph{острый угол одного равен острому углу другого}} или

2) \textbf{\emph{катеты одного пропорциональны катетам другого, то такие треугольники подобны.}}

\paragraph{Признак, требующий особого доказательства.}\label{1938/164}\ 

\smallskip
\so{Теорема}.
\textbf{\emph{Если гипотенуза и катет одного треугольника пропорциональны гипотенузе и катету другого, то такие треугольники подобны.}}

Пусть $ABC$ и $A_1B_1C_1$ — два треугольника (рис.~\ref{1938/ris-174}), у которых углы $B$ и $B_1$ прямые и
\[\frac{AB}{A_1B_1}=\frac{AC}{A_1C_1}.\eqno(1)\]
Требуется доказать, что такие треугольники подобны.

\begin{figure}[!ht]
\centering
\includegraphics{mppics/ris-174}
\caption{}\label{1938/ris-174}
\end{figure}

Для доказательства применим тот же приём, которым мы пользовались ранее.
На $AB$ отложим $BD=A_1B_1$ и проведём $DE\parallel AC$.
Тогда получим вспомогательный $\triangle DBE$, подобный $\triangle ABC$.
Докажем, что он равен $\triangle A_1B_1C_1$.
Из подобия треугольников $ABC$ и $DBE$ следует:
\[\frac{AB}{DB}=\frac{AC}{DE}.\eqno(2)\]

Сравнивая эту пропорцию с данной (1), находим, что первые отношения их одинаковы;
следовательно, равны и вторые отношения, то есть 
\[\frac{AC}{DE}=\frac{AC}{A_1C_1},\]
откуда
\[DE=A_1C_1.\]

Теперь видим, что треугольники $BDE$ и $A_1B_1C_1$ имеют по равной гипотенузе и равному катету, следовательно, они равны;
а так как один из них подобен $\triangle ABC$, то и другой ему подобен.

\paragraph{}\label{1938/165}
\so{Теорема} (об отношении высот).
\textbf{\emph{В подобных треугольниках соответственные стороны пропорциональны соответственным высотам}}, то есть тем высотам, которые опущены на соответственные стороны.

\begin{figure}[!ht]
\centering
\includegraphics{mppics/ris-175}
\caption{}\label{1938/ris-175}
\end{figure}

Действительно, если треугольники $ABC$ и $A_1B_1C_1$ (рис.~\ref{1938/ris-175}) подобны, то прямоугольные треугольники $BAD$ и $B_1A_1D_1$ также подобны ($\angle A = \angle A_1$ и $\angle D=\angle D_1$);
поэтому:
\[\frac{BD}{B_1D_1}=\frac{AB}{A_1B_1}=\frac{BC}{B_1C_1}=\frac{AC}{A_1C_1}.\]

%% file: 2D/podobie-fig.tex
\section{Подобие многоугольников}

\paragraph{}\label{1938/168}
\so{Определение}.
Два $n$-угольника называются \rindex{подобные!многоугольники}\textbf{подобными}, если углы одного равны соответственно углам другого и соответственные стороны пропорциональны.

Во избежание путаницы у подобных многоугольников принято записывать соответственные вершины в том же порядке.
Например «пятиугольник
$ABCDE$ подобен пятиугольнику $A_1B_1C_1D_1E_1$» (рис.~\ref{1938/ris-180}) обычно означает, что  
\[\angle A = \angle A_1, 
\quad
\angle B=\angle B_1,
\quad
\angle C=\angle C_1,
\quad
\angle D=\angle D_1,
\quad
\angle E=\angle E_1\]
и
\[\frac{AB}{A_1B_1}=\frac{BC}{B_1C_1}=\frac{CD}{C_1D_1}=\frac{DE}{D_1E_1}=\frac{EA}{E_1A_1}\]
При этом пары сторон $AB$ и $A_1B_1$, $BC$ и $B_1C_1$, $CD$ и $C_1D_1$, $DE$ и $D_1E_1$, $AE$ и $A_1E_1$
являются соответственными.

То, что такие многоугольники существуют, будет видно из решения следующей задачи:

\paragraph{}\label{1938/169}
\so{Задача}.
\emph{Дан многоугольник $ABCDE$ и отрезок $a$.
Построить другой многоугольник, который был бы подобен данному и у которого сторона, соответственная стороне $AB$ 
равнялась бы $a$} (рис.~\ref{1938/ris-179}).

\begin{figure}[h]
\centering
\includegraphics{mppics/ris-179}
\caption{}\label{1938/ris-179}
\end{figure}

На стороне $AB$ отложим $AB_1=a$ (если $a>AB$, то точка $B_1$ расположится на продолжении $AB$).
Затем, проведя из $A$ все диагонали, построим $B_1C_1 \parallel BC$, $C_1D_1\parallel CD$ и $D_1E_1\parallel DE$.

Тогда получим многоугольник $AB_1C_1D_1E_1$, подобный многоугольнику $ABCDE$.
Действительно, во-первых, углы одного из них соответственно равны углам другого;
так, угол $A$ у них общий, $\angle B_1=\angle B$ и $\angle E_1=\angle E$ как соответственные углы при параллельных прямых;
$\angle C_1=\angle C$ и $\angle D_1=\angle D$, так как углы эти состоят из частей, соответственно равных друг другу.

Во-вторых, мы имеем пропорции:
\begin{align*}
\text{из подобия }&\ \triangle AB_1C_1\sim\triangle ABC:
&
\frac{AB_1}{AB}=\frac{B_1C_1}{BC}=\frac{AC_1}{AC};
\\
\text{—\textquotedbl—\textquotedbl—\quad}&\ \triangle AC_1D_1\sim\triangle ACD:
&
\frac{AC_1}{AC}=\frac{C_1D_1}{CD}=\frac{AD_1}{AD};
\\
\text{—\textquotedbl—\textquotedbl—\quad}&\ \triangle AD_1E_1\sim\triangle ADE:
&
\frac{AD_1}{AD}=\frac{D_1E_1}{DE}=\frac{AE_1}{AE}.
\end{align*}

Так как третье отношение первого ряда равно первому отношению второго ряда и третье отношение второго ряда равно первому отношению третьего ряда, то значит, все 9 отношений равны между собой.
Выбросив из них отношения, в которые входят диагонали, можем написать:
\[\frac{AB_1}{AB}=\frac{B_1C_1}{BC}=\frac{C_1D_1}{CD}=\frac{D_1E_1}{DE}=\frac{AE_1}{AE}.\]

Мы видим, что у $n$-угольников $ABCDE$ и $AB_1C_1D_1E_1$ соответственные углы  равны и соответственные стороны пропорциональны;
значит, многоугольники эти подобны.

{\small

\paragraph{}\label{1938/170}
\so{Замечание}.
Для треугольников, как мы видели (§~\ref{1938/161}), равенство углов влечёт пропорциональность сторон, и обратно:
пропорциональность сторон влечёт равенство углов;
вследствие этого для треугольников одно равенство углов или одна пропорциональность сторон служит достаточным признаком их подобия.
Для многоугольников же одного равенства углов или одной пропорциональности сторон ещё недостаточно для их подобия;
например, у квадрата и прямоугольника углы равны, но стороны не пропорциональны, у квадрата же и ромба стороны пропорциональны, а углы не равны.

}

\paragraph{}\label{1938/171}
\so{Теорема} (о разложении подобных многоугольников на подобные треугольники).
\textbf{\emph{Подобные многоугольники можно разложить на одинаковое число подобных и одинаково расположенных треугольников.}}

Например, подобные многоугольники $ABCDE$ и $AB_1C_1D_1E_1$ (рис. \ref{1938/ris-179}) разделены диагоналями на подобные треугольники, которые одинаково расположены — если пара треугольников в $ABCDE$ имеет общую сторону, то тоже верно и про пару им сходственных треугольников в $AB_1C_1D_1E_1$ и наоборот.

Укажем ещё такой способ разложения.
Возьмём внутри многоугольника $ABCDE$ (рис.~\ref{1938/ris-180}) произвольную точку $O$ и соединим её со всеми вершинами.
Тогда многоугольник $ABCDE$ разобьётся на столько треугольников, сколько в нём сторон.
Возьмём один из них, например $AOE$ (покрытый на рисунке штрихами), и на соответственной стороне $A_1E_1$ другого многоугольника построим углы $O_1A_1E_1$ и $O_1E_1A_1$ соответственно равные углам $OAE$ и $OEA$;
точку пересечения $O_1$ соединим с прочими вершинами многоугольника $A_1B_1C_1D_1E_1$.
Тогда и этот многоугольник разобьётся на то же число треугольников.
Докажем, что треугольники первого многоугольника соответственно подобны треугольникам второго многоугольника.

По построению $\triangle AOE\sim \triangle A_1O_1E_1$. 
Чтобы доказать подобие соседних треугольников $ABO$ и $A_1B_1O_1$, примем во внимание, что из подобия многоугольников следует:
\[\angle BAE=\angle B_1A_1E_1
\quad\text{и}\quad
\frac{BA}{B_1A_1}=\frac{AE}{A_1E_1},
\eqno(1)\]
и из подобия треугольников $AOE$ и $A_1O_1E_1$ выводим:
\[\angle OAE=\angle O_1A_1E_1
\quad\text{и}\quad
\frac{AO}{A_1O_1}=\frac{AE}{A_1E_1}.
\eqno(2)\]
Из равенств (1) и (2) следует:
\[\angle BAO=\angle B_1A_1O_1
\quad\text{и}\quad
\frac{BA}{B_1A_1}=\frac{AO}{A_1O_1}.\]
Теперь видим, что треугольники $ABO$ и $A_1B_1O_1$ имеют по равному углу, заключённому между пропорциональными сторонами;
значит, они подобны.

Совершенно так же докажем, подобие $\triangle BCO\sim \triangle B_1C_1O_1$, затем $\triangle CDO\sim\triangle C_1D_1O_1$, и~т.~д.
При этом заметим, что треугольники первого многоугольника располагаются одинаково с подобными им треугольниками второго многоугольника. 

\begin{figure}[h]
\centering
\includegraphics{mppics/ris-180}
\caption{}\label{1938/ris-180}
\end{figure}

\paragraph{}\label{1938/172}
\so{Теорема}.
\textbf{\emph{Периметры подобных многоугольников относятся как соответственные стороны.}}

Пусть многоугольники $ABCDE$ и $A_1B_1C_1D_1E_1$ (рис.~\ref{1938/ris-180}) подобны;
тогда по определению: 
\[\frac{AB}{A_1B_1}=\frac{BC}{B_1C_1}=\frac{CD}{C_1D_1}=\frac{DE}{D_1E_1}=\frac{EA}{E_1A_1}.\]
Если имеем ряд равных отношений, то сумма всех предыдущих членов относится к сумме всех последующих, как какой-нибудь из предыдущих членов относится к своему последующему, поэтому
\[\frac{AB+BC+CD+DE+EA}{A_1B_1+B_1C_1+C_1D_1+D_1E_1+E_1A_1}=\frac{AB}{A_1B_1}=\frac{BC}{B_1C_1}=\dots\]

\paragraph{Коэффициент подобия.}\label{1938/173}
Отношение соответственных сторон двух подобных многоугольников (в частности  треугольников) называется \rindex{коэффициент подобия}\textbf{коэффициентом подобия} этих многоугольников.

\paragraph{Центрально-подобное преобразование.}\label{1938/174}
Способ построения многоугольника, подобного данному, изложенный в §~\ref{1938/169}, является частным видом так называемого центрально-подобного преобразования.

\begin{wrapfigure}{o}{56mm}
\centering
\includegraphics{mppics/ris-181}
\caption{}\label{1938/ris-181}
\end{wrapfigure}

Общий метод такого преобразования состоит в следующем.
Пусть требуется подобно преобразовать четырёхугольник $ABCD$ (рис.~\ref{1938/ris-181}) при коэффициенте подобия, равном $k$.
Возьмём какую-нибудь точку $O$ на плоскости;
соединив её с вершинами данного четырёхугольника, получим прямые $OA$, $OB$, $OC$ и $OD$.
На прямой $OA$ отложим от точки $O$ в сторону точки $A$ отрезок $OA_1$, равный $k\cdot OA$, так что $OA_1=k\cdot OA$ (на рисунке $k=\tfrac53$).

Продолжим также прямую $OB$ и отложим на ней от точки $O$ в сторону точки $B$ отрезок $OB_1$, равный $k\cdot OB$, так что $OB_1=k\cdot OB$.

Точно так же поступим с прямыми $OC$ и $OD$.
Мы получим на них точки $C_1$ и $D_1$, причём $OC_1=k\cdot OC$ и $OD_1=k\cdot OD$.
Соединив последовательно точки $A_1$, $B_1$, $C_1$ и $D_1$, получим искомый четырёхугольник $A_1B_1C_1D_1$.
В самом деле, из равенств $OA_1=k\cdot OA$, $OB_1=k\cdot OB$, $OC_1=k\cdot OC$ и $OD_1=k\cdot OD$ следует:
\[\frac{OA_1}{OA}=
\frac{OB_1}{OB}=
\frac{OC_1}{OC}=
\frac{OD_1}{OD}=k.\]

Сравним треугольники $OAB$ и $OA_1B_1$.
Они имеют общий угол в вершине $O$ и, кроме того,
\[\frac{OA_1}{OA}=
\frac{OB_1}{OB},\]
следовательно, эти треугольники подобны (§~\ref{1938/161}, 2-й случай).
Из их подобия заключаем:
\[\frac{A_1B_1}{AB}=
\frac{OA_1}{OA}=k
\quad\text{и}\quad
\angle OAB=\angle OA_1B_1,
\eqno(1)
\]
следовательно, $AB\parallel A_1B_1$ (§~\ref{1938/73}).

Совершенно так же покажем, что треугольники $OBC$ и $OB_1C_1$ подобны.
Отсюда следует:
\[\frac{B_1C_1}{BC}=
\frac{OB_1}{OB}=k
\quad\text{и}\quad
\angle OBC=\angle OB_1C_1,
\eqno(2)
\]
и, следовательно, $BC\parallel B_1C_1$.

Таким же образом докажем подобие следующих треугольников:
$OCD$ и $OC_1D_1$, затем треугольников $OAD$ и $OA_1D_1$.

Из подобия $\triangle OCD$ и $\triangle OC_1D_1$ следует:
\[\frac{C_1D_1}{CD}=
\frac{OC_1}{OC}=k
\quad\text{и}\quad
CD\parallel C_1D_1,
\eqno(3)
\]

Из подобия $\triangle OAD$ и $\triangle OA_1D_1$ следует:
\[\frac{A_1D_1}{AD}=
\frac{OD_1}{OD}=k
\quad\text{и}\quad
AD\parallel A_1D_1,
\eqno(4)
\]
Кроме того, $\angle DAB=\angle D_1A_1B_1$ как углы с параллельными сторонами (§~\ref{1938/79}).

По той же причине имеем равенство углов.
\begin{align*}
\angle ABC &= \angle A_1B_1C_1,
\\
\angle BCD &= \angle B_1C_1D_1,
\\
\angle CDA &= \angle C_1D_1A_1.
\end{align*}
Мы видим, что у четырёхугольников $ABCD$ и $A_1B_1C_1D_1$ соответственные углы равны и соответственные стороны пропорциональны;
значит, эти четырёхугольники подобны, причём коэффициент их подобия равен $k$.

\paragraph{Центр подобия.}\label{1938/175}\rindex{центр!подобия}
При центрально-подобном преобразовании многоугольника способом, изложенным в §~\ref{1938/174}, точка $O$ называется центром подобия обоих многоугольников.

Центрально-подобное преобразование многоугольника можно выполнять несколько иначе.
Именно, взяв точку $O$ и соединив её с вершинами четырёхугольника $ABCD$, можно продолжить прямые $OA$, $OB,\dots$
за точку $O$;
затем на прямой $OA$ от точки $O$ в сторону, противоположную точке $A$, отложим отрезок $OA'$, равный $k\cdot OA$.
Точно так же на продолжениях прямых $OB$, $OC,\dots$
от точки $O$ отложим отрезки $OB', OC',\dots$, равные соответственно отрезкам $k\cdot OB, k\cdot OC,\dots$
(рис.~\ref{1938/ris-182});
соединив последовательно точки $A'$, $B'$, $C'$, $D'$ получим четырёхугольник $A'B'C'D'$, очевидно симметричный с $A_1B_1C_1D_1$ относительно точки $O$.
Следовательно, четырёхугольники $A'B'C'D'$ и $A_1B_1C_1D_1$ равны и, значит, четырёхугольники $ABCD$ и $A'B'C'D'$ подобны, причём коэффициент их подобия равен $k$.
При первом способе преобразования точка $O$ называется {}\textbf{внешним центром} подобия многоугольников (рис.~\ref{1938/ris-181});
при втором способе — \rindex{центр!подобия}\textbf{внутренним центром} их подобия (рис.~\ref{1938/ris-182}).

\begin{figure}[h]
\centering
\includegraphics{mppics/ris-182}
\caption{}\label{1938/ris-182}
\end{figure}

{\small

\smallskip
\so{Замечание}.
При выполнении центрально-подобного преобразования можно одинаково пользоваться как внутренним, так и внешним центром подобия.
И тот и другой можно выбирать совершенно произвольно, в частности, если принять одну из вершин многоугольника за внешний центр подобия и выполнить центрально-подобное преобразование, то получим как раз тот способ построения подобного многоугольника, который был изложен в §~\ref{1938/169}.

}

\paragraph{Центрально-подобное расположение многоугольников.}\label{1938/176}
Расположение двух многоугольников $ABCD$ и $A_1B_1C_1D_1$ на рис.~\ref{1938/ris-181}, а также многоугольников $ABCD$ и $A'B'C'D'$ на рис.~\ref{1938/ris-182} имеет следующие свойства:
1) соответственные стороны обоих многоугольников параллельны;
2) прямые, соединяющие соответственные вершины, пересекаются в одной точке.

\begin{figure}[h]
\centering
\includegraphics{mppics/ris-183}
\caption{}\label{1938/ris-183}
\end{figure}

{\sloppy 
Такое расположение многоугольников называется центрально-подобным. 
Докажем, что в такое расположение можно привести любые два подобных многоугольника.

}

Пусть даны два подобных многоугольника $ABCDE$ и $A_1B_1C_1D_1E_1$ (рис.~\ref{1938/ris-183}).
Возьмём какую-либо точку $O$ за центр подобия и построим многоугольник, центрально-подобный  с $ABCDE$, причём коэффициент подобия возьмём равным отношению $\frac{A_1B_1}{AB}$.
Мы получим многоугольник $A'B'C'D'E'$, подобный $ABCDE$ и в то же время равный $A_1B_1C_1D_1E_1$.
В самом деле, так как коэффициент подобия многоугольников $ABCDE$ и $A'B'C'D'E'$ равен $\frac{A_1B_1}{AB}$, то $\frac{A'B'}{AB}=\frac{A_1B_1}{AB}$, отсюда $A'B'=A_1B_1$.
Но многоугольники $A_1B_1C_1D_1E_1$ и $A'B'C'D'E'$ подобны между собой, следовательно:
\[\frac{A'B'}{A_1B_1}=\frac{B'C'}{B_1C_1}=\frac{C'D'}{C_1D_1}=\frac{D'E'}{D_1E_1}=\frac{A'E'}{A_1E_1}.\]
А потому из равенства $A'B'=A_1B_1$ вытекает 
\begin{align*}
B'C'&=B_1C_1,&
C'D'&=C_1D_1,&
D'E'&=D_1E_1,&
A'E'&=A_1E_1.
\end{align*}
Так как, кроме того, углы многоугольника $A_1B_1C_1D_1E_1$ равны соответствующим углам многоугольника $A'B'C'D'E'$, то эти многоугольники равны между собой. 
Если наложить многоугольник $A_1B_1C_1D_1E_1$ на  $A'B'C'D'E'$ так, чтобы они совпадали, то многоугольник $A_1B_1C_1D_1E_1$ примет центрально-подобное расположение с $ABCDE$. 

В случае если $k=1$ то многоугольники равны.
В этом случае следует взять $O$ как внутренний центр подобия, иначе многоугольник $A'B'C'D'E'$ сольётся с $ABCDE$.

\section{Подобие произвольных фигур}

\paragraph{}\label{1938/177}
{\sloppy 
Центрально-подобное преобразование многоугольников даёт возможность обобщить само понятие о подобии на случай, когда фигура образована кривыми линиями.
Именно такой способ построения подобной фигуры можно применить к любой фигуре.
Пусть, например, на плоскости дана фигура $A$ совершенно произвольной формы (рис.~\ref{1938/ris-184}).

}

\begin{figure}[h]
\centering
\includegraphics{mppics/ris-184}
\caption{}\label{1938/ris-184}
\end{figure}

Возьмём произвольную точку $O$ на плоскости этой фигуры и будем соединять её с различными точками $M$, $N$, $P,\dots$
фигуры $A$.
На каждой из проведённых прямых $OM$, $ON$, $OP,\dots$
отложим отрезки $OM_1$, $ON_1$, $OP_1$, такие, что 
\[\frac{OM_1}{OM}=\frac{ON_1}{ON}=\frac{OP_1}{OP}=\dots\quad\text{и так далее}\]
Точки $M_1$, $N_1$, $P_1,\dots$ будут лежать на некоторой новой фигуре $A_1$.
Чем больше точек $M$, $N$, $P,\dots$
мы возьмём на фигуре $A$, тем больше мы получим точек фигуры $A_1$.
Чтобы получить всю фигуру $A_1$, нужно провести прямые из точки $O$ ко всем точкам фигуры $A$ и построить на них соответствующие точки фигуры $A_1$.
Считается, что фигуры $A_1$ и $A$ расположены центрально-подобно,
а любая фигура равная $A_1$ считается подобной~$A$.

В отдельных случаях, чтобы получить фигуру $A_1$, нет необходимости проводить лучи ко всем точкам фигуры $A$, достаточно построить лишь несколько её точек, и затем, пользуясь частными свойствами фигуры $A$, восстановить всю фигуру $A_1$.
Так, в том случае, когда $A$ — многоугольник, достаточно было соединить точку $O$ лишь с вершинами этого многоугольника и построить вершины подобного многоугольника, а затем соединить прямолинейными отрезками полученные вершины между собой.

{\sloppy 
Описанный переход от фигуры $A$ к фигуре $A_1$ называется \rindex{центрально-подобное преобразование}\textbf{центрально-подобным преобразованием}. 
Этот тип преобразований имеет заметное применение на практике.
Показываемая в кино картина на экране подобна изображению, сделанному на плёнке;
технические чертежи планов и фасадов зданий, планов местности, планов городов
получаются в результате подобного преобразования.

}

{\sloppy

\paragraph{Подобие окружностей.}\label{1938/178}
Докажем, что фигура, подобная окружности, есть также окружность.

}

\smallskip

\so{Теорема}.
\textbf{\emph{Геометрическое место точек, делящих в данном отношении отрезки, соединяющие какую-нибудь точку с точками окружности, есть окружность.}}

\begin{figure}[h]
\centering
\includegraphics{mppics/ris-185}
\caption{}\label{1938/ris-185}
\end{figure}

Пусть дана окружность радиуса $R$ с центром в точке $O$ (рис.~\ref{1938/ris-185}).
Возьмём произвольную точку $S$ и, соединив её с точкой $O$, разделим отрезок $SO$ точкой $O_1$ в некотором отношении так, что $\frac{SO_1}{SO}= k$.

Возьмём произвольную точку $M$ на данной окружности и соединим её с точкой $S$.
На отрезке $SM$ найдём точку $M_1$ такую, что $\frac{SM_1}{SM}\z=\frac{SO_1}{SO}= k$.
Для этой цели следует из точки $O_1$ провести прямую, параллельную $OM$, до пересечения с прямой $SM$.
Из подобия треугольников $SOM$ и $SO_1M_1$, следует $\frac{OM_1}{OM}=\frac{SO_1}{SO}$.
Следовательно, $\frac{OM_1}{OM}=k$.
Отсюда найдём длину отрезка $O_1M_1$, именно $O_1M_1=k\cdot OM$, или $O_1M_1=k\cdot R$.

Мы видим, что величина $O_1M_1$ не зависит от положения точки $M$ на данной окружности.
Следовательно, если точка $M$ будет перемещаться по окружности, то точка $M_1$ будет перемещаться по плоскости, описывая окружность с центром $O_1$ и радиусом~$k\cdot R$.

\paragraph{}\label{1938/179}
\so{Теорема}.
\textbf{\emph{Две окружности на плоскости всегда можно рассматривать как центрально-подобные фигуры, причём окружности разных радиусов имеют два центра подобия:
один внешний, другой внутренний.}} 

Пусть даны две окружности с центрами $O_1$ и $O_2$ и радиусами $R_1$ и $R_2$ (рис.~\ref{1938/ris-186}).
Проведём линию центров $O_1O_2$ и построим на ней две точки $I$ и $E$, определяемые равенствами
\[\frac{O_1I}{O_2I}=\frac{R_1}{R_2}\quad\text{и}\quad\frac{O_1E}{O_2E}=\frac{R_1}{R_2}\]
(Точка $I$ на отрезке $O_1O_2$ существует всегда, а точка $E$ на продолжении отрезка существует только если $R_1\ne R_2$.)

Покажем, что точки $I$ и $E$ обладают свойствами центров подобия.
Возьмём какую-либо точку $M_1$ на первой окружности, проведём прямую $IM_1$ и отложим на ней отрезок $IM_2$ так, что 
\[\frac{IM_1}{IM_2}=\frac{R_1}{R_2}.\] 
Тогда $\triangle IO_1M_1\sim\triangle IO_2M_2$, так как 
\begin{align*}
\angle O_1IM_1&=\angle O_2IM_2,
&
\frac{IM_1}{IM_2}&=\frac{R_1}{R_2}
&&\text{и}&
\frac{O_1I}{O_2I}&=\frac{R_1}{R_2}.
\end{align*}
Следовательно:
\[\frac{O_1M_1}{O_2M_2}=\frac{R_1}{R_2},\]
и, так как $O_1M_1=R_1$, получаем
\[O_2M_2=R_2.\]
Это означает, что точка $M_2$ лежит на второй окружности.
Следовательно, точка $I$ есть внутренний центр подобия данных окружностей.
Таким же образом можно доказать, что $E$ есть внешний центр подобия.

\begin{figure}[h]
\centering
\includegraphics{mppics/ris-186}
\caption{}\label{1938/ris-186}
\end{figure}

Построение точек $I$ и $E$ можно выполнить так:
проводим в данных окружностях два каких-либо параллельных радиуса и соединяем их концы, полученная прямая пересечёт линию центров в центре подобия.
При этом если проведённые радиусы направлены в одну сторону (рис.~\ref{1938/ris-186}, $O_1A_1$, и $O_2A_2$), то центр подобия будет внешним;
если они направлены в противоположные стороны (рис.~\ref{1938/ris-186}, $O_1M_1$ и $O_2M_2$), то центр подобия будет внутренним.
В случае равенства $R_1=R_2$ вторая прямая будет параллельна линии центров, в этом случае нет внешнего центра $E$.

\smallskip

Легко далее заметить, что если две окружности касаются, то один из центров подобия совпадает с точкой касания.
При этом если касание окружностей внешнее, то в точке касания находится внутренний центр подобия, если же касание внутреннее, то с точкой касания совпадает внешний центр подобия окружностей.

\smallskip
\so{Упражнения 1}.
Доказать, что если две окружности лежат одна вне другой, то их внешний центр подобия совпадает с точкой пересечения их общих внешних касательных, а внутренний — с точкой пересечения общих внутренних касательных.

2.
Какое положение должны иметь две окружности на плоскости, чтобы их внешний центр подобия совпал с внутренним? 

\smallskip
\so{Ответ}.
Окружности концентричны.

\subsection*{Задачи на построение}

\paragraph{Метод подобия.}\label{1938/181}
Центрально-подобное преобразование фигур можно с успехом применить к решению многих задач на построение.
На этом основан так называемый метод подобия.

Метод подобия состоит в том, что, пользуясь некоторыми данными задачи, строят сначала фигуру, \so{подобную} искомой, а затем переходят к искомой.
Этот метод особенно удобен тогда, когда только одна данная величина есть длина, а все прочие величины — или углы, или отношения линий;
таковы, например, задачи:
построить треугольник по данному углу, стороне и отношению двух других сторон или по двум углам и длине некоторого отрезка (например высоте, медиане или биссектрисе);
построить квадрат по данной сумме или разности между диагональю и стороной.

Решим, например, такую задачу:

\smallskip
\so{Задача 1}.
\emph{Построить треугольник, когда даны один из его углов $C$, отношение сторон $AC:BC$, заключающих этот угол, и высота $h$, опущенная из вершины этого угла на противоположную сторону} (рис.~\ref{1938/ris-189}).

Пусть $AC:BC=m:n$, где $m$ и $n$ — два данных отрезка или два данных числа.
Строим угол $C$, на его сторонах откладываем отрезки $CA_1$ и $CB_1$, пропорциональные $m$ и~$n$.
Если $m$ и $n$ — отрезки, то берём прямую $CA_1=m$ и $CB_1=n$.
Если $m$ и $n$ — числа, то, выбрав произвольный отрезок $\ell$, строим отрезки $CA_1\z=m \cdot \ell$ и $CB_1=n\cdot \ell$.
В обоих случаях имеем
\[\frac{CA_1}{CB_1} = \frac mn.\]

Треугольник $CA_1B_1$, очевидно, подобен искомому треугольнику.

Чтобы получить искомый треугольник, построим в треугольнике $CA_1B_1$ высоту $CD_1$, и обозначим её через $h_1$.
Выбираем теперь произвольный центр подобия и строим треугольник, подобный треугольнику $A_1B_1C$ с коэффициентом подобия, равным отношению $\tfrac{h}{h_1}$,
где $h$ — данная высота искомого треугольника.
Полученный таким путём треугольник и будет искомым.

\begin{wrapfigure}{o}{30mm}
\centering
\includegraphics{mppics/ris-189}
\caption{}\label{1938/ris-189}
\end{wrapfigure}

Удобнее всего выбрать центр подобия прямо в точке $C$.
В таком случае построение искомого треугольника становится особенно простым (рис.~\ref{1938/ris-189}).
Продолжаем высоту $CD_1$ треугольника $A_1B_1C$, откладываем на ней отрезок $CD$, равный $h$, и проводим прямую $AB$, параллельную $A_1B_1$.
Треугольник $ABC$ — искомый. 

\medskip

В задачах этого рода положение искомой фигуры остаётся произвольным;
но во многих вопросах требуется построить фигуру, положение которой относительно данных точек или линий вполне определено.
При этом может случиться, что, отрешившись от какого-нибудь одного из условий положения и оставив все остальные, мы получим бесчисленное множество фигур, \so{подобных} искомой.
В таком случае метод подобия может быть употреблён с пользой.
Приведём примеры.

\smallskip
\so{Задача 2}.
\emph{В данный угол $ABC$ вписать окружность, которая проходила бы через данную точку $M$} (рис.~\ref{1938/ris-190}).

\begin{figure}[h]
\centering
\includegraphics{mppics/ris-190}
\caption{}\label{1938/ris-190}
\end{figure}

Отбросим на время требование, чтобы окружность проходила через точку $M$.
Тогда данному условию удовлетворяет бесчисленное множество окружностей, центры которых лежат на биссектрисе $BD$.
Построим одну из таких окружностей; обозначим её центр буквой $o$. 
Возьмём на ней точку $m$, соответственную точке $M$, то есть лежащую на луче $MB$, и проведём радиус $mo$.
Если построим $MO\parallel mo$, то точка $O$ будет центром искомого круга.
Действительно, проведя к стороне $AB$ перпендикуляры $ON$ и $on$, мы получим две пары подобных треугольников $\triangle MBO\sim \triangle mBo$ и $\triangle NBO\sim \triangle nBo$, из которых будем иметь:
\[\frac{MO}{mo} = \frac{BO}{Bo}, \quad 
 \frac{NO}{no} = \frac{BO}{Bo},
\]
откуда 
\[\frac{MO}{mo} = \frac{NO}{no}.\]
Но $mo =no$;
следовательно, $MO=NO$, то есть окружность, описанная радиусом $OM$ с центром $O$, касается стороны $AB$;
а так как её центр лежит на биссектрисе угла, то она касается и стороны $BC$.

За соответственную точку можно взять и другую точку $m_1$ пересечения окружности с лучом $MB$.
Так мы найдём другой центр $O_1$ искомого круга.
Следовательно, задача допускает два решения.

\medskip

\begin{wrapfigure}{r}{40mm}
\vskip-5mm
\centering
\includegraphics{mppics/ris-191}
\caption{}\label{1938/ris-191}
\end{wrapfigure}

\mbox{\so{Задача 3}.}
\emph{В данный треугольник $ABC$ вписать ромб с данным острым углом так, чтобы одна из его сторон лежала на основании $AB$ треугольника $ABC$, а две его вершины — на боковых сторонах $AC$ и $BC$} (рис.~\ref{1938/ris-191}).

Отбросим на время требование, чтобы одна из вершин ромба лежала на стороне $BC$.
Тогда можно построить бесчисленное множество ромбов, удовлетворяющих остальным условиям задачи.
Построим один из них.

Берём на стороне $AC$ произвольную точку $M$.
Строим угол с вершиной в этой точке, равный данному, одна сторона которого была бы параллельна основанию $AB$, а другая пересекала основание $AB$ в некоторой точке $N$.
На стороне $AB$ от точки $N$ откладываем отрезок $NP$, равный $MN$, и строим ромб со сторонами $MN$ и $NP$.

Пусть $Q$ — его четвёртая вершина.
Далее, выбираем вершину $A$ за центр подобия и строим ромб, подобный ромбу $MNPQ$, выбирая коэффициент подобия так, чтобы вершина нового ромба, соответствующая вершине $Q$, оказалась на стороне $BC$.
Для этой цели продолжаем прямую $AQ$ до пересечения со стороной $BC$ в некоторой точке $x$.
Эта точка $x$ будет одной из вершин искомого ромба.

{\sloppy

Проводя из этой точки прямые, параллельные сторонам ромба $MNPQ$, получаем искомый ромб $xyzu$.

}

Предоставляем самим учащимся решить методом подобия следующие задачи:

1.
Построить треугольник, зная два его угла и радиус описанной окружности.

2.
Построить треугольник, зная отношение высоты к основанию, угол при вершине и медиану боковой стороны.

3.
Дан $\angle AOB$ и внутри него точка $C$.
Найти на стороне $OB$ точку $M$, равно отстоящую от $OA$ и от точки $C$.

%% file: 2D/proportzii.tex
\section{Пропорции} 

\paragraph{Теорема Фалеса.}\label{1938/182}
Следующая теорема названа в честь древнегреческого математика Фалеса Милетского (VII—VI век до нашей эры).

\smallskip
\so{Теорема}.
\textbf{\emph{Стороны угла}} ($ABC$), \textbf{\emph{пересекаемые рядом параллельных прямых}} ($DD_1, EE_1, FF_1, \dots$), \textbf{\emph{рассекаются ими на пропорциональные части}} (рис.~\ref{1938/ris-192}).

\begin{figure}[!ht]
\centering
\includegraphics{mppics/ris-192}
\caption{}\label{1938/ris-192}
\end{figure}

Требуется доказать, что
\[\frac{BD}{BD_1}=\frac{DE}{D_1E_1}=\frac{EF}{E_1F_1},\]
или
\begin{align*}
\frac{BD}{DE}&=\frac{BD_1}{D_1E_1},
\\
\frac{DE}{EF}&=\frac{D_1E_1}{E_1F_1}\quad\text{и так далее}
\end{align*}
Проводя вспомогательные прямые $DM$, $EN$ и~т.~д., параллельные $BA$, мы получим треугольники $BDD_1$, $DEM$, $EFN$ и~т.~д., которые все подобны между собой, так как углы у них соответственно равны (вследствие параллельности прямых).
Из их подобия следует:
\[\frac{BD}{BD_1}=\frac{DE}{DM}=\frac{EF}{EN}\quad\text{и так далее}\]
Заменив в этом ряду равных отношений отрезок $DM$ на $D_1E_1$, отрезок $EN$ на $E_1F_1$ и~т.~д.
(противоположные стороны параллелограммов равны), мы получим то, что требовалось доказать.

\paragraph{}\label{1938/183}
\so{Теорема}.
\textbf{\emph{Две параллельные прямые}} ($MN, M_1N_1$, рис. \ref{1938/ris-193}), \textbf{\emph{пересекаемые рядом прямых}} ($OA, OB, OC, \dots$), \textbf{\emph{исходящих из одной и той же точки}} ($O$), \textbf{\emph{рассекаются ими на пропорциональные части.}}

Требуется доказать, что отрезки $AB$, $BC$, $CD,\dots$
прямой $MN$ пропорциональны отрезкам $A_1B_1$, $B_1C_1$, $C_1D_1,\dots$
прямой $M_1N_1$.

\begin{figure}[!ht]
\centering
\includegraphics{mppics/ris-193}
\caption{}\label{1938/ris-193}
\end{figure}

Из подобия треугольников $OAB$ и $OA_1B_1$ (§~\ref{1938/159}) и треугольников $OBC$ и $OB_1C_1$ выводим:
\[\frac{AB}{A_1B_1}=\frac{BO}{B_1O}
\quad\text{и}\quad
\frac{BO}{B_1O}=\frac{BC}{B_1C_1},
\]
откуда
\[\frac{AB}{A_1B_1}=\frac{BC}{B_1C_1},
\]
Подобным же образом доказывается пропорциональность и прочих отрезков.

\paragraph{}\label{1938/184}
\so{Задача}.
\emph{Разделить отрезок прямой $AB$ \emph{(рис.~\ref{1938/ris-194})} на три части в отношении $m:n:p$, где $m$, $n$, $p$ — данные отрезки или данные числа.}

\begin{figure}[!ht]
\centering
\includegraphics{mppics/ris-194}
\caption{}\label{1938/ris-194}
\end{figure}

Проведя луч $AC$ под произвольным углом к $AB$, отложим на нём от точки $A$ отрезки, равные отрезкам $m$, $n$ и $p$.
Точку $F$ — конец отрезка $p$ — соединяем с $B$ прямой $BF$ и через концы $G$ и $H$ отложенных отрезков проводим прямые $GD$ и $HE$, параллельные $BF$.
Тогда отрезок $AB$ разделится в точках $D$ и $E$ на части в отношении $m:n:p$.

Если $m$, $n$ и $p$ означают какие-нибудь числа, например 2, 5, 3, то построение выполняется так же, с той лишь разницей, что на $AC$ откладываются отрезки, равные 2, 5 и 3 произвольным единицам длины.

Конечно, указанное построение применимо к делению отрезка не только на три части, но на какое угодно иное число частей.

\paragraph{}\label{1938/185}
\so{Задача}.
\emph{К трём данным отрезкам $a$, $b$ и $c$ найти четвёртый пропорциональный} (рис.~\ref{1938/ris-195}), то есть
найти такой отрезок $x$, который удовлетворял бы пропорции $a:b=c:x$.

\begin{figure}[!ht]
\centering
\includegraphics{mppics/ris-195}
\caption{}\label{1938/ris-195}
\end{figure}

На сторонах произвольного угла $ABC$ откладываем отрезки:
$BD\z=a$, $BF=b$, $DE=c$.
Проведя затем через $D$ и $F$ прямую, построим $EG\parallel DF$.
Отрезок $FG$ будет искомый.

\subsection*{Свойство биссектрисы угла треугольника}

\paragraph{}\label{1938/186}
\so{Теорема}.
\textbf{\emph{Биссектриса}} ($BD$, рис.~\ref{1938/ris-196}) \textbf{\emph{любого угла треугольника}} ($ABC$) \textbf{\emph{делит противоположную сторону на части}} ($AD$ и $CD$), \textbf{\emph{пропорциональные прилежащим сторонам треугольника.}}

\begin{wrapfigure}{r}{30mm}
\vskip-4mm
\centering
\includegraphics{mppics/ris-196}
\caption{}\label{1938/ris-196}
\end{wrapfigure}

Требуется доказать, что если $\angle ABD\z=\angle DBC$, то 
\[\frac{AD}{DC}\z=\frac{AB}{BC}.\]

Проведём $CE \parallel BD$ до пересечения в точке $E$ с продолжением стороны $AB$.
Тогда, согласно теореме в §~\ref{1938/182}, мы будем иметь пропорцию:
\[\frac{AD}{DC}=\frac{AB}{BE}.\]
Чтобы от этой пропорции перейти к той, которую требуется доказать, достаточно обнаружить, что $BE=BC$, то есть что $\triangle BCE$ равнобедренный.
В этом треугольнике $\angle E=\angle ABD$ (как углы, соответственные при параллельных прямых) и $\angle BCE \z= \angle DBC$ (как углы, накрест лежащие при тех же параллельных прямых).
Но $\angle ABD=\angle DBC$ по условию;
значит, $\angle E = \angle BCE$, а потому равны и стороны $BC$ и $BE$, лежащие против равных углов.
Заменив в написанной выше пропорции $BE$ на $BC$, получим ту пропорцию, которую требуется доказать.

\medskip

\smallskip
\so{Численный пример}.
Пусть $AB = 10$,
$BC = 7$ и $AC \z= 6$.
Тогда, обозначив $AD$ буквой $x$, можем написать пропорцию:
\[\frac{x}{6 - x} = \frac{10}7;\]
отсюда найдём:
\begin{align*}
7x&=60-10x;
\\
7x+10x&=60;
\\
17x&=60;
\\
x&=\tfrac{60}{17}=3\tfrac9{17}.
\end{align*}
Следовательно,
\[DC=6-x=6-3\tfrac9{17}=2\tfrac8{17}.\]

{\small

\paragraph{}\label{1938/187}
\so{Теорема} (выражающая свойство биссектрисы внешнего угла треугольника).
\textbf{\emph{Если биссектриса}} ($BD$, рис.~\ref{1938/ris-197}) \textbf{\emph{внешнего угла}} ($CBF$) \textbf{\emph{треугольника}} ($ABC$) \textbf{\emph{пересекает продолжение противоположной стороны}} ($AC$) \textbf{\emph{в некоторой точке}} ($D$), \textbf{\emph{тогда расстояния}} ($AD$ и $DC$) \textbf{\emph{от этой точки до концов этой стороны пропорциональны прилежащим сторонам треугольника}} ($AB$ и $BC$).
Требуется доказать, что если $\angle CBD\z=\angle FBD$, то $DA:DC=AB:BC$.

\begin{wrapfigure}[9]{o}{45mm}
\vskip-4mm
\centering
\includegraphics{mppics/ris-197}
\caption{}\label{1938/ris-197}
\end{wrapfigure}

Проведя $CE \parallel BD$, получим пропорцию
\[\frac{DA}{DC}=\frac{BA}{BE}.\]

Так как $\angle BEC=\angle FBD$ (как соответственные), а $\angle BCE\z=\angle CBD$ (как накрест лежащие при параллельных прямых) и углы $FBD$ и $CBD$ равны по условию, то $\angle BEC\z=\angle BCE$;
значит, $\triangle BCE$ равнобедренный, то есть $BE=BC$.
Подставив в пропорцию вместо $BE$ равный отрезок $BC$, получим ту пропорцию, которую требовалось доказать:
\[\frac{DA}{DC}=\frac{AB}{BC}.\]

{\small 

\smallskip
\so{Примечание}.
Особый случай представляет биссектриса внешнего угла при вершине равнобедренного треугольника, которая параллельна основанию.

}

\paragraph{}\label{1914/227}
\so{Обратная теорема.}
\textbf{\emph{Если прямая, исходящая из вершины треугольника, пересекает противоположную сторону (или её продолжение) в точке, расстояния от которой до концов противоположной стороны пропорциональны соответственно двум другим сторонам, то она есть биссектриса угла треугольника (внутреннего или внешнего).}}

Пусть  $E$ есть точка, лежащая на стороне $AC$ треугольника $ABC$ или на её продолжении, такая, что
$AE:EC=AB:BC$.
По доказанному (§~\ref{1938/186}), для основания $D$ биссектрисы угла $B$, а также для основания $D'$ биссектрисы его внешнего угла к $B$ выполнена та же пропорция, то есть
\[\frac{AD}{DC}=\frac{AD'}{D'C}=\frac{AB}{BC}.\]

Но если $AB\ne BC$, то существуют только две точки (§~\ref{extra/proportions}).
Значит $E$ совпадает с $D$ или $D'$ и следовательно прямая $BE$ сливается с биссектрисой угла $B$ или его внешнего угла.

Если $AB=BC$, то точка $E$ есть середина основания $AC$ и треугольник $ABC$ равнобедренный. 
То есть $BE$ есть медиана, проведённая к основанию равнобедренного треугольника,
по доказанному (§~\ref{1938/38}), совпадает с биссектрисой угла $B$.
(Биссектриса внешнего угла в этом случае параллельна основанию.)

\paragraph{}\label{1914/228}
\so{Теорема}.
\textbf{\emph{Геометрическое место точек, до которых расстояния от двух данных точек находятся в постоянном отношении $\bm{m:n}$, есть окружность, когда $\bm{m\ne n}$, и прямая, когда $\bm{m=n}$.}}

Обозначим данные точки буквами $A$ и $B$.

Если $m\z=n$, то искомое место точек есть срединный перпендикуляр к отрезку $AB$ (§~\ref{1938/58}).

\begin{wrapfigure}{o}{50mm}
\centering
\includegraphics{mppics/ris-1914-209}
\caption{}\label{1914/ris-209}
\end{wrapfigure}

Предположим, что $m>n$.
Тогда на прямой $AB$, можно найти две точки $C$ и $C'$ (рис. \ref{1914/ris-209}),
принадлежащие искомому геометрическому месту (§~\ref{extra/proportions}); то есть такие, что
\[\frac{CA}{CB}=\frac{C'A}{C'B}=\frac mn.\]
Точка $C$ лежит на отрезке $AB$, а точка $C'$ на его продолжении.

Пусть ещё другая точка $M$ удовлетворяет пропорции
\[\frac{MA}{MB}=\frac mn.\eqno(1)\]
Проведя $MC$ и $MC'$ мы должны заключить (§~\ref{1914/227}), что первая из этих прямых есть биссектриса угла $AMB$, а вторая — биссектриса угла $BMN$;
вследствие этого угол $CMC'$, составленный из двух половин смежных углов, должен быть прямой.
Поэтому вершина $M$ лежит на окружности, описанной на $CC'$ как на диаметре (§~\ref{1938/125}).

Таким образом, мы доказали, что всякая точка $M$ удовлетворяющая пропорции (1), лежит на окружности с диаметром $CC'$.
Докажем обратное предложение, то есть, что всякая точка $M$, лежащая на окружности с диаметром $CC'$, удовлетворяет пропорции (1).
Обозначим центр этой окружности буквой $O$, это середина отрезка $CC'$.

Найдём сначала пропорции $\frac{OA}{OC}$ и $\frac{OC}{OB}$.
Поскольку
\[\frac{AC}{BC}=\frac{AC'}{BC'}=\frac{m}{n}\]
верны также следующие пропорции
\[\frac{AC'+AC}{BC'+BC}=\frac{AC'-AC}{BC'-BC}=\frac{m}{n}. \eqno(2)\]

Поскольку $OC=OC'$,
\begin{align*}
2\cdot AO&=AC'-OC'+AC+OC=
\\
&=AC'+AC.
\\
2\cdot BO&=BC'-OC'+OC-BC=
\\
&=BC'-BC.
\\
2\cdot CO&=CC'=AC'-AC=
\\
&=BC'+BC.
\end{align*}
Подставляя эти значения в пропорции (2) получаем
\[\frac{OA}{OC}=\frac{OC}{OB}=\frac mn.\]

Пусть $M$ произвольная точка окружности отличная от $C$ и $C'$.
Тогда $OM=OC$, следовательно
\[\frac{OA}{OM}=\frac{OM}{OB}=\frac mn.\]
Значит $\triangle OBM\sim\triangle OMA$ и в частности
\[\frac{AM}{BM}=\frac{OA}{OM}=\frac mn,\]
то есть выполняется пропорция (1).

Случай $m<n$ аналогичен.

\so{Замечание}.
Окружность, о которой говорится в этой теореме, известна под названием \rindex{окружность!Аполлония}\textbf{окружность Аполлония} (Аполлоний Пергский — греческий геометр, живший во II вeке до нашей эры).

}

%% file: 2D/teorema-pifagora.tex
\section{Теорема Пифагора}

\paragraph{}\label{1938/188}
\so{Теорема}.
\textbf{\emph{В прямоугольном треугольнике перпендикуляр, опущенный из вершины прямого угла на гипотенузу, есть средняя пропорциональная величина между отрезками, на которые основание перпендикуляра делит гипотенузу, а каждый катет есть средняя пропорциональная между гипотенузой и прилежащим к этому катету отрезком гипотенузы.}}

Пусть $AD$ (рис.~\ref{1938/ris-198}) есть перпендикуляр, опущенный из вершины прямого угла $A$ на гипотенузу $BC$.
Требуется доказать следующие три пропорции:
\begin{align*}
1)\ \frac{BD}{AD}&=\frac{AD}{DC};
&
2)\ \frac{BC}{AB}&=\frac{AB}{BD};
&
3)\ \frac{BC}{AC}&=\frac{AC}{DC}.
\end{align*}

\begin{wrapfigure}{o}{45mm}
\centering
\includegraphics{mppics/ris-198}
\caption{}\label{1938/ris-198}
\end{wrapfigure}

Первую пропорцию мы докажем из подобия $\triangle ABD\sim \triangle CAD$.
Эти треугольники подобны, потому что
\[\angle 1 = \angle 4\quad \text{и} \quad\angle 2 = \angle 3.\]
вследствие перпендикулярности их сторон (§~\ref{1938/80}).
Возьмём в $\triangle ABD$ те стороны $BD$ и $AD$, которые составляют первое отношение доказываемой пропорции;
соответственными сторонами в $\triangle CAD$ будут $AD$ и $DC$%
\footnote{Чтобы безошибочно определить, какие стороны взятых треугольников соответственны между собой, полезно держаться такого пути:

1) указать углы, против которых лежат взятые стороны одного треугольника.

2) найти равные им углы в другом треугольнике.

3) взять противолежащие им стороны.

Например, для треугольников $ABD$ и $CAD$ рассуждаем так:
в треугольнике $ABD$ стороны $BD$ и $AD$ лежат против углов 1 и 3;
в треугольнике $CAD$ этим углам равны 4 и 2;
против них лежат стороны $AD$ и $DC$.
Значит, стороны $AD$ и $DC$ соответственны сторонам $BD$ и $AD$.
}%
, поэтому
\[\frac{BD}{AD}=\frac{AD}{DC}.\]

Вторую пропорцию докажем из подобия $\triangle ABC\sim \triangle DBA$.
Эти треугольники подобны, потому что они прямоугольные и острый угол $B$ у них общий.
В $\triangle ABC$ возьмём те стороны $BC$ и $AB$, которые составляют первое отношение доказываемой пропорции;
соответственными сторонами в $\triangle DBA$ будут $AB$ и $BD$;
поэтому
\[\frac{BC}{AB}=\frac{AB}{BD}.\]

Третью пропорцию докажем из подобия $\triangle ABC\sim \triangle DAC$.
Эти треугольники подобны, потому что они оба прямоугольные и имеют общий острый угол $C$.
В $\triangle ABC$ возьмём стороны $BC$ и $AC$;
соответственными сторонами в $\triangle DAC$ будут $AC$ и $DC$;
поэтому
\[\frac{BC}{AC}=\frac{AC}{DC}.\]

\begin{figure}
\begin{minipage}{.32\textwidth}
\centering
\includegraphics{mppics/ris-199}
\end{minipage}\hfill
\begin{minipage}{.32\textwidth}
\vskip6mm
\centering
\includegraphics{mppics/ris-200}
\end{minipage}\hfill
\begin{minipage}{.32\textwidth}
\vskip4mm
\centering
\includegraphics{mppics/ris-201}
\end{minipage}

\medskip

\begin{minipage}{.32\textwidth}
\centering
\caption{}\label{1938/ris-199}
\end{minipage}
\begin{minipage}{.32\textwidth}
\vfill
\centering
\caption{}\label{1938/ris-200}
\end{minipage}
\begin{minipage}{.32\textwidth}
\vfill
\centering
\caption{}\label{1938/ris-201}
\end{minipage}
\vskip-4mm
\end{figure} 

\paragraph{}\label{1938/189}
\mbox{\so{Следствие}.}
Пусть $A$ (рис.~\ref{1938/ris-199}) есть произвольная точка окружности, описанной на диаметре $BC$.
Соединив концы диаметра с этой точкой, мы получим прямоугольный треугольник $ABC$, у которого гипотенуза есть диаметр, а катеты образованы хордами (§~\ref{1938/125}, \ref{1938/2}).
Применяя доказанную выше теорему к этому треугольнику, приходим к следующему заключению:

\emph{Перпендикуляр, опущенный из какой-либо точки окружности на диаметр, есть средняя пропорциональная величина между отрезками, на которые основание перпендикуляра делит диаметр, а хорда, соединяющая эту точку с концом диаметра, есть средняя пропорциональная между диаметром и прилежащим к хорде отрезком диаметра.}

{\sloppy 

\paragraph{}\label{1938/190}
\mbox{\so{Задача}.}
\emph{Построить отрезок, средний пропорциональный между двумя отрезками $a$ и $b$,}
то есть построить отрезок $x$, такой, что $a:x= x:b$.

}

Приведём два решения.

1) на произвольной прямой (рис.~\ref{1938/ris-200}) откладываем отрезки $AB=a$ и $BC\z=b$;
на $AC$, как на диаметре, описываем полуокружность;
из $B$ восстанавливаем перпендикуляр $BD$ до пересечения с окружностью.
Этот перпендикуляр и есть искомая средняя пропорциональная между $AB$ и $BC$.

2) на произвольной прямой (рис.~\ref{1938/ris-201}) откладываем от точки $A$ отрезки $a$ и $b$.
На большем из этих отрезков описываем полуокружность.
Проведя из конца меньшего отрезка перпендикуляр к $AB$ до пересечения его с окружностью в точке $D$, соединяем $A$ с $D$.
Хорда $AD$ есть средняя пропорциональная между $a$ и $b$.

\paragraph{Теорема Пифагора.}\label{1938/191}
Доказанные выше теоремы позволяют обнаружить замечательное соотношение между сторонами любого прямоугольного треугольника.
Это соотношение названо в честь греческого геометра Пифагора (VI век до нашей эры).

\textbf{\emph{Квадрат длины гипотенузы равен сумме квадратов длин катетов.}} 

Пусть $ABC$ (рис.~\ref{1938/ris-202}) есть прямоугольный треугольник, $AD$ — перпендикуляр, опущенный на гипотенузу из вершины прямого угла.
Положим, что стороны и отрезки гипотенузы измерены одной и той же единицей, причём получились числа $a$, $b$, $c$, $c'$ и $b'$ (принято длины сторон треугольника обозначать малыми буквами, соответствующими большим буквам, которыми обозначены противолежащие углы).

\begin{wrapfigure}{o}{47mm}
\vskip-2mm
\centering
\includegraphics{mppics/ris-202}
\caption{}\label{1938/ris-202}
\end{wrapfigure}

Применяя теорему §~\ref{1938/188}, можем написать пропорции:
\[\frac ac=\frac c{c'}
\quad\text{и}\quad
\frac ab=\frac b{b'}\]
откуда 
\[a\cdot c'=c^2
\quad\text{и}\quad
a\cdot b'=b^2.\]
Сложив почленно эти два равенства, найдём:
\[a\cdot c'+a\cdot b'=c^2+b^2
\quad\text{или}\quad
a\cdot (c'+b')=c^2+b^2.\]
Но $c'+b'=a$, следовательно,
\[a^2=c^2+b^2.\]

\smallskip
\so{Пример}.
Если катеты, измеренные какой-нибудь линейной единицей, равны 3 и 4; 
тогда гипотенуза  $x$, удовлетворяет уравнению:
\[x^2=3^2+4^2=9+16=25,\]
откуда $x = \sqrt{25} = 5$.

{\small
\smallskip
\so{Замечание}.
Прямоугольный треугольник со сторонами 3, 4 и 5 называется \rindex{египетский треугольник}\textbf{египетским треугольником}. 
Им можно пользоваться для построения прямого угла на земной поверхности:
бечёвку посредством узлов разделяем на 12 равных частей;
затем, связав концы, натягиваем её на земле (посредством кольев) в виде треугольника со сторонами в 3, 4 и 5 делений;
тогда угол между сторонами, равными 3 и 4, оказывался прямым.%
\footnote{Прямоугольные треугольники, у которых стороны измеряются целыми числами, носят название \so{пифагоровых треугольников}.
Прямые вычисления показывают, что если $a$ и $b$ — произвольные целые числа такие, что $a>b$ то для
\[x=2ab,
\quad
y=a^2-b^2,
\quad
z=a^2+b^2\]
выполнено тождество $x^2+y^2=z^2$.
В частности треугольник со сторонами $x$, $y$ и $z$ является пифагоровым.
Можно доказать, также что любой пифагоров треугольник имеет стороны $k{\cdot}x$, $k{\cdot}y$ и $k{\cdot}z$ для некоторых целых $a$, $b$ и $k$. 

}
}

Теорема Пифагора имеет ещё другую формулировку, именно ту, которая была для неё получена самим Пифагором.
С этой формулировкой мы познакомимся позднее (§~\ref{1938/257}).

\paragraph{}\label{1938/192}
\so{Следствие}.
\emph{Квадраты катетов относятся между собой как прилежащие отрезки гипотенузы.}
Действительно, из уравнений предыдущего параграфа находим:
\[\frac{c^2}{b^2}=\frac{a\cdot c'}{a\cdot b'}=\frac{c'}{b'}.\]

{\small

\paragraph{}\label{1938/193}
\so{Замечание 1}.
К трём равенствам, которые мы вывели выше:
\[1)\ a\cdot c'=c^2;
\qquad
2)\ a\cdot b'=b^2
\quad
\text{и}
\quad
3)\ a^2=b^2+c^2,
\]
можно присоединить ещё следующие два:
\[4)\ b'+c'=a
\quad
\text{и}
\quad
5)\ h^2=b'\cdot c',
\]
(если буквой $h$ обозначим длину высоты $AD$).
Из этих равенств третье, как мы видели, составляет следствие первых двух и четвёртого, так что из пяти равенств только четыре независимы;
вследствие этого можно по данным двум из шести чисел находить остальные четыре.

Для примера положим, что нам даны отрезки гипотенузы $b' = 5 \text{м}$ и $c' \z= 7\text{м}$;
тогда
\begin{align*}
a&=b'+c'=12;
\\
c&=\sqrt{a\cdot c'}=
\sqrt{12\cdot 7}=
\sqrt{84}=9{,}165\dots;
\\
b&=\sqrt{a\cdot b'}=\sqrt{12\cdot 5}=\sqrt{60}=7{,}745\dots;
\\
h&=\sqrt{c'\cdot b'}=\sqrt{5\cdot 7}=\sqrt{35} = 5{,}916\dots
\end{align*}

}

%% file: 2D/zadach-na-vych.tex
\section{Задачи на вычисление}

\paragraph{}\label{1938/194} 
\so{Теорема}.
\textbf{\emph{Во всяком треугольнике квадрат стороны, лежащей против острого угла, равен сумме квадратов двух других сторон без удвоенного произведения какой-нибудь из этих сторон на отрезок её от вершины острого угла до высоты.}}

Пусть $BC$ — сторона $\triangle ABC$ (рис.~\ref{1938/ris-203} и \ref{1938/ris-204}), лежащая против острого угла $A$, и $BD$ — высота, опущенная на сторону $AC$ (или её продолжение).
Требуется доказать, что
\[BC^2=AB^2+AC^2-2\cdot AC\cdot  AD.\]
или, обозначая длины линий малыми буквами, как указано на рисунке, надо доказать равенство:
\[a^2=b^2+c^2-2\cdot b\cdot c'\]

Из прямоугольного $\triangle BDC$ находим:
\[a^2=h^2+(a')^2.
\eqno(1)\]

Найдём каждый из квадратов $h^2$ и $(a')^2$. 
Из прямоугольного $\triangle BAD$ находим:
\[h^2=c^2-(c')^2.
\eqno(2)\]

\begin{figure}[!ht]
\begin{minipage}{.48\textwidth}
\centering
\includegraphics{mppics/ris-203}
\end{minipage}
\hfill
\begin{minipage}{.48\textwidth}
\centering
\includegraphics{mppics/ris-204}
\end{minipage}

\medskip

\begin{minipage}{.48\textwidth}
\centering
\caption{}\label{1938/ris-203}
\end{minipage}
\hfill
\begin{minipage}{.48\textwidth}
\centering
\caption{}\label{1938/ris-204}
\end{minipage}
\vskip-4mm
\end{figure}

С другой стороны, $a'=b-c'$ (рис.~\ref{1938/ris-203})) или $a'=c'-b$ (рис.~\ref{1938/ris-204}).
В обоих случаях для $(a')^2$ получаем одно и то же выражение:
\[
\begin{aligned}
(a')^2&=(b-c')^2=b^2-2\cdot b\cdot c'+(c')^2;
\\
(a')^2&=(c'-b)^2=(c')^2-2\cdot b\cdot c'+b^2.
\end{aligned}
\eqno(3)
\]

Равенство (1) можно переписать так:
\[a^2=c^2-(c')^2+b^2-2\cdot b\cdot c'+(c')^2=c^2+b^2-2\cdot b\cdot c'.\]

{\sloppy
\paragraph{}\label{1938/195}
\mbox{\so{Теорема}.}
\textbf{\emph{В тупоугольном треугольнике квадрат стороны, лежащей против тупого угла, равен сумме квадратов двух других сторон, сложенной с удвоенным произведением какой-нибудь из этих сторон на отрезок её продолжения от вершины тупого угла до высоты.}}

}

Пусть $AB$ — сторона $\triangle ABC$ (рис. \ref{1938/ris-205}), лежащая против тупого угла $C$, и $BD$ — высота, опущенная на продолжение стороны $AC$;
требуется доказать, что
\[AB^2=AC^2+BC^2+2\cdot AC \cdot  CD,\]
или, применяя сокращённые обозначения, согласно указанию на рисунке:
\[c^2=a^2+b^2+2\cdot b\cdot a'.\]

\begin{wrapfigure}{o}{45mm}
\vskip-8mm
\centering
\includegraphics{mppics/ris-205}
\caption{}\label{1938/ris-205}
\vskip0mm
\end{wrapfigure}

Из треугольников $ABD$ и $CBD$ находим:
\begin{align*}
c^2&=h^2+(c')^2=
\\
&=a^2-(a')^2+(a'+b)^2=
\\
&=a^2-(a')^2+(a')^2+2\cdot b\cdot a'+b^2=
\\
&=a^2+b^2+2\cdot b\cdot a',
\end{align*}
что и требовалось доказать.

\paragraph{}\label{1938/196}
\so{Следствие}.
Из трёх последних теорем выводим, что \emph{квадрат стороны треугольника равен, меньше или больше суммы квадратов двух других сторон, смотря по тому, будет ли противолежащий угол прямой, острый или тупой.}
Отсюда следует обратное предложение:
\emph{угол треугольника окажется прямым, острым или тупым, смотря по тому, будет ли квадрат противолежащей этому углу стороны равен, меньше или больше суммы квадратов двух других сторон.}

\paragraph{}\label{1938/197}
\so{Теорема}.
\textbf{\emph{Сумма квадратов диагоналей параллелограмма равна сумме квадратов его сторон}} (рис.~\ref{1938/ris-206}).

\begin{wrapfigure}{O}{48mm}
\centering
\includegraphics{mppics/ris-206}
\caption{}\label{1938/ris-206}
\end{wrapfigure}

Из вершин $B$ и $C$ параллелограмма $ABCD$ опустим на основание $AD$ перпендикуляры $BE$ и $CF$.
Тогда из треугольников $ABD$ и $ACD$ находим:
\[BD^2=AB^2+AD^2-2\cdot AD\cdot AE\]
\[AC^2=AD^2+CD^2+2\cdot AD\cdot  DF.\]

Прямоугольные треугольники $ABE$ и $DCF$ равны, так как они имеют по равной гипотенузе и равному острому углу;
поэтому $AE=DF$.
Заметив это, сложим почленно два выведенных выше равенства;
тогда $2AD\cdot  AE$ и $2AD\cdot  DF$ взаимно уничтожаются, и мы получим:
\begin{align*}
BD^2+AC^2&=AB^2+AD^2+AD^2+CD^2=
\\
&=AB^2+BC^2+CD^2+AD^2.
\end{align*}

\paragraph{Вычисление медианы треугольника.}\label{1914/241}
Медиана треугольника обыкновенно обозначается буквой $m$ (от латинского слова \emph{mediāna} — средняя), сопровождаемою (внизу) одною из маленьких букв $a$, $b$ или $c$ в зависимости от стороны треугольника, к которой проведена обозначаемая медиана.

\begin{wrapfigure}{r}{35mm}
\centering
\includegraphics{mppics/ris-1914-221}
\caption{}\label{1914/ris-221}
\end{wrapfigure}

Определим длину $m_a$ медианы, проведённой к стороне $a$ (рис.~\ref{1914/ris-221}).
Для этого продолжим медиану на расстояние $DE=AD$ и соединим точку $E$ с $B$ и с~$C$.
Мы получим параллелограмм $ABEC$ (§~\ref{1938/90}).
Применив к нему теорему о сумме квадратов диагоналей (§~\ref{1938/197}), получим:
\[a^2+(2m_a)^2=2b^2+2c^2;\]
откуда: 
\[4m_a^2=2b^2+2c^2-a^2\]
и, следовательно:
\[m_a=\tfrac 12\sqrt{2b^2+2c^2-a^2}\]

Подобным же образом можем найти $m_b$ и $m_c$.

\begin{figure}[!ht]
\begin{minipage}{.48\textwidth}
\centering
\includegraphics{mppics/ris-207}
\end{minipage}
\hfill
\begin{minipage}{.48\textwidth}
\centering
\includegraphics{mppics/ris-208}
\end{minipage}

\medskip

\begin{minipage}{.48\textwidth}
\centering
\caption{}\label{1938/ris-207}
\end{minipage}
\hfill
\begin{minipage}{.48\textwidth}
\centering
\caption{}\label{1938/ris-208}
\end{minipage}
\vskip-4mm
\end{figure}

\paragraph{Вычисление высот треугольника по его сторонам.}\label{1938/198}
Определим высоту $h_a$ треугольника $ABC$, опущенную на сторону $BC=a$ (рис.~\ref{1938/ris-207} и \ref{1938/ris-208}).
Обозначим отрезки стороны $a$ (продолженной в случае тупого угла $C$, рис.~\ref{1938/ris-208}) таким образом:
отрезок $BD$, прилежащий к стороне $c$, через $c'$, а отрезок $DC$, прилежащий к стороне $b$, через $b'$.
Пользуясь теоремой о квадрате стороны треугольника, лежащей против острого угла (§~\ref{1938/194}), можем написать:
\[b^2=a^2+c^2-2\cdot a\cdot c'.\]
Из этого уравнения находим отрезок $c'$:
\[c'=\frac{a^2+c^2-b^2}{2\cdot a}\]
после чего из треугольника $ABD$ определяем высоту как катет:
\[h_a=\sqrt{c^2-\left(\tfrac{a^2+c^2-b^2}{2\cdot a}\right)^2}\]
Таким же путём можно определить в зависимости от сторон треугольника длины $h_b$, и $h_c$ высот, опущенных на стороны $b$ и $c$.

%% file: 2D/proportzii-v-kruge.tex
\section{Пропорции в круге}

\paragraph{}\label{1938/199}
Некоторые пропорциональные линии в круге мы указали ранее (§~\ref{1938/189});
теперь укажем ещё другие.

\smallskip
\mbox{\so{Теорема}.}
\textbf{\emph{Если через точку}} ($M$, рис.~\ref{1938/ris-209}), \textbf{\emph{взятую внутри круга, проведены какая-нибудь хорда}} ($AB$) \textbf{\emph{и диаметр}} ($CD$), \textbf{\emph{то произведение отрезков хорды}} ($AM\cdot  MB$) \textbf{\emph{равно произведению отрезков диаметра}} ($MD\cdot  MC$).

Проведя две вспомогательные хорды $AC$ и $BD$, мы получим два треугольника $AMC$ и $DMB$ (покрытые на рисунке штрихами), которые подобны, так как у них углы $A$ и $D$ равны как вписанные, опирающиеся на одну и ту же дугу $BC$, и углы $C$ и $B$ равны как вписанные, опирающиеся на одну и ту же дугу $AD$.

\begin{wrapfigure}{r}{45mm}
\vskip-0mm
\centering
\includegraphics{mppics/ris-209}
\caption{}\label{1938/ris-209}
\end{wrapfigure}

Из подобия треугольников выводим:
\[\frac{AM}{MD}=\frac{MC}{MB}.\]
откуда
\[AM\cdot  MB=MD\cdot  MC.\]

\paragraph{}\label{1938/200}
\mbox{\so{Следствие}.}
\emph{Если через точку \emph{($M$, рис.~\ref{1938/ris-209}),} взятую внутри круга, проведено сколько угодно хорд ($AB$, $EF$, $KL,\dots$), то произведение отрезков каждой хорды есть число постоянное для всех хорд,} так как для каждой хорды это произведение равно произведению отрезков диаметра $CD$, проходящего через взятую точку $M$.

\begin{wrapfigure}{r}{50mm}
\vskip-6mm
\centering
\includegraphics{mppics/ris-210}
\caption{}\label{1938/ris-210}
\end{wrapfigure}

\paragraph{}\label{1938/201}
\mbox{\so{Теорема}.}
\textbf{\emph{Если из точки}} ($M$, рис.~\ref{1938/ris-210}), \textbf{\emph{взятой вне круга, проведены к нему какая-нибудь секущая}} ($MA$) \textbf{\emph{и касательная}} ($MC$), \textbf{\emph{то произведение секущей на её внешнюю часть равно квадрату касательной}} (предполагается, что секущая ограничена второй точкой пересечения, а касательная — точкой касания).

Проведём вспомогательные хорды $AC$ и $BC$;
тогда получим два треугольника $MCA$ и $MBC$ (покрытые на чертеже штрихами), которые подобны, потому что у них угол $M$ общий и углы $MCB$ и $CAB$ равны, так как каждый из них измеряется половиной дуги $BC$.

Возьмём в $\triangle MAC$ стороны $MA$ и $MC$;
соответственными сторонами в $\triangle MBC$ будут $MC$ и $MB$;
поэтому
\[\frac{MA}{MC} = \frac{MC}{MB},
\qquad\text{откуда}\qquad
MA\cdot MB=MC^2.\]

{\small

\paragraph{}\label{1938/202}
\mbox{\so{Следствие}.}
\emph{Если из точки \emph{($M$, рис.~\ref{1938/ris-210}),} взятой вне круга, проведены к нему сколько угодно секущих ($MA$, $MD$, $ME,\dots$), то произведение каждой секущей на её внешнюю часть есть число постоянное для всех секущих, так как для каждой секущей это произведение равно квадрату касательной ($MC^2$), проведённой из точки $M$.}


Величина $d^2- R^2$, где $d$ — расстояние от точки до центра окружности, а $R$ — её радиус называется \rindex{степень точки}\textbf{степенью точки} относительно окружности.
Заметим, что степень отрицательна для точек внутри круга,
положительна вне его и обращается в ноль на самой окружности.

Используя понятие степени точки, следствия, приведённые в §§ \ref{1938/200} и \ref{1938/202}, можно переформулировать следующим образом:
\emph{Для любой хорды $AB$ и любой точки $M$ на ней или её продолжении, произведение
$MA\cdot MB$
равно абсолютной величине степени $M$ относительно окружности}.

\begin{figure}[!ht]
\begin{minipage}{.48\textwidth}
\centering
\includegraphics{mppics/ris-1914-228}
\end{minipage}
\hfill
\begin{minipage}{.48\textwidth}
\centering
\includegraphics{mppics/ris-1914-229}
\end{minipage}

\medskip

\begin{minipage}{.48\textwidth}
\centering
\caption{}\label{1914/ris-228}
\end{minipage}
\hfill
\begin{minipage}{.48\textwidth}
\centering
\caption{}\label{1914/ris-229}
\end{minipage}
\vskip-4mm
\end{figure}

\paragraph{}\label{1914/250}\so{Теорема}. \textbf{\emph{Произведение двух сторон треугольника равно произведению диаметра круга, описанного около этого треугольника, на высоту его, опущенную на третью сторону.}}

Обозначив буквою $R$ радиус круга, описанного около $\triangle ABC$ (рис. \ref{1914/ris-228} и \ref{1914/ris-229}), докажем, что
\[b\cdot c=2R\cdot h_a.\]

Проведём диаметр $AD$ и соединим $D$ с $B$.
Треугольники $ABD$ и $AEC$ подобны, потому что углы $B$ и $E$ прямые и $\angle D=\angle C$, как углы вписанные, опирающиеся на одну и ту же дугу.
Из подобия выводим:
\begin{align*}
\frac{c}{h_a}&=\frac{2R}{b};
&
&\text{откуда:}
&
b\cdot c&=2R\cdot h_a.
\end{align*}

}

%% file: 2D/algebra-k-geomtrii.tex
{\small

\section{О приложении алгебры к геометрии}

\paragraph{}\label{1938/209}
\so{Задача}.
\emph{Данный отрезок разделить в среднем и крайнем отношении.}

Эту задачу надо понимать так:
разделить данный отрезок на такие две части, чтобы б\'{о}льшая из них была средней пропорциональной между всей линией и меньшей её частью.

Задача будет решена, если мы найдём одну из двух частей, на которые требуется разделить данный отрезок.
Будем находить б\'{о}льшую часть, то есть ту, которая должна быть \so{средней пропорциональной} между всем отрезком и меньшей его частью.
Предположим сначала, что речь идёт не о построении этой части, а только о \so{вычислении} её длины.
Тогда задача решается \so{алгебраически} так:
если длину данного отрезка обозначим $a$, а длину б\'{о}льшей части $x$, то длина другой части будет равна $a-x$ и, согласно требованию задачи, мы будем иметь пропорцию:
\[\frac ax=\frac x{(a-x)}.\]
Откуда
\[x^2=a(a-x),
\quad\text{или}\quad
x^2+ax-a^2=0.\]
Решив это квадратное уравнение, находим:
\[
x_1=-\frac a2 + \sqrt{\left(\frac a2\right)^2+a^2};
\quad
x_2=-\frac a2 - \sqrt{\left(\frac a2\right)^2+a^2}.
\]
Отбросив второе решение, как отрицательное, возьмём только первое положительное решение, которое удобнее представить так:
\begin{align*}
x_1&= \sqrt{\left(\frac a2\right)^2+a^2}-\frac a2 =
\\
&= \sqrt{\frac {a^2}4+a^2}-\frac a2 =
\\
&= \sqrt{\frac {5a^2}4}-\frac a2 =
\\
&= \frac a2\sqrt{5}-\frac a2 =
\\
&= \frac {a(\sqrt{5}-1)}2 \approx
\\
&\approx a\cdot 0{,}61803.
\end{align*}
Таким образом, задача всегда имеет решение и притом только одно.

Если бы нам удалось построить такой отрезок, длина которого выражается найденной выше формулой, то нанеся этот отрезок на данный отрезок, мы разделили бы его в среднем и крайнем отношении.
Итак, вопрос сводится к построению найденной формулы.
Построить эту формулу всего удобнее, если мы её возьмём в том виде, в каком она была до упрощения, то есть возьмём:
\[x_1= \sqrt{\left(\frac a2\right)^2+a^2}-\frac a2.\]

\begin{wrapfigure}{o}{45mm}
\centering
\includegraphics{mppics/ris-217}
\caption{}\label{1938/ris-217}
\end{wrapfigure}

\noindent
Выражение $\sqrt{(\frac a2)^2+a^2}$ представляет собой длину гипотенузы такого прямоугольного треугольника, у которого один катет равен $a$, а другой $\frac a2$.
Построив такой треугольник, мы найдём отрезок, выражаемый формулой $\sqrt{(\frac a2)^2+a^2}$.
Чтобы получить затем отрезок $x_1$, достаточно из гипотенузы построенного треугольника вычесть $\frac a2$.
Таким образом, построение можно выполнить так:

Делим (рис.~\ref{1938/ris-217}) данный отрезок $AB=a$ пополам в точке $C$.
Из конца $B$ восстанавливаем перпендикуляр и откладываем на нём $BD = BC$.
Соединив $A$ с $D$ прямой, получим прямоугольный $\triangle ABD$, у которого один катет $AB=a$, а другой катет $BD=\frac a2$.
Следовательно, его гипотенуза $AD$ равна $\sqrt{(\frac a2)^2+a^2}$.
Чтобы вычесть из гипотенузы длину $\frac a2$, опишем дугу радиусом $BD=\frac a2$ с центром в точке $D$.
Тогда отрезок $AE$ будет равен $\sqrt{(\frac a2)^2+a^2}-\frac a2$, 
то есть будет равен $x_1$.
Отложив $AE$ на $AB$ (от $A$ до $G$), получим точку $G$, в которой отрезок $AB$ делится в среднем и крайнем отношении%
\footnote{Деление отрезка в среднем и крайнем отношении известно под названием \rindex{золотое сечение}\textbf{золотого сечения}.}%
.

{\small
\smallskip
\so{Замечание}.
Деление отрезка в среднем и крайнем отношении нужно в геометрии для построения правильного 10-угольника, вписанного в данный круг.
}

\paragraph{Алгебраический способ решения геометрических задач.}\label{1938/210}
Мы решили предложенную задачу путём приложения алгебры к геометрии.
Этот приём состоит в следующем:
сперва определяют, какой отрезок прямой нужно отыскать, чтобы решить задачу.
Затем, обозначив длины данных отрезков буквами $a, b, c,\dots$, а искомого буквой $x$, составляют из условий задачи и известных теорем уравнение, связывающее длину искомого отрезка с длинами данных, и полученное уравнение решают по правилам алгебры.
Найденную формулу исследуют, то есть определяют, при всяких ли данных $a, b, c,\dots$ эта формула определяет решение или только при некоторых и получается ли одно решение или несколько.
Затем строят формулу, то есть находят построением такой отрезок, длина которого выражается этой формулой.

Таким образом, алгебраический приём решения геометрических задач состоит в общем из следующих четырёх частей:
1) составление уравнения;
2) решение его;
3) исследование полученной формулы и 4) построение её.

Иногда задача сводится к отысканию нескольких отрезков.
Тогда, обозначив их длины буквами $x,y,\dots$, стремятся составить столько уравнений, сколько неизвестных.

\paragraph{Построение простейших формул.}\label{1938/211}
Укажем простейшие формулы, которые можно построить посредством циркуля и линейки;
при этом будем предполагать, что буквы $a$, $b$, $c,\dots$
означают длины данных отрезков, а $x$ — длину искомого.
Не останавливаясь на таких формулах:
\[x=a+b+c,
 \quad
 x=a-b,
 \quad
 x=2\cdot a,\ 3\cdot a,\dots,
\]
построение которых выполняется весьма просто, перейдём к более сложным:

1) Формулы $x=\frac a2, \frac a3,\dots$ и~т.~д. строятся посредством деления отрезка $a$ на равные части, затем, если нужно, повторением одной части слагаемых $2, 3,\dots$ и так далее раз.

2) Формула $x=\frac{ab}c$ выражает \so{четвёртую пропорциональную} к отрезкам $c$, $a$ и $b$.
Действительно, из этого равенства выводим:
\[c\cdot x=a\cdot b,
\quad\text{откуда}\quad
\frac ca=\frac bx.\]

Следовательно, $x$ найдётся способом, указанным нами для построения четвёртой пропорциональной (§~\ref{1938/185}).

3) Формула $x=\frac{a^2}b$ выражает четвёртую пропорциональную к отрезкам $b$, $a$ и $a$, или, как говорят, \so{третью пропорциональную} к отрезкам $b$ и $a$.
Действительно, из данного равенства выводим:
\[b\cdot x=a^2,
\quad\text{откуда}\quad
\frac ba=\frac ax.\]
Следовательно, $x$ найдётся тем же способом, каким отыскивается четвёртая пропорциональная (только отрезок $a$ придётся откладывать два раза).

4) Формула $x=\sqrt{a\cdot b}$ выражает \so{среднюю пропорциональную} между $a$ и $b$.
Действительно, из неё выводим:
\[x^2=a\cdot b,
\quad\text{откуда}\quad\frac ax=\frac xb.\]

Следовательно, $x$ найдётся способом, указанным раньше для построения средней пропорциональной (§~\ref{1938/190}).

5) Формула $x=\sqrt{a^2+b^2}$ выражает \so{гипотенузу} прямоугольного треугольника с катетами $a$ и $b$.

6) Формула $x=\sqrt{a^2-b^2}$ представляет катет прямоугольного треугольника с гипотенузой $a$ и другим катетом $b$.

Построение всего удобнее выполнить так, как указано в §~\ref{1938/126}.
Указанные формулы можно считать основными.
При помощи их строятся более сложные формулы.
Например:

7) $x=a\sqrt{\frac23}$.
Подведя $a$ под знак радикала, получим:
\[x=\sqrt{\tfrac23a^2}=\sqrt{a\cdot\tfrac23a}.\]
Отсюда видно, что $x$ есть средняя пропорциональная между отрезками $a$ и~$\tfrac23a$.

8) $x=\sqrt{a^2 + b^2 - c^2 + d^2}$.
Положим, что  $a^2+b^2=k^2$.
Тогда $k$ найдётся как гипотенуза прямоугольного треугольника с известными катетами $a$ и $b$.
Построив $k$, положим, что $k^2+d^2=\ell^2$.
Тогда $\ell$ найдётся как гипотенуза прямоугольного треугольника с катетами $k$ и $d$.
Построив $\ell$, будем иметь $x=\sqrt{\ell^2-c^2}$.
Следовательно, $x$ есть катет такого треугольника, у которого гипотенуза $\ell$, а другой катет~$c$.

Ограничимся этими примерами.
Заметим, что подробное рассмотрение способов построения алгебраических формул приводит к следующему важному выводу.

\emph{При помощи линейки и циркуля возможно строить только такие алгебраические выражения, которые могут быть получены из длин известных отрезков с помощью конечного числа арифметических операций \emph{(то есть сложения, вычитания, умножения и деления)} 
и извлечения квадратных корней}.

}

{\small

\subsection*{Упражнения}

\begin{center}
\so{Доказать теоремы}
\end{center}

\begin{enumerate}[noitemsep]

\item
Прямая, проведённая через середины оснований трапеции, проходит через точку пересечения продолжений боковых сторон и через точку пересечения диагоналей.

\item
Если в треугольнике из вершины угла, лежащего между неравными сторонами, проведены биссектриса и медиана, то первая меньше второй.

\item
Если два круга касаются извне, то часть внешней общей касательной, ограниченная точками касания, есть средняя пропорциональная между диаметрами кругов.

\item
Если на сторонах угла отложим от вершины пропорциональные отрезки, то прямые, соединяющие соответственные концы их, параллельны.

\item
Если в прямоугольный треугольник $ABC$ вписать квадрат $DEFG$ так, чтобы сторона $DE$ лежала на гипотенузе $BC$, то эта сторона есть средняя пропорциональная между отрезками гипотенузы $BD$ и $EC$ (точки на гипотенузе следуют в порядке $B$, $D$, $E$, $C$).

\item
Если два отрезка $AB$ и $CD$ пересекаются в точке $E$ так, что $BE\cdot  EA\z=EC\cdot  ED$, то точки $A$, $B$, $C$ и $D$ лежат на одной окружности (теорема, обратная изложенной в §~\ref{1938/200}).

\item
Дана окружность с центром $O$ и две точки $A$ и $B$.
Через эти точки проведено несколько окружностей, пересекающих окружность $O$ или касающихся её.
Доказать, что все хорды, соединяющие точки пересечения каждой из этих трёх окружностей с окружностью $O$, а также и общие касательные либо параллельны, либо сходятся (при продолжении) в одной точке, лежащей на прямой $AB$. 

\item
Основываясь на свойстве, изложенном в предыдущей задаче, вывести способ построения такой окружности, которая проходит через две данные точки $A$ и $B$ и касается данной окружности.

\item
Даны два каких-нибудь круга на плоскости.
Если два радиуса этих кругов движутся, оставаясь постоянно параллельными, то прямая, проходящая через концы их, пересекает линию центров всегда в одной и той же точке (в \so{центре подобия} двух кругов).

\item
Медиана треугольника делит пополам все прямые, проведённые внутри треугольника параллельно той стороне, относительно которой взята медиана.

\item
Даны три прямые, исходящие из одной точки.
Если по одной из них движется какая-нибудь точка, то расстояния от неё до двух других прямых сохраняют всегда одно и то же отношение.

\item
Если две окружности концентрические, то сумма квадратов расстояний от всякой точки одной из них до концов какого угодно диаметра другой есть величина постоянная (§~\ref{1938/197}).

\item
Если соединим прямыми основания трёх высот какого-нибудь треугольника, то образовавшиеся при этом три треугольника у вершин данного подобны ему.

Вывести отсюда, что для треугольника, имеющего сторонами прямые, соединяющие основания высот данного треугольника, эти высоты служат биссектрисами.

\item
Диаметр $AB$ данной окружности продолжен за точку $B$.
Через какую-нибудь точку $C$ этого продолжения проведена прямая $CD\z\perp AB$.
Если произвольную точку $M$ этого перпендикуляра соединим с $A$, (обозначив через $A_1$ вторую точку пересечения с окружностью этой прямой), то произведение $AM\cdot  AA_1$ есть величина постоянная для всякой точки $M$.

\end{enumerate}

\begin{center}
\so{Найти геометрические места}
\end{center}

\begin{enumerate}[resume,noitemsep]

\item
Середин всех хорд, проходящих через данную точку окружности.

\item
Точек, делящих в одном и том же отношении $m:n$ все хорды, проходящие через данную точку окружности.

\item \label{1938/upr-3-17}
Точек, расстояние от которых до сторон данного угла имеет одно и то же отношение $m:n$.

\item
Точек, для которых сумма квадратов расстояний до двух данных точек есть величина постоянная (§~\ref{1938/197}).

\item
Точек, для которых разность квадратов расстояний до двух данных точек есть величина постоянная.

\item
Точек, делящих в данном отношении $m:n$ все прямые, соединяющие точки окружности с данной точкой $O$ (лежащей вне или внутри круга).

\end{enumerate}

\begin{center}
\so{Задачи на построение}
\end{center}

\begin{enumerate}[resume,noitemsep]

\item
Через точку, данную внутри или вне угла, провести прямую так, чтобы отрезки её, заключённые между этой точкой и сторонами угла, имели данное отношение $m:n$.

\item
Найти в треугольнике такую точку, чтобы перпендикуляры, опущенные из неё на стороны треугольника, находились в данном отношении $m:n:p$ (смотри упражнение \ref{1938/upr-3-17}).

\item
Построить треугольник по углу, одной из сторон, прилежащих к нему, по отношению этой стороны к третьей стороне (сколько решений?).

\item
То же — по углу при вершине, основанию и отношению его к одной из боковых сторон.

\item
То же — по высоте, углу при вершине и отношению отрезков основания.

\item
То же — по углу при вершине, основанию и данной на основании точке, через которую проходит биссектриса угла при вершине.

\item
То же — по двум углам и сумме или разности основания с высотой.

\item
Построить равнобедренный треугольник по углу при вершине и сумме основания с высотой.

\item
На бесконечной прямой $MN$ даны две точки $A$ и $B$.
Найти на этой прямой третью точку $C$, чтобы $CA:CB=m:n$, где $m$ и $n$ — данные отрезки или данные числа (если $m\ne n$, то таких точек существует две:
одна между $A$ и $B$, другая вне отрезка $AB$).

\item
Вписать в данный круг треугольник, у которого даны:
основание и отношение двух сторон.

\item
Вписать в данный круг треугольник, у которого даны:
основание и медиана к одной из неизвестных сторон. 

\item
Вписать квадрат в данный сегмент так, чтобы одна его сторона лежала на хорде, а вершины противолежащих углов — на дуге.

\smallskip
\so{Указание}.
Задачи решаются методом подобия (§~\ref{1938/181}).

\item
Вписать квадрат в данный треугольник так, чтобы одна сторона его лежала на основании треугольника, а вершины противолежащих углов — на боковых сторонах треугольника.

\item
В данный треугольник вписать прямоугольник (смотри предыдущую задачу), у которого стороны относились бы как $m:n$.

\item
Около данного квадрата описать треугольник, подобный данному треугольнику.

\item
Дана окружность и на ней две точки $A$ и $B$.
Найти на этой окружности третью точку $C$, чтобы расстояния от неё до $A$ и $B$ находились в данном отношении.

\item 
Построить треугольник по двум сторонам и биссектрисе угла между ними (смотри рис.~\ref{1938/ris-196};
сначала находим прямую $CE$ из пропорции $CE\z:BD\z=AE:AB$, затем строим $\triangle BCE$ и~т.~д.).

\item
Построить отрезок $x$, который относился бы к данному отрезку $m$ как $a^2:b^2$ ($a$ и $b$ — данные отрезки).

\item
Найти вне данного круга такую точку, чтобы касательная, проведённая из неё к этой окружности, была вдвое менее секущей, проведённой из той же точки через центр (приложением алгебры к геометрии).

\item
Через данную вне круга точку провести такую секущую, которая разделилась бы этой окружностью в данном отношении (приложением алгебры к геометрии).

\item
Построить треугольник по трём его высотам $h_1$, $h_2$, $h_3$

\smallskip

\so{Решение}.
Предварительно из подобия прямоугольных треугольников надо доказать, что высоты \so{обратно пропорциональны} соответствующим сторонам.
Если стороны, на которые опущены высоты $h_1$, $h_2$, $h_3$ обозначим соответственно через $x_1$, $x_2$, $x_3$, то
\[x_1:x_2=h_2:h_1,\]
\[x_2:x_3=h_3:h_2=1:\frac{h_2}{h_3}=h_1:\frac{h_1h_2}{h_3}\]
откуда 
\[x_1:x_2:x_3=h_2:h_1:\frac{h_1h_2}{h_3}\]
Выражение $\frac{h_1h_2}{h_3}$ есть четвёртая пропорциональная к $h_1$, $h_2$ и $h_3$.
Построив её (пусть это будет $k$), мы будем иметь три отрезка:
$h_2$, $h_1$ и $k$, которым искомые стороны пропорциональны;
значит, треугольник, имеющий эти отрезки сторонами, подобен искомому, и потому вопрос сводится к построению такого треугольника, который, будучи подобен данному, имел бы данную высоту.
Задача не имеет решения, если по трём прямым $h_1$, $h_2$ и $k$ нельзя построить треугольник.

\item
Построить отрезки, выражаемые формулами: 
\[1)\ x=\frac{abc}{de}=\frac{ab}{d}\cdot \frac{c}{e}\]
(придётся два раза построить четвёртую пропорциональную).
\[2)\ x=\sqrt{a^2+bc}.\]
(предварительно построить отрезки $k=\sqrt{bc}$ и $x=\sqrt{a^2+k^2}$).

\end{enumerate}

\begin{center}
\so{Задачи на вычисление}
\end{center}

\begin{enumerate}[resume,noitemsep]

\item
По данному основанию $a$ и высоте $h$ остроугольного треугольника вычислить сторону $x$ квадрата, вписанного в этот треугольник так, что одна сторона квадрата лежит на основании треугольника, а две вершины квадрата — на боковых сторонах треугольника.

\item
Стороны треугольника имеют длины 10, 12 и 17 м.
Вычислить отрезки стороны, равной 17 м, на которые она делится биссектрисой противолежащего угла.

\item
Перпендикуляр, опущенный из вершины прямого угла на гипотенузу, делит её на два отрезка $m$ и~$n$.
Вычислить катеты.

\smallskip
\so{Указание}.
Продолжив $AD$ на расстояние $DE=AD$ и соединив точку $E$ с $B$ и $C$, получим параллелограмм, к которому применим теорему §~\ref{1938/197}.

\item
В $\triangle ABC$ стороны равны:
$AB$ = 7, $BC$ = 15 и $AC$ = 10.
Определить, какого вида угол $A$, и вычислить высоту, опущенную из вершины $B$.

\item
Из точки вне круга проведена касательная $a$ и секущая.
Вычислить длину секущей, зная, что отношение внешней её части к внутренней равно $m:n$.

\item
К двум кругам, радиусы которых $R$ и $r$, а расстояние между центрами $d$, проведена общая касательная.
Определить вычислением положение точки пересечения этой касательной с линией центров, во-первых, когда эта точка лежит по одну сторону от центров, во-вторых, когда она расположена между ними.

\end{enumerate}

}

%% file: 2D/pravilnye-mnogougi.tex
\section{Правильные многоугольники}

\paragraph{}\label{1938/212}
\so{Определения}.
Ломаная линия называется \rindex{правильная ломаная}\textbf{правильной}, если она удовлетворяет следующим трём условиям:
1) отрезки прямых, составляющие её, равны;
2) углы, составленные каждыми двумя соседними отрезками, равны;
3) из каждых трёх последовательных отрезков первый и третий расположены по одну сторону от прямой, на которой лежит второй.

\begin{figure}[h]
\centering
\includegraphics{mppics/ris-218}
\caption{}\label{1938/ris-218}
\end{figure}

Таковы, например, линии $ABCDE$ и $FGHKL$ (рис.~\ref{1938/ris-218});
но ломаную $MNPQR$ нельзя назвать правильной, потому что она не удовлетворяет третьему условию.

Правильная ломаная может быть \rindex{выпуклая ломаная}\textbf{выпуклой}, как например линия $ABCDE$.

Многоугольник называется \rindex{правильный многоугольник}\textbf{правильным}, если он ограничен правильной ломаной линией, то есть если он имеет равные стороны и равные углы.
Таковы, например, квадрат, равносторонний треугольник.

\begin{wrapfigure}{o}{55mm}
\centering
\includegraphics{mppics/ris-ru-219}
\caption{}\label{1938/ris-219}
\end{wrapfigure}

Многоугольник, изображённый на рис.~\ref{1938/ris-219}а, есть выпуклый правильный пятиугольник;
многоугольник на рис.~\ref{1938/ris-219}б — также правильный пятиугольник, но не выпуклый (так называемый звёздчатый).
В нашем курсе геометрии мы будем рассматривать \so{только выпуклые} правильные многоугольники и поэтому, когда мы скажем «правильный многоугольник», мы будем подразумевать слово «выпуклый».

{\sloppy 

Последующие теоремы показывают, что построение правильных многоугольников тесно связано с разделением окружности на равные части.

}

\paragraph{}\label{1938/213}
\so{Теорема}.
\textbf{\emph{Если окружность разделена на $n$ равных частей ($n\ge 3$), то:}}

1) \textbf{\emph{соединив хордами каждые две соседние точки деления, получим правильный $n$-угольник (вписанный).}}

2) \textbf{\emph{проведя через все точки деления касательные и продолжив каждую из них до взаимного пересечения с касательными соседних точек деления, получим правильный $n$-угольник (описанный).}}

\begin{wrapfigure}{r}{36mm}
\centering
\includegraphics{mppics/ris-220}
\caption{}\label{1938/ris-220}
\end{wrapfigure}

1) Пусть окружность (рис.~\ref{1938/ris-220}) разделена на несколько равных частей в точках $A, B, C$ и~т.~д.
и через эти точки проведены хорды $AB, BC,\dots$
и касательные $MBN$, $NCP$ и~т.~д.
Тогда 1) вписанный многоугольник $ABCDEF$ — правильный, потому что все его стороны равны (как хорды, стягивающие равные дуги) и все углы равны (как вписанные, опирающиеся на равные дуги).

2) Чтобы доказать правильность описанного многоугольника $MNPQRS$, рассмотрим треугольники $AMB$, $BNC$ и~т.~д.
У них основания $AB, BC$ и~т.~д.
равны;
углы, прилежащие к этим основаниям, также равны, потому что каждый из них имеет одинаковую меру (угол, составленный касательной и хордой, измеряется половиной дуги, заключённой внутри него).
Значит, все эти треугольники равнобедренные и равны между собой, а потому 
\begin{align*}
MN&=NP=\dots,
\\
\angle M&=\angle N=\dots;
\end{align*}
то есть многоугольник $MNPQRS$ правильный.

\paragraph{}\label{1938/215}
\so{Теорема}.
\textbf{\emph{Если многоугольник правильный, то:}}

1) \textbf{\emph{около него можно описать окружность.}}

2) \textbf{\emph{в него можно вписать окружность.}}

\begin{wrapfigure}{o}{35mm}
\centering
\includegraphics{mppics/ris-222}
\caption{}\label{1938/ris-222}
\end{wrapfigure}

1) Проведём окружность через какие-нибудь три соседние вершины $A, B$ и $C$ (рис.~\ref{1938/ris-222}) правильного многоугольника $ABCDE$ и докажем, что она пройдёт через следующую, четвёртую вершину $D$.
Опустим из центра $O$ перпендикуляр $OK$ на хорду $BC$ и соединим $O$ с $A$ и $D$.
Повернём четырёхугольник $ABKO$ вокруг стороны $OK$ так, чтобы он упал на четырёхугольник $ODCK$.
Тогда $KB$ пойдёт по $KC$ (вследствие равенства прямых углов при точке $K$), 
$B$ совпадёт с $C$ (так как хорда $BC$ делится в точке $K$ пополам), 
сторона $BA$ пойдёт по $CD$ (вследствие равенства углов $B$ и $C$)
и, наконец, точка $A$ совпадёт с $D$ (вследствие равенства сторон $BA$ и $CD$).
Из этого следует, что $OA$ совместится с $OD$, и, значит, точки $A$ и $D$ одинаково удалены от центра;
поэтому вершина $D$ должна лежать на окружности, проходящей через $A, B$ и $C$.
Точно так же докажем, что эта окружность, проходя через три соседние вершины $B, C$ и $D$, пройдёт через следующую вершину $E$, и~т.~д.;
значит, она пройдёт через все вершины многоугольника.

2) Из доказанного следует, что стороны правильного многоугольника всегда можно рассматривать как равные хорды одной окружности;
но такие хорды одинаково удалены от центра;
значит, все перпендикуляры $OM$, $ON$ и~т.~д., опущенные из $O$ на стороны многоугольника, равны между собой, и потому окружность, описанная радиусом $OM$ с центром в точке $O$, будет вписанной в многоугольник $ABCDE$.

\paragraph{}\label{1938/216}
\so{Следствие}.
Из предыдущего видно, что две окружности, описанная около правильного многоугольника и вписанная в него, имеют один и тот же центр.
Так как этот общий центр одинаково удалён от всех вершин многоугольника, то он должен лежать на срединном перпендикуляре, восстановленном к любой стороне многоугольника, а будучи одинаково удалён от сторон каждого угла, он должен находиться на его биссектрисе.
Поэтому, чтобы найти центр окружности, описанной около правильного многоугольника или вписанной в него, достаточно определить точку пересечения двух срединных перпендикуляров, восстановленных к сторонам многоугольника, или двух биссектрис углов, или одного срединного перпендикуляра с биссектрисой.

Легко заметить, что срединные перпендикуляры, восстановленные к сторонам правильного многоугольника, а также биссектрисы всех углов правильного многоугольника являются его осями симметрии.

\paragraph{}\label{1938/217}
\so{Определения}.
Общий центр окружностей, описанной около правильного многоугольника и вписанной в него, называется \rindex{центр!правильного многоугольника}\textbf{центром} этого многоугольника, радиус вписанной окружности — \rindex{апофема}\textbf{апофемой} его.

{\sloppy 

Угол, составленный двумя радиусами, проведёнными к концам какой-нибудь стороны правильного многоугольника, называется \rindex{центральный угол}\textbf{центральным углом}.
Центральных углов в многоугольнике столько, сколько сторон;
все они равны, как измеряющиеся равными дугами.

}

Так как сумма всех центральных углов равна $360\degree$, то в правильном $n$-угольнике каждый из них равен $\tfrac{360\degree}n$;
так, центральный угол правильного шестиугольника равен $\tfrac{360\degree}6\z=60\degree$, правильного восьмиугольника равен $\tfrac{360\degree}8 = 45\degree$ и так далее.

Так как сумма всех внутренних углов $n$-угольника (§~\ref{1938/82}) равна $180\degree\cdot(n-2)$, то каждый внутренний угол правильного $n$-угольника равен
\[\frac{180\degree\cdot(n-2)}{n}\]

Например, у правильного восьмиугольника внутренний угол равен
\[\frac{180\degree\cdot(8-2)}{8}=\frac{6}{8}\cdot 180\degree=\frac{3}{4}\cdot 180\degree=135\degree.\]

\paragraph{}\label{1938/218}
\so{Теорема}.
\textbf{\emph{Правильные $\bm{n}$-угольники подобны, а стороны их относятся как радиусы или апофемы.}}

\begin{figure}[!ht]
\centering
\includegraphics{mppics/ris-223}
\caption{}\label{1938/ris-223}
\end{figure}

{\sloppy 
1) Чтобы доказать подобие (рис.~\ref{1938/ris-223}) правильных $n$-угольников $ABCDEF$ и $A_1B_1C_1D_1E_1F_1$, достаточно обнаружить, что у них углы равны и стороны пропорциональны.

}

Углы $n$-угольников равны, так как каждый из них содержит одно и то же число градусов, а именно: $\frac{180\degree\cdot(n-2)}{n}$.
Так как 
\[AB=BC = CD=\dots,
\quad\text{и}\quad A_1B_1=B_1C_1 = C_1D_1=\dots,\]
то, очевидно, что
\[\frac{AB}{A_1B_1}=\frac{BC}{B_1C_1} = \frac{CD}{C_1D_1}=\dots;\]
то есть у таких $n$-угольников стороны пропорциональны.

2) Пусть $O$ и $O_1$ (рис.~\ref{1938/ris-223}) будут центры данных $n$-угольников, $OA$ и $O_1A_1$ — их радиусы, $OM$ и $O_1M_1$ — апофемы.
Треугольники $OAB$ и $O_1A_1B_1$ подобны, так как углы одного соответственно равны углам другого.

Из подобия их следует:
\[\frac{AB}{A_1B_1}=\frac{OA}{O_1A_1} = \frac{OM}{O_1M_1}\]

\smallskip
\so{Следствие}.
Так как периметры подобных многоугольников относятся как соответственные стороны (§~\ref{1938/172}), то \emph{периметры правильных $n$-угольников относятся как радиусы или как апофемы.}


\begin{wrapfigure}{r}{38mm}
\vskip-4mm
\centering
\includegraphics{mppics/ris-224}
\caption{}\label{1938/ris-224}
\end{wrapfigure}

\paragraph{}\label{1938/219}
\mbox{\so{Задача}.}
\emph{Вычислить сторону вписанного в круг:
1) квадрата;
2) правильного шестиугольника;
3) правильного треугольника.}

Условимся обозначать длину стороны правильного $n$-угольника буквой $a_n$, а его периметр — буквой $p_n$.

{\sloppy

Формулы для сторон вписанного квадрата, шестиугольника и треугольника можно легко получить из рассмотрения рис.~\ref{1938/ris-224}—\ref{1938/ris-226}.

}

1) На рис.~\ref{1938/ris-224} проведены два взаимно перпендикулярных диаметра $AC$ и $BD$ и последовательные концы их соединены хордами;
от этого получился вписанный квадрат $ABCD$.

Из прямоугольного $\triangle AOB$ находим:
\[AB^2=AO^2+OB^2=2R^2.\]
откуда
\[a_4=R\sqrt2.\]

\paragraph{}\label{1938/220}
2) На рис.~\ref{1938/ris-225} построена хорда, соответствующая центральному углу в $60\degree$ (сторона правильного вписанного шестиугольника).

\begin{wrapfigure}{o}{38mm}
\vskip-4mm
\centering
\includegraphics{mppics/ris-225}
\caption{}\label{1938/ris-225}
\end{wrapfigure}

Так как у равнобедренного $\triangle AOB$ каждый из углов $A$ и $B$ равен \[\frac{180\degree-60\degree} 2 = 60\degree,\] то $\triangle AOB$ есть равноугольный и, следовательно, равносторонний;
значит:
\[AB=AO,\quad\text{то есть}\quad a_6 = R.\]
Отсюда мы получаем простой способ деления окружности на шесть равных частей.

\paragraph{}\label{1938/221}
3) На рис.~\ref{1938/ris-226} окружность разделена на шесть равных частей, и точки деления через одну последовательно соединены хордами, отчего образовался вписанный равносторонний треугольник $ABC$.

\begin{wrapfigure}{o}{38mm}
\centering
\includegraphics{mppics/ris-226}
\caption{}\label{1938/ris-226}
\end{wrapfigure}

Проведя хорду $AD$, получаем прямоугольный треугольник $ABD$ (угол $BAD$ прямой как вписанный и опирающийся на диаметр).
Из $\triangle ABD$ находим:
\[AB=\sqrt{BD^2-AD^2},\]
то есть
\[a_3=\sqrt{(2R)^2-R^2}\]
и, значит,
\[a_3=R\cdot \sqrt3.\]

\paragraph{}\label{1938/222}
\so{Задача}.
\emph{Вписать в данный круг правильный десятиугольник и определить его сторону в зависимости от радиуса.}

\begin{wrapfigure}{r}{45mm}
\centering
\includegraphics{mppics/ris-227}
\caption{}\label{1938/ris-227}
\end{wrapfigure}

Предварительно докажем одно важное свойство правильного 10-угольника.
Пусть хорда $AB$ (рис.~\ref{1938/ris-227}) есть сторона правильного 10-угольника.
Тогда угол $AOB$ равен $36\degree$ , а каждый из углов $A$ и $B$ содержит по $\tfrac12(180\degree-36\degree)$, то есть
по $72\degree$.
Разделим угол $A$ пополам прямой $AC$.
Каждый из углов, образовавшихся при точке $A$, равен $36\degree$;
следовательно, $\triangle ACO$, имея два равных угла, есть равнобедренный, то есть $AC=CO$, $\triangle ABC$ также равнобедренный, потому что $\angle B = 72\degree$ и $\angle ACB = 180\degree- 72\degree\z- 36\degree = 72\degree$;
следовательно, $AB=AC\z=CO$.
По свойству биссектрисы угла треугольника (§~\ref{1938/186}) можно написать:
\[\frac{AO}{AB}=\frac{OC}{CB}.
\eqno(1)\]
Заменив $AO$ и $AB$ равными им отрезками $OB$ и $OC$, получим:
\[\frac{OB}{OC}=\frac{OC}{CB}, \eqno(2)\]
то есть радиус $OB$ разделён в точке $C$ в среднем и крайнем отношении (§~\ref{1938/209}), причём $OC$ есть его б\'{о}льшая часть.
Но $OC$ равна стороне правильного вписанного 10-угольника;
значит, \emph{сторона правильного вписанного 10-угольника равна б\'{о}льшей части радиуса, разделённого в среднем и крайнем отношении.}

\begin{wrapfigure}{O}{45mm}
\centering
\includegraphics{mppics/ris-228}
\caption{}\label{1938/ris-228}
\end{wrapfigure}

Теперь задача решается легко:

1) Делят радиус круга (например, $OA$, рис.~\ref{1938/ris-228}) в среднем и крайнем отношении
;
затем, дав циркулю раствор, равный б\'{о}льшей части радиуса, откладывают им по окружности дуги, одна за другой, и точки деления последовательно соединяют хордами.

2) Обозначив буквой $x$ длину стороны правильного вписанного 10-угольника. Пропорцию (2) можно переписать так:
\[\frac Rx=\frac x{(R-x)},\]
откуда
\[x^2+Rx-R^2=0.\]

Решив это квадратное уравнение, найдём:
\[x=a_{10}=R\cdot\tfrac{\sqrt5-1}{2}\approx R \cdot  0{,}61803\dots\]

{\small

\paragraph{}\label{1938/223}
\mbox{\so{Замечания}.}
1) Чтобы вписать в данный круг правильный пятиугольник, делят окружность на 10 равных частей (как указано выше) и точки деления соединяют через одну хордами.

2) Из равенства
\[\tfrac16-\tfrac1{10}=\tfrac5{30}-\tfrac3{30}=\tfrac2{30}=\tfrac1{15}\]
видно, что если из $\tfrac16$ части окружности вычесть $\tfrac1{10}$ её часть, то остаток будет равен $\tfrac1{15}$ окружности.
Это даёт нам простой способ вписать в окружность правильный 15-угольник, так как делить окружность на 6 и на 10 равных частей мы умеем.

\begin{wrapfigure}{r}{35mm}
\centering
\includegraphics{mppics/ris-229}
\caption{}\label{1938/ris-229}
\end{wrapfigure}

3) Чтобы построить пятиконечную звезду (рис.~\ref{1938/ris-229}), делят окружность на 10 равных частей и какую-нибудь из точек деления соединяют хордами с другими точками деления через три (как указано на рисунке).

\paragraph{}\label{1938/224}
\mbox{\so{Задача}.}
\emph{Удвоить число сторон правильного вписанного $n$-угольника.}

В этом сокращённом выражении подразумеваются две задачи:

1) по данному правильному вписанному $n$-угольнику \so{построить} правильный $2n$-угольник, вписанный в ту же окружность;

2) \so{вычислить сторону} этого $2n$-угольника по данной стороне данного $n$-угольника и данному радиусу круга.

\begin{wrapfigure}{O}{45mm}
\centering
\includegraphics{mppics/ris-230}
\caption{}\label{1938/ris-230}
\end{wrapfigure}

1) Пусть $AB$ (рис.~\ref{1938/ris-230}) есть сторона правильного вписанного $n$-угольника и $O$ — центр круга.
Проведём $OC\perp AB$ и соединим $A$ с $C$.
Дуга $AB$ делится в точке $C$ пополам, следовательно, хорда $AC$ есть сторона правильного вписанного $2n$-угольника.

2) В $\triangle ACO$ угол $O$ всегда острый (так как дуга $ACB$ всегда меньше полуокружности, и, следовательно, половина её, дуга $AC$, меньше четверти окружности);
поэтому (§~\ref{1938/194})
\[AC^2=OA^2+OC^2-2OC\cdot OD,\]
то есть
\[a_{2n}^2=R^2+R^2-2R\cdot OD=2R^2-2R\cdot OD.\]

Из прямоугольного $\triangle AOD$ определим катет $OD$:
\begin{align*}
OD&=\sqrt{AO^2-AD^2}=
\\
&=\sqrt{R^2-(\tfrac{a_n}2)^2}=
\\
&=\sqrt{R^2-\tfrac{a_n^2}4}.
\end{align*}
Следовательно,
\[a_{2n}^2=2R^2-2R\sqrt{R^2-\frac{a_n^2}4}.\]
Такова формула удвоения числа сторон правильного вписанного многоугольника (из неё сторону $a_{2n}$ получим посредством извлечения квадратного корня).

}

\smallskip
\so{Пример}.
Вычислим сторону правильного 12-угольника.
Для простоты примем, что $R=1$ (и, следовательно, $a_6 = 1$).
\begin{align*}
a_{12}^2&=2-2\sqrt{1-\tfrac14}=
\\
&=2-2\sqrt{\tfrac34}=
\\
&=2-\sqrt{3},
\end{align*}
откуда
\[a_{12}=\sqrt{2-\sqrt3}\approx 0{,}517\dots\]
Так как стороны правильных $n$-угольников пропорциональны их радиусам (§~\ref{1938/218}), то при радиусе, равном не единице, а какому-нибудь числу $R$, для стороны правильного 12-угольника получим такую формулу:
\[a_{12}=R\cdot \sqrt{2-\sqrt3}\approx R\cdot 0{,}517\dots\]

{\small

\paragraph{На сколько равных частей можно делить окружность с помощью циркуля и линейки?}\label{1938/225}
Применяя указанные в предыдущих задачах способы, мы можем с помощью циркуля и линейки делить окружность на такое число равных частей (и, следовательно, вписывать в окружность правильные многоугольники с таким числом сторон), которое заключается в следующей таблице:
\begin{align*}
3,&&3\cdot 2,&&3\cdot2\cdot2,&&\text{вообще}&&3\cdot 2^n,&
\\
4,&&4\cdot 2,&&4\cdot2\cdot2,&&\text{—\textquotedbl—\ \ }&&2^n,&
\\
5,&&5\cdot 2,&&5\cdot2\cdot2,&&\text{—\textquotedbl—\ \ }&&5\cdot 2^n,&
\\
15,&&15\cdot 2,&&15\cdot2\cdot2,&&\text{—\textquotedbl—\ \ }&&15\cdot 2^n.&
\end{align*}

Немецкий математик Гаусс (живший в 1777—1855 годах) доказал, что посредством циркуля и линейки можно делить окружность на такое число равных частей, которое, будучи простым, выражается формулой $2^{2^n} + 1$.
Например, можно разделить окружность на 17 равных частей и на 257 равных частей, так как 17 и 257 простые числа вида $2^{2^n} + 1$ 
($17 = 2^{2^2} + 1$;
$257 = 2^{2^3} + 1$).
Доказательство Гаусса выходит за пределы элементарной математики.

Доказано также, что с помощью линейки и циркуля окружность можно делить на такое составное число равных частей, в состав которого не входят никакие иные простые множители, кроме:
1) множителей вида $2^{2^n} + 1$ при условии, что все эти множители различны 
и 
2) множителя 2 в какой угодно степени.

Например, в окружность с помощью циркуля и линейки можно вписать правильный 170-угольник ($170 = 2 \cdot 5 \cdot 17 = 2 \cdot (2^2 + 1) \cdot (2^{2^2} +1)$),
но нельзя вписать правильный 9-угольник (хотя множитель 3 имеет вид $2^{2^n} + 1$, но в составе 9 он повторяется).

На всякое иное число равных частей окружность может быть разделена \so{приближённо}.
Пусть, например, требуется разделить окружность на 7 равных частей (или вписать правильный семиугольник).
Тогда предварительно вычислим величину центрального угла,
$\frac{360\degree}7\z=51\tfrac37\degree$.
Построить точно такой угол мы не можем, но по транспортиру приблизительно можем отложить при центре угол в $51\degree$ и тогда получим приблизительно $\tfrac17$ часть окружности.

}

{\small

\subsection*{Упражнения}

\begin{enumerate}[noitemsep]

\item
Составить формулу для стороны правильного вписанного 24-угольника.

\item
Составить формулу для сторон правильных вписанных восьмиугольника и 16-угольника.

\item
Составить формулу для сторон правильных описанных треугольника и шестиугольника.

\item
Пусть $AB, BC$ и $CD$ будут три последовательные стороны правильного многоугольника, имеющего центр в $O$.
Доказать, что если продолжим стороны $AB$ и $CD$ до взаимного пересечения в точке $E$, то четырёхугольник $OAEC$ может быть вписан в окружность.

\item
Доказать, что:
1) всякий вписанный равносторонний многоугольник — правильный;
2) всякий описанный равноугольный многоугольник — правильный.

\item
Доказать, что:
1) каждый правильный $n$-угольник имеет $n$ осей симметрии, причём все эти оси симметрии проходят через его центр;
2) для многоугольника с чётным числом сторон центр многоугольника является центром его симметрии.

\item
Доказать, что две диагонали правильного пятиугольника, не исходящие из одной вершины, пересекаясь, делятся в среднем и крайнем отношении.

\smallskip
\so{Указание}.
Пусть $ABCDE$ — правильный пятиугольник, $AC$ и $BE$ — его диагонали, $F$ — точка их пересечения.
$\triangle ABC \sim \triangle ABF$ и~т.~д.

\item
На данной стороне построить:
1) правильный восьмиугольник;
2) правильный 10-угольник.

\item
Срезать от данного квадрата углы так, чтобы образовался правильный восьмиугольник.

\item
В данный квадрат вписать равносторонний треугольник, помещая одну из его вершин или в вершине квадрата, или в середине какой-либо стороны.

\item
Вписать в равносторонний треугольник другой равносторонний треугольник, стороны которого были бы перпендикулярны к сторонам данного.

\item
Построить углы в $18\degree$, $30\degree$, $72\degree$, $75\degree$.

\item
Около окружности описан какой-нибудь правильный многоугольник.
Пользуясь им, вписать в эту окружность правильный многоугольник, имеющий вдвое более сторон, чем описанный.

\end{enumerate}

}

%% file: 2D/dlina-okr.tex
\section{Длина окружности и её частей}

\paragraph{Предварительное разъяснение.}\label{1938/226}
Отрезок можно сравнить с другим отрезком, принятым за единицу, так как прямые линии при наложении совмещаются.
Действительно, только по этой причине мы можем установить, какие отрезки считать равными и неравными;
что такое сумма отрезков или какой отрезок более другого в $2, 3, 4,\dots$ раза.
Точно так же дуги окружностей одинакового радиуса можно сравнить между собой вследствие того, что такие дуги при наложении совмещаются.
Но так как никакая часть окружности (или другой кривой) не может совместиться с прямой, то нельзя путём наложения установить, какой криволинейный отрезок надлежит считать равным данному прямолинейному отрезку, а следовательно, и то, какой криволинейный отрезок больше данного прямолинейного в $2,3,4,\dots$
раза.
Таким образом, является необходимость особо определить, что мы будем подразумевать под длиной окружности (или части её), когда сравниваем её с прямолинейным отрезком.

Для этой цели мы должны ввести новое понятие, имеющее исключительно большое значение во всей математике, именно понятие о пределе.

\subsection*{Пределы}

\paragraph{Предел числовой последовательности.}\label{1938/227}
Во многих вопросах алгебры и геометрии приходится встречаться с последовательностями чисел, написанных одно за другим по определённому закону.
Например, натуральный ряд чисел:
\[1, 2, 3, 4, 5,\dots,\]
арифметическая и геометрическая прогрессии, продолженные неограниченно:
\[a,a+d,a+2d,a+3d,\dots,\]
\[a,aq,aq^2,aq^3,\dots,\]
представляют собой бесконечные последовательности чисел или бесконечные числовые последовательности.

Для каждой такой последовательности можно указать правило, по которому составляются её члены.
Так, в арифметической прогрессии каждый член разнится от предыдущего на одно и то же число, в геометрической прогрессии каждый член равен предшествующему, умноженному на некоторое определённое число (знаменатель прогрессии).

Многие последовательности составляются по более сложным правилам.
Так, например, вычисляя $\sqrt{2}$ с недостатком, сначала с точностью до $\tfrac1{10}$, затем с точностью до $\tfrac1{100}$, затем до $\tfrac1{1000}$ и продолжая это вычисление неограниченно, мы получим бесконечную числовую последовательность:
\[1{,}4;
1{,}41;
1{,}414;
1{,}4142,\dots,\]
дающую приближённое значение  $\sqrt{2}$  с возрастающей степенью точности.

Для этой последовательности нельзя указать простого правила, по которому можно было бы получить новые её члены, зная предыдущие, но все же можно определить любой член этой последовательности.
Так, чтобы получить 4-й её член, нужно вычислить с точностью до $0{,}0001$, для получения 5-го члена нужно вычислить  $\sqrt{2}$  с точностью до $0{,}00001$ и~т.~д.

Допустим, что члены данной бесконечной последовательности $a_1$, $a_2$, $a_3, \dots, a_n,\dots$ по мере повышения их номера неограниченно приближаются к некоторому числу $A$.
Это значит следующее:
существует некоторое число $A$, такое, что какое бы малое положительное число $q$ мы ни взяли, в данной последовательности можно отыскать член, начиная с которого все члены последовательности по абсолютной величине отличаются от $A$ меньше, чем на $q$.
Мы будем это свойство коротко выражать так:
 абсолютная величина разности $a_n-A$
неограниченно убывает с возрастанием номера~$n$.

В этом случае число $A$ называется \rindex{предел}\textbf{пределом} данной бесконечной числовой последовательности.
Приведём пример такой последовательности.
Составим последовательность десятичных дробей.
\[0{,}9;
0{,}99;
0{,}999;\dots,\]
Здесь каждый член получается из предыдущего приписыванием нового десятичного знака 9.

Легко заметить, что члены этой последовательности неограниченно приближаются к единице.

Именно, первый член отличается от единицы на $\tfrac1{10}$, 
второй на $\tfrac1{100}$, третий на $\tfrac1{1000}$ и, если достаточно продолжить эту последовательность, то можно найти в ней член, начиная с которого все последующие члены будут отличаться от единицы на сколь угодно малую, заранее указанную величину.
Следовательно, мы можем сказать, что взятая нами бесконечная числовая последовательность имеет пределом единицу.

Другим примером числовой последовательности, имеющей предел, служит последовательность приближённых значений длины отрезка, несоизмеримого с единицей длины (§~\ref{1938/150}), вычисленных с недостатком, сначала с точностью до $\tfrac1{10}$, затем $\tfrac1{100}$, затем — до $\tfrac1{1000}$ и~т.~д.

Пределом этой последовательности служит бесконечная десятичная дробь, представляющая точную меру длины данного отрезка.
В самом деле, величина бесконечной десятичной дроби заключена между двумя её приближёнными значениями, вычисленными с одинаковой точностью — одно с недостатком, другое с избытком. 
Как было показано выше, эта разность неограниченно убывает по мере повышения степени точности приближённых значений.
Следовательно, должна неограниченно убывать и разность между самой бесконечной десятичной дробью и её приближёнными значениями по мере повышения степени точности этих значений.
Значит, бесконечная десятичная дробь служит пределом последовательности всех её приближённых значений, взятых с недостатком (или всех приближённых значений, взятых с избытком).

Легко заметить, что не всякая бесконечная последовательность имеет предел;
например, натуральный ряд чисел.
\[1, 2, 3, 4, 5,\dots,\]
очевидно, никакого предела не имеет, так как его члены неограниченно возрастают и ни к какому числу не приближаются.

\paragraph{}\label{1938/228}
\so{Теорема}.
\textbf{\emph{Всякая бесконечная числовая последовательность может иметь только один предел.}}

В справедливости этой теоремы легко убедиться доказательством от противного.
В самом деле, предположим, что дана последовательность
\[a_1,a_2,a_3,\dots,a_n,\dots,\]
которая имеет два различных предела $A$ и $B$.
В таком случае, в силу того, что $A$ есть предел данной последовательности, абсолютная величина разности $a_n-A$ должна неограниченно убывать с возрастанием~$n$.
В силу того, что $B$ есть тоже предел данной последовательности, абсолютная величина разности $a_n-B$ также должна неограниченно убывать с возрастанием~$n$.
Но в таком случае абсолютная величина разности
\[(a_n-A)-(a_n-B)\]
должна также или неограниченно убывать, или быть равной нулю.
Но эта последняя разность равна разности чисел $B-A$ и, следовательно, есть некоторое вполне определённое, отличное от нуля число.
Это число не зависит от номера $n$ и при возрастании $n$ вовсе не изменяется.
Таким образом, предположение, что существуют два предела числовой последовательности, привело нас к противоречию.

{\sloppy 
\paragraph{Предел возрастающей бесконечной числовой последовательности.}\label{1938/229}
Рассмотрим такую последовательность 
\[a_1, a_2, a_3,\dots,a_n,\dots,\]
в которой каждый следующий член больше предыдущего, то есть $a_{n+1} \z> a_n$, и в то же время все члены последовательности меньше некоторого определённого числа $M$, то есть $a_n < M$ для любого номера~$n$.

}

В этом случае последовательность имеет определённый предел (теорема Вейерштрасса).

\paragraph{}\label{1938/230}
\so{Доказательство}.
Пусть дана бесконечная числовая последовательность
\[a_1,a_2,a_3,\dots,a_n,\dots,
\eqno(1)\]
в которой каждый член больше предыдущего или равен ему (то есть $a_{n+1}\ge a_n$ для всякого $n$), при этом среди членов последовательности нет числа, большего данного числа $M$, например, нет числа, большего, чем $10$.

Возьмём число $9$ и смотрим, нет ли среди членов последовательности (1) чисел больше, чем $9$.
Допустим, что таких нет.
Возьмём число $8$ и смотрим, имеются ли в последовательности (1) числа больше, чем $8$.
Допустим, что такие есть.
Тогда записываем число $8$, затем делим промежуток от $8$ до $9$ на $10$ частей и испытываем последовательно числа:
$8{,}1; 8{,}2; 8{,}3;\dots$, то есть смотрим, имеются ли среди членов последовательности (1) числа б\'{о}льшие, чем $8{,}1$.
Если есть, то ставим тот же вопрос для числа $8{,}2$ и~т.~д.
Допустим, что в последовательности (1) есть числа, б\'{о}льшие, чем $8{,}6$, но нет чисел, б\'{о}льших, чем $8{,}7$.
Тогда делаем вторую запись:
пишем число $8{,}6$, затем разбиваем промежуток от $8{,}6$ до $8{,}7$ на $10$ частей и испытываем таким же образом последовательно числа: $8{,}61; 8{,}62; 8{,}63;\dots$
Допустим, что в последовательности (1) есть числа, б\'{о}льшие, чем $8{,}64$, но нет чисел, б\'{о}льших, чем $8{,}65$.
Тогда делаем третью запись $8{,}64$ и поступаем таким же образом для промежутка от $8{,}64$ до $8{,}65$.

Продолжая этот процесс неограниченно, мы придём к бесконечной десятичной дроби:
$8{,}64\dots$, то есть
к некоторому вещественному числу.
Назовём его $\alpha$ и возьмём его приближённые значения с $n$ десятичными знаками с недостатком и с избытком (§~\ref{1938/150}).
Первое назовём $\alpha_n$, второе —  $\alpha_n'$.
Эти приближённые значения можно выбрать так, что.
\[\alpha_n< \alpha\le\alpha_n'
\quad\text{и}\quad
\alpha_n'-\alpha_n=\tfrac1{10^n}.\] 
Из способа образования вещественного числа $\alpha$ следует, что среди членов последовательности (1) нет чисел, б\'{о}льших $\alpha_n'$, но имеются числа, б\'{о}льшие $\alpha_n$.
Пусть $a_k$, — одно из таких чисел
\[\alpha_n< a_k\le\alpha_n'.\]
В силу возрастания последовательности (1) и отсутствия в ней членов, б\'{о}льших $\alpha_n'$, заключаем, что все следующие члены последовательности $a_{k+1}$, $a_{k+2}, \dots$ также заключены между $\alpha_n$ и  $\alpha_n'$, то есть
если $m > k$, то 
\[\alpha_n< a_m<\alpha_n'.\]

Так как вещественное число $\alpha$ также заключено между $\alpha_n$ и $\alpha_n'$, то абсолютная величина разности $a_m-\alpha$ меньше разности чисел $\alpha_n'$ и $\alpha_n$.
Но $\alpha_n'-\alpha_n=\tfrac1{10^n}$, следовательно,
\[|a_m-\alpha|<\tfrac1{10^n}.
\eqno(2)\]
Таким образом, для любого значения $n$ можно указать такое число $k$, что при $m \ge k$ и имеет место неравенство (2).
Так как при неограниченном возрастании $n$ дробь $\tfrac1{10^n}$ неограниченно убывает, то из неравенства (2) следует, что вещественное число $\alpha$ есть предел последовательности (1).
Таким образом, числовая последовательность (1) имеет определённый предел.

\subsection*{Длина окружности}

Понятие о пределе даст возможность точно определить, что мы подразумеваем под длиной окружности.
Предварительно докажем следующие леммы.

\paragraph{}\label{1938/232}
\so{Лемма} 1.
\textbf{\emph{Выпуклая ломаная}} ($ABCD$, рис.~\ref{1938/ris-231}) \textbf{\emph{меньше всякой другой ломаной}} ($AEFGD$), \textbf{\emph{объемлющей первую.}}

\begin{wrapfigure}{O}{45mm}
\centering
\includegraphics{mppics/ris-231}
\caption{}\label{1938/ris-231}
\end{wrapfigure}

Выражения «объемлющая ломаная», «объемлемая ломаная» имеют следующий смысл.

Пусть две ломаные (как те, которые изображены у нас на рисунке) имеют одни и те же концы $A$ и $D$ и расположены таким образом, что одна ломаная ($ABCD$) вся лежит внутри многоугольника, образованного другой ломаной и отрезком $AD$, соединяющим концы $A$ и $D$;
тогда внешняя ломаная называется \rindex{объемлющая ломаная}\textbf{объемлющей}, а внутренняя ломаная — {}\textbf{объемлемой}.

Предстоит доказать, что объемлемая ломаная $ABCD$ (если она выпуклая) короче всякой объемлющей линии $AEFGD$ (всё равно — выпуклой или невыпуклой), то есть
что
\[AB+BC+CD<AE+EF+FG+GD.\]
Продолжим стороны выпуклой ломаной так, как указано на рисунке.
Тогда, приняв во внимание, что отрезок меньше всякой ломаной, соединяющей концы отрезка, мы можем написать следующие неравенства.
\[AB+BH<AE+EH;\]
\[BC+CK<BH+HF+FG+GK;\]
\[CD<CK+KD.\]

Сложим почленно все эти неравенства и затем от обеих частей полученного неравенства отнимем вспомогательные отрезки $BH$ и $CK$;
далее, заменив сумму $EH+HF$ отрезком $EF$ и сумму $GK\z+KD$ отрезком $GD$, получим то неравенство, которое требовалось доказать.

\begin{figure}[!ht]
\begin{minipage}{.48\textwidth}
\centering
\includegraphics{mppics/ris-232}
\end{minipage}\hfill
\begin{minipage}{.48\textwidth}
\centering
\includegraphics{mppics/ris-233}
\end{minipage}

\medskip

\begin{minipage}{.48\textwidth}
\centering
\caption{}\label{1938/ris-232}
\end{minipage}\hfill
\begin{minipage}{.48\textwidth}
\centering
\caption{}\label{1938/ris-233}
\end{minipage}
\vskip-4mm
\end{figure}

{\small

\smallskip
\mbox{\so{Замечание}.}
Если бы объемлемая линия не была выпуклой (рис. \ref{1938/ris-232}), то изложенное доказательство нельзя было бы применить.
В этом случае объемлемая ломаная может оказаться и больше объемлющей.

}

\paragraph{}\label{1938/233}
\so{Лемма} 2.
\textbf{\emph{Периметр любого выпуклого многоугольника}} ($ABCD$) \textbf{\emph{меньше периметра всякого другого многоугольника}} ($MNPQRL$), \textbf{\emph{объемлющего первый}} (рис.~\ref{1938/ris-233}).

Требуется доказать, что
\[AB+BC+CD+DA < LM+MN+NP+PQ+QR+RL.\]
Продолжив в обоих направлениях сторону $AB$ выпуклого многоугольника, применим к ломаным линиям $ABCD$ и $ATMNPQRSD$, соединяющим точки $A$ и $B$, лемму предыдущего параграфа;
получим неравенство:
\[AB+BC+CD<AT+TM+MN+NP+PQ+QR+RS+SD.\]
С другой стороны, так как отрезок $ST$ меньше ломаной $SLT$, то можем написать:
\[TA+AD+DS<TL+SL.\]
Сложим почленно эти два неравенства и отнимем от обеих частей вспомогательные отрезки $AT$ и $DS$;
далее, заменив сумму $TL+TM$ отрезком $LM$ и сумму $LS+RS$ отрезком $LR$, получим то, что требовалось доказать.

\begin{wrapfigure}{o}{30mm}
\centering
\includegraphics{mppics/ris-234}
\caption{}\label{1938/ris-234}
\end{wrapfigure}

\paragraph{Определение длины окружности.}\label{1938/234}
Впишем в данную окружность (рис.~\ref{1938/ris-234}) правильный многоугольник, например треугольник, и на какой-нибудь прямой $MN$ (рис.~\ref{1938/ris-23}5) отложим отрезок $OP_1$, равный его периметру.

Удвоим теперь число сторон вписанного треугольника, то есть
вместо треугольника возьмём правильный вписанный шестиугольник.
Найдём также его периметр и отложим его на той же прямой $MN$ от той же точки $O$;
пусть тогда получится отрезок $OP_2$, который должен быть больше $OP_1$, так как вместо каждой стороны треугольника мы теперь берём ломаную (из двух сторон шестиугольника), которая длиннее прямой.

\begin{figure}[!ht]
\centering
\includegraphics{mppics/ris-235}
\caption{}\label{1938/ris-235}
\end{figure}

Удвоим снова число сторон вписанного шестиугольника, то есть возьмём теперь правильный 12-угольник, найдём его периметр и отложим его на $MN$ от той же точки $O$; 
мы получим тогда отрезок $OP_3$, который будет больше $OP_2$ по той же причине, по какой $OP_2$ больше $OP_1$.

Вообразим, что такой процесс удвоения и откладывания периметров продолжается всё далее и далее.
Тогда мы получим неограниченную последовательность периметров $OP_1, OP_2, OP_3,$ и так далее, которая является возрастающей последовательностью.
Однако возрастание это не может быть неограниченным, так как периметр всякого вписанного многоугольника (выпуклого), каково бы ни было число его сторон, всегда остаётся меньше периметра любого описанного многоугольника (как его объемлющего).
Вследствие этого полученная последовательность периметров правильных вписанных многоугольников имеет определённый предел (§~\ref{1938/229}), обозначенный на чертеже $Ox$.
Этот предел и принимают за длину окружности.

Таким образом, мы принимаем следующее \so{определение}:
\emph{за длину окружности принимается тот предел, к которому стремится (приближается) периметр правильного многоугольника, вписанного в эту окружность, когда число сторон его неограниченно удваивается.}

{\small
\smallskip
\so{Замечание}.
Можно доказать (мы опускаем это доказательство), что предел этот не зависит от того, с какого многоугольника мы начинаем удвоение.
Более того, если даже вписанные многоугольники и не будут правильные, всё же периметры их стремятся к тому же самому пределу, как и периметры правильных многоугольников, лишь бы только стороны их неограниченно уменьшались (и, следовательно, число сторон их неограниченно увеличивалось) путём ли удвоения, как мы это предполагали для правильных многоугольников, или по какому-нибудь иному закону (мы опускаем это доказательство).

}

{

\begin{wrapfigure}[7]{o}{41mm}
\centering
\includegraphics{mppics/ris-236}
\caption{}\label{1938/ris-236}
\end{wrapfigure}

Таким образом, для каждой окружности существует свой единственный предел, к которому стремится периметр вписанного выпуклого многоугольника, когда стороны его неограниченно уменьшаются;
предел этот и принимается за длину окружности.

Равным образом за длину какой-нибудь дуги окружности $AB$ (рис.~\ref{1938/ris-236}) принимается предел, к которому стремятся длины ломаных линий, вписанных в эту дугу и имеющих с ней одни и те же концы, 
когда их стороны неограниченно уменьшаются.

}

\paragraph{}\label{1938/235}
\so{Лемма}.
\textbf{\emph{Длина дуги окружности
1) больше стягивающей её хорды, но 2) меньше периметра всякой ломаной линии, описанной около этой дуги и имеющей с ней одни и те же концы}} (рис.~\ref{1938/ris-237}).

1) Пусть $ACB$ (рис.~\ref{1938/ris-236}) — дуга окружности и $AB$ — стягивающая её хорда;
требуется доказать, что дуга больше этой хорды.

Предположим, что в дугу мы вписываем правильные ломаные таким образом:
первая ломаная будет составлена из двух хорд $AC$ и $CD$;
вторую ломаную получим путём удвоения числа сторон первой ломаной;
это будет ломаная $ADCEB$, состоящая из четырёх хорд;
третью ломаную получим удвоением числа сторон второй ломаной;
она будет состоять из восьми хорд.
Вообразим, что этот процесс удвоения продолжается неограниченно.
Тогда с каждым удвоением периметр ломаной будет возрастать;
например:
\[AD+DC+CE+EB>AC+CB,\]
так как
\[AD+DC>AC\quad\text{и}\quad CE+EB>CB.\]
Вследствие этого предел, к которому стремятся длины ломаных, должен быть больше длины первой ломаной, то есть больше суммы $AC+CB$, и, значит, должен быть и подавно больше хорды $AB$.
Но предел этот принимается за длину дуги $ACB$, значит, эта дуга больше хорды $AB$.

\begin{wrapfigure}[6]{o}{45mm}
\vskip-3mm
\centering
\includegraphics{mppics/ris-237}
\caption{}\label{1938/ris-237}
\end{wrapfigure}

2) Пусть около дуги описана какая-нибудь ломаная линия (правильная или неправильная — всё равно) (рис.~\ref{1938/ris-237}).
Если концы ломаной совпадают с концами дуги, то эту дугу можно рассматривать как сумму нескольких дуг, из которых каждая объемлется ломаной, состоящей только из двух отрезков.
Пусть одна из таких частей будет дуга $AB$ (рис.~\ref{1938/ris-238}).
Докажем, что длина этой дуги меньше суммы $AC+CB$, которую мы для краткости обозначим одной буквой $S$.
Для доказательства возьмём вспомогательную ломаную $AmnB$, которая получится, если мы срежем угол $C$ каким-нибудь отрезком $mn$, \so{не пересекающимся с дугой} $AB$ (что всегда возможно, если ломаная описана, то есть составлена из касательных).
Обозначим длину этой вспомогательной ломаной $AmnB$ буквой $S_1$.
Так как $mn<mC+Cn$, то $S_1<S$.

Докажем, что предел, к которому стремятся длины правильных ломаных, вписанной в дугу $AB$, при неограниченном удвоении числа их сторон не может быть больше $S_1$.

\begin{wrapfigure}{o}{45mm}
\centering
\includegraphics{mppics/ris-238}
\caption{}\label{1938/ris-238}
\end{wrapfigure}

Обозначим этот предел буквой $L$ и допустим, что $L>S_1$.
Так как длины ломаных приближается к своему пределу \so{как угодно близко}, то разность между $L$ и такой длиной сделаться меньше разности $L-S_1$; тогда длина такой ломаной сделается больше $S_1$.
Но это невозможно, так как всякая выпуклая ломаная линия, вписанная в дугу $AB$, есть \so{объемлемая} по отношению к объемлющей ломаной $AmnB$ и потому она меньше $S_1$.
Следовательно, нельзя допустить, что $L>S_1$.
Но тогда $L$ должно быть или меньше $S_1$, или в крайнем случае равно $S_1$.
Но так как $S_1<S$, то и в этом и в другом случае должно быть:
$L<S$, что и требуется доказать.

\paragraph{Нахождение длины окружности.}\label{1938/237}
Для этой цели можно пользоваться \so{формулой удвоения}, которую мы вывели раньше (§~\ref{1938/224}), то есть формулой:
\[a_{2n}^2=2R^2-2R\sqrt{R^2-\frac{a_n^2}4}.\]

Если радиус $R$ примем за 1, то формула эта примет более простой вид:
\[a_{2n}^2=2-2\sqrt{1-\frac{a_n^2}4}.\]

Обозначая, по принятому, через $a_n$ сторону правильного вписанного $n$-угольника будем иметь:
$a_6=R=1$.
Применяя формулу удвоения, находим:
\begin{align*}
a_{12}^2&=2-2\sqrt{1-\tfrac14}=2-\sqrt3;
\\
a_{24}^2&=2-2\sqrt{1-\tfrac{a_{12}^2}{4}};
\\
a_{48}^2&=2-2\sqrt{1-\tfrac{a_{24}^2}{4}}\quad\text{и так далее.}
\end{align*}
Положим, что мы прекратили удвоение на 96-угольнике.
Чтобы получить его периметр, надо сторону умножить на 96.
Этот периметр можно принять за приближённое значение длины окружности.
Обозначив его через $p_{96}$ и выполнив вычисления, найдём:
\[p_{96} = 6{,}2820638\dots\]
При радиусе, равном $R$, получим:
\[p_{96}=R\cdot6{,}2820638\dots,
\quad\text{или}\quad
p_{96}=2R\cdot3{,}1410319\dots\]
Обозначая длину окружности буквой $C$, мы получим для неё приближённую формулу:
\[C\approx 2R\cdot3{,}1410319\dots\]
Если бы мы прекратили процесс удвоения на 192-угольнике, то получили бы для длины окружности более точное значение, именно:
\[C\approx2R\cdot3{,}14145247\dots\]
Продолжая процесс удвоения, можно получать для длины окружности всё более и более точные значения.

\paragraph{Отношение длины окружности к диаметру.}\label{1938/238}
Рассматривая процесс нахождения длины окружности, можно заметить, что число, на которое нужно умножить диаметр, чтобы получить длину окружности, не зависит от величины самого диаметра, так что если мы нашли, что длина какой-нибудь окружности равна её диаметру, умноженному на некоторое число, то и длина всякой другой окружности будет равна её диаметру, умноженному на то же самое число.

В самом деле, возьмём две окружности:
одну радиуса $R$, другую радиуса $r$.
Длину первой окружности обозначим через $C$, длину второй — через $c$.
Впишем в каждую из них правильный многоугольник с одним и тем же числом сторон и будем удваивать число сторон каждого из этих многоугольников.

Обозначим через $P_n$ периметр правильного $n$-угольника, вписанного в первую окружность, и через $p_n$ периметр правильного $n$-угольника, вписанного во вторую окружность.

В силу теорем, доказанных в §~\ref{1938/218}, мы можем написать:
\[\frac{P_n}{R}=\frac{p_n}{r}
\quad\text{или}\quad
\frac{P_n}{2R}=\frac{p_n}{2r}.\]

Периметры $P_n$ имеют пределом длину $C$ первой окружности.
Периметры $p_n$ имеют пределом длину $c$ второй окружности.
Значит, из равенства $\frac{P_n}{2R}=\frac{p_n}{2r}$ вытекает
\[\frac{C}{2R}=\frac{c}{2r}.\]
Действительно, поскольку $P_n$ стремится к $C$,
разница $C-P_n$ неограниченно убывает с возрастанием $n$;
а значит и разница
\[\frac{C}{2R}-\frac{P_n}{2R}=\frac{1}{2R}(C-P_n)\]
также неограниченно убывает, то есть последовательность $\frac{P_n}{2R}$ стремится к $\frac{C}{2R}$.
Далее, из равенств $\frac{P_n}{2R}=\frac{p_n}{2r}$,
мы получаем, что последовательность $\frac{p_n}{2r}$ также стремится к $\frac{C}{2R}$.
С другой стороны, эта последовательность должна стремиться к $\frac{c}{2r}$ и значит $\frac{C}{2R}=\frac{c}{2r}$.

Таким образом, мы можем сказать, что \emph{отношение длины окружности к её диаметру есть число постоянное для всех окружностей}. 
Это постоянное число принято обозначать греческой буквой $\pi$ (читается «пи»).%
\footnote{Буква $\pi$ есть начальная буква греческого слова \textgreek{περιφέρεια} (окружность).
Обозначение это введено, по всей вероятности, в XVII веке.
}

Мы можем, таким образом, для длины $C$ окружности написать такую формулу:
\[C=2R\cdot\pi
\quad\text{или}\quad
C=2\pi R.
\]
Доказано, что число $\pi$ является числом иррациональным, и, значит, оно не может быть выражено точно никаким рациональным числом.
Но его приближённые значения можно находить различными способами с какой угодно точностью.
Приняв периметр вписанного 96-угольника за приближённую длину окружности, мы получим для $\pi$ приближённое значение $3{,}14$ с недостатком и с точностью до $0{,}01$.
Эта точность для практических целей почти всегда достаточна.
В случаях особенной точности можно довольствоваться таким приближённым значением (с избытком):
$\pi = 3{,}1416$.\footnote{Для запоминания довольно длинного ряда цифр, выражающих число $\pi$, можно пользоваться следующим  двустишием (придуманным покойным преподавателем средней школы Шенроком):
\begin{verse}
Кто и шутя и скоро пожелает(ъ).\\
Пи узнать число уж(ъ) знает(ъ).
\end{verse}
Если выписать в ряд числа букв, заключающихся в каждом слове этих фраз (написанных по старой орфографии), то получим для $\pi$:
приближённое число (с избытком), $3{,}141596536$, верное до одной половины десятибиллионной.}%

Ещё в III веке до нашей эры знаменитый сиракузский геометр Архимед нашёл для $\pi$ очень простое приближение $\tfrac{22}7$ то есть $3\tfrac17$.
Это число несколько более $\pi$ и разнится от него менее чем на 2 тысячных;
оно равно второй подходящей дроби для разложения $\pi$ в цепную дробь
\[\pi=3+\frac{1}{7+\frac{1}{15+\frac{1}{1+\cdots}}}.\]
В V веке, китайский математик Цзу Чунчжи нашёл приближение $\tfrac{355}{113}$, это четвёртая подходящая дробь. 
Сейчас вычислены многие триллионы знаков в десятичном разложении $\pi$, что далеко превосходит всякие практические требования.
(Например 40 знаков числа $\pi$ более чем достаточно чтобы вычислить длину окружности радиуса сравнимого с расстоянием до самой далёкой видимой звезды и с точностью превышающей размер атома.) 

При решении геометрических задач часто встречается число, обратное числу $\pi$, то есть равное дроби $\tfrac1\pi$.
Полезно запомнить несколько цифр этого числа.
\[\tfrac1\pi = 0{,}3183098\dots\]

{\small
\sloppy

\smallskip
\so{Замечание}.
Число $\pi$ не является рациональным;
более того, оно \so{трансцендентно}, то есть оно не может служить корнем никакого алгебраического уравнения с рациональными коэффициентами.
Это было доказано в 1882 году немецким математиком Фердинандом фон Линдеманом.
В~частности, с помощью циркуля и линейки нельзя решить построением задачу о выпрямлении окружности (§~\ref{1938/211}), то есть нельзя построить такой отрезок, длина которого в точности равнялась бы длине данной окружности.

}

\paragraph{Длина дуги, содержащей $\bm{n}$ градусов.}\label{1938/239}
Длина окружности есть $2\pi R$, значит, длина дуги в $1\degree$ равна $\frac{2\pi R}{360}=\frac{\pi R}{180}$; следовательно, длина $s$ дуги, содержащей $n\degree$, выразится так:
\[s=\frac{\pi R n}{180}.\]
Если дуга выражена в минутах ($n'$) или в секундах ($n''$), то длина её определяется соответственно формулами:
\[s=\frac{\pi R n}{180\cdot 60}
\quad\text{или}\quad
s=\frac{\pi R n}{180\cdot 60\cdot 60},\]
где $n$ — число минут или секунд.

\paragraph{}\label{1938/240}
\so{Задача}.
\emph{Вычислить с точностью до 1 мм радиус такой окружности, дуга которой, содержащая $81\degree 21'36''$, равна $0{,}452$ м.}

Обратив $81\degree 21'36''$ в секунды, получим число 292896.
Из уравнения
\[0{,}452 = \frac{\pi R\cdot  292896}{180\cdot 60\cdot 60}\]
находим:
\[R=\frac{0{,}452\cdot 180\cdot 60\cdot 60}{292896\cdot \pi}=\frac1\pi\approx0{,}318 (\text{м}).\]

\paragraph{}\label{1938/241}
\so{Задача}.
\emph{Определить число градусов дуги, длина которой равна радиусу.}

Заменив в формуле, определяющей длину дуги в $n\degree$, величину $s$ на $R$, получим уравнение:
\[R=\frac{\pi R n}{180}
\quad\text{или}\quad
1=\frac{\pi n}{180}\]
откуда
\[n\degree = 180\degree \cdot \frac1\pi \approx 180\degree \cdot 0{,}3183098 = 57{,}295764\degree \approx 57\degree17'44''\]

\paragraph{Радианы.}\label{extra/radians}
Дуга, равная радиусу, называется \rindex{радиан}\textbf{радианом}.
В радианах принято мерить и углы — за угол в один радиан берут центральный угол, который вырезает из окружности дугу равную её радиусу.

Радианы считаются за основную единицу измерения углов — если единица измерения углов не указывается, то всегда имеется в виду радиан.
Например, прямой угол это «угол равный $\tfrac\pi2$», развёрнутый угол — «угол равный $\pi$».
Градус это второстепенная единица удобная для работы с дробными частями прямого угла;
$1\degree =\tfrac\pi{180}\approx 0{,}01745$.

Заметим, что если $\alpha$ есть угловая величина дуги $AB$ окружности радиуса $R$ выраженная в радианах то её длину можно выразить формулой
\[s=\alpha\cdot R.\]
Аналогичная формула для угла измеренного в градусах (§~\ref{1938/239}) существенно сложнее.

{\small

\subsection*{Упражнения}

\begin{enumerate}[noitemsep]

\item
Доказать, что в двух кругах отношение центральных углов, соответствующих дугам, имеющим одинаковую длину, равно обратному отношению радиусов.

\item
На окружности взята точка $A$ и через неё проведены:
диаметр $AB$, сторона правильного вписанного шестиугольника $AC$ и касательная $MN$.
Из центра $O$ опущен на $AC$ перпендикуляр и продолжен до пересечения с касательной в точке $D$.
От этой точки отложен по касательной (через точку $A$) отрезок $BE$, равный трём радиусам.
Точка $E$ соединена с концом диаметра $B$.
Определить, как велика погрешность, если прямую $BE$ возьмём за длину полуокружности.

\item
На диаметре данной полуокружности построены две равные полуокружности, и в ту часть плоскости, которая заключена между тремя полуокружностями, вписан круг.
Доказать, что диаметр этого круга относится к диаметру равных полуокружностей, как $2:3$.

\item
Вычислить в градусах, минутах и секундах дугу, равную стороне квадрата, вписанного в эту окружность.

\item
Вычислить длину $1\degree$ земного экватора, принимая радиус Земли равным 6400 км.

\end{enumerate}

}

%% file: 2D/ponyatie-o-ploschadi.tex
\section{Понятие о площади}

Каждый из нас имеет некоторое представление о площади из повседневной жизни.

Мы займёмся уточнением понятия о площади фигуры и установлением способов её измерения.

\paragraph{Основные допущения о площадях.}\label{1938/243}
Площадь многоугольника, это положительное число, удовлетворяющее следующим условиям:

1) площади двух равных многоугольников должны быть равны между собой;

\begin{wrapfigure}{o}{42mm}
\centering
\includegraphics{mppics/ris-239}
\caption{}\label{1938/ris-239}
\end{wrapfigure}

2) если данный многоугольник разбит на несколько многоугольников ($M$, $N$, $P$, рис.~\ref{1938/ris-239}), то  площадь всего многоугольника равна сумме площадей отдельных его частей.

3) площадь квадрата со стороной равной линейной единице длины считается равной квадратной единице.

{\small
\smallskip
\mbox{\so{Замечания}.}
1) Единица площади зависит от выбора единицы длины и в каждой конкретной задаче мы вольны выбрать удобную единицу длины и соответственно площади, важно только, чтобы все измерения были проведены, используя только эту единицу.
Например, если за линейную единицу взят отрезок длины 1 метр, то квадрат со стороной 1 метр имеет площадь равную 1 квадратный метр.

2) Многоугольники, имеющие равные площади, принято называть \rindex{равновеликие многоугольники}\textbf{равновеликими}.
В силу условия 1, равные многоугольники всегда и равновелики, но равновеликие многоугольники могут быть неравными (как те, которые изображены на рис.~\ref{1938/ris-240}).

{\sloppy
3) Поскольку площади измеряются положительными числами, а сумма двух положительных чисел всегда больше каждого из слагаемых, то из условия 2 мы получаем следующее заключение:
\emph{Площадь любого многоугольника больше площади любого другого многоугольника целиком в нём лежащего.}
Это утверждение иногда формулируется сокращённо как \emph{целое больше части}.
}

}

\begin{wrapfigure}{o}{60mm}
\centering
\includegraphics{mppics/ris-240}
\caption{}\label{1938/ris-240}
\end{wrapfigure}

4) Положим, что, разбив данный многоугольник $M$ на несколько многоугольников, мы будем переставлять эти части и получать таким образом новые многоугольники (подобно тому, как на рис.~\ref{1938/ris-240} перемещены части $A$ и $B$).
Спрашивается:
нельзя ли путём этих перестановок получить многоугольник $M'$, который мог бы целиком уместиться внутри $M$?
Если бы это оказалось возможным, то поскольку целое больше части, мы получили бы, что
\[\text{площадь}~M'<\text{площади}~M,\]
а при этом, в силу условий 1 и 2, мы получили бы, что
\[\text{площадь}~M'=\text{площади}~M.\]

Эти два утверждения противоречат друг другу.
Значит площадь многоугольников $M$ и $M'$ невозможно было бы определить так, чтобы удовлетворялись все условия.

Таким образом возможность определить площадь для многоугольников вовсе не очевидна.
Впервые обратил внимание на этот вопрос итальянский математик Антонио Де Цольт (1881).
Невозможность указанной выше перестановки частей многоугольника принималась вначале как некоторый постулат, но позднее эта невозможность была строго доказана.\footnote{Это доказательство довольно сложное; оно приводится например в статье «Площадь и объём» В. А. Рохлина (Энциклопедия элементарной математики, книга пятая, Геометрия).}
Используя  это утверждение можно доказать, что каждый многоугольник (а также фигуры из более широкого класса) имеет определённую площадь, удовлетворяющую трём условиям выше (§~\ref{extra/kvad-fig}).

Начиная с §~\ref{1938/245}, нас будет интересовать как измерить площадь данного многоугольника принимая без доказательства то, что каждый многоугольник имеет определённую площадь удовлетворяющую трём условиям, указанным выше.

\begin{wrapfigure}{r}{22mm}
\centering
\includegraphics{mppics/ris-extra-2}
\caption{}\label{extra/ris-2}
\end{wrapfigure}

{\small

\paragraph{Об измерении площади.}\label{1938/244}
Заметим, что единичный квадрат можно разбить на $n^2$ равных квадратов со стороной $\tfrac1n$. 
Поскольку площади равных фигур равны, мы получаем, что все эти квадраты имеют одну и ту же площадь, обозначим её за $s$.
Далее из свойства 2 (§~\ref{1938/243}), заключаем, что площадь единичного квадрата равна $n^2\cdot s$.
Поскольку площадь единичного квадрата равна единице (условие 3 в §~\ref{1938/243}) получаем, что площадь квадрата со стороной $\tfrac1n$ равна $\tfrac1{n^2}$.

\begin{wrapfigure}{o}{55mm}
\centering
\includegraphics{mppics/ris-1931-250}
\caption{}\label{1931/ris-250}
\end{wrapfigure}

Допустим, что на многоугольник $M$, площадь которого надо измерить, наложена сеть из единичных квадратов.
По отношению к данному многоугольнику $M$, все квадраты сети можно разбить на
три рода: 
1) внешние квадраты, 
2) внутренние квадраты и 
3) оставшиеся квадраты, то есть те через которые проходит контур многоугольника.
Оставив без внимания внешние квадраты, сосчитаем отдельно квадраты внутренние и квадраты 3-го рода.
Пусть первых окажется $m$, а вторых~$n$.
Тогда, очевидно, площадь $M$ будет больше $m$, но меньше $m+n$ квадратных единиц.
Числа $m$ и $m+n$ будут в этом случае приближённые значения измеряемой площади, первое число с недостатком, второе с избытком, причём погрешность этого измерения меньше $n$ квадратных единиц.

Чтобы получить более точный результат, уплотним сеть квадратов, подразделив каждый из них на более
мелкие квадраты.
Например, мы можем разбить каждый единичный квадрат на 100 квадратов со стороной $\tfrac1{10}$.
Тогда мы получим другие приближённые меры площади, причём погрешность будет не больше прежней (так как все квадраты 3-го рода в уплотнённой сети лежат внутри квадратов 3-го рода в изначальной сети).

\paragraph{Квадрируемые фигуры.}\label{extra/kvad-fig}
Построение в предыдущем параграфе, применимо к произвольным ограниченным фигурам, необязательно многоугольникам.
Если при последовательном уплотнении сети квадратов погрешность измерения стремится к нулю,
то фигура называется \textbf{квадрируемой}.
Для квадрируемой фигуры, общий предел приближённых значений площадей с недостатком и избытком принимается равным её площади.
При этом площадь некоторых квадрируемых фигур может равняться нулю;
например несложно видеть, что фигура состоящая из одной точки квадрируема и её площадь  равна нулю.

Как пример неквадрируемой фигуры, представьте себе фигуру $F$ состоящую из точек единичного квадрата, отстоящих от одной из его сторон на рациональное расстояние (такую фигуру невозможно нарисовать).
При попытке приближённо измерить площадь $F$ с помощью любой сетки, мы увидим только квадраты рода 1 и 3 — любой квадрат содержащий точки из $F$ будет также содержать точки вне $F$ (лежащие на иррациональном расстоянии от стороны).
При этом погрешность измерения, то есть общая площадь квадратов рода 3 не может быть меньше единицы. 
То есть $F$ не является квадрируемой.

Следующие три утверждения являются ключевыми, при формальном введении понятия площади.
Мы приводим только идею их доказательства; полные доказательства не сложные, но довольно громоздкие.

1) \emph{Любой многоугольник является квадрируемой фигурой};
то есть при последовательном уплотнении сети квадратов погрешность измерения его площади стремится к нулю. \emph{То же верно для любой фигуры ограниченной конечным набором дуг окружностей и отрезков прямых.}
В частности, круги, сектора и сегменты являются кварируемыми фигурами.

Действительно, заметим, что погрешность измерения равна общей площади квадратов сетки прорезаемых контуром фигуры.
При этом сам контур состоит из нескольких отрезков и дуг.
Значит достаточно доказать, что общая площадь $s$ всех квадратов сетки прорезаемых одним отрезком (или дугой) стремится к нулю при неограниченном уплотнении сетки.
Последнее утверждение следует из неравенства, которое мы предлагаем проверить его самостоятельно:
\[s\le 4\cdot (\ell+a)\cdot a,\]
где $\ell$ есть длина отрезка (или дуги), а $a$ — сторона квадрата сетки.

2) \emph{Если фигура квадрируема в данной сетке то она квадрируема и в любой другой сетке.}
Рассмотрим две сети из квадратов равного размера.
Не трудно видеть, что любой квадрат первой сети может пересекать не более чем девять квадратов второй сети.

Отсюда следует, что общая площадь квадратов третьего рода для первой сети  отличается не более чем в девять раз от общей площади квадратов третьего рода для второй сети.
Действительно, любой квадрат третьего рода содержит точки контура фигуры.
Значит любой квадрат третьего рода второй сети должен пересекаться с каким-то квадратом третьего рода первой сети.
То есть на один квадрат первой сети приходится не более девяти квадратов второй.

Погрешность измерения площади для каждой сети равна общей площади квадратов третьего рода.
Значит погрешности измерений в обеих сетках отличаются друг от друга не более чем в 9 раз.
То же верно и для одинаковых уплотнений сетей.
Следовательно, если погрешность измерения для уплотнений первой сети стремится к нулю, то тоже верно и для второй сети. 

3) \emph{Полученное значение площади квадрируемой фигуры не зависит от выбора сетки.}
В противном случае при неограниченном уплотнении обоих сеток мы получили бы, что
приближённое значение площади взятое с избытком одной сетке меньше
приближённого значения площади взятое с недостатком в другой.
Но эти приближённые значения равны площадям некоторых многоугольников — многоугольника $M_1$, составленного из квадратов 2-го и 3-го рода первой сетки,
и многоугольника $M_2$, составленного из квадратов 2-го рода второй сетки.
Заметим, что многоугольник $M_1$ содержит $F$, а многоугольник $M_2$ содержится в $F$.
В частности $M_2$ есть часть $M_1$, а целое должно быть больше части (§~\ref{1938/243}) — противоречие.

}


%% file: 2D/ploschadi-mnogugov.tex
\section{Подсчёт площадей многоугольников}

\paragraph{Основание и высота.}\label{1938/245}
Условимся одну из сторон треугольника или параллелограмма называть \rindex{основание}\textbf{основанием} этих многоугольников, а перпендикуляр, опущенный на эту сторону из вершины треугольника или из какой-нибудь точки противоположной стороны параллелограмма, будем называть \rindex{высота}\textbf{высотой}.

В прямоугольнике за высоту можно взять сторону, перпендикулярную к той, которая принята за основание.

В трапеции основаниями называют обе параллельные стороны, а высотой — общий перпендикуляр между ними.

Основание и высота прямоугольника называются его \rindex{измерения!прямоугольника}\textbf{измерениями}.

\paragraph{}\label{1938/246}
\so{Теорема}.
\textbf{\emph{Площадь прямоугольника равна произведению его основания на высоту.}}

При доказательстве могут представиться три случая:

1) Длины основания и высоты (измеренных одной и той же единицей) выражаются \so{целыми числами}.

\begin{wrapfigure}{o}{29mm}
\vskip-4mm
\centering
\includegraphics{mppics/ris-242}
\caption{}\label{1938/ris-242}
\end{wrapfigure}

Пусть у данного прямоугольника (рис.~\ref{1938/ris-242}) основание равно целому числу $b$ линейных единиц, а высота — целому числу $h$ тех же единиц.
Разделим основание на $b$ и высоту на $h$ равных частей, проведём через точки деления ряд прямых, параллельных высоте, и другой ряд прямых, параллельных основанию.
От взаимного пересечения этих прямых образуются некоторые четырёхугольники.

Возьмём какой-нибудь один из них, например четырёхугольник $K$ (покрытый на рисунке штрихами).
Так как стороны этого четырёхугольника, по построению, параллельны соответствующим сторонам данного прямоугольника, то все углы его прямые;
значит, четырёхугольник $K$ есть прямоугольник.
С другой стороны, каждая сторона этого прямоугольника равна расстоянию между соседними параллельными прямыми, то есть равна одной и той же линейной единице. 
Значит, прямоугольник $K$ представляет собой квадрат, а именно, ту квадратную единицу, которая соответствует взятой линейной единице (если, например, за линейную единицу принят линейный сантиметр, то площадь квадрата $K$ есть квадратный сантиметр).

Так как сказанное об одном четырёхугольнике справедливо и для всякого другого, то, значит, проведением указанных параллельных прямых мы разбиваем всю площадь данного прямоугольника на квадратные единицы.
Найдём их число.

Очевидно, что ряд прямых, параллельных основанию, разделяет прямоугольник на столько равных горизонтальных полос, сколько в высоте содержится линейных единиц, то есть на $h$ равных полос.
С другой стороны, ряд прямых, параллельных высоте, разбивает каждую горизонтальную полосу на столько квадратных единиц, сколько в основании содержится линейных единиц, то есть на $b$ квадратных единиц.
Значит, всех квадратных единиц окажется $b\cdot h$.
Таким образом,
\[\text{площадь прямоугольника} = b\cdot h.\]
то есть она равна произведению основания на высоту.

2) Длины основания и высоты (измеренных одной и той же единицей) выражаются \so{дробными числами}.

Пусть, например, у данного прямоугольника
\[\text{основание} = 3\tfrac12=\tfrac72~\text{линейных единиц.}\]
\[\text{высота} = 4\tfrac35 = \tfrac{23}5~\text{той же единицы.}\]
Приведя дроби к одинаковому знаменателю, получим:
\[\text{основание} = \tfrac{35}{10};
\quad
\text{высота} = \tfrac{46}{10}.
\]
Примем $\tfrac1{10}$ долю линейной единицы за новую единицу длины.

Тогда мы можем сказать, что основание содержит 35 этих новых единиц, а высота 46 тех же единиц.
Значит, по доказанному, в случае 1-м, площадь прямоугольника равна $35 \cdot 46$ таких квадратных единиц, которые соответствуют $\tfrac1{100}$ новой единице длины.
Но эта квадратная единица составляет часть квадратной единицы, соответствующей прежней линейной единице;
значит, площадь прямоугольника в прежних квадратных единицах равна:
\[\frac{35\cdot 46}{100}=\frac{35}{10}\cdot\frac{46}{10}=3\tfrac12\cdot4\tfrac35(\text{квадратных единиц}).\]

3) Основание и высота (или только одно из этих измерений) несоизмеримы с единицей длины, и, следовательно, их длины выражаются \so{иррациональными числами}.

В этом случае можно довольствоваться приближённым результатом измерения площади с желаемой степенью точности.
Но можно и в этом случае найти точную меру площади прямоугольника.

\begin{wrapfigure}{O}{45mm}
\centering
\includegraphics{mppics/ris-243}
\caption{}\label{1938/ris-243}
\end{wrapfigure}

Пусть длина основания $AB$ прямоугольника $ABCD$ (рис.~\ref{1938/ris-243}) выражается  числом $\alpha$, а для высоты $AD$ — числом $\beta$.
Каждое из этих чисел может быть представлено в виде бесконечной десятичной дроби (§~\ref{1938/150}).
Возьмём приближённые значения этих чисел в виде десятичных дробей с $n$ десятичными знаками сначала с недостатком, затем с избытком.
Приближённые значения с недостатком обозначим через $\alpha_n$ (для первого числа) и $\beta_n$ (для второго числа), а приближённые значения с избытком соответственно через  $\alpha_n'$ и $\beta_n'$.

Отложим на основании $AB$ от точки $A$ сначала отрезок $AB_1$ длины $\alpha_n$, затем отрезок $AB_2$ длины  $\alpha_n'$.
Очевидно, 
\[AB_1\le AB\le AB_2.\]
Отложим, далее, на высоте $AD$ от точки $A$ отрезки $AD_1$ и $AD_2$, длины которых равны соответственно $\beta_n$ и $\beta_n'$.
Очевидно, 
\[AD_1\le AD\le AD_2.\]

Построим вспомогательные прямоугольники $AB_1C_1D_1$ и $AB_2C_2D_2$, у каждого из них основание и высота выражаются рациональными числами:
\begin{align*}
AB_1&=\alpha_n,
&
AD_1&=\beta_n,
\\
AB_2&=\alpha_n',
&
AD_2&=\beta_n'.
\end{align*}
Поэтому, согласно доказанному в случае 2-м.
\begin{align*}
\text{площадь}~AB_1C_1D_1 &= \alpha_n\cdot \beta_n, 
\\
\text{площадь}~AB_2C_2D_2 &= \alpha_n'\cdot \beta_n'. 
\end{align*}

Поскольку прямоугольник $AB_1C_1D_1$ составляет часть прямоугольника $ABCD$, 
а прямоугольник $ABCD$ составляет часть прямоугольника $AB_2C_2D_2$;
имеем
\[\alpha_n\cdot \beta_n\le\text{площадь}~ABCD\le\alpha_n'\cdot \beta_n'.\]
Это двойное неравенство остаётся верным при всяком значении $n$, то есть оно остаётся верным, с какой бы точностью мы ни находили приближённые значения чисел $\alpha$ и $\beta$.
Значит, мы можем сказать, что площадь $ABCD$ должна выражаться таким числом, которое больше произведения любых приближённых значений чисел $\alpha$ и $\beta$, если эти значения взяты с недостатком, но меньше произведения любых их приближённых значений, если эти значения взяты с избытком.
Но существует единственное такое число, это произведение чисел $\alpha$ и $\beta$ (§~\ref{1938/154}).
Следовательно 
\[\text{площадь}~ABCD=\alpha\cdot \beta.\]
Таким образом, и в этом случае площадь прямоугольника равна произведению основания на высоту.

\paragraph{}\label{1938/247}
\so{Теорема}.
\textbf{\emph{Площадь параллелограмма}} ($ABCD$, рис.~\ref{1938/ris-244} и \ref{1938/ris-245}) \textbf{\emph{равна произведению основания на высоту.}}

На основании $AD$ (на том и другом рисунке) построим прямоугольник $AEFD$, у которого сторона $EF$ составляет продолжение стороны $BC$.

При этом могут представиться два случая:

\begin{figure}[!ht]
\begin{minipage}{.58\textwidth}
\centering
\includegraphics{mppics/ris-244}
\end{minipage}
\hfill
\begin{minipage}{.38\textwidth}
\centering
\includegraphics{mppics/ris-245}
\end{minipage}

\medskip

\begin{minipage}{.58\textwidth}
\centering
\caption{}\label{1938/ris-244}
\end{minipage}
\hfill
\begin{minipage}{.38\textwidth}
\centering
\caption{}\label{1938/ris-245}
\end{minipage}
\vskip-4mm
\end{figure}

1) сторона $BC$ лежит вне стороны $EF$ и 2) сторона $BC$ частью совпадает с $EF$ (первый случай изображён на рис.~\ref{1938/ris-244}, второй — на рис.~\ref{1938/ris-245}).
Докажем, что и в том, и другом случае
\[\text{площадь}~ABCD = \text{площади}~AEFD.\]
Если параллелограмм дополним треугольником $AEB$, а прямоугольник дополним треугольником $DFC$, то мы получим одну и ту же трапецию $AECD$.
Так как дополняющие треугольники равны (они имеют по две стороны и углу, заключённому между ними, соответственно равными), то параллелограмм и прямоугольник должны быть равновелики.
Но площадь $AEFD=b\cdot h$;
следовательно, и площадь $ABCD=b\cdot h$, причём $b$ можно рассматривать как основание параллелограмма и $h$ — как его высоту.


\paragraph{}\label{1938/248}
\so{Теорема}.
\textbf{\emph{Площадь треугольника}} ($ABC$, рис.~\ref{1938/ris-246}) \textbf{\emph{равна половине произведения основания на высоту.}}

Проведём $BE \parallel  AC$ и $AE \parallel BC$.
Тогда получим параллелограмм $AEBC$, площадь которого, по доказанному, равна $b\cdot h$.
Поскольку треугольники $ABC$ и $BAE$ равны, площадь $\triangle ABC$ составляет половину площади $AEBC$;
то есть
\[\text{площадь}~ \triangle ABC=\tfrac12b\cdot h.\]

\begin{figure}[!ht]
\begin{minipage}{.48\textwidth}
\centering
\includegraphics{mppics/ris-246}
\end{minipage}
\hfill
\begin{minipage}{.48\textwidth}
\centering
\includegraphics{mppics/ris-247}
\end{minipage}

\medskip

\begin{minipage}{.48\textwidth}
\centering
\caption{}\label{1938/ris-246}
\end{minipage}
\hfill
\begin{minipage}{.48\textwidth}
\centering
\caption{}\label{1938/ris-247}
\end{minipage}
\vskip-4mm
\end{figure}

{\small
\smallskip
\mbox{\so{Замечание}.}
Легко убедиться, что всякий треугольник разлагается на части, перемещением которых можно образовать прямоугольник, имеющий одинаковое с треугольником основание и высоту, вдвое меньшую высоты треугольника (рис.~\ref{1938/ris-247}).

}

\paragraph{}\label{1938/249}
\so{Следствия}.
1) \emph{Треугольники с равными основаниями и равными высотами равновелики.}

Если, например, вершину $B$ треугольника $ABC$ (рис.~\ref{1938/ris-248}) будем перемещать по прямой, параллельной основанию $AC$, а основание оставим то же самое, то площадь треугольника не будет изменяться.

\begin{figure}[!ht]
\begin{minipage}{.48\textwidth}
\centering
\includegraphics{mppics/ris-248}
\end{minipage}
\hfill
\begin{minipage}{.48\textwidth}
\centering
\includegraphics{mppics/ris-249}
\end{minipage}

\medskip

\begin{minipage}{.48\textwidth}
\centering
\caption{}\label{1938/ris-248}
\end{minipage}
\hfill
\begin{minipage}{.48\textwidth}
\centering
\caption{}\label{1938/ris-249}
\end{minipage}
\vskip-4mm
\end{figure}

2) \emph{Площадь прямоугольного треугольника равна половине произведения его катетов,} потому что один катет можно взять за основание, а другой — за высоту.


3) \emph{Площадь ромба равна половине произведения его диагоналей.}
Действительно, если $ABCD$ (рис.~\ref{1938/ris-249}) есть ромб, то его диагонали взаимно перпендикулярны.
Поэтому

\begin{center}
\begin{tabular}{rrcl}
 &$\text{площадь}~\triangle ABC$
 &$=$ 
 &$\tfrac12AC\cdot  OB$
\\
\raisebox{3mm}[3mm][3mm]{$+$} 
& 
$\text{площадь}~\triangle ACD$
&$=$ &$\tfrac12AC\cdot  OD$
\\
\hline
&$\text{площадь}~ABCD$ & $=$&$\tfrac12AC\cdot (OB+OD)=\tfrac12AC\cdot BD.$
\end{tabular}
\end{center}

4) \emph{Площади двух треугольников относятся как произведения их оснований на высоты} (множитель $\tfrac12$ сокращается).

\paragraph{}\label{1938/250}
\so{Теорема}.
\textbf{\emph{Площадь трапеции равна произведению полусуммы оснований на высоту.}}

Проведя в трапеции $ABCD$ (рис.~\ref{1938/ris-250}) диагональ $AC$, мы можем рассматривать её площадь как сумму площадей двух треугольников $CAD$ и $ABC$.
Поэтому
\[\text{площадь трапеции}~ABCD
=\tfrac12AD\cdot  h+\tfrac12BC\cdot  h= 
\tfrac12(AD+BC) \cdot  h.\]

\begin{figure}[!ht]
\begin{minipage}{.48\textwidth}
\centering
\includegraphics{mppics/ris-250}
\end{minipage}
\hfill
\begin{minipage}{.48\textwidth}
\centering
\includegraphics{mppics/ris-251}
\end{minipage}

\medskip

\begin{minipage}{.48\textwidth}
\centering
\caption{}\label{1938/ris-250}
\end{minipage}
\hfill
\begin{minipage}{.48\textwidth}
\centering
\caption{}\label{1938/ris-251}
\end{minipage}
\vskip-4mm
\end{figure}

\paragraph{}\label{1938/251}
\so{Следствие}.
Если $MN$ (рис.~\ref{1938/ris-251}) есть средняя линия трапеции, то, как известно (§~\ref{1938/99}).
\[MN = \tfrac12(AD+BC).\]
Поэтому
\[\text{площадь трапеции}~ABCD=MN\cdot  h.\]
то есть \emph{площадь трапеции равна произведению средней линии на высоту.}
Это же можно видеть и непосредственно из рис.~\ref{1938/ris-251}.

\paragraph{}\label{1938/252}
\so{Теорема}.
\textbf{\emph{Площадь всякого описанного многоугольника равна произведению периметра на половину радиуса.}}

\begin{wrapfigure}{o}{35mm}
\centering
\includegraphics{mppics/ris-252}
\caption{}\label{1938/ris-252}
\end{wrapfigure}

Соединив центр $O$ (рис.~\ref{1938/ris-252}) со всеми вершинами описанного многоугольника, мы разделим его на треугольники, в которых за основания можно взять стороны многоугольника, а за высоты — радиус круга.

Обозначив этот радиус через $R$, будем иметь:
\begin{align*}
\text{площадь}~\triangle AOB&=AB \cdot  \tfrac12R,
\\
\text{площадь}~\triangle BOC &= BC \cdot  \tfrac12R\quad\text{и~т.~д.}
\end{align*}
Следовательно,
\[\text{площадь}~ABCDE = (AB+BC+CD+DE+EA) \cdot  \tfrac12R= P \cdot \tfrac12R,\]
где буквой $P$ обозначен периметр многоугольника.

\smallskip
\so{Следствие}.
\emph{Площадь правильного многоугольника равна произведению периметра на половину апофемы,} потому что всякий правильный многоугольник можно рассматривать как описанный около круга, у которого радиус есть апофема.

\paragraph{Площадь неправильного многоугольника.}\label{1938/253}
Для нахождения площади какого-нибудь неправильного многоугольника можно его разбить на треугольники (например, диагоналями), вычислить площадь каждого треугольника в отдельности и результаты сложить.

{\sloppy

\paragraph{}\label{1938/254}
\mbox{\so{Задача}.}
\emph{Построить треугольник, равновеликий данному многоугольнику} ($ABCDE$, рис.~\ref{1938/ris-253}).

}

\begin{wrapfigure}{o}{38mm}
\vskip-2mm
\centering
\includegraphics{mppics/ris-253}
\caption{}\label{1938/ris-253}
\end{wrapfigure}

Какой-нибудь диагональю $AC$ отсекаем от данного многоугольника треугольник $ABC$.
Через ту вершину $B$ этого треугольника, которая лежит против взятой диагонали, проводим прямую $MN\parallel AC$.
Затем продолжим одну из сторон $EA$ или $DC$, прилежащих к отсечённому треугольнику, до пересечения с прямой $MN$ (на рисунке продолжена сторона $EA$).
Точку пересечения $F$ соединим прямой с $C$.
Треугольники $CBA$ и $CFA$ равновелики, так как у них общее основание $AC$, а вершины $B$ и $F$ лежат на прямой, параллельной основанию.

Если от данного многоугольника отделим $\triangle CBA$ и вместо него приложим равновеликий ему $\triangle CFA$, то величина площади не изменится;
следовательно, данный многоугольник равновелик многоугольнику $FCDE$, у которого, очевидно, число углов на единицу меньше, чем у данного многоугольника.

Таким же приёмом можно число углов полученного многоугольника уменьшить ещё на единицу и продолжать такое последовательное уменьшение до тех пор, пока не получится треугольник ($FCG$ на нашем рисунке).

\paragraph{}\label{1938/255}
\so{Задача}.
\emph{Построить квадрат, равновеликий данному многоугольнику.}

Сначала преобразовывают многоугольник в равновеликий треугольник, а затем этот треугольник — в квадрат.
Пусть основание и высота треугольника $b$ и $h$, а сторона искомого квадрата $x$.
Тогда площадь первого равна $\tfrac12 b\cdot h$, а второго — $x^2$;
следовательно,
\begin{align*}
\tfrac12 b\cdot h&=x^2,
\end{align*}
откуда
$\frac12 b:{x}=x:h$,
то есть $x$ есть средняя пропорциональная между $\tfrac12 b$ и $h$.
Значит, сторону квадрата можно построить способом, указанным раньше (§~\ref{1938/190}) для нахождения средней пропорциональной.

{\small
\smallskip
\so{Замечание}.
Преобразование данного многоугольника в треугольник не всегда необходимо.
Например, если речь идёт о преобразовании в квадрат данной трапеции, то достаточно найти среднюю пропорциональную между высотой трапеции и её средней линией и на полученном отрезке построить квадрат.

}

{\small

\paragraph{}\label{1938/256}
\mbox{\so{Задача}.}
\emph{Вычислить площадь $S$ треугольника, зная длины $a$, $b$ и $c$ его сторон.}

\begin{wrapfigure}{r}{40mm}
\vskip-0mm
\centering
\includegraphics{mppics/ris-254}
\caption{}\label{1938/ris-254}
\end{wrapfigure}

Пусть высота $\triangle ABC$ (рис.~\ref{1938/ris-254}), опущенная на сторону $a$, есть $h_a$.
Тогда
\[S=\tfrac12ah_a.\]

Чтобы найти высоту $h_a$, возьмём равенство (§~\ref{1938/194}):
\[b^2 =  a^2+c^2 -2ac'\]
и определим из него отрезок $c'$.
\[c'=\frac{a^2+c^2-b^2}{2a}.\]

Из $\triangle ABD$ находим:
\begin{align*}
h_a&=\sqrt{c^2-(c')^2}=
\\
&=\sqrt{c^2-\left(\frac{a^2+c^2-b^2}{2a}\right)^2}=
\\
&=\frac{1}{2a}\sqrt{4a^2c^2-(a^2+c^2-b^2)^2}.
\end{align*}
Преобразуем подкоренное выражение так:
\begin{align*}
4a^2c^2&-(a^2+c^2-b^2)^2=
\\
&=[2ac+(a^2+c^2-b^2)]\cdot[2ac-(a^2+c^2-b^2)]=
\\
&=[(a^2+2ac+c^2)-b^2]\cdot[b^2-(a^2-2ac+c^2)]=
\\
&=[(a+c)^2-b^2]\cdot[b^2-(a-c)^2]=
\\
&=(a+c+b)(a+c-b)(b+a-c)(b-a+c).
\end{align*}
Следовательно,
\[S=\tfrac12 ah_a=\tfrac14\sqrt{(a+c+b)(a+c-b)(b+a-c)(b-a+c)}.\]
(Так как в треугольнике сумма любых двух сторон больше третьей, то все разности $a+b-c$, $a+c-b$ и $b+c-a$ — числа положительные.)

Если положим, что $a+b+c=2p$, то
\[a+c-b=(a+b+c)-2b=2p-2b=2(p-b).\]
Подобно этому
\begin{align*}
 b+a-c&=2(p-c);
 \\
 b+c-a&=2(p-a).
\end{align*}
Тогда
\begin{align*}
S&=\tfrac14\sqrt{2p\cdot2(p-a)\cdot2(p-b)\cdot2(p-c)},
\intertext{то есть}
S&=\sqrt{p(p-a)(p-b)(p-c)}.
\end{align*}
Это выражение известно под названием \rindex{формула Герона}\textbf{формулы Герона} (по имени математика Герона из Александрии, жившего приблизительно в III—II веках до нашей эры).

\smallskip
\so{Частный случай.}
Площадь равностороннего треугольника со стороной $a$ выражается следующей формулой:
\begin{align*}S&=\sqrt{\frac{3a}2\cdot \frac a2\cdot \frac a2\cdot \frac a2}=
 \\
 &=\frac{a^2\sqrt{3}}{4}.
\end{align*}

}

\subsection*{Теорема Пифагора и основанные на ней задачи}

\paragraph{}\label{1938/257}
\so{Теорема}.
\textbf{\emph{Сумма площадей квадратов, построенных на катетах прямоугольного треугольника, равна площади квадрата, построенного на гипотенузе этого треугольника.}}

Это предложение является другой формой теоремы Пифагора, доказанной ранее (§~\ref{1938/191}):
квадрат числа, измеряющего длину гипотенузы, равен сумме квадратов чисел, измеряющих длины катетов.
Действительно, квадрат числа, измеряющего длину отрезка, и является мерой площади квадрата, построенного на этом отрезке.
Поэтому теорема §~\ref{1938/191} равносильна указанной теореме Пифагора.

Приведём другое доказательство теоремы Пифагора, основанное не на вычислении площадей, а на непосредственном их сравнении между собой.

{

\begin{wrapfigure}{o}{47mm}
\vskip-6mm
\centering
\includegraphics{mppics/ris-255}
\caption{}\label{1938/ris-255}
\end{wrapfigure}

\smallskip
\mbox{\so{Доказательство}} (Евклида).
Пусть $ABC$ (рис.~\ref{1938/ris-255}) — прямоугольный треугольник, а $BDEA$, $AFGC$ и $BCKH$ — квадраты, построенные на его катетах и гипотенузе;
требуется доказать, что сумма площадей двух первых квадратов равна площади третьего квадрата.

Проведём $AM\perp BC$.
Тогда квадрат $BCKH$ разделится на два прямоугольника.
Докажем, что прямоугольник $BLMH$ равновелик квадрату $BDEA$, а прямоугольник $LCKM$ равновелик квадрату $AFGC$.

}

Проведём вспомогательные прямые $DC$ и $AH$.
Рассмотрим два треугольника, покрытые на рисунке штрихами.
Треугольник $BCD$, имеющий основание $BD$, общее с квадратом $BDEA$, а высоту $CN$, 
равную высоте $AB$ этого квадрата, равновелик половине квадрата.
Треугольник $ABH$, имеющий основание $BH$, общее с прямоугольником $BLMH$, и высоту $AP$, 
равную высоте $BL$ этого прямоугольника, равновелик половине его.
Сравнивая эти два треугольника между собой, находим, что у них $BD = BA$ и $BC=BH$ (как стороны квадрата);
сверх того $\angle DBC=\angle ABH$, так как каждый из этих углов состоит из общей части $ABC$ и прямого угла.
Значит, треугольники $ABH$ и $DBC$ равны.
Отсюда следует, что прямоугольник $BLMH$ равновелик квадрату $BDEA$.

Соединив $G$ с $B$ и $A$ с $K$, мы совершенно так же докажем, что прямоугольник $LCKM$ равновелик квадрату $AFGC$.
Отсюда следует, что квадрат $BCKH$ равновелик сумме квадратов $BDEA$ и $AFGC$.


\paragraph{}\label{1938/258}
\so{Задачи}.
1) \emph{Построить квадрат, площадь которого равна сумме площадей двух данных квадратов.}

Строим прямоугольный треугольник, у которого катетами были бы стороны данных квадратов.
Квадрат, построенный на гипотенузе этого треугольника, имеет площадь, равную сумме площадей данных квадратов.

2) \emph{Построить квадрат, площадь которого равна разности площадей двух данных квадратов.}

{

\begin{wrapfigure}{o}{50mm}
\vskip-0mm
\centering
\includegraphics{mppics/ris-256}
\caption{}\label{1938/ris-256}
\end{wrapfigure}

Строим прямоугольный треугольник, у которого гипотенузой была бы сторона большего из данных квадратов, а катетом — сторона меньшего квадрата.
Квадрат, построенный на другом катете этого треугольника, является искомым.

3) \emph{Построить квадрат, площадь которого относится к площади данного квадрата как $m:n$.}

}

На произвольной прямой (рис. \ref{1938/ris-256}) откладываем отрезки $AB=m$ и $BC\z=n$, и на $AC$ как на диаметре описываем полуокружность.
Из точки $B$ восстанавливаем перпендикуляр $BD$ до пересечения с окружностью.
Проведя хорды $AD$ и $DC$, получим прямоугольный треугольник, у которого (§~\ref{1938/192})
\[\frac{AD^2}{DC^2}=\frac{AB}{BC}=\frac mn.\]
На катете $DC$ этого треугольника отложим отрезок $DE$, равный стороне данного квадрата,%
\footnote{Если сторона данного квадрата больше $DC$, то точки $E$ и $E$ будут лежать на продолжениях катетов $DC$ и $DA$.}
и проведём $EF\parallel CA$.
Отрезок $DF$ есть сторона искомого квадрата, потому что
\[
\frac{DF}{DE}=\frac{AD}{DC},
\quad\text{откуда}\quad\left(\frac{DF}{DE}\right)^2=\left(\frac{AD}{DC}\right)^2;\]
следовательно,
\[\frac{DF^2}{DE^2}=\frac{AD^2}{DC^2}=\frac mn.\]

\subsection*{Отношение площадей подобных многоугольников}

\paragraph{}\label{1938/259}
\so{Теорема}.
\textbf{\emph{Площади двух треугольников, имеющих по равному углу, относятся как произведения сторон, заключающих эти углы.}}

Пусть в треугольниках $ABC$ и $A_1B_1C_1$ (рис.~\ref{1938/ris-257}) углы $A$ и $A_1$ равны.

\begin{figure}[!ht]
\centering
\includegraphics{mppics/ris-257}
\caption{}\label{1938/ris-257}
\end{figure}

Проведя высоты $BD$ и $B_1D_1$, будем иметь:
\[\frac{\text{площадь}~ABC}{\text{площадь}~A_1B_1C_1}=\frac{AC\cdot  BD}{A_1C_1\cdot  B_1D_1}=\frac{AC}{A_1C_1}\cdot\frac{BD}{B_1D_1}.\]
Треугольники $ABD$ и $A_1B_1D_1$ подобны ($\angle A = \angle A_1$ и $\angle D\z=\angle D_1$), поэтому 
\[\frac{BD}{B_1D_1}=\frac{AB}{A_1B_1};\]
заменив первое отношение вторым, получим:
\[\frac{\text{площадь}~ABC}{\text{площадь}~A_1B_1C_1}=\frac{AC}{A_1C_1}\cdot\frac{BD}{B_1D_1}=\frac{AC}{A_1C_1}\cdot\frac{AB}{A_1B_1}.\]

{\sloppy

\paragraph{}\label{1938/260}
\so{Теорема}.
\textbf{\emph{Площади подобных треугольников или многоугольников относятся как квадраты соответственных сторон.}}

}

1) Если $ABC$ и $A_1B_1C_1$ — два подобных треугольника, то углы одного равны соответственно углам другого;
пусть $\angle A \z= \angle A_1$, $\angle B\z=\angle B_1$ и $\angle C = \angle C_1$.
Применим к ним предыдущую теорему:
\[\frac{\text{площадь}~ABC}{\text{площадь}~A_1B_1C_1}=
 \frac{AB\cdot  AC}{A_1B_1\cdot  A_1C_1}=
 \frac{AB}{A_1B_1}\cdot
 \frac{AC}{A_1C_1}.
 \eqno(1)
\]
Но из подобия треугольников следует:
\[\frac{AB}{A_1B_1}=\frac{AC}{A_1C_1}=\frac{BC}{B_1C_1} \eqno(2)\]
Поэтому в равенстве (1) мы можем каждое из отношений $\frac{AB}{A_1B_1}$ и $\frac{AC}{A_1C_1}$ заменить любым отношением ряда (2);
следовательно:
\begin{align*}
\frac{\text{площадь}~ABC}{\text{площадь}~A_1B_1C_1}&=
\left(\frac{AB}{A_1B_1}\right)^2=
\left(\frac{AC}{A_1C_1}\right)^2=
\left(\frac{BC}{B_1C_1}\right)^2=
\\
&=\frac{AB^2}{A_1B_1^2}=
\frac{AC^2}{A_1C_1^2}=
\frac{BC^2}{B_1C_1^2}.
\end{align*}

2) Если $ABCDE$ и $A_1B_1C_1D_1E_1$ (рис.~\ref{1938/ris-258}) — два подобных многоугольника, то их можно, как мы видели (§~\ref{1938/171}), разложить на одинаковое число подобных и одинаково расположенных треугольников.

\begin{figure}[!ht]
\centering
\includegraphics{mppics/ris-258}
\caption{}\label{1938/ris-258}
\end{figure}

Пусть эти треугольники будут:
$AOB$ и $A_1O_1B_1$, $BOC$ и $B_1O_1C_1$ и~т.~д.
Согласно доказанному в первой части этой теоремы, мы получим пропорции:
\begin{align*}
\frac{\text{площадь}~AOB}{\text{площадь}~A_1O_1B_1}&=\left(\frac{AB}{A_1B_1}\right)^2;
\\
\frac{\text{площадь}~BOC}{\text{площадь}~B_1O_1C_1}&=\left(\frac{BC}{B_1C_1}\right)^2\quad\text{и так далее}
\end{align*}
Но из подобия многоугольников следует:
\[\frac{AB}{A_1B_1}=\frac{BC}{B_1C_1}=\frac{CD}{C_1D_1}=\dots\]
и потому
\[\left(\frac{AB}{A_1B_1}\right)^2=\left(\frac{BC}{B_1C_1}\right)^2=\left(\frac{CD}{C_1D_1}\right)^2=\dots\]

Значит,
\begin{align*}
\frac{\text{площадь}~AOB}{\text{площадь}~A_1O_1B_1}
&=
\frac{\text{площадь}~BOC}{\text{площадь}~B_1O_1C_1}
=
\\
&=
\frac{\text{площадь}~COD}{\text{площадь}~C_1O_1D_1}=
\\
&=\dots,
\end{align*}
откуда
\begin{align*}
&\frac{\text{площадь}~AOB+\text{площадь}~BOC+\dots+\text{площадь}~EOA}{\text{площадь}\,A_1O_1B_1+\text{площадь}\,B_1O_1C_1+\dots+\text{площадь}\,E_1O_1A_1}
=
\\
&\quad=
\frac{\text{площадь}~ABCDE}{\text{площадь}~A_1B_1C_1D_1E_1}
=\frac{AB^2}{A_1B_1^2}.
\end{align*}

\smallskip
\so{Следствие}.
\emph{Площади правильных  $n$-угольников относятся как квадраты сторон, или квадраты радиусов описанных окружностей, или квадраты апофем.}

\paragraph{}\label{extra/pifagor} Покажем, что теорему Пифагора (§~\ref{1938/191}) можно получить как следствие из теоремы в §~\ref{1938/260}.

\begin{wrapfigure}{r}{47mm}
\vskip-0mm
\centering
\includegraphics{mppics/ris-extra-9}
\caption{}\label{extra/ris-9}
\end{wrapfigure} 

Пусть $ABC$ (рис.~\ref{extra/ris-9}) есть прямоугольный треугольник, $AD$ — перпендикуляр, опущенный на гипотенузу из вершины прямого угла.
Треугольники $ABC$, $DBA$ подобны, так как имеют по прямому углу и ещё общий угол $B$.
Значит 
\[\frac{\text{площадь}~DBA}{\text{площадь}~ABC}=
 \frac{c^2}{a^2}.
\]

Точно также получаем, что
\[\frac{\text{площадь}~DAC}{\text{площадь}~ABC}=
 \frac{b^2}{a^2}.
\]

Поскольку треугольники $DBA$ и $DAC$ вместе составляют $\triangle ABC$, сумма площадей
$\triangle DBA$ и $\triangle DAC$ равна площади $\triangle ABC$.
Отсюда
\[\frac{b^2+c^2}{a^2}=1
\quad\text{или}\quad 
b^2+c^2=a^2.\]

\paragraph{}\label{1938/261}
\so{Задача}.
\emph{Разделить данный треугольник на $m$ равновеликих частей прямыми, параллельными его стороне.}

Пусть, например, требуется разделить $\triangle ABC$ (рис.~\ref{1938/ris-259}) на три равновеликие части отрезками, параллельными основанию $AC$.

Предположим, что искомые отрезки будут $DE$ и $FG$.
Очевидно, что если мы найдём отрезки $BE$ и $BG$, то определятся и отрезки $DE$ и $FG$.
Треугольники $BDE$, $BFG$ и $BAC$ подобны;
поэтому
\begin{align*}
\frac{\text{площадь}~BDE}{\text{площадь}~BAC}&=\frac{BE^2}{BC^2}
&&\text{и}
&\frac{\text{площадь}~BFG}{\text{площадь}~BAC}&=\frac{BG^2}{BC^2}.
\intertext{Но}
\frac{\text{площадь}~BDE}{\text{площадь}~BAC}&=\frac{1}{3}
&&\text{и}
&\frac{\text{площадь}~BFG}{\text{площадь}~BAC}&=\frac{2}{3}.
\intertext{Следовательно,}
\frac{BE^2}{BC^2}&=\frac13
&&\text{и}
&\frac{BG^2}{BC^2}&=\frac23,
\end{align*}
откуда
\begin{align*}
BE&=\sqrt{\tfrac13BC^2}=\sqrt{\tfrac13BC\cdot BC}
\intertext{и}
BG&=\sqrt{\tfrac23BC^2}=\sqrt{\tfrac23BC\cdot BC}
\end{align*}
Из этих выражений видим, что $BE$ есть средняя пропорциональная между $BC$ и $\tfrac13BC$, а $BG$ есть средняя пропорциональная между $BC$ и $\tfrac23BC$.

\begin{wrapfigure}{o}{45mm}
\centering
\includegraphics{mppics/ris-259}
\caption{}\label{1938/ris-259}
\end{wrapfigure}

Поэтому построение можно выполнить так:
разделим $BC$ на три равные части в точках $a$ и $b$;
опишем на $BC$ полуокружность;
из $a$ и $b$ восстановим к $BC$ перпендикуляры $aH$ и $bK$.
Хорды $HB$ и $KB$ будут искомыми средними пропорциональными;
первая — между всем диаметром $BC$ и его третьей частью $Ba$, вторая — между $BC$ и В$b$, то есть между $BC$ и $\tfrac23BC$.
Остаётся отложить эти хорды на $BC$ от точки $B$, тогда получим искомые точки $E$ и~$G$.

Подобным образом можно разделить треугольник на какое угодно число равновеликих частей.

%% file: 2D/ploschad-kruga.tex
\section{Площадь круга и его частей}

\paragraph{}\label{1938/262}
\so{Лемма}.
\textbf{\emph{При неограниченном увеличении $\bm{n}$, сторона правильного вписанного $\bm{n}$-угольника делается как угодно малой.}}

Пусть $p_n$ есть периметр правильного вписанного $n$-угольника;
тогда длина одной его стороны выразится дробью $\frac {p_n}n$.
При неограниченном увеличении $n$ знаменатель этой дроби будет, очевидно, возрастать неограниченно, а числитель, то есть $p_n$, хотя и может возрастать, но не беспредельно,
так как периметр всякого вписанного выпуклого многоугольника всегда остаётся меньшим периметра любого описанного многоугольника.

Если же в какой-нибудь дроби знаменатель неограниченно возрастает, а числитель остаётся меньше некоторой постоянной величины, то дробь эта может сделаться как угодно малой.
Значит, то же самое можно сказать о стороне правильного вписанного $n$-угольника:
при неограниченном $n$ она может сделаться как угодно малой.

\begin{wrapfigure}[14]{r}{36mm}
\vskip-0mm
\centering
\includegraphics{mppics/ris-260}
\caption{}\label{1938/ris-260}
\end{wrapfigure}

\paragraph{}\label{1938/263}
\mbox{\so{Следствие}.}
Пусть $AB$ (рис. \ref{1938/ris-260}) есть сторона правильного вписанного многоугольника, $OA$ — радиус и $OC$ — апофема.
Из $\triangle OAC$ находим (§~\ref{1938/50}):
\[OA-OC<AC,\]
то есть
\[AO-OC<\tfrac12 AB.\]
Но при увеличении $n$, сторона правильного вписанного $n$-угольника, как мы сейчас доказали, может сделаться как угодно малой;
значит, то же самое можно сказать и о разности $OA-OC$.

Таким образом, \emph{при неограниченном увеличении $n$, разность между радиусом и апофемой правильного вписанного $n$-угольника может сделаться как угодно малой}.
Это же можно высказать другими словами так: 
\emph{предел, к которому стремится апофема $a_n$ правильного вписанного $n$-угольника, есть радиус.}

\paragraph{}\label{1914/230}
\so{Лемма}.
\textbf{\emph{Разность между площадью правильного $\bm{n}$-угольника, описанного около круга, и площадью правильного $\bm{n}$-угольника, вписанного в тот же
круг, стремится к нулю при неограниченном увеличении $\bm{n}$.}}

Впишем в окружность (рис. \ref{1914/ris-292}) и опишем около неё по правильному $n$-угольнику (на чертеже изображены $6$-угольники).
Пусть $R$ будет радиус круга, $a_n$ — апофема вписанного $n$-угольника, $q_n$ его площадь и $Q_n$ — площадь описанного $n$-угольника.
Тогда (§~\ref{1938/260})
\[\frac {Q_n}{q_n}=\frac{R^2}{a_n^2}\]

\begin{wrapfigure}[12]{o}{30mm}
\vskip-3mm
\centering
\includegraphics{mppics/ris-1914-292}
\caption{}\label{1914/ris-292}
\end{wrapfigure}
\noindent
и следовательно 
\[\frac {Q_n-q_n}{q_n}=\frac{R^2-a_n^2}{a_n^2}.\]
Откуда
\[(Q_n-q_n)a_n^2=q_n(R^2-a_n^2)\]
или
\[(Q_n-q_n)a_n^2=q_n(R+a_n)(R-a_n).\]

При неограниченном увеличении $n$, разность $R-a_n$ стремится к нулю (§~\ref{1938/263}).
Сомножитель $q_n$ остаётся ограниченным, так как площадь вписанного многоугольника меньше площади любого описанного.  
Сомножитель $R+a_n$ всегда остаётся меньше $R+R$.
Вследствие этого правая часть последнего равенства (а значит и левая его часть),
при неограниченном увеличении $n$, стремится к нулю.

Заметим, что сомножитель $a_n^2$ в левой части возрастает при увеличении $n$ (§~\ref{1938/110}).
Значит левая часть может стремиться к нулю только тогда, когда другой сомножитель $Q_n-q_n$ стремится к нулю.
То есть разница площадей описанного и вписанного $n$-угольников стремится к нулю.

{\small
\paragraph{}\label{1914/231}
\so{Замечание}.
Таким же путём мы можем доказать, что \emph{разность между периметром описанного и
периметром вписанного правильного $n$-угольника, при неограниченном увеличении $n$, стремится к нулю.}

Действительно, если периметры правильных $n$-угольников, описанного и вписанного, обозначены буквами $P_n$ и $p_n$, то (Следствие в §~\ref{1938/218}):
\[\frac {P_n}{p_n}=\frac R{a_n}.\]
Откуда
\[\frac{P_n-p_n}{p_n}=\frac{R-a_n}{a_n},\]
и следовательно
\[(P_n-p_n)a_n=(R-a_n)p_n.\]

При неограниченном увеличении $n$ правая часть последнего равенства (а следовательно, и его левая часть) стремится к нулю, так как множитель $R-a_n$, по доказанному, стремится к нулю, а множитель $p_n$ остаётся ограниченным (так как периметр вписанного многоугольника не превосходит периметр любого описанного).
Но левая часть равенства, представляя собою произведение, в котором сомножитель $a_n$ увеличивается с ростом $n$, может стремиться к нулю только тогда, когда его другой сомножитель  стремится к нулю;
а этот сомножитель и есть разность периметров $P_n$ и $p_n$.

}

\paragraph{}\label{1914/232} \so{Теорема}.
\textbf{\emph{Площадь круга есть oбщий предел площадей правильных вписанных в этот круг и описанных около него $\bm{n}$-угольников при неограниченном возрастании $\bm{n}$.}}

Пусть около круга, площадь которого мы обозначим $K$, описан правильный $n$-угольник и в него вписан правильный $n$-угольник (рис. \ref{1914/ris-292}).
Обозначим площадь первого $Q_n$, а площадь второго $q_n$.
При увеличении $n$, величины $Q_n$ и $q_n$ меняются, тогда как величина $K$ остаётся неизменной.
Требуется доказать, что при неограниченном увеличении $n$ величины $Q_n$ и $q_n$ стремятся к одному
и тому же пределу, а именно к площади $K$.

Заметим, что для любого $n$, выполняется неравенство 
\[Q_n>K>q_n\]
и потому каждая из двух разностей: 
$Q_n-K$ и $K-q_n$ всегда меньше разности $Q_n-q_n$.
Но при неограниченном увеличении $n$ разность $Q_n-q_n$ стремится к нулю (§~\ref{1914/230}).
Следовательно, при этом каждая из меньших разностей: $Q_n-K$ и $K-q_n$ и подавно стремится к нулю;
а это, согласно определению предела, означает, что $K$ является общим пределом для $Q_n$ и $q_n$.

{\small
\smallskip
\paragraph{}\label{1914/233} \so{Замечание}.
Можно также утверждать, что \emph{длина окружности есть общий предел периметров правильных вписанных в эту окружность и описанных около неё $n$-угольников при неограниченном увеличении $n$}.
Действительно, из того обстоятельства, что разность $P_n-p_n$ стремится к нулю (§~\ref{1914/231}), надо заключить, что периметры $P_n$ и $p_n$ могут стремиться только к одному и тому же пределу;
но предел $p_n$ есть то, что принимается за длину окружности;
значит, и предел $P_n$ равен длине окружности.
}

\paragraph{Площадь круга.}\label{1938/264}\ 
Впишем в круг, радиус которого обозначим $R$, какой-нибудь правильный $n$-угольник.
Пусть
\begin{align*}
\text{площадь}\quad&\text{этого многоугольника будет}&&q_n,
\\
\text{периметр}\quad&\text{—\textquotedbl—\qquad\quad —\textquotedbl—\qquad\quad —\textquotedbl—}&&p_n,
\\
\text{апофема}\quad&\text{—\textquotedbl—\qquad\quad —\textquotedbl—\qquad\quad —\textquotedbl—}&&a_n.
\end{align*}

Мы видели (§~\ref{1938/252}, следствие), что между этими величинами есть такая зависимость.
\[q_n=\tfrac12p_n\cdot a_n.\]

Вообразим теперь, что $n$ неограниченно увеличивается.
Тогда периметр $p_n$ и апофема  $a_n$ (следовательно, и площадь $q_n$) будут меняться, причём периметр будет стремиться к пределу, принимаемому за длину окружности $C$, апофема будет стремиться к пределу, равному радиусу круга $R$.
Из этого следует, что площадь $n$-угольника будет стремиться к пределу, равному $\tfrac12 C\cdot R$.
Но согласно §~\ref{1914/232} этот предел и есть площадь круга.

{\sloppy
Таким образом, обозначив площадь круга буквой $K$, можем написать:
\[K=\tfrac12 C\cdot R\]
то есть \emph{\textbf{площадь круга равна половине произведения длины окружности на радиус.}}

}

Так как $C=2\pi R$, то
\[K=\tfrac12\cdot 2\pi  R\cdot R=\pi R^2,\]
то есть \emph{\textbf{площадь круга равна квадрату радиуса, умноженному на отношение длины окружности к диаметру.}}

\paragraph{}\label{1938/265}
\so{Следствие}.
\emph{Площади кругов относятся, как квадраты радиусов или диаметров.}

Действительно, если $K$ и $K_1$ будут площади двух кругов, а $R$ и $R_1$ — их радиусы, то
\[K=\pi R^2\]
и
\[K_1=\pi R_1^2,\]
откуда
\begin{align*}
\frac{K}{K_1}&=\frac{\pi R^2}{\pi R^2_1}=\frac{R^2}{R^2_1}=
\\
&=\frac{4R^2}{4R^2_1}=\frac{(2R)^2}{(2R_1)^2}.
\end{align*}

\paragraph{}\label{1938/266}
\so{Задача 1}.
\emph{Вычислить площадь круга, длина окружности которого равна 2 м.}

Для этого предварительно находим радиус $R$ из равенства.
\[2\pi R= 2.\]
откуда
\[R=\tfrac1\pi=0{,}3183\dots\]

Затем определим площадь круга.
\[K=\pi R^2=\pi(\tfrac1\pi)^2=\tfrac1\pi=0{,}3183\dots\text{м}^2.\]

\paragraph{}\label{1938/267}
\so{Задача 2}.
\emph{Построить квадрат, равновеликий данному кругу.}

Эта задача, известная под названием \so{квадратуры круга}, не может быть решена при помощи циркуля и линейки.
Действительно, если обозначим буквой $x$ сторону искомого квадрата, а буквой $R$ радиус круга, то получим уравнение:
\[x^2=\pi R^2.\]
откуда
\[\frac{\pi R}{x}=\frac{x}{R}\]
то есть
$x$ есть средняя пропорциональная между полуокружностью и радиусом.
Следовательно, если известен отрезок, длина которого равна длине полуокружности, то легко построить квадрат, равновеликий данному кругу, и обратно, если известна сторона квадрата, равновеликого кругу, то можно построить отрезок, равный по длине полуокружности.
Но с помощью циркуля и линейки нельзя построить отрезок, длина которого равнялась бы длине полуокружности (§~\ref{1938/238});
следовательно, нельзя в точности решить задачу о построении квадрата, равновеликого кругу.
Приближённое решение можно выполнить, если предварительно найти приближённую длину полуокружности и затем построить среднюю пропорциональную между отрезком этой длины и радиусом.

\paragraph{}\label{1938/268}
\mbox{\so{Теорема}.}
\textbf{\emph{Площадь сектора равна произведению длины его дуги на половину радиуса.}}

Пусть дуга $AB$ (рис.~\ref{1938/ris-261}) сектора $AOB$ содержит $n\degree$.
Очевидно, что площадь сектора, дуга которого содержит $1\degree$, составляет $\tfrac1{360}$ часть площади круга, то есть она равна $\frac{\pi R^2}{360}$.
Следовательно, площадь $S$ сектора, дуга которого содержит $n\degree$, равна:
\[S=\frac{\pi R^2}{360}\cdot n\degree=\frac{\pi R\cdot  n\degree}{180}\cdot \frac R2.\]

{

\begin{wrapfigure}{o}{35mm}
\vskip-4mm
\centering
\includegraphics{mppics/ris-261}
\caption{}\label{1938/ris-261}
\end{wrapfigure}

\noindent
Так как дробь $\frac{\pi R n\degree}{180}$ выражает длину дуги $AB$ (§~\ref{1938/239}), то, обозначив
её буквой $s$, получим:
\[S=s\cdot \frac R2.\]

Заметим также, что если $\alpha$ угловая мера угла $AOB$, выраженная в радианах, то $s=\alpha\cdot R$ (§~\ref{extra/radians}) и значит
\[S=\frac{\alpha\cdot R^2}2.\]

}

{\small
\medskip
\mbox{\so{Замечание}.}
Строго говоря, предложенное рассуждение доказывает формулу
\[S=s\cdot \frac R2\]
только при целом значении~$n$.
Аналогичным способом можно установить равенства в случае если угол $AOB$ имеет общую меру с полным углом.

Верность этого равенства для произвольных секторов доказывается аналогично лемме в §~\ref{1938/159};
то есть требуется предъявить произвольно близкие приближения для $S$ с недостатком и избытком между которыми лежит значение $s\cdot \frac R2$.


}

{\sloppy

\paragraph{Площадь сегмента.}\label{1938/269}
Для нахождения площади сегмента, ограниченного дугой $s$ и хордой $AB$ (рис.~\ref{1938/ris-261}), надо отдельно вычислить площадь сектора $AOB$ и площадь треугольника $AOB$ и из первой вычесть вторую.

}

Впрочем, когда градусное измерение дуги $s$ невелико, площадь сегмента можно вычислять по следующей \so{приближённой} формуле (мы её приводим без доказательства):
\[\text{площадь сегмента}~\approx\tfrac23bh.\eqno(1)\]
где $b$ есть основание сегмента (рис.~\ref{1938/ris-262}), а $h$ — его высота (обыкновенно называемая \rindex{стрелка сегмента}\textbf{стрелкой сегмента}). 
Иначе говоря, площадь сегмента на рис.~\ref{1938/ris-262} близка к площади прямоугольника на том же рисунке. 

\begin{figure}[h]
\centering
\includegraphics{mppics/ris-262}
\caption{}\label{1938/ris-262}
\end{figure}

Доказано, что погрешность результата вычисления, получаемого по этой приближённой формуле, тем меньше, чем меньше отношение $\tfrac hb$;
так, если $h<\tfrac19b$ (что бывает тогда, когда дуга $s$ содержит меньше $50\degree$), то погрешность оказывается меньше 1\% площади.
Более точные результаты даёт более сложная формула.
\[\text{площадь сегмента}~\approx\tfrac23bh+\frac{h^3}{2b}.\eqno(2)\]

{\small

\subsection*{Упражнения}

\begin{center}
\so{Доказать теоремы}
\end{center}

\begin{enumerate}[noitemsep]
\item
В параллелограмме расстояния от любой точки диагонали до двух прилежащих сторон обратно пропорциональны этим сторонам.

\item
Площадь трапеции равна произведению одной из непараллельных сторон на перпендикуляр, опущенный из середины другой непараллельной стороны на первую.

\item
Два четырёхугольника равновелики, если у них равны соответственно диагонали и угол между ними.

\item
Если площади двух треугольников, прилежащих к основаниям трапеции и образуемых пересечением её диагоналей, равны соответственно $p^2$ и $q^2$, то площадь всей трапеции равна $(p+q)^2$.

\item
Площадь правильного вписанного шестиугольника равна $\tfrac34$ площади правильного описанного шестиугольника.

\item
В четырёхугольнике $ABCD$ через середину диагонали $BD$ проведена прямая, параллельная другой диагонали $AC$.
Эта прямая пересекает сторону $AD$ в точке $E$.
Доказать, что отрезок $CE$ делит четырёхугольник пополам.

\item
Если медианы треугольника взять за стороны другого треугольника, то площадь последнего равна площади первого.

\item
В круге с центром $O$ проведена хорда $AB$.
На радиусе $OA$, как на диаметре, описана окружность.
Доказать, что площади двух сегментов, отсекаемых хордой $AB$ от обоих кругов, относятся, как 4:1.

\item Точка $X$ внутри треугольника $ABC$ лежит на медиане $AM$ тогда и только тогда, когда площадь треугольника $ABX$ равна площади треугольника $ACX$.

\end{enumerate}

\begin{center}
\so{Задачи на вычисление}
\end{center}

\begin{enumerate}[resume,noitemsep]

\item
Вычислить площадь прямоугольной трапеции, у которой один из углов равен $60\degree$, зная или оба основания, или одно основание и высоту, или одно основание и боковую сторону, наклонную к основанию.

\item
Даны основания трапеции $B$ и $b$ и её высота Н.
Вычислить высоту треугольника, образованного продолжением непараллельных сторон трапеции до взаимного пересечения.

\item
В треугольник вписан другой треугольник, вершины которого делят пополам стороны первого треугольника;
в другой треугольник вписан подобным же образом третий;
в третий — четвёртый и~т.~д.
неограниченно.
Найти предел суммы площадей этих треугольников.

\item
По трём данным сторонам $a$, $b$ и $c$ треугольника вычислить радиус $r$ круга, вписанного в этот треугольник.

\smallskip
\so{Указание}.
Если $S$ есть площадь треугольника, то легко усмотреть, что
\[S=\tfrac12 ar+\tfrac12 br+\tfrac12 cr=pr,\]
где $p$ означает полупериметр треугольника.
С другой стороны, площадь $S$ выражается формулой Герона (§~\ref{1938/256}).
Отсюда можно получить формулу для $r$.

\item
Вычислить стрелку (высоту) и площадь сегмента в зависимости от радиуса $r$ круга, если центральный угол, соответствующий сегменту, содержит $60\degree$.
Вычисление это произвести тремя способами:
1) посредством вычитания из площади сектора площади треугольника;
2) по первой сокращённой формуле, указанной в §~\ref{1938/269};
3) по второй сокращённой формуле, указанной там же.
Сравните результаты вычисления друг с другом с целью определить абсолютную и относительную погрешности приближённых результатов.

\smallskip
\so{Решение}.
\[b=r;\]
\[h=r-\tfrac12 r\sqrt3=\tfrac12 r(2-\sqrt3)\approx 0{,}1340\cdot r;\]

1) площадь $p_1=\frac{\pi r^2}{6}-\frac{r^2\sqrt3}{4}\approx 0{,}0906\cdot r^2$;

2) площадь $p_2\approx \tfrac23 bh=\tfrac23\cdot r\cdot 0{,}1340\cdot r=0{,}0893\cdot r^2$;

3) площадь $p_3\approx \tfrac23 bh+\frac{h^3}{2b}=0{,}0893r^2+0{,}0012r^2=0{,}0905\cdot r^2$.

Абсолютная погрешность.
\[\text{для площади}~p_2 \approx 0{,}0906\cdot r^2 - 0{,}0893\cdot r^2 = 0{,}0013\cdot r^2;\]
\[\text{для площади}~p_3 \approx 0{,}0906\cdot r^2 - 0{,}0905\cdot r^2 = 0{,}0001\cdot r^2.\]
Относительная погрешность (то есть отношение абсолютной погрешности к измеряемой величине):
\[\text{для площади}~p_2 = \frac{p_1-p_2}{p_1}\approx\frac{0{,}0013r^2}{0{,}0906r^2}\approx0{,}014 = 1{,}4\%;\]
\[\text{для площади}~p_3 = \frac{p_1-p_3}{p_1}\approx\frac{0{,}0001r^2}{0{,}0906r^2}\approx0{,}001 = 0{,}1\%;\]
Таким образом, результат, вычисленный по первой приближённой формуле, меньше истинного результата (приблизительно) на $1{,}4\%$, а результат, вычисленный по второй приближённой формуле, меньше истинного на $0{,}1\%$.

\item
1) Зная основание $b$ сегмента и высоту его (стрелку) $k$, вычислить радиус $r$ круга.

\smallskip
\so{Указание}.
Из прямоугольного треугольника, у которого гипотенуза есть $r$, один катет $\frac b2$, а другой $r-h$, находим уравнение:
\[(\tfrac b2)^2+(r-h)^2=r^2\]
из которого легко определить $r$.

2) Вычислить диаметр круга, если известно, что при основании сегмента, равном $67{,}2$ см, стрелка его есть $12{,}8$ см (смотри предыдущее указание).

\smallskip
\so{Задачи на построение}

\item
Разделить треугольники прямыми, проходящими через его вершину, на три части, площади которых относятся, как $m\z:n:p$.

\item
Разделить треугольник на две равновеликие части прямой, проходящей через данную точку его стороны.

\item
Найти внутри треугольника такую точку, чтобы прямые, соединяющие её с вершинами треугольника, делили его на три равновеликие части.

\smallskip
\so{Указание}.
Разделим сторону $AC$ на $3$ равные части в точках $E$ и $F$.
Через $E$ проводим прямую, параллельную $AB$, и через $F$ — прямую, параллельную $BC$.
Точка пересечения этих прямых — искомая.

\item
То же — на три части в отношении $2:3:4$ (или вообще $m:n\z:p$).

\item
Разделить параллелограмм на три равновеликие части прямыми, исходящими из вершины его.

\item
Разделить параллелограмм на две части прямой, проходящей через данную точку на продолжении его стороны, так, чтобы их площади относились, как $m:n$. 

\smallskip
\so{Указание}.
Среднюю линию параллелограмма разделить в отношении $m:n$ и точку деления соединить прямой с данной точкой.  

\item
Разделить параллелограмм на три равновеликие части прямыми, параллельными диагонали.

\item
Разделить площадь треугольника в среднем и крайнем отношении прямой, параллельной основанию.

\smallskip
\so{Указание}.
Решается приложением алгебры к геометрии.

\item
Разделить треугольник на три равновеликие части прямыми, перпендикулярными к основанию.

\item
Разделить круг на $2$, на $3,\dots$ равновеликие части концентрическими окружностями.

\item
Разделить трапецию на две равновеликие части прямой, параллельной основаниям.

\smallskip
\so{Указание}.
Продолжив непараллельные стороны до взаимного пересечения, взять за неизвестную величину расстояние от конца искомой линии до вершины треугольника;
составить пропорции, исходя из площадей подобных треугольников.

\item
Данный треугольник преобразовать в другой равновеликий прямоугольник с данным основанием.

\item
Построить квадрат, равновеликий $\tfrac23$ данного квадрата.

\item
Преобразовать квадрат в равновеликий прямоугольник, у которого сумма или разность $d$ двух смежных сторон дана.

\item
Построить круг, равновеликий кольцу, заключённому между двумя данными концентрическими окружностями.

\item
Построить треугольник, подобный одному и равновеликий другому из двух данных треугольников.

\item
Данный треугольник преобразовать в равновеликий равносторонний (посредством приложения алгебры к геометрии).

\item
В данный круг вписать прямоугольник с данной площадью $m$ (посредством приложения алгебры к геометрии).

\item
В данный треугольник вписать прямоугольник с данной площадью $m$ (приложением алгебры к геометрии; исследовать).

\end{enumerate}

}

%% file: 3D/pryamye_i_ploskosti.tex
\chapter{Прямые и плоскости}

\section{Определение положения плоскости}

\begin{wrapfigure}{r}{34 mm}
\centering
\includegraphics{mppics/s-ris-1}
\caption{}\label{1938/s-ris-1}
\end{wrapfigure}

\paragraph{Изображение плоскости.}\label{1938/s2}
В обыденной жизни многие предметы, поверхность которых напоминает геометрическую плоскость, имеют форму прямоугольника: переплёт книги, оконное стекло, поверхность письменного стола и тому подобное.
При этом если смотреть на эти предметы под углом и с большего расстояния, то они представляются нам имеющими форму параллелограмма.
Поэтому принято изображать плоскость на рисунке в виде параллелограмма.

Плоскость обычно обозначают одной буквой, например «плоскость $M$» (рис.~\ref{1938/s-ris-1}).

\paragraph{Основные свойства плоскости.}\label{1938/s3}
Укажем следующие свойства плоскости, которые принимаются без доказательства, то есть являются аксиомами:

1) \emph{Если две точки прямой принадлежат плоскости, то и каждая точка этой прямой принадлежит плоскости.}

2) \emph{Если две плоскости имеют общую точку, то они пересекаются по прямой, проходящей через эту точку.}

3) \emph{Через всякие три точки, не лежащие на одной прямой, можно провести плоскость, и притом только одну.}

\paragraph{} \label{1938/s4}
\so{Следствия}. Из последнего предложения можно вывести следствия:

1) \emph{Через прямую и точку вне её можно провести плоскость (и только одну).} Действительно, точка вне прямой вместе с какими-нибудь двумя точками этой прямой составляют три точки, через которые можно провести плоскость (и притом одну).

2) \emph{Через две пересекающиеся прямые можно провести плоскость (и только одну).} Действительно, взяв точку пересечения и ещё по одной точке на каждой прямой, мы будем иметь три точки, через которые можно провести плоскость (и притом одну).

3) \emph{Через две параллельные прямые можно провести только одну плоскость.} Действительно, параллельные прямые, по определению, лежат в одной плоскости;
эта плоскость единственная, так как через одну из параллельных и какую-нибудь точку другой можно провести не более одной плоскости.

\paragraph{Вращение плоскости вокруг прямой.}\label{1938/s5} 
\textbf{\emph{Через каждую прямую в пространстве можно провести бесчисленное множество плоскостей.}}

\begin{wrapfigure}{o}{31 mm}
\centering
\includegraphics{mppics/s-ris-2}
\caption{}\label{1938/s-ris-2}
\end{wrapfigure}

В самом деле, пусть дана прямая $a$ (рис.~\ref{1938/s-ris-2}).
Возьмём какую-нибудь точку $A$ вне её.
Через точку $A$ и прямую $a$ проходит единственная плоскость (§~\ref{1938/s4}).
Назовём её плоскостью $M$.
Возьмём новую точку $B$ вне плоскости $M$.
Через точку $B$ и прямую $a$ в свою очередь проходит плоскость.
Назовём её плоскостью $N$.
Она не может совпадать с $M$, так как в ней лежит точка $B$, которая не принадлежит плоскости $M$.
Мы можем далее взять в пространстве ещё новую точку $C$ вне плоскостей $M$ и $N$.
Через точку $C$ и прямую $a$ проходит новая плоскость.
Назовём её~$P$.
Она не совпадает ни с $M$, ни с $N$, так как в ней находится точка $C$, не принадлежащая ни плоскости $M$, ни плоскости $N$.
Продолжая брать в пространстве всё новые и новые точки, мы будем таким путём получать всё новые и новые плоскости, проходящие через данную прямую $a$.
Таких плоскостей будет бесчисленное множество.

Все эти плоскости можно рассматривать как различные положения одной и той же плоскости, которая вращается вокруг прямой $a$.

\paragraph{Задачи на построение в пространстве.}\label{1938/s6}
Все построения, которые делались в планиметрии, выполнялись в одной плоскости при помощи чертёжных инструментов.
Для построений в пространстве чертёжные инструменты становятся уже непригодными, так как чертить фигуры в пространстве невозможно.
Кроме того, при построениях в пространстве появляется ещё новый элемент — \rindex{плоскость}\textbf{плоскость}, построение которой в пространстве нельзя выполнять столь простыми средствами, как построение прямой на плоскости.

Поэтому при построениях в пространстве необходимо точно {определить, что значит выполнить то или иное построение} и, в частности, что значит построить плоскость в пространстве.
Во всех построениях в пространстве мы будем предполагать, что:

1) плоскость может быть построена, если найдены элементы, определяющие её положение в пространстве (§§~\ref{1938/s3} и \ref{1938/s4}), то есть что мы умеем построить плоскость, проходящую через три данные точки, через прямую и точку вне её, через две пересекающиеся или две параллельные прямые;

2) если даны две пересекающиеся плоскости, то дана и линия их пересечения, то есть мы умеем найти линию пересечения двух плоскостей;

3) если в пространстве дана плоскость, то мы можем выполнять в ней все построения, которые выполнялись в планиметрии.

\emph{Выполнить какое-либо построение в пространстве — это значит свести его к конечному числу только что указанных основных построений.}
При помощи этих основных задач можно решать и задачи более сложные.

В этих предложениях и решаются задачи на построение в стереометрии.

\paragraph{Пример задачи на построение в пространстве.}\label{1938/s7}\ 

\so{Задача}.
\emph{Найти точку пересечения данной прямой $a$ \emph{(рис.~\ref{1938/s-ris-3})} с данной плоскостью~$P$.}

\begin{wrapfigure}{o}{55 mm}
\centering
\includegraphics{mppics/s-ris-3}
\caption{}\label{1938/s-ris-3}
\end{wrapfigure}

Возьмём на плоскости $P$ какую-либо точку $A$.
Через точку $A$ и прямую $a$ проводим плоскость $Q$.
Она пересекает плоскость $P$ по некоторой прямой $b$.
В плоскости $Q$ находим точку $C$ пересечения прямых $a$ и $b$.
Эта точка и будет искомой.

Если прямые $a$ и $b$ окажутся параллельными, то задача не будет иметь решения.

\section{Параллельные прямые и плоскости}

\subsection*{Параллельные прямые}

\paragraph{Предварительное замечание.}\label{1938/s8}
Две прямые могут быть расположены в пространстве так, что через них нельзя провести плоскость.
Возьмём, например (рис.~\ref{1938/s-ris-4}), две такие прямые $AB$ и $DE$, из которых одна пересекает некоторую плоскость $P$, а другая лежит на ней, но не проходит через точку ($C$) пересечения первой прямой и плоскости~$P$.
Через такие две прямые нельзя провести плоскость, потому что в противном случае через прямую $DE$ и точку $C$ проходили бы две различные плоскости:
одна $P$, пересекающая прямую $AB$, и другая, содержащая её, а это невозможно (§~\ref{1938/s3}).

\begin{wrapfigure}{o}{44 mm}
\centering
\includegraphics{mppics/s-ris-4}
\caption{}\label{1938/s-ris-4}
\end{wrapfigure}

Две прямые, не лежащие в одной плоскости, конечно, не пересекаются, сколько бы их ни продолжали;
однако их не называют параллельными, оставляя это название для таких прямых, которые, \so{находясь в одной плоскости}, не пересекаются, сколько бы их ни продолжали.

Две прямые, не лежащие в одной плоскости, называются \rindex{скрещивающиеся прямые}\textbf{скрещивающимися}.

\subsection*{Прямая и плоскость, параллельные между собой}

\paragraph{}\label{1938/s9}
\so{Определение}. Плоскость и прямая, не лежащая в этой плоскости, называются \rindex{параллельность}\textbf{параллельными}, если они не пересекаются, сколько бы их ни продолжали.
При этом плоскость считается параллельной любой прямой в ней лежащей.

\begin{wrapfigure}{r}{55 mm}
\vskip-6mm
\centering
\includegraphics{mppics/s-ris-5}
\caption{}\label{1938/s-ris-5}
\end{wrapfigure}

\paragraph{}\label{1938/s10}
\mbox{\so{Теорема}.} \textbf{\emph{Если прямая}} ($AB$, рис.~\ref{1938/s-ris-5}) \textbf{\emph{параллельна какой-нибудь прямой}} ($CD$), \textbf{\emph{расположенной в плоскости}} ($P$), \textbf{\emph{то она параллельна самой плоскости.}}

Если $AB$ лежит в $P$, то $AB\z\parallel P$ по определению.
Значит можно предположить, что $AB$ не лежит в~$P$.

Проведём через $AB$ и $CD$ плоскость $R$ и предположим, что прямая $AB$ где-нибудь пересекается с плоскостью~$P$.
Тогда точка пересечения, находясь на прямой $AB$, должна принадлежать также и плоскости $R$, на которой лежит прямая $AB$, в то же время точка пересечения, конечно, должна принадлежать и плоскости~$P$.
Значит, точка пересечения, находясь одновременно и на плоскости $R$, и на плоскости $P$, должна лежать на прямой $CD$, по которой пересекаются эти плоскости;
следовательно, прямая $AB$ пересекается с прямой $CD$.
Но это невозможно, так как по условию $AB\parallel CD$.
Значит, нельзя допустить, что прямая $AB$ пересекалась с плоскостью $P$, и потому $AB\parallel P$.

\paragraph{}\label{1938/s11}
\so{Теорема}. \textbf{\emph{Если плоскость}} ($R$, рис.~\ref{1938/s-ris-5}) \textbf{\emph{проходит через прямую}} ($AB$), \textbf{\emph{параллельную другой плоскости}} ($P$), \textbf{\emph{и пересекает эту плоскость, то линия пересечения}} ($CD$) \textbf{\emph{параллельна первой прямой}} ($AB$).

Действительно, во-первых, прямая $CD$ лежит в одной плоскости с прямой $AB$, во-вторых, эта прямая не может пересечься с прямой $AB$, потому что в противном случае прямая $AB$ пересекалась бы с плоскостью $P$, что невозможно.

{

\begin{wrapfigure}[21]{o}{34 mm}
\vskip-6mm
\centering
\includegraphics{mppics/s-ris-6}
\caption{}\label{1938/s-ris-6}
\bigskip
\includegraphics{mppics/s-ris-7}
\caption{}\label{1938/s-ris-7}
\end{wrapfigure}

\paragraph{}\label{1938/s12}
\mbox{\so{Следствие} 1.}
\emph{Если прямая} ($AB$, рис.~\ref{1938/s-ris-6}) \emph{параллельна каждой из двух пересекающихся плоскостей} ($P$ и $Q$), \emph{то она параллельна линии их пересечения} ($CD$).

Проведём плоскость через $AB$ и какую-нибудь точку $M$ прямой $CD$.
Эта плоскость должна пересечься с плоскостями $P$ и $Q$ по прямым, параллельным $AB$ и проходящим через точку $M$.
Но через точку $M$ можно провести только одну прямую, параллельную $AB$;
значит, две линии пересечения проведённой плоскости с плоскостями $P$ и $Q$ должны слиться в одну прямую.
Эта прямая, находясь одновременно на плоскости $P$ и на плоскости $Q$, должна совпадать с прямой $CD$, по которой плоскости $P$ и $Q$ пересекаются;
значит, $CD \parallel AB$.

\paragraph{}\label{1938/s13}
\mbox{\so{Следствие} 2.}
\emph{Если две прямые} ($AB$ и $CD$, рис.~\ref{1938/s-ris-7}) \emph{параллельны третьей прямой} ($EF$), \emph{то они параллельны между собой.}

Если прямые $AB$, $CD$ и $EF$ лежат в одной плоскости, то $AB\parallel CD$.
Остаётся рассмотреть случай когда прямые $AB$, $CD$ и $EF$ не лежат в одной плоскости.

Проведём плоскость $M$ через параллельные прямые $AB$ и $EF$.
Так как $CD \parallel EF$, то $CD \parallel M$ (§~\ref{1938/s10}).

}

Проведём также плоскость $N$ через $CD$ и некоторую точку $A$ прямой $AB$.
Так как $EF\z\parallel CD$, то $EF\parallel N$.
Значит, плоскость $N$ должна пересечься с плоскостью $M$ по прямой, параллельной $EF$ (§~\ref{1938/s11}) и в то же время проходящей через точку $A$.
Но в плоскости $M$ через $A$ проходит единственная прямая, параллельная $EF$, а именно прямая $AB$.
Следовательно, плоскость $N$ пересекается с $M$ по прямой $AB$, значит, $CD \parallel AB$.

\subsection*{Параллельные плоскости}

\paragraph{}\label{1938/s14}
\so{Определение}.
Две плоскости называются \rindex{параллельность}\textbf{параллельными}, если они не пересекаются, сколько бы их ни продолжали.
При этом плоскость считается параллельной самой себе.

\begin{figure}[!ht]
\centering
\includegraphics{mppics/s-ris-8}
\caption{}\label{1938/s-ris-8}
\end{figure}

\paragraph{}\label{1938/s15}
\so{Теорема}. \textbf{\emph{Если две пересекающиеся прямые}} ($AB$ и $AC$, рис.~\ref{1938/s-ris-8}) \textbf{\emph{одной плоскости}} ($P$) \textbf{\emph{соответственно параллельны двум прямым}} ($A_1B_1$ и $A_1C_1$) \textbf{\emph{другой плоскости}} ($Q$), \textbf{\emph{то эти плоскости параллельны}}.

Прямые $AB$ и $AC$ параллельны плоскости $Q$ (§~\ref{1938/s10}).

\begin{wrapfigure}{r}{34 mm}
\centering
\includegraphics{mppics/s-ris-9}
\caption{}\label{1938/s-ris-9}
\end{wrapfigure}

Допустим, что плоскости $P$ и $Q$ пересекаются по некоторой прямой $DE$ (рис.~\ref{1938/s-ris-8}).
В таком случае $AB \parallel DE$ и $AC \parallel DE$ (§~\ref{1938/s11}).
Таким образом, в плоскости $P$ через точку $A$ проходят две прямые $AB$ и $AC$, параллельные прямой $DE$, что невозможно.
Значит, плоскости $P$ и $Q$ не пересекаются.

\paragraph{}\label{1938/s16}
\mbox{\so{Теорема}.} \textbf{\emph{Если две параллельные плоскости}} ($P$ и $Q$, рис.~\ref{1938/s-ris-9}) \textbf{\emph{пересекаются третьей плоскостью}} ($R$), \textbf{\emph{то линии пересечения}} ($AB$ и $CD$) \textbf{\emph{параллельны.}}

Действительно, во-первых, прямые $AB$ и $CD$ находятся в одной плоскости ($R$);
во-вторых, они не могут пересечься, так как в противном случае пересекались бы плоскости $P$ и $Q$, что противоречит условию.

\paragraph{}\label{1938/s17} \mbox{\so{Теорема}.}
\textbf{\emph{Отрезки параллельных прямых}} ($AC$ и $BD$, рис.~\ref{1938/s-ris-9}), \textbf{\emph{заключённые между параллельными плоскостями}} ($P$ и $Q$), \textbf{\emph{равны.}}

Через параллельные прямые $AC$ и $BD$ проведём плоскость $R$;
она пересечёт плоскости $P$ и $Q$ по параллельным прямым $AB$ и $CD$;
следовательно, фигура $ABDC$ есть параллелограмм, и потому $AC = BD$.

\paragraph{}\label{1938/s18}
\mbox{\so{Теорема}.} \textbf{\emph{Два угла}} ($BAC$ и $B_1A_1C_1$, рис.~\ref{1938/s-ris-10}) \textbf{\emph{с соответственно параллельными и одинаково направленными сторонами равны и лежат в параллельных плоскостях}} ($P$ и $Q$).

\begin{wrapfigure}{o}{44 mm}
\centering
\includegraphics{mppics/s-ris-10}
\caption{}\label{1938/s-ris-10}
\end{wrapfigure}

Что плоскости $P$ и $Q$ параллельны, было доказано выше (§~\ref{1938/s15});
остаётся доказать, что углы $A$ и $A_1$ равны.

Отложим на сторонах углов произвольные, но равные отрезки $AB\z=A_1B_1$, $AC=A_1C_1$ и проведём прямые $AA_1$, $BB_1$, $CC_1$, $BC$ и $B_1C_1$.
Так как отрезки $AB$ и $A_1B_1$ равны и параллельны, то фигура $ABB_1A_1$ есть параллелограмм;
поэтому отрезки $AA_1$ и $BB_1$ равны и параллельны.
По той же причине равны и параллельны отрезки $AA_1$ и $CC_1$.
Следовательно, $BB_1\parallel CC_1$ и $BB_1=CC_1$ (§§~\ref{1938/s13} и \ref{1938/s17}).
Поэтому $BC=B_1C_1$ и $\triangle ABC=\triangle A_1B_1C_1$ (по трём сторонам);
значит, $\angle A=\angle A_1$.

\subsection*{Задачи на построение}

\begin{wrapfigure}{r}{30 mm}
\centering
\vskip-0mm
\includegraphics{mppics/s-ris-11}
\caption{}\label{1938/s-ris-11}
\end{wrapfigure}

\paragraph{}\label{1938/s19}
\emph{Через точку} ($A$, рис.~\ref{1938/s-ris-11}), \emph{расположенную вне данной прямой} ($a$), \emph{в пространстве провести прямую, параллельную данной прямой}~($a$).

\medskip

\mbox{\so{Решение}.}
Через прямую $a$ и точку $A$ проводим плоскость $M$.
В этой плоскости строим прямую $b$, параллельную прямой $a$.

Задача имеет единственное решение.
В самом деле, искомая прямая должна лежать с прямой $a$ в одной плоскости.
В этой же плоскости должна находиться точка $A$, через которую проходит искомая прямая.
Значит, эта плоскость должна совпадать с $M$.
Но в плоскости $M$ через точку $A$ можно провести только одну прямую, параллельную прямой~$a$.

\paragraph{}\label{1938/s20}
\emph{Через данную точку} ($A$, рис.~\ref{1938/s-ris-12}) \emph{провести плоскость, параллельную данной плоскости} ($P$), \emph{не проходящей через точку} $A$.

\begin{wrapfigure}{o}{44 mm}
\centering
\includegraphics{mppics/s-ris-12}
\caption{}\label{1938/s-ris-12}
\end{wrapfigure}

\mbox{\so{Решение}.}
Проводим на плоскости $P$ через какую-либо точку $B$ две какие-либо прямые $BC$ и $BD$.
Построим две вспомогательные плоскости: плоскость $M$ — через точку $A$ и прямую $BC$ и плоскость $N$ — через точку $A$ и прямую $BD$.
Искомая плоскость, параллельная плоскости $P$, должна пересечь плоскость $M$ по прямой, параллельной $BC$, а плоскость $N$ — по прямой, параллельной $BD$ (§~\ref{1938/s16}).

Отсюда вытекает такое построение: через точку $A$ проводим в плоскости $M$ прямую $AC_1\parallel BC$, а в плоскости $N$ прямую $AD_1\parallel BD$.
Через прямые $AC_1$ и $AD_1$ проводим плоскость $Q$.
Она и будет искомой.
В самом деле, стороны угла $D_1AC_1$, расположенного в плоскости $Q$, параллельны сторонам угла $DBC$, расположенного в плоскости~$P$.
Следовательно, $Q\parallel P$.

Так как в плоскости $M$ через точку $A$ можно провести лишь одну прямую, параллельную $BC$, а в плоскости $N$ через точку $A$ лишь одну прямую, параллельную $BD$, то задача имеет единственное решение.
Следовательно, через каждую точку пространства можно провести единственную плоскость, параллельную данной плоскости.

\paragraph{}\label{1938/s21}
\emph{Через данную прямую} ($a$, рис.~\ref{1938/s-ris-13}) \emph{провести плоскость, параллельную другой данной прямой} ($b$).

\begin{wrapfigure}{o}{54 mm}
\centering
\includegraphics{mppics/s-ris-13}
\caption{}\label{1938/s-ris-13}
\end{wrapfigure}

\mbox{\so{Решение}.}
1-й \so{случай}.
Прямые $a$ и $b$ не параллельны.
Через какую-нибудь точку $A$ прямой $a$ проводим прямую $b_1$, параллельную $b$;
через прямые $a$ и $b_1$ проводим плоскость.
Она и будет искомой (§~\ref{1938/s10}).
Задача имеет в этом случае единственное решение.

2-й \so{случай}.
Прямые $a$ и $b$ параллельны.
В этом случае задача имеет много решений: всякая плоскость, проходящая через прямую $a$, будет параллельна прямой $b$.

\paragraph{Пример более сложной задачи на построение.}\label{1938/s22}
Даны две скрещивающиеся прямые ($a$ и $b$, рис.~\ref{1938/s-ris-14}) и точка $A$, не лежащая ни на одной из данных прямых.
Провести через точку $A$ прямую, пересекающую обе данные прямые ($a$ и $b$).

\begin{wrapfigure}{o}{55 mm}
\centering
\includegraphics{mppics/s-ris-14}
\caption{}\label{1938/s-ris-14}
\end{wrapfigure}

\mbox{\so{Решение}.}
Так как искомая прямая должна проходить через точку $A$ и пересекать прямую $a$, то она должна лежать в плоскости, проходящей через прямую $a$ и точку $A$ (так как две её точки должны лежать в этой плоскости: точка $A$ и точка пересечения с прямой $a$).
Совершенно так же убеждаемся, что искомая прямая должна лежать в плоскости, проходящей через точку $A$ и прямую $b$.
Следовательно, она должна служить линией пересечения этих двух плоскостей.

Отсюда такое построение.
Через точку $A$ и прямую $a$ проводим плоскость $M$;
через точку $A$ и прямую $b$ проводим плоскость $N$.
Берём прямую $c$ пересечения плоскостей $M$ и $N$.
Если прямая $c$ не параллельна ни одной из данных прямых, то она пересечётся с каждой из данных прямых (так как с каждой из них она лежит в одной плоскости: $a$ и $c$ лежат в плоскости $M$, $b$ и $c$ — в плоскости $N$).
Прямая $c$ будет в этом случае искомой.
Если же $a\parallel c$ или $b\parallel c$, то задача не имеет решения.

\section{Перпендикуляр и наклонные к плоскости}

\paragraph{}\label{1938/s23}
Поставим задачу определить, в каком случае прямая может считаться перпендикулярной к плоскости.
Докажем предварительно следующее предложение:

\medskip

\so{Теорема}.
\textbf{\emph{Если прямая}} ($AA_1$, рис.~\ref{1938/s-ris-15}), \textbf{\emph{пересекающаяся с плоскостью}} ($M$), \textbf{\emph{перпендикулярна к каким-нибудь двум прямым}} ($OB$ и $OC$), \textbf{\emph{проведённым на этой плоскости через точку пересечения}} ($O$) \textbf{\emph{данной прямой и плоскости, то она перпендикулярна и ко всякой третьей прямой}} ($OD$), \textbf{\emph{проведённой на плоскости через ту же точку пересечения}} ($O$).

\begin{wrapfigure}{o}{40 mm}
\centering
\includegraphics{mppics/s-ris-15}
\caption{}\label{1938/s-ris-15}
\end{wrapfigure}

Отложим на прямой $AA_1$ произвольной длины, но равные отрезки $OA$ и $OA_1$ и проведём на плоскости какую-нибудь прямую, которая пересекала бы три прямые, исходящие из точки $O$, в каких-нибудь точках $C$, $D$ и $B$.
Эти точки соединим с точками $A$ и $A_1$.
Мы получим тогда несколько треугольников.
Рассмотрим их в такой последовательности.

Сначала возьмём треугольники $ACB$ и $A_1CB$;
они равны, так как у них $CB$ — общая сторона, $AC=A_1C$ как наклонные к прямой $AA_1$, одинаково удалённые от основания $O$ перпендикуляра $OC$;
по той же причине $AB=A_1B$.
Из равенства этих треугольников следует, что $\angle ABC\z=\angle A_1BC$.

После этого перейдём к треугольникам $ADB$ и $A_1DB$;
они равны, так как у них $DB$ — общая сторона, $AB=A_1B$ и $\angle ABD=\angle A_1BD$.
Из равенства этих треугольников выводим, что $AD=A_1D$.

Теперь возьмём треугольники $AOD$ и $A_1OD$;
они равны, так как имеют соответственно равные стороны.
Из их равенства выводим, что $\angle AOD=\angle A_1OD$;
а так как эти углы смежные, то, следовательно, $AA_1\perp OD$.

\paragraph{}\label{1938/s24}
\so{Определение}.
Прямая называется \rindex{перпендикулярность}\textbf{перпендикулярной к плоскости}, если она, пересекаясь с этой плоскостью, образует прямой угол с каждой прямой, проведённой на плоскости через точку пересечения.
В этом случае говорят также, что плоскость перпендикулярна к прямой.

Из предыдущей теоремы (§~\ref{1938/s23}) следует, что прямая перпендикулярна к плоскости, если она перпендикулярна к двум прямым, лежащим в данной плоскости и проходящим через точку пересечения данной прямой и плоскости.

Прямая, пересекающая плоскость, но не перпендикулярная к ней, называется \rindex{наклонная}\textbf{наклонной} к этой плоскости.
Точка пересечения прямой с плоскостью называется \rindex{основание!наклонной}\rindex{основание!перпендикуляра}\textbf{основанием} перпендикуляра или наклонной.

\paragraph{Сравнительная длина перпендикуляра и наклонных.}\label{1938/s25}%
\footnote{Для краткости термины «перпендикуляр» и «наклонная» употребляются вместо «отрезок перпендикуляра, ограниченный данной точкой и основанием перпендикуляра», и «отрезок наклонной, ограниченный данной точкой и основанием наклонной».}
Когда из одной точки $A$ (рис.~\ref{1938/s-ris-16}) проведены к плоскости перпендикуляр $AB$ и наклонная $AC$, условимся называть \rindex{проекция наклонной}\textbf{проекцией} наклонной на плоскость $P$ отрезок $BC$, соединяющий основание перпендикуляра и основание наклонной.
Таким образом, отрезок $BC$ есть проекция наклонной $AC$, отрезок $BD$ есть проекция наклонной $AD$ и так далее.

\medskip

\so{Теорема}.
\textbf{\emph{Если из одной и той же точки}} ($A$, рис.~\ref{1938/s-ris-16}), \textbf{\emph{взятой вне плоскости}} ($P$), \textbf{\emph{проведены к этой плоскости перпендикуляр}} ($AB$) \textbf{\emph{и какие-нибудь наклонные}} ($AC$, $AD$, $AE,\dots$), \textbf{\emph{то:}}

1) \textbf{\emph{две наклонные, имеющие равные проекции, равны;}}

2) \textbf{\emph{из двух наклонных та больше, проекция которой больше.}} 

\begin{wrapfigure}{o}{40 mm}
\centering
\includegraphics{mppics/s-ris-16}
\caption{}\label{1938/s-ris-16}
\end{wrapfigure}

Вращая прямоугольные треугольники $ABC$ и $ABD$ вокруг катета $AB$, мы можем совместить их плоскости с плоскостью $\triangle ABE$.
Тогда все наклонные будут лежать в одной плоскости с перпендикуляром, а все проекции расположатся на одной прямой.
Таким образом, доказываемые теоремы приводятся к аналогичным теоремам планиметрии.

{\small
\medskip

\mbox{\so{Замечание}.}
Так как $AB$ есть катет прямоугольного треугольника, а каждая из наклонных $AC, AD, AE,\dots$ есть гипотенуза, то перпендикуляр $AB$ меньше всякой наклонной;
значит, перпендикуляр, опущенный из точки на плоскость, есть наименьший из всех отрезков, соединяющих данную точку с любой точкой плоскости, и потому он принимается за меру расстояния точки $A$ от плоскости~$P$.

}

\paragraph{}\label{1938/s27}
\so{Обратные теоремы}.
\textbf{\emph{Если из одной и той же точки, взятой вне плоскости, проведены перпендикуляр и какие-нибудь наклонные, то: }}

1) \textbf{\emph{равные наклонные имеют равные проекции;}}

2) \textbf{\emph{из двух проекций та больше, которая соответствует большей наклонной.}}

Доказательство (от противного) предоставляем самим учащимся.

\paragraph{}\label{1938/s28}\so{Теорема}.
\textbf{\emph{Прямая}} ($DE$, рис.~\ref{1938/s-ris-17}), \textbf{\emph{проведённая на плоскости}} ($P$) \textbf{\emph{через основание наклонной}} ($AC$) \textbf{\emph{перпендикулярно к её проекции}} ($BC$), \textbf{\emph{перпендикулярна и к самой наклонной.}}

\begin{wrapfigure}{o}{50 mm}
\centering
\includegraphics{mppics/s-ris-17}
\caption{}\label{1938/s-ris-17}
\end{wrapfigure}

Отложим произвольные, но равные отрезки $CD$ и $CE$ и соединим прямолинейными отрезками точки $A$ и $B$ с точками $D$ и $E$.
Тогда $BD=BE$ как наклонные к прямой $DE$, одинаково удалённые от основания $C$ перпендикуляра $BC$.
Далее, $AD=AE$ как наклонные к плоскости $P$, имеющие равные проекции $BD$ и $BE$.
Вследствие этого $\triangle ADE$ равнобедренный, и потому его медиана $AC$ перпендикулярна к основанию $DE$.

Эта теорема носит название \so{теоремы о трёх перпендикулярах}.
Действительно, в ней говорится о связи, соединяющей следующие три перпендикуляра: 1) $AB$ к плоскости $P$, 2) $BC$ к прямой $DE$ и 3) $AC$ к той же прямой $DE$.

\paragraph{}\label{1938/s29}
\so{Обратная теорема}.
\textbf{\emph{Прямая}} ($DE$, рис.~\ref{1938/s-ris-17}), \textbf{\emph{проведённая на плоскости}} ($P$) \textbf{\emph{через основание наклонной}} ($AC$) \textbf{\emph{перпендикулярно к этой наклонной, перпендикулярна и к её проекции.}}

Проделаем те же построения, что и при доказательстве прямой теоремы.
Отложим произвольные, но равные отрезки $CD$ и $CE$ и соединим прямолинейными отрезками точки $A$ и $B$ с точками $D$ и $E$, тогда будем иметь: $AD=AE$ как наклонные к прямой $DE$, одинаково удалённые от основания $C$ перпендикуляра $AC$;
$BD=BE$ как проекции равных наклонных $AD$ и $AE$.
Вследствие этого $\triangle BDE$ равнобедренный, и потому его медиана $BC$ перпендикулярна к основанию $DE$.

\section[Параллельность и перпендикулярность]{Параллельность и перпендикулярность\\ прямых и плоскостей}

\paragraph{Предварительное замечание.}\label{1938/s30}
Параллельность прямых и плоскостей в пространстве и перпендикулярность прямой к плоскости находятся в некоторой зависимости.
Именно наличие параллельности одних элементов влечёт перпендикулярность других, и, обратно, из перпендикулярности одних элементов можно сделать заключение о параллельности других.
Эта связь между параллельностью и перпендикулярностью прямых и плоскостей в пространстве выражается следующими теоремами.

\paragraph{}\label{1938/s31}
\so{Теорема}.
\textbf{\emph{Если плоскость}} ($P$, рис.~\ref{1938/s-ris-18}) \textbf{\emph{перпендикулярна к одной из параллельных прямых}} ($AB$), \textbf{\emph{то она перпендикулярна и к другой}} ($CD$).

Проведём через точку $B$ на плоскости $P$ две какие-нибудь прямые $BE$ и $BF$, а через точку $D$ проведём прямые $DG$ и $DH$, соответственно параллельные и одинаково направленные прямым $BE$ и $BF$.
Тогда будем иметь: $\angle ABE=\angle CDG$ и $\angle ABF = \angle CDH$ как углы с параллельными сторонами.
Но углы $ABE$ и $ABF$ прямые, так как $AB\perp P$, значит, углы $CDG$ и $CDH$ также прямые (§~\ref{1938/s18}).
Следовательно, $CD \perp P$ (§~\ref{1938/s24}).

\begin{figure}[!ht]
\begin{minipage}{.48\textwidth}
\centering
\includegraphics{mppics/s-ris-18}
\end{minipage}
\hfill
\begin{minipage}{.48\textwidth}
\centering
\includegraphics{mppics/s-ris-19}
\end{minipage}

\medskip

\begin{minipage}{.48\textwidth}
\centering
\caption{}\label{1938/s-ris-18}
\end{minipage}
\hfill
\begin{minipage}{.48\textwidth}
\centering
\caption{}\label{1938/s-ris-19}
\end{minipage}
\vskip-4mm
\end{figure}

\paragraph{}\label{1938/s32}
\mbox{\so{Обратная теорема}.}
\textbf{\emph{Если две прямые}} ($AB$ и $CD$, рис. \ref{1938/s-ris-19}) \textbf{\emph{перпендикулярны к одной и той же плоскости}} ($P$), \textbf{\emph{то они параллельны.}}

Предположим противное, то есть что прямые $AB$ и $CD$ не параллельны.
Проведём тогда через точку $D$ прямую, параллельную $AB$. 
При нашем предположении это будет какая-нибудь прямая $DC_1$, не сливающаяся с $DC$.
Согласно прямой теореме прямая $DC_1$ будет перпендикулярна к плоскости~$P$.
Проведём через $CD$ и $C_1D$ плоскость $Q$ и возьмём линию её пересечения $DE$ с плоскостью~$P$.
Так как (на основании предыдущей теоремы) $C_1D \perp P$, то $\angle C_1DE$ прямой, а так как по условию $CD \perp P$, то $\angle CDE$ также прямой.
Таким образом, окажется, что в плоскости $Q$ к прямой $DE$ из одной её точки $D$ восстановлены два перпендикуляра $DC$ и $DC_1$.
Так как это невозможно, 
то нельзя допустить, что прямые $AB$ и $CD$ были не параллельны.

\paragraph{}\label{1938/s33}
\so{Теорема}.
\textbf{\emph{Если прямая}} ($BB_1$, рис.~\ref{1938/s-ris-20}) \textbf{\emph{перпендикулярна к одной из параллельных плоскостей}} ($P$), \textbf{\emph{то она перпендикулярна и к другой}} ($Q$).

\begin{wrapfigure}{r}{50 mm}
\vskip-4mm
\centering
\includegraphics{mppics/s-ris-20}
\caption{}\label{1938/s-ris-20}
\end{wrapfigure}

Проведём через прямую $BB_1$ какие-нибудь две плоскости $M$ и $N$, каждая из которых пересекается с $P$ и $Q$ по параллельным прямым: одна — по параллельным прямым $BC$ и $B_1C_1$, другая — по параллельным прямым $BD$ и $B_1D_1$.
Согласно условию прямая $BB_1$ перпендикулярна к прямым $BC$ и $BD$;
следовательно, она также перпендикулярна к параллельным им прямым $B_1C_1$ и $B_1D_1$, а потому перпендикулярна и к плоскости $Q$, на которой лежат прямые $B_1C_1$ и $B_1D_1$.

\paragraph{}\label{1938/s34}
\so{Обратная теорема}.
\textbf{\emph{Если две плоскости}} ($P$ и $Q$, рис.~\ref{1938/s-ris-21}) \textbf{\emph{перпендикулярны к одной и той же прямой}} ($AB$), \textbf{\emph{то они параллельны.}}

\begin{figure}[!ht]
\centering
\includegraphics{mppics/s-ris-21}
\caption{}\label{1938/s-ris-21}
\end{figure}

Предположим противное, то есть что две различные плоскости $P$ и $Q$ пересекаются.
Возьмём на линии их пересечения какую-нибудь точку $C$ и проведём плоскость $R$ через $C$ и прямую $AB$.
Плоскость $R$ пересечёт плоскости $P$ и $Q$ соответственно по прямым $AC$ и $BC$.
Так как $AB\perp P$, то $AB\perp AC$, и так как $AB\perp Q$, то $AB\perp BC$.
Таким образом, в плоскости $R$ мы будем иметь два перпендикуляра к прямой $AB$, проходящих через одну и ту же точку $C$, перпендикуляры $AC$ и $BC$.
Так как это невозможно, то предположение, что плоскости $P$ и $Q$ пересекаются, неверно.
Значит, они параллельны.

\subsection*{Задачи на построение}

\paragraph{}\label{1938/s35}
\emph{Через данную точку в пространстве провести плоскость, перпендикулярную к данной прямой} ($AB$).

\medskip

\so{Решение}.
1-й \so{случай}.
Данная точка $C$ лежит на прямой $AB$ (рис. \ref{1938/s-ris-22}).

\begin{wrapfigure}{o}{50 mm}
\centering
\includegraphics{mppics/s-ris-22}
\caption{}\label{1938/s-ris-22}
\end{wrapfigure}

Проведём через прямую $AB$ какие-нибудь две плоскости $P$ и $Q$.
Искомая плоскость должна пересекать эти плоскости по прямым, перпендикулярным к прямой $AB$ (§~\ref{1938/s24}).

Отсюда построение: через $AB$ проводим две произвольные плоскости $P$ и $Q$.
В каждой из этих плоскостей восстановим перпендикуляр к прямой $AB$ в точке $C$ (в плоскости $P$ — перпендикуляр $CD$, в плоскости $Q$ — перпендикуляр $CE$).
Плоскость, проходящая через прямые $CD$ и $CE$, есть искомая.

2-й \so{случай}.
Данная точка $D$ не лежит на прямой $AB$ (рис.~\ref{1938/s-ris-22}).
Через точку $D$ и прямую $AB$ проводим плоскость $P$ и в этой плоскости строим прямую $DC$, перпендикулярную к $AB$.
Через прямую $AB$ проводим произвольно вторую плоскость $Q$ и в этой плоскости строим прямую $CE$, перпендикулярную к $AB$.
Искомая плоскость должна пересечь плоскости $P$ и $Q$ по прямым, перпендикулярным к $AB$.

Отсюда построение: через точку $D$ проводим в плоскости $P$ прямую $DC$, перпендикулярную к $AB$.
Прямая $DC$ пересечёт прямую $AB$ в некоторой точке $C$.
Через точку $C$ проводим в плоскости $Q$ прямую $CE$ перпендикулярно к $AB$.
Плоскость, проходящая через прямые $CD$ и $CE$, — искомая.

Так как в каждой из плоскостей $P$ и $Q$ через данную точку можно провести лишь одну прямую, перпендикулярную к данной, то задача в обоих случаях имеет одно решение, то есть через каждую точку в пространстве можно провести лишь одну плоскость, перпендикулярную к данной прямой.

\paragraph{}\label{1938/s36}
\emph{Через данную точку} ($O$) \emph{пространства провести прямую, перпендикулярную к данной плоскости} ($P$).

1-й \so{случай}.
Точка $O$ лежит на плоскости $P$ (рис.~\ref{1938/s-ris-23}).
Проведём на плоскости $P$ через точку $O$ две какие-либо взаимно перпендикулярные прямые $OA$ и $OB$.
Проведём, далее, через прямую $OA$ какую-либо новую плоскость $Q$ и на этой плоскости $Q$ построим прямую $OC$, перпендикулярную к $OA$.
Через прямые $OB$ и $OC$ проведём новую плоскость $R$ и построим в ней прямую $OM$, перпендикулярную к $OB$.
Прямая $OM$ и будет искомым перпендикуляром к плоскости~$P$.

Действительно, так как $OA \perp OB$ и $OA \perp OC$, то прямая $OA$ перпендикулярна к плоскости $R$ и, следовательно, $OA \perp OM$.
Таким образом, мы видим, что $OM\perp OA$ и $OM\perp OB$;
следовательно, $OM$ перпендикулярна к плоскости~$P$.

\begin{figure}[!ht]
\begin{minipage}{.48\textwidth}
\centering
\includegraphics{mppics/s-ris-23}
\end{minipage}
\hfill
\begin{minipage}{.48\textwidth}
\centering
\includegraphics{mppics/s-ris-24}
\end{minipage}

\medskip

\begin{minipage}{.48\textwidth}
\centering
\caption{}\label{1938/s-ris-23}
\end{minipage}
\hfill
\begin{minipage}{.48\textwidth}
\centering
\caption{}\label{1938/s-ris-24}
\end{minipage}
\vskip-4mm
\end{figure}

2-й \mbox{\so{случай}.}
Точка $O$ не лежит на плоскости $P$ (рис.~\ref{1938/s-ris-24}).
Возьмём на плоскости $P$ какую-нибудь точку $A$ и выполним для неё предыдущее построение.
Мы получим тогда прямую $AB$, перпендикулярную к плоскости~$P$.
После этого через точку $O$ проводим прямую, параллельную $AB$.
Эта прямая и будет искомой (§~\ref{1938/s31}).

\medskip

Задача в обоих случаях имеет одно решение.
В самом деле, так как два перпендикуляра к одной и той же плоскости параллельны, то через одну и ту же точку $O$ нельзя провести двух перпендикуляров к плоскости~$P$.
Следовательно, через каждую точку в пространстве можно провести одну и только одну прямую, перпендикулярную к данной плоскости.

\paragraph{Пример более сложной задачи.}\label{1938/s37}
\emph{Даны две скрещивающиеся прямые} ($a$ и $b$, рис.~\ref{1938/s-ris-25}).
\emph{Построить прямую, пересекающую обе данные прямые и перпендикулярную к ним обеим.}

\medskip

\mbox{\so{Решение}.}
Проведём через прямую $a$ плоскость $M$, параллельную прямой $b$ (§~\ref{1938/s21}).
Из двух каких-нибудь точек прямой $b$ опустим перпендикуляры $AA_1$ и $BB_1$ на плоскость $M$.
Соединим точки $A_1$ и $B_1$ отрезком прямой и найдём точку $C_1$ пересечения прямых $A_1B_1$ и $a$.
Через точку $C_1$ проведём прямую, перпендикулярную к плоскости $M$.

Предоставляем самим учащимся доказать, что эта прямая 1) пересечётся с прямой $b$ в некоторой точке $C$ и 2) будет перпендикулярна как к прямой $a$, так и к прямой $b$.
Прямая $CC_1$ будет, следовательно, искомой прямой.

\begin{wrapfigure}{o}{55 mm}
\vskip-0mm
\centering
\includegraphics{mppics/s-ris-25}
\caption{}\label{1938/s-ris-25}
\end{wrapfigure}

Заметим, что отрезок $CC_1$ меньше всех других отрезков, которые можно получить, соединяя точки прямой $a$ с точками прямой $b$.
В самом деле, возьмём на прямой а какую-нибудь точку $E$ и на прямой $b$ какую-нибудь точку $F$, соединим эти точки отрезком прямой и докажем, что $EF>CC_1$.
Опустим из точки $F$ перпендикуляр $FF_1$ на плоскость $M$.
Тогда будем иметь: $EF>FF_1$ (§~\ref{1938/s25}).
Но $FF_1=CC_1$, следовательно $EF>CC_1$.
На этом основании длина отрезка $CC_1$ называется \rindex{кратчайшее расстояние}\textbf{кратчайшим расстоянием} между данными прямыми $a$ и~$b$.

\section[Углы]{Двугранные углы,
угол прямой с плоскостью,
угол двух скрещивающихся прямых,
многогранные углы}

\subsection*{Двугранные углы}

\begin{wrapfigure}{r}{34 mm}
\vskip-0mm
\centering
\includegraphics{mppics/s-ris-26}
\caption{}\label{1938/s-ris-26}
\bigskip
\includegraphics{mppics/s-ris-27}
\caption{}\label{1938/s-ris-27}
\bigskip
\includegraphics{mppics/s-ris-28}
\caption{}\label{1938/s-ris-28}
\end{wrapfigure}

\paragraph{}\label{1938/s38}
\mbox{\so{Определения}.}
Часть плоскости, лежащая по одну сторону от какой-либо прямой, лежащей в этой плоскости, называется \rindex{полуплоскость}\textbf{полуплоскостью}.
Фигура, образованная двумя полуплоскостями ($P$ и $Q$, рис.~\ref{1938/s-ris-26}), исходящими из одной прямой ($AB$), называется \rindex{двугранный угол}\textbf{двугранным углом}.
Прямая $AB$ называется \rindex{ребро!двугранного угла}\textbf{ребром}, а полуплоскости $P$ и $Q$ — \rindex{грань!двугранного угла}\textbf{гранями} двугранного угла.

Такой угол обозначается обыкновенно двумя буквами, поставленными у его ребра (двугранный угол $AB$).
Но если при одном ребре лежат несколько двугранных углов, то каждый из них обозначают четырьмя буквами, из которых две средние стоят при ребре, а две крайние — у граней (например, двугранный угол $SCDR$, рис.~\ref{1938/s-ris-27}).

Если из произвольной точки $D$ ребра $AB$ (рис.~\ref{1938/s-ris-28}) проведём на каждой грани по перпендикуляру к ребру, то образованный ими угол $CDE$ называется \rindex{линейный угол}\textbf{линейным углом} двугранного угла.

Величина линейного угла не зависит от положения его вершины на ребре.
Так, линейные углы $CDE$ и $C_1D_1E_1$ равны, потому что их стороны соответственно параллельны и одинаково направлены.

Плоскость линейного угла перпендикулярна к ребру, так как она содержит две прямые, перпендикулярные к нему.
Поэтому для получения линейного угла достаточно грани данного двугранного угла пересечь плоскостью, перпендикулярной к ребру, и рассмотреть получившийся в этой плоскости угол.

\paragraph{Равенство и неравенство двугранных углов.}\label{1938/s39}
Два двугранных угла считаются \rindex{равные двугранные углы}\textbf{равными}, если они при вложении могут совместиться;
в противном случае тот из двугранных углов считается меньшим, который составит часть другого угла.

Подобно углам в планиметрии двугранные углы могут быть \rindex{смежные углы}\textbf{смежные}, \rindex{вертикальные углы}\textbf{вертикальные} и прочие.

Если два смежных двугранных угла равны между собой, то каждый из них называется прямым двугранным углом.

\medskip

\so{Теоремы}.
1) \textbf{\emph{Равным двугранным углам соответствуют равные линейные углы.}}

2) \textbf{\emph{Б\'{о}льшему двугранному углу соответствует больший линейный угол.}}

\begin{wrapfigure}{o}{64 mm}
\centering
\includegraphics{mppics/s-ris-29}
\caption{}\label{1938/s-ris-29}
\end{wrapfigure}

Пусть $PABQ$ и $P_1A_1B_1Q_1$ (рис.~\ref{1938/s-ris-29}) — два двугранных угла.
Вложим угол $A_1B_1$ в угол $AB$ так, чтобы ребро $A_1B_1$ совпало с ребром $AB$ и грань $P_1$ с гранью~$P$.
Тогда если эти двугранные углы равны, то грань $Q_1$ совпадёт с гранью $Q$;
если же угол $A_1B_1$ меньше угла $AB$, то грань $Q_1$ займёт некоторое положение внутри двугранного угла, например~$Q_2$.

Заметив это, возьмём на общем ребре какую-нибудь точку $B$ и проведём через неё плоскость $R$, перпендикулярную к ребру.
От пересечения этой плоскости с гранями двугранных углов получатся линейные углы.
Ясно, что если двугранные углы совпадут, то у них окажется один и тот же линейный угол $CBD$;
если же двугранные углы не совпадут, если, например, грань $Q_1$ займёт положение $Q_2$, то у б\'{о}льшего двугранного угла окажется больший линейный угол, а именно $\angle CBD > \angle C_2BD$.

\paragraph{}\label{1938/s40}
\so{Обратные теоремы}.
1) \textbf{\emph{Равным линейным углам соответствуют равные двугранные углы.}}

2) \textbf{\emph{Б\'{о}льшему линейному углу соответствует больший двугранный угол.}}

Эти теоремы легко доказываются от противного.

\paragraph{}\label{1938/s41}
\so{Следствия}. 1) \emph{Прямому двугранному углу соответствует прямой линейный угол, и обратно.}

Пусть (рис.~\ref{1938/s-ris-30}) двугранный угол $PABQ$ прямой.
Это значит, что он равен смежному углу $QABP_1$.
Но в таком случае линейные углы $CDE$ и $CDE_1$ также равны;
а так как они смежные, то каждый из них должен быть прямой.

\begin{wrapfigure}{o}{60 mm}
\centering
\includegraphics{mppics/s-ris-30}
\caption{}\label{1938/s-ris-30}
\end{wrapfigure}

Обратно, если равны смежные линейные углы $CDE$ и $CDE_1$, то равны и смежные двугранные углы, то есть каждый из них должен быть прямой.

2) \emph{Все прямые двугранные углы равны}, потому что у них равны линейные углы.

Подобным же образом легко доказать, что:

3) \emph{Вертикальные двугранные углы равны.}

4) \emph{Двугранные углы с соответственно параллельными и одинаково (или противоположно) направленными гранями равны.}

5) Если за единицу двугранных углов возьмём такой двугранный угол, который соответствует единице линейных углов, то можно сказать, что \emph{двугранный угол измеряется его линейным углом}.

\subsection*{Перпендикулярные плоскости}

\paragraph{}\label{1938/s42}
\mbox{\so{Определение}.}
Две плоскости называются \rindex{перпендикулярность}\textbf{взаимно перпендикулярными}, если, пересекаясь, они образуют прямые двугранные углы.

\begin{wrapfigure}{r}{38 mm}
\vskip-5mm
\centering
\includegraphics{mppics/s-ris-31}
\caption{}\label{1938/s-ris-31}
\end{wrapfigure}

\paragraph{Признак перпендикулярности.}\label{1938/s43}\ 

\mbox{\so{Теорема}.}
\textbf{\emph{Если плоскость}} ($P$, рис. \ref{1938/s-ris-31}) \textbf{\emph{проходит через перпендикуляр}} ($AB$) \textbf{\emph{к другой плоскости}} ($Q$), \textbf{\emph{то она перпендикулярна к этой плоскости.}}

Пусть $DE$ будет линия пересечения плоскостей $P$ и $Q$.
На плоскости $Q$ проведём $BC \z\perp DE$.
Тогда угол $ABC$ будет линейным углом двугранного угла $PDEQ$.
Так как прямая $AB$ по условию перпендикулярна к $Q$, то $AB\perp BC$;
значит, угол $ABC$ прямой, а потому и двугранный угол прямой, то есть плоскость $P$ перпендикулярна к плоскости~$Q$.

\paragraph{}\label{1938/s44}
\so{Теорема}.
\textbf{\emph{Если две плоскости}} ($P$ и $Q$, рис.~\ref{1938/s-ris-31}) \textbf{\emph{взаимно перпендикулярны и к одной из них}} ($Q$) \textbf{\emph{проведён перпендикуляр}} ($AB$), \textbf{\emph{имеющий общую точку}} ($A$) \textbf{\emph{с другой плоскостью}} ($P$), \textbf{\emph{то этот перпендикуляр весь лежит в этой плоскости}}~($P$).

\begin{wrapfigure}{o}{45 mm}
\centering
\includegraphics{mppics/s-ris-32}
\caption{}\label{1938/s-ris-32}
\bigskip
\includegraphics{mppics/s-ris-33}
\caption{}\label{1938/s-ris-33}
\end{wrapfigure}

Предположим, что перпендикуляр $AB$ не лежит в плоскости $P$ (как изображено на рис.~\ref{1938/s-ris-32}).
Пусть $DE$ будет линия пересечения плоскостей $P$ и $Q$.
На плоскости $P$ проведём прямую $AC \z\perp DE$, а на плоскости $Q$ проведём прямую $CF \z\perp DE$.
Тогда угол $ACF$, как линейный угол прямого двугранного угла, будет прямой.
Поэтому линия $AC$, образуя прямые углы с $DE$ и $CF$, будет перпендикуляром к плоскости $Q$.
Мы будем иметь тогда два перпендикуляра, опущенные из одной и той же точки $A$ на плоскость $Q$, а именно $AB$ и $AC$.
Так как это невозможно (§~\ref{1938/s36}), то допущение неверно;
значит, перпендикуляр $AB$ лежит в плоскости~$P$.

\paragraph{}\label{1938/s45}
\mbox{\so{Следствие}.}
\emph{Линия пересечения} ($AB$, рис.~\ref{1938/s-ris-33}) \emph{двух плоскостей} ($P$ и $Q$), \emph{перпендикулярных к третьей плоскости} ($R$), \emph{есть перпендикуляр к этой плоскости}.

Действительно, если через какую-нибудь точку $A$ линии пересечения плоскостей $P$ и $Q$ проведём перпендикуляр к плоскости $R$, то этот перпендикуляр согласно предыдущей теореме должен лежать и в плоскости $Q$, и в плоскости $P$; значит, он сольётся с $AB$.

\subsection*{Угол двух скрещивающихся прямых}

\paragraph{}\label{1938/s46}
\mbox{\so{Определение}.}
Углом двух скрещивающихся прямых ($AB$ и $CD$, рис.~\ref{1938/s-ris-34}), для которых дано положение и направление, называется угол ($MON$), 
\begin{figure}[!ht]
\vskip-0mm
\centering
\includegraphics{mppics/s-ris-34}
\caption{}\label{1938/s-ris-34}
\end{figure}
который получится, если из произвольной точки пространства ($O$) проведём полупрямые ($OM$ и $ON$), соответственно параллельные данным прямым ($AB$ и $CD$) и одинаково с ними направленные.

Величина этого угла не зависит от положения точки $O$, так как если построим указанным путём угол $M_1O_1N_1$ с вершиной в какой-нибудь другой точке $O_1$, то $\angle MON = \angle M_1O_1N_1$, потому что эти углы имеют соответственно параллельные и одинаково направленные стороны.

\subsection*{Угол, образуемый прямой с плоскостью}

\begin{wrapfigure}{r}{34 mm}
\vskip-6mm
\centering
\includegraphics{mppics/s-ris-35}
\caption{}\label{1938/s-ris-35}
\bigskip
\includegraphics{mppics/s-ris-36}
\caption{}\label{1938/s-ris-36}
\end{wrapfigure}

\paragraph{Проекция точки и прямой на плоскость.}\label{1938/s47}
Мы говорили ранее (§~\ref{1938/s25}), что когда из одной точки проведены к плоскости перпендикуляр и наклонная, то проекцией этой наклонной на плоскость называется отрезок, соединяющий основание перпендикуляра с основанием наклонной.
Дадим теперь более общее определение проекции.

1) \rindex{ортогональная проекция}\textbf{Ортогональной} (или \textbf{прямоугольной}) \textbf{проекцией} какой-нибудь \textbf{точки} на данную плоскость (например, точки $M$ на плоскость $P$, рис.~\ref{1938/s-ris-35}) называется основание ($m$) перпендикуляра, опущенного на эту плоскость из взятой точки.

2) \textbf{Ортогональной проекцией} какой-нибудь \textbf{линии} на плоскость называется геометрическое место проекций всех точек этой линии.

В частности, если проектируемая линия есть прямая (например, $AB$, рис.~\ref{1938/s-ris-35}), не перпендикулярная к плоскости ($P$), то проекция её на эту плоскость есть также прямая.
В самом деле, если мы через прямую $AB$ и перпендикуляр $Mm$, опущенный на плоскость проекций из какой-нибудь одной точки $M$ этой прямой, проведём плоскость $Q$, то эта плоскость должна быть перпендикулярна к плоскости $P$;
поэтому перпендикуляр, опущенный на плоскость $P$ из любой точки прямой $AB$ (например, из точки $N$), должен лежать в этой плоскости $Q$ (§~\ref{1938/s44}) и, следовательно, проекции всех точек прямой $AB$ должны лежать на прямой $ab$, по которой пересекаются плоскости $P$ и $Q$.

Обратно, всякая точка этой прямой $ab$ есть проекция какой-нибудь точки прямой $AB$, так как перпендикуляр, восстановленный из любой точки прямой $ab$, лежит на плоскости $Q$ и, следовательно, пересекается с $AB$ в некоторой точке.
Таким образом, прямая $ab$ представляет собой геометрическое место проекций всех точек данной прямой $AB$ и, следовательно, есть её проекция.

Для краткости вместо «ортогональная проекция» мы будем говорить просто «проекция».

\paragraph{Угол прямой с плоскостью.}\label{1938/s48}
Углом прямой ($AB$, рис.~\ref{1938/s-ris-36}) с плоскостью ($P$) в том случае, когда прямая наклонна к плоскости, называется острый угол ($ABC$), составленный этой прямой с её проекцией на плоскость.
Угол прямой с плоскостью считается прямым, если прямая перпендикулярна к плоскости.

Угол этот обладает тем свойством, что он есть наименьший из всех углов, которые данная прямая образует с прямыми, проведёнными на плоскости $P$ через её основание $B$.
Докажем, например, что угол $ABC$ меньше угла $ABD$.
Для этого отложим отрезок $BD=BC$ и соединим $D$ с $A$.
У треугольников $ABC$ и $ABD$ две стороны одного равны соответственно двум сторонам другого, но третьи стороны не равны, а именно $AD>AC$ (§~\ref{1938/s25}).
Вследствие этого угол $ABD$ больше угла $ABC$.

\subsection*{Многогранные углы}

\paragraph{}\label{1938/s49}
\mbox{\so{Определения}.}
Возьмём несколько углов (рис.~\ref{1938/s-ris-37}): $ASB$, $BSD$, $CSD$, которые, примыкая последовательно один к другому, расположены в одной плоскости вокруг общей вершины $S$.
Повернём плоскость угла $ASB$ вокруг общей стороны $SB$ так, чтобы эта плоскость составила некоторый двугранный угол с плоскостью $BSD$.
Затем, не изменяя получившегося двугранного угла, повернём его вокруг прямой $SD$ так, чтобы плоскость $BSD$ составила некоторый двугранный угол с плоскостью $CSD$.
Продолжим такое последовательное вращение вокруг каждой общей стороны.
Если при этом последняя сторона $SF$ совместится с первой стороной $SA$, то образуется фигура (рис.~\ref{1938/s-ris-38}), 
\begin{figure}[!ht]
\begin{minipage}{.32\textwidth}
\centering
\includegraphics{mppics/s-ris-37}
\end{minipage}
\hfill
\begin{minipage}{.32\textwidth}
\centering
\includegraphics{mppics/s-ris-38}
\end{minipage}
\hfill
\begin{minipage}{.32\textwidth}
\centering
\includegraphics{mppics/s-ris-39}
\end{minipage}

\medskip

\begin{minipage}{.32\textwidth}
\centering
\caption{}\label{1938/s-ris-37}
\end{minipage}
\hfill
\begin{minipage}{.32\textwidth}
\centering
\caption{}\label{1938/s-ris-38}
\end{minipage}
\hfill
\begin{minipage}{.32\textwidth}
\centering
\caption{}\label{1938/s-ris-39}
\end{minipage}
\vskip-4mm
\end{figure}
которая называется \rindex{многогранный угол}\textbf{многогранным углом}.
Углы $ASB, BSD,\dots$ называются плоскими углами или \rindex{грань!многогранного угла}\textbf{гранями}, стороны их $SA, SB,\dots$ называются \rindex{ребро!многогранного угла}\textbf{рёбрами}, а общая вершина $S$ — \rindex{вершина!многогранного угла}\textbf{вершиной} многогранного угла.
Каждое ребро многогранного угла является также ребром его двугранного угла;
поэтому в многогранном угле столько двугранных углов и столько плоских, сколько в нём всех рёбер.
Наименьшее число граней в многогранном угле — три;
такой угол называется \rindex{трёхгранный угол}\textbf{трёхгранным}.
Могут быть углы четырёхгранные, пятигранные и так далее.

Многогранный угол обозначается или одной буквой $S$, поставленной у вершины, или же рядом букв $SABCDE$, из которых первая обозначает вершину, а прочие — рёбра по порядку их расположения.

Многогранный угол называется \rindex{выпуклый многогранный угол}\textbf{выпуклым}, если он весь расположен по одну сторону от плоскости каждой из его граней, неограниченно продолженной.
Таков, например, угол, изображённый на рис.~\ref{1938/s-ris-38}.
Наоборот, угол на рис.~\ref{1938/s-ris-39} нельзя назвать выпуклым, так как он расположен по обе стороны от грани $ASB$ или от грани $BSD$.
Если все грани многогранного угла пересечём плоскостью, то в сечении образуется многоугольник ($abcde$).
В выпуклом многогранном угле этот многоугольник тоже выпуклый.

Мы будем рассматривать \so{только выпуклые многогранные углы.}

{

\begin{wrapfigure}{r}{40 mm}
\vskip-9mm
\centering
\includegraphics{mppics/s-ris-40}
\caption{}\label{1938/s-ris-40}
\end{wrapfigure}

\paragraph{}\label{1938/s50}
\mbox{\so{Теорема}.}
\textbf{\emph{В трёхгранном угле каждый плоский угол меньше суммы двух других плоских углов.}}

Пусть в трёхгранном угле $SABC$ (рис. \ref{1938/s-ris-40}) наибольший из плоских углов есть угол $ASC$.
Отложим на этом угле угол $ASD$, равный углу $ASB$, и проведём какую-нибудь прямую $AC$, пересекающую $SD$ в некоторой точке $D$.
Отложим $SB=SD$.

}

Соединив $B$ с $A$ и $C$, получим $\triangle ABC$, в котором
\[AD + DC < AB + BC.\]
Треугольники $ASD$ и $ASB$ равны, так как они содержат по равному углу, заключённому между равными сторонами;
следовательно, $AD \z=AB$.
Поэтому если в выведенном неравенстве отбросить равные слагаемые $AD$ и $AB$, получим, что $DC<BC$.
Теперь замечаем, что у треугольников $SCD$ и $SCB$ две стороны одного равны двум сторонам другого, а третьи стороны не равны;
в таком случае против большей из этих сторон лежит больший угол;
значит,
\[\angle CSD < \angle CSB.\]

Прибавив к левой части этого неравенства угол $ASD$, а к правой равный ему угол $ASB$, получим то неравенство, которое требовалось доказать:
\[\angle ASC < \angle CSB + \angle ASB.\]

Мы доказали, что даже наибольший плоский угол меньше суммы двух других углов.
Значит, теорема доказана.

\medskip

\so{Следствие}.
Отнимем от обеих частей последнего неравенства по углу $ASB$ или по углу $CSB$;
получим:
\[\angle ASC - \angle ASB < \angle CSB;\]
\[\angle ASC - \angle CSB < \angle ASB.\]
Рассматривая эти неравенства справа налево и приняв во внимание, что угол $ASC$ как наибольший из трёх углов больше разности двух других углов, мы приходим к заключению, что \so{в трёхгранном угле каждый плоский угол больше разности двух других углов}.

\begin{wrapfigure}[9]{o}{40 mm}
\vskip-8mm
\centering
\includegraphics{mppics/s-ris-41}
\caption{}\label{1938/s-ris-41}
\end{wrapfigure}

\paragraph{}\label{1938/s51}
\mbox{\so{Теорема}.}
\textbf{\emph{В выпуклом, многогранном угле сумма всех плоских углов меньше $\bm{360\degree}$.}}

Пересечём грани (рис.~\ref{1938/s-ris-41}) выпуклого угла $SABCDE$ такой плоскостью, что в сечении получится выпуклый $n$-угольник $ABCDE$.
Такую плоскость можно получить немного повернув плоскость грани угла $ASB$ вокруг прямой $AB$.
Применяя теорему предыдущего параграфа к каждому из трёхгранных углов, вершины которых находятся в точках $A$, $B$, $C$, $D$ и $E$, находим:
\begin{align*}
\angle ABC &< \angle ABS + \angle SBC,
\\
\angle BCD &< \angle BCS + \angle SCD
\\
&\text{\dots}
\end{align*}
Сложим почленно все эти неравенства.
Тогда в левой части получим сумму всех углов многоугольника $ABCDE$, которая равна $180\degree\cdot n \z- 360\degree$, а в правой — сумму углов треугольников $ABS$, $SBC$ и так далее, кроме тех углов, которые лежат при вершине $S$.
Обозначив сумму этих последних углов буквой $x$, мы получим после сложения:
\[180\degree\cdot n - 360\degree< 180\degree\cdot n - x.\]

Так как в разностях $180\degree\cdot n-360\degree$ и $180\degree\cdot n-x$, уменьшаемые одинаковы, то, чтобы первая разность была меньше второй, необходимо, чтобы вычитаемое $360\degree$ было больше вычитаемого $x$;
значит, $360\degree> x$, то есть $x < 360\degree$.

\paragraph{Симметричные многогранные углы.}\label{1938/s53}
Как известно, вертикальные углы равны, если речь идёт об углах, образованных прямыми или плоскостями.
Посмотрим, справедливо ли это утверждение применительно к углам многогранным.

\begin{wrapfigure}[15]{r}{50 mm}
\vskip-0mm
\centering
\includegraphics{mppics/s-ris-43}
\caption{}\label{1938/s-ris-43}
\end{wrapfigure}

Продолжим (рис.~\ref{1938/s-ris-43}) все рёбра угла $SABCDE$ за вершину $S$, тогда образуется другой многогранный угол $SA_1B_1D_1E_1$, который можно назвать \rindex{вертикальные углы}\textbf{вертикальным} по отношению к первому углу.
Нетрудно видеть, что у обоих углов равны соответственно и плоские углы, и двугранные, но те и другие расположены в \rindex{зеркальный порядок}\textbf{зеркальном порядке}.
Действительно, если мы вообразим наблюдателя, который смотрит извне многогранного угла на его вершину, то рёбра $SA$, $SB$, $SC$, $SD$, $SE$ будут казаться ему расположенными в направлении против движения часовой стрелки, тогда как смотря на угол $SA_1B_1C_1D_1E_1$ он видит рёбра $SA_1, SB_1,\dots$ расположенными по движению часовой стрелки;
то есть если смотреть на один из углов в зеркало, то порядок их рёбер покажется тем же.

Очевидно, что равные многогранные углы должны иметь соответственно равные плоские и двугранные углы, которые \rindex{одинаково расположенные грани}\textbf{одинаково расположены}, то есть расположены в том же порядке.
Многогранные углы с соответственно равными плоскими и двугранными углами, но расположенными в зеркальном порядке, вообще не могут совместиться при вложении;
значит, они не равны.
Такие углы называются центрально симметричными относительно вершины $S$.
Подробнее о симметрии фигур в пространстве будет сказано ниже.

\subsection*{Равенство трёхгранных углов}

\paragraph{}\label{1938/s52} В следующей теореме приведены два простейших признака равенства трёхгранных углов; пара более сложных признаков приведена в §~\ref{1914/402}.

\begin{wrapfigure}[10]{r}{50 mm}
\vskip-6mm
\centering
\includegraphics{mppics/s-ris-42}
\caption{}\label{1938/s-ris-42}
\end{wrapfigure}

\medskip

\mbox{\so{Теоремы}}.
\textbf{\emph{Трёхгранные углы равны, если они имеют:}}

1) \textbf{\emph{по равному двугранному углу, заключённому между двумя соответственно равными и одинаково расположенными плоскими углами;}} или

2) \textbf{\emph{по равному плоскому углу, заключённому между двумя соответственно равными и одинаково расположенными двугранными углами;}} 

1) Пусть $S$ и $S_1$ — два трёхгранных угла (рис.~\ref{1938/s-ris-42}), у которых $\angle ASB\z=\angle A_1S_1B_1$,
$\angle ASC= \angle A_1S_1C_1$ (и эти равные углы одинаково расположены) и двугранный угол $AS$ равен двугранному углу $A_1S_1$.
Вложим угол $S_1$ в угол $S$ так, чтобы у них совпали точки $S_1$ и $S$, прямые $S_1A_1$ и $SA$ и плоскости $A_1S_1B_1$ и $ASB$.
Тогда ребро $S_1B_1$ пойдёт по $SB$ (в силу равенства углов $A_1S_1B_1$ и $ASB$), плоскость $A_1S_1C_1$ пойдёт по $ASC$ (по равенству двугранных углов) и ребро $S_1C_1$ пойдёт по ребру $SC$ (в силу равенства углов $A_1S_1C_1$ и $ASC$).
Таким образом, трёхгранные углы совместятся всеми своими рёбрами, то есть они будут равны.

2) Второй признак, подобно первому, доказывается вложением.

\begin{wrapfigure}{r}{42 mm}
\vskip0mm
\centering
\includegraphics{mppics/s-ris-348}
\caption{}\label{1914/s-ris-348}
\end{wrapfigure}

{\small

\paragraph{Дополнительный угол.}\label{1914/399}
Из вершины $S$ (рис.~\ref{1914/s-ris-348}) трёхгранного угла $SABC$ восстановим к грани $ASB$ перпендикуляр $SC_1$ направляя его в ту сторону от этой грани, в которой расположено противоположное ребро $SC$.
Подобно этому проведём перпендикуляр $SA_1$ к грани $BSC$ и $SB_1$ к грани $ASC$.
Трёхгранный угол, у которого рёбрами служат полупрямые $SA_1$, $SB_1$ и $SC_1$, называется \so{дополнительным} для угла $SABC$.

Заметим, что \emph{если для угла $SABC$ дополнительным углом служит угол $SA_1B_1C_1$, то и наоборот: для угла $SA_1B_1C_1$ дополнительным углом будет $SABC$.}
Действительно, плоскость $SA_1B_1$, проходя через перпендикуляры к плоскостям $BSC$ и $ASC$, перпендикулярна к ним обеим, а следовательно, и к линии их пересечения $SC$; значит, ребро $SC$ есть перпендикуляр к грани $SA_1B_1$ и, кроме того, оно расположено по ту же сторону от этой грани, по которую лежит противоположное ребро $SC_1$.
Подобно этому убедимся, что рёбра $SB$ и $SA$ соответственно перпендикулярны к граням $SA_1C_1$ и $SB_1C_1$ и расположены по ту сторону от них, по которую лежат рёбра $SB_1$ и $SA_1$.
Значит, углы $SABC$ и $SA_1B_1C_1$ взаимно дополнительны.

\begin{wrapfigure}{r}{55 mm}
\vskip-4mm
\centering
\includegraphics{mppics/s-ris-349}
\caption{}\label{1914/s-ris-349}
\end{wrapfigure}

\paragraph{}\label{1914/400}
\mbox{\so{Лемма} 1.}
\textbf{\emph{Если два трёхгранные угла взаимно дополнительны, то плоские углы одного служат дополнением до $180\degree$ к противоположным двугранным углам другого.}}

Каждый плоский угол одного из взаимно дополнительных трёхгранных углов образован двумя перпендикулярами, восстановленными к граням противоположного двугранного угла другого трёхгранного, из одной точки его ребра.

Заметив это, возьмём какой-нибудь двугранный угол $AB$ (рис. \ref{1914/s-ris-349}) и из произвольной точки $B$ его ребра построим перпендикуляры: $BE$ к грани $AD$ и $BF$ к грани $AC$.
Затем через $BE$ и $BF$ построим плоскость, перпендикулярную к ребру $AB$ (§§~\ref{1938/s43}, \ref{1938/s45}).
Пусть пересечения этой плоскости с гранями угла $AB$ будут прямые $BC$ и $BD$.
Тогда угол $CBD$ должен быть линейным углом двугранного угла $AB$.
Заметим, что стороны угла $EBF$ соответственно перпендикулярны к сторонам угла $CBD$.
Приняв во внимание направления перпендикуляров, получим, что сумма углов $EBF$ и $CBD$ равна $180\degree$ (рис.~\ref{1914/s-ris-349}); что и требовалось доказать.

\paragraph{}\label{1914/401} 
\so{Лемма} 2. \textbf{\emph{Равным трёхгранным углам соответствуют равные дополнительные углы и обратно.}}

Равные трёхгранные углы при вложении совмещаются; поэтому совмещаются и те перпендикуляры, которые образуют рёбра дополнительных углов; значит, дополнительные углы также совмещаются.
Обратно: если совмещаются дополнительные углы, то совмещаются и данные углы.

\paragraph{}\label{1914/402}
Следующая теорема даёт два признака равенства трёхгранных углов в дополнение к признакам приведённым в §~\ref{1938/s52}.

\medskip

\so{Теоремы}.
\textbf{\emph{Трёхгранные углы равны, если они имеют:}}

1) \textbf{\emph{по три соответственно равных и одинаково расположенных плоских угла;}}

2) \textbf{\emph{по три соответственно равных и одинаково расположенных двугранных угла.}} 

\begin{wrapfigure}{o}{60 mm}
\vskip-0mm
\centering
\includegraphics{mppics/s-ris-351}
\caption{}\label{1914/s-ris-351}
\end{wrapfigure}

1) Пусть $S$ и $S_1$ (рис.~\ref{1914/s-ris-351}) два трёхгранные угла, у которых плоские углы одного равны соответственно плоским углам другого и, кроме того, равные углы одинаково расположены.

Отложим на всех рёбрах произвольные, но равные, отрезки 
\begin{align*}SA &= SB = SC =
\\
 =S_1A_1 &=S_1B_1 = S_1C_1. 
\end{align*}
и построим треугольники $ABC$ и $A_1B_1C_1$.
Из равенства треугольников $ABS$ и $A_1B_1S_1$ находим: $AB \z= A_1B_1$.
Подобно этому из равенства других боковых треугольников выводим: $AC \z= A_1C_1$ и $BC\z=B_1C_1$.
Следовательно, $\triangle ABC\z=\triangle A_1B_1C_1$.

Опустим на плоскости этих треугольников перпендикуляры $SO$ и $S_1O_1$.
Так как наклонные $SA$, $SB$ и $SC$ равны, то должны быть равны их проекции $OA$, $OB$ и $OC$ (§~\ref{1938/s25}).
Значит, точка $O$ есть центр круга, описанного около треугольника $ABC$.
Точно так же точка $O_1$ есть центр круга, описанного около треугольника $A_1B_1C_1$.
У равных треугольников радиусы описанных кругов равны; значит, $OB = O_1B_1$.
Поэтому $\triangle SBO=\triangle S_1B_1O_1$ (по гипотенузе и катету), и, следовательно, $OS = O_1S_1$.

Вложим теперь фигуру $S_1A_1B_1C_1$ в фигуру $SABC$ так, чтобы равные треугольники $A_1B_1C_1$ и $ABC$ совместились.
Тогда совместятся описанные окружности, и, следовательно, их центры $O_1$ и $O$.
Поскольку плоские углы одинаково расположены, перпендикуляр $O_1S_1$ пойдёт по $OS$ и точка $S_1$ совпадёт с $S$.
Таким образом, трёхгранные углы совместятся всеми своими рёбрами, значит, они равны.

2) Четвёртый признак легко доказывается при помощи дополнительных углов.
Если у двух трёхгранных углов соответственно равны и одинаково расположены двугранные углы, то у дополнительных углов соответственно равны и одинаково расположены плоские углы (§~\ref{1914/400});
следовательно, дополнительные углы равны; а если равны дополнительные, то равны и данные углы (§~\ref{1914/401}).

}

{\small

\subsection*{Упражнения}

\so{Доказать теоремы:}

\begin{enumerate}[noitemsep]

\item
Две плоскости, параллельные третьей, параллельны между собой.

\item
Все прямые, параллельные данной плоскости и проходящие через одну точку, лежат в одной плоскости, параллельной данной.

\item Дана плоскость $P$ и параллельная ей прямая $a$.
Доказать, что все точки прямой $a$ находятся на одинаковом расстоянии от плоскости~$P$.

\item
Доказать, что все точки одной из двух параллельных плоскостей находятся на одинаковом расстоянии от другой плоскости.

\item
Две плоскости, проходящие через две данные параллельные прямые и не параллельные между собой, пересекаются по прямой, параллельной данным прямым.

\item
Если прямая $a$ параллельна какой-либо прямой $b$, лежащей на плоскости $M$, то всякая плоскость, проходящая через $a$, пересекает плоскость $M$ по прямой, параллельной $b$ (возможно сливающейся с $b$).

\item
Если прямая $a$ параллельна плоскости $M$, то всякая прямая, проходящая через точку, лежащую в плоскости $M$, и параллельная прямой $a$, лежит в плоскости $M$.

\item
Если даны две скрещивающиеся прямые $a$ и $b$ и через первую проведена плоскость, параллельная второй, а через вторую — плоскость, параллельная первой, то эти две плоскости параллельны.

\item
Все прямые, проходящие через какую-нибудь точку на прямой $a$ и перпендикулярные к этой прямой, лежат в одной плоскости, перпендикулярной к $a$.

\item
Если плоскость и прямая перпендикулярны к одной прямой, то они параллельны.

\item
Если прямая $a$, параллельная плоскости $M$, пересекает прямую $b$, перпендикулярную этой плоскости, то прямые $a$ и $b$ перпендикулярны.
\end{enumerate}

\so{Задачи на построение}

\begin{enumerate}[resume,noitemsep]
\item
Через данную точку провести плоскость, параллельную двум данным прямым $a$ и $b$.

\item
Через данную точку провести прямую, параллельную данной плоскости и пересекающую данную прямую.

\item
Построить прямую, пересекающую две данные прямые и параллельную третьей данной прямой.

\item
Построить какую-либо прямую, пересекающую две данные прямые и параллельную данной плоскости (задача может иметь много решений).

\item
Построить какую-либо прямую, пересекающую три данные прямые (задача может иметь много решений).

\item
Через данную точку провести прямую, перпендикулярную двум данным скрещивающимся прямым.

\item
Через данную прямую провести плоскость, перпендикулярную к данной плоскости.

\item
Даны плоскость $M$ и прямая $a\parallel M$.
Через прямую $a$ провести плоскость, пересекающую плоскость $M$ под данным углом.

\item
Дана плоскость $M$ и две точки $A$ и $B$ по одну сторону от неё.
Найти на плоскости $M$ такую точку $C$, что сумма $AC$ + $BC$ была наименьшей.

\end{enumerate}

}

%% file: 3D/mnogogranniki.tex
\chapter{Многогранники}

\section{Параллелепипед и пирамида}

\paragraph{Многогранник.}\label{1938/s67}
\so{Многогранником} называется тело, ограниченное плоскими многоугольниками.
Общие стороны смежных многоугольников называются \rindex{ребро!многогранника}\textbf{рёбрами} многогранника.
Многоугольники, которые ограничивают многогранник, называются его \rindex{грань!многогранника}\textbf{гранями}.
Грани многогранника, сходящиеся в одной точке, образуют многогранный угол;
вершины таких многогранных углов называются \rindex{вершина!многогранника}\textbf{вершинами многогранника}.
Прямые, соединяющие две какие-нибудь вершины, не лежащие на одной грани, называются диагоналями многогранника.

Мы будем рассматривать только \so{выпуклые} многогранники, то есть такие, которые расположены по одну сторону от плоскости каждой из его граней.

Наименьшее число граней в многограннике — четыре;
такой многогранник получается от пересечения трёхгранного угла какой-нибудь плоскостью.

\paragraph{Призма.}\label{1938/s68}
\so{Призмой} называется многогранник, у которого две грани — равные многоугольники с соответственно параллельными сторонами, а все остальные грани — параллелограммы.

Чтобы убедиться в существовании такого многогранника, возьмём (рис.~\ref{1938/s-ris-73}) какой-нибудь многоугольник $ABCDE$ и через его вершины проведём ряд параллельных прямых, не лежащих в его плоскости.
Взяв затем на одной из этих прямых произвольную точку $A_1$, проведём через неё плоскость, параллельную плоскости $ABCDE$;
через каждые две соседние параллельные прямые также проведём по плоскости.
Эти плоскости образуют продолжение граней многогранника $ABCDEA_1B_1C_1D_1E_1$, удовлетворяющего определению призмы.
Действительно, параллельные плоскости $ABCDE$ и $A_1B_1C_1D_1E_1$ пересекаются боковыми плоскостями по параллельным прямым (§~\ref{1938/s16});
поэтому фигуры $AA_1E_1E$, $EE_1D_1D$ и так далее — параллелограммы.
С другой стороны, у многоугольников $ABCDE$ и $A_1B_1C_1D_1E_1$ равны соответственно стороны (как противоположные стороны параллелограммов) и углы (как углы с параллельными и одинаково направленными сторонами);
следовательно, эти многоугольники равны.

{\sloppy

Многоугольники $ABCDE$ и $A_1B_1C_1D_1E_1$, лежащие в параллельных плоскостях, называются \rindex{основание!призмы}\textbf{основаниями} призмы; перпендикуляр $OO_1$, опущенный из какой-нибудь точки одного основания на плоскость другого, называется \rindex{высота!призмы}\textbf{высотой} призмы.
Параллелограммы $AA_1B_1B$, $BB_1C_1C$ и так далее называются \rindex{боковая грань!призмы}\textbf{боковыми гранями} призмы, а их стороны $AA_1$, $BB_1$ и так далее, соединяющие соответственные вершины оснований, — \rindex{боковое ребро призмы}\textbf{боковыми рёбрами}.
У призмы все боковые рёбра равны как отрезки параллельных прямых, заключённые между параллельными плоскостями.

}

\begin{wrapfigure}{o}{35 mm}
\vskip-0mm
\centering
\includegraphics{mppics/s-ris-73}
\caption{}\label{1938/s-ris-73}
\end{wrapfigure}

Отрезок прямой, соединяющий какие-нибудь две вершины, не прилежащие к одной грани, называется диагональю призмы.
Таков, например, отрезок $AD_1$ (рис.~\ref{1938/s-ris-73}).

Плоскость, проведённая через какие-нибудь два боковых ребра, не прилежащих к одной боковой грани призмы (например, через рёбра $AA_1$ и $CC_1$, рис.~\ref{1938/s-ris-73}), называется диагональной плоскостью (на рисунке не показанной).

Призма называется \rindex{прямая призма}\textbf{прямой} или \rindex{наклонная призма}\textbf{наклонной}, смотря по тому, будут ли её боковые рёбра перпендикулярны или наклонны к основаниям.
У прямой призмы боковые грани — прямоугольники.
За высоту такой призмы можно принять боковое ребро.

Прямая призма называется \rindex{правильная призма}\textbf{правильной}, если её основания — правильные многоугольники.
У такой призмы все боковые грани — равные прямоугольники.

Призмы бывают треугольные, четырёхугольные и так далее, смотря по тому, что является основанием: треугольник, четырёхугольник и так далее.

\paragraph{Параллелепипед.}\label{1938/s69}\rindex{параллелепипед}
Параллелепипедом называют призму, у которой основаниями служат параллелограммы (рис.~\ref{1938/s-ris-74}).
Параллелепипеды, как и всякие призмы, могут быть прямые и наклонные.
Прямой параллелепипед называется \rindex{прямоугольный параллелепипед}\textbf{прямоугольным}, если его основание — прямоугольник (рис.~\ref{1938/s-ris-75}).
Из этих определений следует:

1) у параллелепипеда все шесть граней — параллелограммы;

2) у прямого параллелепипеда четыре боковые грани — прямоугольники, а два основания — параллелограммы;

\begin{wrapfigure}[39]{r}{35 mm}
\vskip-0mm
\centering
\includegraphics{mppics/s-ris-74}
\caption{}\label{1938/s-ris-74}
\bigskip
\includegraphics{mppics/s-ris-75}
\caption{}\label{1938/s-ris-75}
\bigskip
\includegraphics{mppics/s-ris-76}
\caption{}\label{1938/s-ris-76}
\bigskip
\includegraphics{mppics/s-ris-77}
\caption{}\label{1938/s-ris-77}
\end{wrapfigure}

3) у прямоугольного параллелепипеда все шесть граней — прямоугольники.

Три ребра прямоугольного параллелепипеда, сходящиеся к одной вершине, называются его \rindex{измерения!параллелепипеда}\textbf{измерениями};
одно из них можно рассматривать как длину, другое — как ширину, а третье — как высоту.

Прямоугольный параллелепипед, имеющий равные измерения, называется кубом.
У куба все грани — квадраты.

\paragraph{Пирамида.}\label{1938/s70}
Пирамидой называется многогранник, у которого одна грань, называемая основанием, есть какой-нибудь многоугольник, а все остальные грани, называемые боковыми, — треугольники, имеющие общую вершину.

Чтобы получить пирамиду, достаточно какой-нибудь многогранный угол $S$ (рис.~\ref{1938/s-ris-76}) пересечь произвольной плоскостью $ABCD$ и взять отсечённую часть $SABCD$.

Общая вершина $S$ боковых треугольников называется \rindex{вершина!пирамиды}\textbf{вершиной} пирамиды, а перпендикуляр $SO$, опущенный из вершины на плоскость основания, — \rindex{высота!пирамиды}\textbf{высотой}.

Обыкновенно, обозначая пирамиду буквами, пишут сначала ту, которой обозначена вершина, например $SABCD$ (рис.~\ref{1938/s-ris-76}).

Плоскость, проведённая через вершину пирамиды и какую-нибудь диагональ основания (например, через диагональ $BD$, рис.~\ref{1938/s-ris-78}), называется {}\textbf{диагональной плоскостью}.

Пирамиды бывают треугольные (рис.~\ref{1938/s-ris-77}), четырёхугольные и так далее, смотря по тому, что является основанием — треугольник, четырёхугольник и так далее.

Пирамида называется \rindex{правильная пирамида}\textbf{правильной} (рис. \ref{1938/s-ris-78}), если, во-первых, её основание есть правильный многоугольник и, во-вторых, высота проходит через центр этого многоугольника.
В правильной пирамиде все боковые рёбра равны между собой (как наклонные с равными проекциями).

Поэтому все боковые грани правильной пирамиды — равные равнобедренные треугольники.
Высота $SM$ (рис.~\ref{1938/s-ris-78}) каждого из этих треугольников называется \rindex{апофема!пирамиды}\textbf{апофемой}.
Все апофемы в правильной пирамиде равны.

\begin{figure}[!ht]
\begin{minipage}{.48\textwidth}
\centering
\includegraphics{mppics/s-ris-78}
\end{minipage}
\hfill
\begin{minipage}{.48\textwidth}
\centering
\includegraphics{mppics/s-ris-79}
\end{minipage}

\medskip

\begin{minipage}{.48\textwidth}
\centering
\caption{}\label{1938/s-ris-78}
\end{minipage}
\hfill
\begin{minipage}{.48\textwidth}
\centering
\caption{}\label{1938/s-ris-79}
\end{minipage}
\vskip-4mm
\end{figure}

\paragraph{Усечённая пирамида.}\label{1938/s71}
Часть пирамиды (рис.~\ref{1938/s-ris-79}), заключённая между основанием ($ABCDE$) и секущей плоскостью ($A_1B_1C_1D_1E_1$), параллельной основанию, называется усечённой пирамидой.
Параллельные грани называются основаниями, а отрезок перпендикуляра $OO_1$, опущенного из какой-нибудь точки $O_1$ основания $A_1B_1C_1D_1E_1$ на другое основание, — высотой усечённой пирамиды.
Усечённая пирамида называется правильной, если она составляет часть правильной пирамиды.

\subsection*{Свойства граней и диагоналей параллелепипеда}

\paragraph{}\label{1938/s72}
\mbox{\so{Теорема}.}
\textbf{\emph{В параллелепипеде:}}

1) \textbf{\emph{противоположные грани равны и параллельны;}}

2) \textbf{\emph{все четыре диагонали пересекаются в одной точке и делятся в ней пополам.}}

1) Грани (рис.~\ref{1938/s-ris-80}) $BB_1C_1C$ и $AA_1D_1D$ параллельны, потому что две пересекающиеся прямые $BB_1$ и $B_1C_1$ одной грани параллельны двум пересекающимся прямым $AA_1$ и $A_1D_1$ другой (§~\ref{1938/s15});
эти грани и равны, так как $B_1C_1=A_1D_1$, $B_1B=A_1A$ (как противоположные стороны параллелограммов) и $\angle BB_1C_1=\angle AA_1D_1$.

\begin{figure}[!ht]
\begin{minipage}{.32\textwidth}
\centering
\includegraphics{mppics/s-ris-80}
\end{minipage}\hfill
\begin{minipage}{.32\textwidth}
\centering
\includegraphics{mppics/s-ris-81}
\end{minipage}\hfill
\begin{minipage}{.32\textwidth}
\centering
\includegraphics{mppics/s-ris-82}
\end{minipage}

\medskip

\begin{minipage}{.32\textwidth}
\centering
\caption{}\label{1938/s-ris-80}
\end{minipage}\hfill
\begin{minipage}{.32\textwidth}
\centering
\caption{}\label{1938/s-ris-81}
\end{minipage}\hfill
\begin{minipage}{.32\textwidth}
\centering
\caption{}\label{1938/s-ris-82}
\end{minipage}
\vskip-4mm
\end{figure} 

2) Возьмём (рис.~\ref{1938/s-ris-81}) какие-нибудь две диагонали, например $AC_1$ и $BD_1$, и проведём вспомогательные прямые $AD_1$ и $BC_1$.
Так как рёбра $AB$ и $D_1C_1$ соответственно равны и параллельны ребру $DC$, то они равны и параллельны между собой;
вследствие этого фигура $AD_1C_1B$ есть параллелограмм, в котором прямые $C_1A$ и $BD_1$ — диагонали, а в параллелограмме диагонали делятся в точке пересечения пополам.

Возьмём теперь одну из этих диагоналей, например $AC_1$, с третьей диагональю, положим, с $B_1D$.
Совершенно так же мы можем доказать, что они делятся в точке пересечения пополам.
Следовательно, диагонали $B_1D$ и $AC_1$ и диагонали $AC_1$ и $BD_1$ (которые мы раньше брали) пересекаются в одной и той же точке, а именно в середине диагонали $AC_1$.
Наконец, взяв эту же диагональ $AC_1$ с четвёртой диагональю $A_1C$, мы также докажем, что они делятся пополам.
Значит, точка пересечения и этой пары диагоналей лежит в середине диагонали $AC_1$.
Таким образом, все четыре диагонали параллелепипеда пересекаются в одной и той же точке и делятся этой точкой пополам.

\paragraph{}\label{1938/s73}
\mbox{\so{Теорема}.}
\textbf{\emph{В прямоугольном параллелепипеде квадрат любой диагонали}} ($AC_1$, рис.~\ref{1938/s-ris-82}) \textbf{\emph{равен сумме квадратов трёх его измерений.}}

Проведя диагональ основания $AC$, получим треугольники $AC_1C$ и $ACB$.
Оба они прямоугольные: первый потому, что параллелепипед прямой и, следовательно, ребро $CC_1$ перпендикулярно к основанию;
второй потому, что параллелепипед прямоугольный и, значит, в основании его лежит прямоугольник.
Из этих треугольников находим:
\begin{align*}
AC_1^2 &= AC^2 + CC_1^2
&&
\text{и}
&
AC^2 &= AB^2 + BC^2.
\end{align*}
Следовательно,
\begin{align*}
AC_1^2 &= AB^2 + BC^2 + CC_1^2 = 
\\
&=AB^2 + AD^2 + AA_1^2.
\end{align*}

\medskip

\so{Следствие}.
\emph{В прямоугольном параллелепипеде все диагонали равны.}

\subsection*{Свойства параллельных сечений в пирамиде}

\begin{wrapfigure}{r}{40 mm}
\vskip-0mm
\centering
\includegraphics{mppics/s-ris-83}
\caption{}\label{1938/s-ris-83}
\vskip-0mm
\end{wrapfigure}

\paragraph{}\label{1938/s74}
\mbox{\so{Теоремы}.}
\textbf{\emph{Если пирамида}} (рис. \ref{1938/s-ris-83}) \textbf{\emph{пересечена плоскостью, параллельной основанию, то:}}

1) \textbf{\emph{боковые рёбра и высота делятся этой плоскостью на пропорциональные части;}}

2) \textbf{\emph{в сечении получается многоугольник}} ($abcde$), \textbf{\emph{подобный основанию}} ($ABCDE$);

3) \textbf{\emph{площади сечения и основания относятся, как квадраты их расстояний от вершины.}}

1) Прямые $ab$ и $AB$ можно рассматривать как линии пересечения двух параллельных плоскостей (основания и секущей) третьей плоскостью $ASB$;
поэтому $ab\parallel AB$ (§~\ref{1938/s16}).
По этой же причине $bc\parallel BC$, $cd\z\parallel CD\dots$ и $am\parallel AM$;
вследствие этого
\[\frac{Sa}{aA}=\frac{Sb}{bB}=\frac{Sc}{cC}=\dots=\frac{Sm}{mM}.\]

2) Из подобия треугольников $ASB$ и $aSb$, затем $BSC$ и $bSc$ и так далее выводим:
\[\frac{AB}{ab}=\frac{BS}{bS};\quad \frac{BS}{bS}=\frac{BC}{bc},\]
откуда
\[\frac{AB}{ab}=\frac{BC}{bc}.\]
Также
\[\frac{BC}{bc}=\frac{CS}{cS};\quad \frac{CS}{cS}=\frac{CD}{cd},\]
откуда
\[\frac{BC}{bc}=\frac{CD}{cd}.\]
Так же докажем пропорциональность остальных сторон многоугольников $ABCDE$ и $abcde$.
Так как, сверх того, у этих многоугольников равны соответственные углы (как образованные параллельными и одинаково направленными сторонами), то они подобны.
Площади подобных многоугольников относятся, как квадраты сходственных сторон (§~\ref{1938/260});
поэтому
\[\frac{\text{площадь\,}ABCDE}{\text{площадь\,} abcde}=\frac{AB^2}{ab^2}=\left(\frac{AB}{ab}\right)^2.\]
но
\[\frac{AB}{ab}=\frac{AS}{aS}=\frac{MS}{mS},\]
значит
\[\frac{\text{площадь\,}ABCDE}{\text{площадь\,} abcde}=\frac{MS^2}{mS^2}=\left(\frac{MS}{mS}\right)^2.\]

\paragraph{}\label{1938/s75}
\so{Следствие}.
\emph{У правильной усечённой пирамиды верхнее основание есть правильный многоугольник, подобный нижнему основанию, а боковые грани — равные и равнобочные трапеции} (рис.~\ref{1938/s-ris-83}).
Высота любой из этих трапеций называется \rindex{апофема!усечённой пирамиды}\textbf{апофемой} правильной усечённой пирамиды.

\paragraph{}\label{1938/s76}
\so{Теорема}.
\textbf{\emph{Если две пирамиды с равными высотами рассечены на одинаковом расстоянии от вершины плоскостями, параллельными основаниям, то площади сечений пропорциональны площадям оснований.}}

\begin{figure}[!ht]
\vskip-0mm
\centering
\includegraphics{mppics/s-ris-84}
\caption{}\label{1938/s-ris-84}
\vskip-0mm
\end{figure}

Пусть (рис.~\ref{1938/s-ris-84}) $B$ и $B_1$ — площади оснований двух пирамид, $H$ — высота каждой из них, $b$ и $b_1$ — площади сечений плоскостями, параллельными основаниям и удалёнными от вершин на одно и то же расстояние $h$.

Согласно предыдущей теореме мы будем иметь:
\[\frac{b}{B}=\frac{h^2}{H^2}\quad\text{и}\quad\frac{b_1}{B_1}=\frac{h^2}{H^2},\]
откуда
\[\frac{b}{B}=\frac{b_1}{B_1}\quad\text{и}\quad\frac{b}{b_1}=\frac{B}{B_1}.\]

\paragraph{}\label{1938/s77}
\so{Следствие}.
\emph{Если $B=B_1$, то и $b=b_1$, то есть если у двух пирамид с равными высотами основания равновелики, то равновелики и сечения, равноотстоящие от вершины.}

\subsection*{Боковая поверхность призмы и пирамиды}

\paragraph{}\label{1938/s78} \rindex{перпендикулярное сечение}\textbf{Перпендикулярным сечением} (рис.~\ref{1938/s-ris-85}) называется многоугольник $abcde$, получаемый от пересечения боковых граней призмы (или их продолжений за основания) с плоскостью, перпендикулярной к боковому ребру.

Заметим, что стороны этого многоугольника перпендикулярны к рёбрам (§§~\ref{1938/s31}, \ref{1938/s24}).

\medskip

\mbox{\so{Теорема}.}
\textbf{\emph{Площадь\footnote{В дальнейшем ради краткости термин «поверхность» может употребляться вместо «площадь поверхности».} боковой поверхности призмы равна произведению периметра перпендикулярного сечения на боковое ребро.}}

\begin{wrapfigure}[14]{r}{53 mm}
\vskip-2mm
\centering
\includegraphics{mppics/s-ris-85}
\caption{}\label{1938/s-ris-85}
\vskip-0mm
\end{wrapfigure}

Боковая поверхность призмы состоит из параллелограммов;
в каждом из них за основание можно взять боковое ребро, а за высоту — сторону перпендикулярного сечения.
Поэтому площадь боковой поверхности призмы равна:
\begin{align*}
AA_1\cdot ab
&+ BB_1\cdot bc+CC_1\cdot cd+
\\
&+DD_1\cdot de+EE_1\cdot ea=
\\
=(ab&+ bc+ cd+ de+ ea)\cdot AA_1.
\end{align*}

\paragraph{}\label{1938/s79}
\mbox{\so{Следствие}.}
\emph{Площадь боковой поверхности прямой призмы равна произведению периметра основания на высоту} потому, что в такой призме за перпендикулярное сечение можно взять само основание, а боковое ребро её равно высоте.

\paragraph{}\label{1938/s80}
\so{Теорема}.
\textbf{\emph{Площадь боковой поверхности правильной пирамиды равна произведению периметра основания на половину апофемы.}}

Пусть (рис.~\ref{1938/s-ris-86}) $SABCDE$ — правильная пирамида и $SM$ — её апофема.
Площадь боковой поверхности этой пирамиды есть сумма площадей равных равнобедренных треугольников.
Площадь одного из них, например $ASB$, равна $AB\cdot\tfrac12SM$.
Если всего треугольников $n$, то площадь боковой поверхности равна 
\[AB\cdot\tfrac12SM\cdot n= AB\cdot n\cdot\tfrac12SM,\]
где $AB\cdot n$ есть периметр основания, а $SM$ — апофема.

{

\begin{wrapfigure}{r}{47 mm}
\vskip-4mm
\centering
\includegraphics{mppics/s-ris-86}
\caption{}\label{1938/s-ris-86}
\vskip-0mm
\end{wrapfigure}

\paragraph{}\label{1938/s81}
\mbox{\so{Теорема}.}
\textbf{\emph{Площадь боковой поверхности правильной усечённой пирамиды равна произведению полусуммы периметров обоих оснований на апофему.}}

Боковая поверхность правильной усечённой пирамиды состоит из равных трапеций.
Площадь одной трапеции, например $AabB$ (рис.~\ref{1938/s-ris-86}), равна $\tfrac12(AB \z+ ab)\cdot Mm$.
Если число всех трапеций есть $n$, то площадь боковой поверхности равна:
\[\frac{AB+ab}{2}\cdot Mm\cdot n=\frac{AB\cdot n+ab\cdot n}{2}\cdot Mm\]
где $AB\cdot n$ и $ab\cdot n$ — периметры нижнего и верхнего оснований.

}

{\small

\subsection*{Упражнения}

\begin{enumerate}[noitemsep]

\item
Высота прямой призмы, основание которой есть правильный треугольник, равна 12 м, сторона основания 3 м.
Вычислить полную поверхность призмы.

\item
Площадь полной поверхности прямоугольного параллелепипеда равна 1714~м$^2$, а неравные стороны основания равны 25~м и 14~м.
Вычислить площадь боковой поверхности и боковое ребро.

\item
В прямоугольном параллелепипеде с квадратным основанием и высотой $h$ проведена секущая плоскость через два противоположных боковых ребра.
Вычислить площадь полной поверхности параллелепипеда, зная, что площадь сечения равна $S$.

\item
Правильная шестиугольная пирамида имеет сторону основания $a$ и высоту $h$.
Вычислить боковое ребро, апофему, площади боковой поверхности и полной поверхности.

\item
Вычислить площадь полной поверхности и высоту треугольной пирамиды, у которой каждое ребро равно $a$.

\item
Правильная шестиугольная пирамида, у которой высота 25~см, а сторона основания 5~см, рассечена плоскостью, параллельной основанию.
Вычислить расстояние этой плоскости от вершины пирамиды, зная, что площадь сечения равна $\tfrac23\sqrt{3}$~см$^2$.

\item
Высота усечённой пирамиды с квадратным основанием равна $h$, сторона нижнего основания $a$, а верхнего $b$.
Найти площадь полной поверхности усечённой пирамиды.

\item
Высота усечённой пирамиды равна 6, а площади оснований 18 и~8.
Пирамида рассечена плоскостью, параллельной основаниям и делящей высоту пополам.
Вычислить площадь сечения.
\end{enumerate}

}

\section{Объём призмы и пирамиды}

\paragraph{Основные допущения в объёмах.}\label{1938/s82}
Величина части пространства, занимаемого геометрическим телом, называется \rindex{объём}\textbf{объёмом} этого тела.

Мы ставим задачу — выразить эту величину положительным числом.
При этом мы будем руководствоваться следующими исходными положениями.

1) \emph{Равные тела имеют равные объёмы.}

2) \emph{Объём какого-нибудь тела}
(например, каждого параллелепипеда, изображённого на рис.~\ref{1938/s-ris-87}),
\emph{состоящего из частей}
($P$ и $Q$),
\emph{равен сумме объёмов этих частей.}

\begin{wrapfigure}{o}{50 mm}
\vskip-0mm
\centering
\includegraphics{mppics/s-ris-87}
\caption{}\label{1938/s-ris-87}
\vskip-0mm
\end{wrapfigure}

Так как объём измеряется положительным числом, получаем, что 
объём тела (например, одного из тел состоящего из частей $P$ и $Q$ на рис.~\ref{1938/s-ris-87})
больше объёма любой его части ($P$ и $Q$).
Это свойство можно выразить как \emph{объём части меньше объёма целого}.

Два тела, имеющие одинаковые объёмы, называются \rindex{равновеликие тела}\textbf{равновеликими}.
Заметим, что равновеликие тела (например, два тела на рис.~\ref{1938/s-ris-87}) могут быть неравны, то есть одно из них невозможно совместить с другим.

\paragraph{Единица объёма.}\label{1938/s83}
За единицу объёма берут объём такого куба, у которого каждое ребро равно линейной единице.
Так, употребительны кубические метры (м$^3$), кубические сантиметры (см$^3$) и так далее.

\begin{wrapfigure}{o}{40 mm}
\vskip-0mm
\centering
\includegraphics{mppics/s-ris-91}
\caption{}\label{1938/s-ris-91}
\vskip-0mm
\end{wrapfigure}

Отношение двух кубических единиц разных названий равно третьей степени отношения тех линейных единиц, которые служат рёбрами для этих кубических единиц.
Так, отношение кубического метра к кубическому дециметру равно $10^3$, то есть $1000$.
Поэтому, например, если мы имеем куб с ребром длиной $a$ линейных единиц и другой куб с ребром длиной $3a$ линейных единиц, то отношение их объёмов будет равно $3^3$, то есть $27$, что ясно видно из рис.~\ref{1938/s-ris-91}.

{\small

\paragraph{Замечание о числе, измеряющем объём.}\label{1930/366}
Относительно числа, измеряющего данный объём в кубических единицах, можно сделать разъяснение, аналогичное тому, какое было нами приведено в §~\ref{1938/244}
относительно числа, измеряющего данную площадь в квадратных единицах.
Повторим вкратце это разъяснение в применении к объёмам.%
\footnote{Этот материал подробно изложен в статье «Площадь и объём» В. А. Рохлина (Энциклопедия элементарной математики, книга пятая, Геометрия).}

\begin{wrapfigure}{o}{40 mm}
\vskip-0mm
\centering
\includegraphics{mppics/s-ris-392}
\caption{}\label{1914/s-ris-392}
\vskip-0mm
\end{wrapfigure}

Возьмём три взаимно перпендикулярные прямые (рис. \ref{1914/s-ris-392}) $OA$, $OB$ и $OC$ и через каждые две из них проведём плоскость.
Мы получим тогда три взаимно перпендикулярные плоскости $AOC$, $COB$ и $BOA$.
Вообразим теперь три ряда параллельных плоскостей: ряд плоскостей, параллельных плоскости $AOC$, другой ряд плоскостей, параллельных плоскости $BOA$, и третий ряд плоскостей, параллельных плоскости $BOC$.
Допустим, кроме того, что соседние плоскости каждого ряда отстоят одна от другой на одно и то же расстояние, равное какой-нибудь $\tfrac1k$ доле линейной единицы.
Тогда от взаимного пересечения этих трёх рядов плоскостей образуется пространственная сеть кубов, из которых каждый представляет собой $(\tfrac1k)^3$ часть кубической единицы.

Вообразим, что в эту сеть мы поместили то тело, объём которого желаем измерить.
Тогда все кубы сети мы можем подразделить на три рода:
1) кубы, которые расположены полностью внутри тела, 
2) кубы, которые некоторой частью выступают из тела (которые, другими словами, пересекаются поверхностью тела), и  
3) кубы, расположенные полностью вне тела.
Если кубов 1-го рода будет $m$, а 2-го рода $n$, то объём данного тела не может быть меньше $\tfrac{m}{k^3}$ и не может быть больше $\tfrac{m+n}{k^3}$ кубических единиц.
То есть, эти два числа будут приближённые меры данного объёма с точностью до $\tfrac{n}{k^3}$ кубических единиц, первое число с недостатком, второе — с избытком.
Уменьшая всё более и более расстояние между параллельными плоскостями, мы будем заполнять пространство всё меньшими и меньшими кубами и можем получать приближённые результаты измерения всё с большей и большей точностью.
Если точность измерения $\tfrac{n}{k^3}$ стремится к нулю при неограниченном возрастании $k$,
то число $V$ равное общему пределу приближённых мер (с недостатком и с избытком) принимается за точную меру данного объёма. 

Доказано, что такое число $V$ существует для всякого многогранника и что оно не зависит от выбора тех трёх прямых $OA$, $OB$ и $OC$ (рис.~\ref{1914/s-ris-392}), которые были взяты для построения пространственной сети кубов.
(То же выполняется для тел, ограниченных кривыми поверхностями — цилиндров, конусов, шаров, шаровых сегментов и шаровых слоёв, о которых речь пойдёт далее.)
Более того, число это обладает свойствами, указанными в §~\ref{1938/s82}.

}

\subsection*{Объём параллелепипеда}

\paragraph{}\label{1938/s84}
\so{Теорема}.
\textbf{\emph{Объём прямоугольного параллелепипеда равен произведению трёх его измерений.}}

\begin{wrapfigure}{o}{35 mm}
\vskip-0mm
\centering
\includegraphics{mppics/s-ris-88}
\caption{}\label{1938/s-ris-88}
\vskip-0mm
\end{wrapfigure}

В таком кратком выражении теорему эту надо понимать так: число, выражающее объём прямоугольного параллелепипеда в кубической единице, равно произведению чисел, выражающих три его измерения в соответствующей линейной единице, то есть в единице, являющейся ребром куба, объём которого принят за кубическую единицу.
Так, если $x$ есть число, выражающее объём прямоугольного параллелепипеда в кубических сантиметрах, и $a$, $b$ и $c$ — числа, выражающие три его измерения в линейных сантиметрах, то теорема утверждает, что $x=abc$.
При доказательстве рассмотрим особо следующие три случая:

1) Измерения выражаются \so{целыми числами}.

Пусть, например, измерения, будут (рис.~\ref{1938/s-ris-88}) $AB=a$, $BC=b$ и $BD=c$, где $a$, $b$ и $c$ — какие-нибудь целые числа (например, как изображено у нас на рисунке: $a=4$, $b=2$ и $c=5$).
Тогда основание параллелепипеда содержит $ab$ таких квадратов, из которых каждый представляет собой соответствующую квадратную единицу.
На каждом из этих квадратов, очевидно, можно поместить по одной кубической единице.
Тогда получится слой (изображённый на рисунке), состоящий из $ab$ кубических единиц.
Так как высота этого слоя равна одной линейной единице, а высота всего параллелепипеда содержит $c$ таких единиц, то внутри параллелепипеда можно поместить $c$ таких слоёв.
Следовательно, объём этого параллелепипеда равен $abc$ кубических единиц.

2) Измерения выражаются \so{дробными числами}.

Пусть измерения параллелепипеда будут:
\[\frac mn, \frac pq, \frac rs\]
(некоторые из этих дробей могут равняться целому числу).
Приведя дроби к одинаковому знаменателю, будем иметь:
\[\frac {mqs}{nqs}, \frac {pns}{qns}, \frac {rnq}{snq}.\]

Примем $\frac 1{nqs}$ долю линейной единицы за новую (вспомогательную) единицу длины.
Тогда в этой новой единице измерения данного параллелепипеда выразятся целыми числами, а именно: $mqs$, $pns$ и $rnq$
и потому по доказанному (в случае 1) объём параллелепипеда равен произведению
\[(mqs)\cdot (pns)\cdot (rnq),\]
 если измерять этот объём новой кубической единицей, соответствующей новой линейной единице.
Таких кубических единиц в одной кубической единице, соответствующей прежней линейной единице, содержится $(nqs)^3$; значит, новая кубическая единица составляет $\tfrac1{(nqs)^3}$ прежней.
Поэтому объём параллелепипеда, выраженный в прежних единицах равен
\begin{align*}\frac1{(nqs)^3}\cdot(mqs)\cdot (pns)\cdot (rnq)&=\frac{mqs}{nqs}\cdot \frac{pns}{nqs}\cdot \frac{rnq}{nqs}= 
\\&=\frac mn\cdot \frac pq\cdot \frac rs.
\end{align*}

3) Измерения выражаются \so{иррациональными числами}.

Пусть у данного параллелепипеда (рис.~\ref{1938/s-ris-89}), который для краткости мы обозначим одной буквой $Q$, измерения будут:
\[AB=\alpha,\quad AC=\beta,\quad AD=\gamma,\]
где все числа $\alpha$, $\beta$ и $\gamma$ или только некоторые из них иррациональные.

\begin{wrapfigure}{o}{60 mm}
\vskip-0mm
\centering
\includegraphics{mppics/s-ris-89}
\caption{}\label{1938/s-ris-89}
\vskip-0mm
\end{wrapfigure}

Каждое из чисел $\alpha$, $\beta$ и $\gamma$ может быть представлено в виде бесконечной десятичной дроби.
Возьмём приближённые значения этих дробей с $n$ десятичными знаками сначала с недостатком, а затем с избытком.
Значения с недостатком обозначим $\alpha_n$, $\beta_n$, $\gamma_n$, значения с избытком $\alpha_n'$, $\beta_n'$, $\gamma_n'$.
Отложим на ребре $AB$, начиная от точки $A$, два отрезка $AB_1 = \alpha_n$ и $AB_2=\alpha_n'$.
На ребре $AC$ от той же точки $A$ отложим отрезки $AC_1=\beta_n$ и $AC_2=\beta_n'$ и на ребре $AD$ от той же точки — отрезки $AD_1=\gamma_n$ и $AD_2=\gamma_n'$.
При этом мы будем иметь
\begin{align*}
AB_1&\le AB<AB_2;
\\
AC_1&\le AC<AC_2;
\\ 
AD_1&\le AD<AD_2.
\end{align*}

Построим теперь два вспомогательных параллелепипеда: один (обозначим его $Q_1$) с измерениями $AB_1$, $AC_1$ и $AD_1$ и другой (обозначим его $Q_2$) с измерениями $AB_2$, $AC_2$ и $AD_2$.
Параллелепипед $Q_1$ будет весь помещаться внутри параллелепипеда $Q$, а параллелепипед $Q_2$ будет содержать внутри себя параллелепипед $Q$.

По доказанному (в случае 2) будем иметь:
\[\text{объём}\, Q_1 = \alpha_n\beta_n\gamma_n, \]
\[\text{объём}\, Q_2 = \alpha_n'\beta_n'\gamma_n'.\]

Поскольку $Q_1$ лежит в $Q$, а $Q$ лежит в $Q_2$, выполняется следующее двойное неравенство
\[\text{объём}\, Q_1 < \text{объём}\, Q <\text{объём}\, Q_2,\]
или 
\[\alpha_n\beta_n\gamma_n < \text{объём}\, Q <\alpha_n'\beta_n'\gamma_n'.\]
Это двойное неравенство остаётся верным при всякой степени точности, с которою мы находим приближённые значения чисел $\alpha$, $\beta$ и $\gamma$.
Значит, неравенство это мы можем высказать так: \emph{число, измеряющее объём данного параллелепипеда, должно быть больше произведения любых приближенных значений чисел $\alpha$, $\beta$ и $\gamma$, если эти значения взяты с недостатком, но меньше произведения любых приближенных значений тех же чисел, если эти значения взяты с избытком.}
Такое число, как известно из алгебры, равно произведению $\alpha\beta\gamma$.
Значит, и в этом случае
\[\text{объём}\, Q=\alpha\beta\gamma.\]

\paragraph{}\label{1938/s85}
\so{Следствие}.
Пусть измерения прямоугольного параллелепипеда, служащие сторонами его основания, выражаются числами $a$ и $b$, а третье измерение (высота) — числом $c$.
Тогда, обозначая объём его в соответствующих кубических единицах буквой $V$, можем написать
\[V = abc.\]
Так как произведение $ab$ выражает площадь основания, то можно сказать, что \emph{объём прямоугольного параллелепипеда равен произведению площади основания на высоту.}

\paragraph{}\label{1938/s86}
\so{Лемма}.
\textbf{\emph{Наклонная призма равновелика такой прямой призме, основание которой равно перпендикулярному сечению наклонной призмы, а высота — её боковому ребру.}}

Пусть дана наклонная призма $ABCDEA_1B_1C_1D_1E_1$ (рис.~\ref{1938/s-ris-92}).
Продолжим все её боковые рёбра и боковые грани в одном направлении.

Возьмём на продолжении одного какого-нибудь ребра произвольную точку $a$ и проведём через неё перпендикулярное сечение $abcde$.
Затем, отложив $aa_1=AA_1$, проведём через $a_1$ перпендикулярное сечение $a_1b_1c_1d_1e_1$.
Так как плоскости обоих сечений параллельны, то 
\[bb_1=cc_1=dd_1=ee_1=aa_1=AA_1\]
(§~\ref{1938/s17}).
Вследствие этого многогранник $abcdea_1b_1c_1d_1e_1$, у которого за основания приняты проведённые нами сечения, есть прямая призма, о которой говорится в теореме.

\begin{wrapfigure}[18]{o}{45 mm}
\vskip-0mm
\centering
\includegraphics{mppics/s-ris-92}
\caption{}\label{1938/s-ris-92}
\vskip-0mm
\end{wrapfigure}

Докажем, что данная наклонная призма равновелика этой прямой.
Для этого предварительно убедимся, что многогранники $aD$ и $a_1D_1$ равны.
Основания их $abcde$ и $a_1b_1c_1d_1e_1$ равны как основание призмы $ad_1$;
с другой стороны, прибавив к обеим частям равенства 
\[A_1A=a_1a\] по одному и тому же отрезку прямой $A_1a$, получим 
\[aA=a_1A_1.\] 
Подобно этому получим $bB=b_1B_1$, $cC\z=c_1C_1$ и так далее.
Вообразим теперь, что многогранник $aD$ вложен в многогранник $a_1D_1$ так, что основания их совпали;
тогда боковые рёбра, будучи перпендикулярны к основаниям и соответственно равны, также совпадут;
поэтому многогранник $aD$ совместится с многогранником $a_1D_1$;
значит, эти тела равны.

Теперь заметим, что если к прямой призме $a_1d$ добавим многогранник $aD$,
a к наклонной призме $A_1D$ добавим многогранник $a_1D_1$, равный $aD$, то получим один и тот же многогранник $a_1D$.
Из этого следует, что две призмы $A_1D$ и $a_1d$ равновелики.

\paragraph{}\label{1938/s87}
\so{Теорема}.
\textbf{\emph{Объём, параллелепипеда равен произведению площади основания на высоту.}}

Ранее мы доказали эту теорему для параллелепипеда прямоугольного, теперь докажем её для параллелепипеда прямого, а потом и наклонного.

\begin{wrapfigure}{O}{60 mm}
\vskip-0mm
\centering
\includegraphics{mppics/s-ris-93}
\caption{}\label{1938/s-ris-93}
\vskip-0mm
\end{wrapfigure}

1) Пусть (рис.~\ref{1938/s-ris-93}) $AC_1$ — прямой параллелепипед, то есть такой, у которого основание $ABCD$ — какой-нибудь параллелограмм, а все боковые грани — прямоугольники.
Возьмём в нём за основание боковую грань $AA_1B_1B$;
тогда параллелепипед будет наклонный.
Рассматривая его как частный случай наклонной призмы, мы на основании леммы предыдущего параграфа можем утверждать, что этот параллелепипед равновелик такому прямому параллелепипеду, у которого основание есть перпендикулярное сечение $MNPQ$, а высота $BC$.

Четырёхугольник $MNPQ$ — прямоугольник, потому что его углы служат линейными углами прямых двугранных углов;
поэтому прямой параллелепипед, имеющий основанием прямоугольник $MNPQ$, должен быть прямоугольным и, следовательно, его объём равен произведению трёх его измерений, за которые можно принять отрезки $MN$, $MQ$ и $BC$.
Таким образом,
\[\text{объём}\, AC_1 = MN\cdot MQ\cdot BC = MN\cdot (MQ\cdot BC).\]
Но произведение $MQ\cdot BC$ выражает площадь параллелограмма $ABCD$, поэтому
\[
\text{объём}\, AC_1
= (\text{площади}\, ABCD)\cdot MN
= (\text{площади}\, ABCD)\cdot BB_1.
\]

\begin{figure}[!ht]
\centering
\includegraphics{mppics/s-ris-94}
\caption{}\label{1938/s-ris-94}
\vskip-0mm
\end{figure}

2) Пусть (рис.~\ref{1938/s-ris-94}) $AC_1$ — наклонный параллелепипед.
Он равновелик такому прямому, у которого основанием служит перпендикулярное сечение $MNPQ$ (то есть перпендикулярное к рёбрам $AD$, $BC,\dots$), а высотой — ребро $BC$.
Но, по доказанному, объём прямого параллелепипеда равен произведению площади основания на высоту;
значит,
\[\text{объём}\, AC_1 = (\text{площади}\, MNPQ)\cdot BC.\]

Если $RS$ есть высота сечения $MNPQ$, то площадь $MNPQ$ равна $MQ\cdot RS$, поэтому
\[\text{объём}\, AC_1 = MQ\cdot RS\cdot BC = (BC\cdot MQ)\cdot RS.\]
Произведение $BC\cdot MQ$ выражает площадь параллелограмма $ABCD$;
следовательно, 
\[\text{объём}\, AC_1 =(\text{площади}\, ABCD) \cdot RS.\]

Остаётся доказать, что отрезок $RS$ представляет собой высоту параллелепипеда.
Действительно, сечение $MNPQ$, будучи перпендикулярно к рёбрам $BC$, $B_1C_1,\dots$, должно быть перпендикулярно к граням $ABCD$, $BB_1C_1C,\dots$, проходящим через эти рёбра (§~\ref{1938/s43}).
Поэтому если мы из точки $S$ восстановим перпендикуляр к плоскости $ABCD$, то он должен лежать весь в плоскости $MNPQ$ (§~\ref{1938/s44}) и, следовательно, должен слиться с прямой $SR$, лежащей в этой плоскости и перпендикулярной к $MQ$.
Значит, отрезок $SR$ есть высота параллелепипеда.
Таким образом, объём и наклонного параллелепипеда равен произведению площади основания на высоту.

\medskip

\so{Следствие}.
Если $V$, $B$ и $H$ — числа, выражающие в соответствующих единицах объём, площадь основания и высоту параллелепипеда, то можно написать:
\[V = B\cdot H.\]

\subsection*{Объём призмы}

\paragraph{}\label{1938/s88}
\so{Теорема}.
\textbf{\emph{Объём, призмы равен произведению площади основания на высоту.}}

Сначала докажем эту теорему для треугольной призмы, а потом и для многоугольной.

{\sloppy

1) Проведём через ребро $AA_1$ треугольной призмы $ABCA_1B_1C_1$ (рис.~\ref{1938/s-ris-95}) плоскость, параллельную грани $BB_1C_1C$, а через ребро $CC_1$ — плоскость, параллельную грани $AA_1B_1B$;
затем продолжим плоскости обоих оснований призмы до пересечения с проведёнными плоскостями.
Тогда мы получим параллелепипед $BD_1$, который диагональной плоскостью $AA_1C_1C$ делится на две треугольные призмы (одна из них данная).

}

\begin{wrapfigure}{o}{45 mm}
\vskip-4mm
\centering
\includegraphics{mppics/s-ris-95}
\caption{}\label{1938/s-ris-95}
\vskip-0mm
\end{wrapfigure}

Докажем, что эти призмы равновелики.
Для этого проведём перпендикулярное сечение $abcd$.
В сечении получится параллелограмм, который диагональю $ac$ делится на два равных треугольника.
Данная призма равновелика такой прямой призме, у которой основание есть $\triangle abc$, а высота — ребро $AA_1$ (§~\ref{1938/s86}).
Другая треугольная призма равновелика такой прямой, у которой основание есть $\triangle adc$, а высота — ребро $AA_1$.
Но две прямые призмы с равными основаниями и равными высотами равны (потому что при вложении они совмещаются), значит, призмы $ABCA_1B_1C_1$ и $ADCA_1D_1C_1$ равновелики.
Из этого следует, что объём данной призмы составляет половину объёма параллелепипеда $BD_1$;
поэтому, обозначив высоту призмы через $H$, получим:
\begin{align*}
\text{объём треугольной призмы}
&=
\frac{(\text{площади}\,ABCD)\cdot H}2=
\\
&=
\frac{(\text{площади}\,ABCD)}2\cdot H=
\\
&=
(\text{площади}\,ABC)\cdot H.
\end{align*}

{

\begin{wrapfigure}[12]{o}{30 mm}
\vskip-6mm
\centering
\includegraphics{mppics/s-ris-96}
\caption{}\label{1938/s-ris-96}
\vskip-0mm
\end{wrapfigure}

2) Проведём через ребро $AA_1$ многоугольной призмы (рис.~\ref{1938/s-ris-96}) диагональные плоскости $AA_1C_1C$ и $AA_1D_1D$.
Тогда данная призма рассечётся на несколько треугольных призм.
Сумма объёмов этих призм составляет искомый объём.
Если обозначим площади их оснований через $b_1,b_2,b_3$, а общую высоту через $H$, то получим: 
\begin{align*}
&\text{объём многоугольной призмы}=
\\
&= b_1\cdot H + b_2 \cdot H+ b_3\cdot H=
\\
&=(b_1+ b_2+b_3)\cdot H = 
\\
&=(\text{площади}\, ABCDE)\cdot H.
\end{align*}

}

\mbox{\so{Следствие}.}
Если $V$, $B$ и $H$ будут числа, выражающие в соответствующих единицах объём, площадь основания и высоту призмы, то, по доказанному, можно написать:
\[V = B\cdot H.\]

\subsection*{Объём пирамиды}

{\small

\paragraph{}\label{1938/s90}
\so{Лемма}.
\textbf{\emph{Треугольные пирамиды с равновеликими основаниями и равными высотами равновелики.}}

Доказательство наше будет состоять из трёх частей.
В первой части мы докажем равновеликость не самих пирамид, а вспомогательных тел, составленных из ряда треугольных призм, поставленных друг на друга.
Во второй части мы докажем, что объёмы этих вспомогательных тел при увеличении числа составляющих их призм приближаются к объёмам пирамид как угодно близко.
Наконец, в третьей части мы убедимся, что сами пирамиды должны быть равновелики.

\begin{figure}[!ht]
\vskip-0mm
\centering
\includegraphics{mppics/s-ris-99}
\caption{}\label{1938/s-ris-99}
\vskip-0mm
\end{figure}

I.
Вообразим, что пирамиды поставлены основаниями на некоторую плоскость (как изображено на рис.~\ref{1938/s-ris-99}), тогда их вершины будут находиться на одной прямой, параллельной плоскости оснований, и высота пирамид может быть изображена одним и тем же отрезком прямой $H$.
Разделим эту высоту на какое-нибудь целое число $n$ равных частей (например, на 4, как это указано на рис.~\ref{1938/s-ris-99}) и через точки деления проведём ряд плоскостей, параллельных плоскости оснований.
Плоскости эти, пересекаясь с пирамидами, дают в сечениях ряд треугольников, причём треугольники пирамиды $S$ будут равновелики соответствующим треугольникам пирамиды $S_1$ (§~\ref{1938/s77}).
Поставим внутри каждой пирамиды ряд таких призм, чтобы верхними основаниями у них были треугольники сечений, боковые рёбра были параллельны ребру $SA$ в одной пирамиде и ребру $S_1A_1$ в другой, а высота каждой призмы равнялась бы $H/n$.
Таких призм в каждой пирамиде окажется $n-1$;
они образуют некоторое ступенчатое тело, объём которого, очевидно, меньше объёма той пирамиды, в которой призмы построены.
Обозначим объёмы призм пирамиды $S$ по порядку, начиная от вершины буквами $p_1,p_2,p_3,\dots,p_{n-1}$, а объёмы призм пирамиды $S_1$ — также по порядку от вершины буквами $q_1, q_2, q_3, \dots ,q_{n-1}$.
Тогда, принимая во внимание, что у каждой пары соответствующих призм (у $p_1$ и $q_1$, у $p_2$ и $q_2$ и так далее) основания равновелики и высоты равны, мы можем написать ряд равенств на их объёмы: 
\[p_1=q_1,
\quad
p_2=q_2,
\quad
p_3=q_3,
\quad \dots,\quad p_{n-1}=q_{n-1}.
\]
Сложив все равенства почленно, получим:
\[p_1+p_2+p_3+\dots+p_{n-1}=q_1+q_2+q_3+\dots+q_{n-1}\eqno(1)\]
Мы доказали, таким образом, что объёмы построенных нами вспомогательных ступенчатых тел равны между собой (при всяком числе $n$, на которое мы делим высоту $H$).

II. Обозначив объёмы пирамид $S$ и $S_1$ соответственно буквами $V$
и $V_1$, положим, что
\[V - (p_1+p_2+p_3 + \dots+ p_{n-1}) = x\] 
и
\[V_1 - (q_1+q_2+q_3 + \dots+ q_{n-1}) = y,\]
откуда
\[p_1+p_2+p_3 + \dots+ p_{n-1} = V - x\] 
и
\[q_1+q_2+q_3 + \dots+ q_{n-1} = V_1 - y,\]

Тогда равенство (1) мы можем записать так:
\[V-x=V_1-y.\eqno(2)\]

Предположим теперь, что число $n$ равных частей, на которое мы делим высоту $H$, неограниченно возрастает;
например, предположим, что, вместо того чтобы делить высоту на 4 равные части, мы разделим её на 8 равных частей, потом на 16, на 32 и так далее, и пусть каждый раз мы строим указанным образом ступенчатые тела в обеих пирамидах.
Как бы ни возросло число призм, составляющих ступенчатые тела, равенство (1), а следовательно, и равенство (2) остаются в силе.
При этом объёмы $V$ и $V_1$, конечно, не будут изменяться, тогда как величины $x$ и $y$, показывающие, на сколько объёмы пирамид превосходят объёмы соответствующих ступенчатых тел, будут, очевидно, всё более и более уменьшаться.
Докажем, что величины $x$ и $y$ могут сделаться как угодно малы (другими словами, что они стремятся к нулю).
Это достаточно доказать для какой-нибудь одной из двух величин $x$ и $y$, например для~$x$.

\begin{wrapfigure}{o}{43 mm}
\vskip-0mm
\centering
\includegraphics{mppics/s-ris-100}
\caption{}\label{1938/s-ris-100}
\vskip-0mm
\end{wrapfigure}

С этой целью построим для пирамиды $S$ (рис.~\ref{1938/s-ris-100}) ещё другой ряд призм, который составит тоже ступенчатое тело, но по объёму большее пирамиды.
Призмы эти мы построим так же, как строили внутренние призмы, с той только разницей, что треугольники сечении мы теперь примем не за верхние основания призм, а за нижние.
Вследствие этого мы получим ряд призм, которые некоторой своей частью будут выступать из пирамид наружу, и потому они образуют новое ступенчатое тело с объёмом, б\'{о}льшим, чем объём пирамиды.
Таких призм будет теперь не $n-1$, как внутренних призм, а~$n$.
Обозначим их объёмы по порядку, начиная от вершины, буквами: $p_1',p_2',p_3',\dots,p_n'$.
Рассматривая чертёж, мы легко заметим, что
\[p_1'=p_1,\quad p_2'=p_2,\quad p_3'=p_3,\quad\dots, \quad p_{n-1}'=p_{n-1}.\]
Поэтому
\[(p_1'+p_2'+p_3'+\dots+p_{n-1}'+p_n')-(p_1+p_2+p_3+\dots+p_{n-1})=p_n'\]
Так как
\[p_1'+p_2'+p_3'+\dots+p_{n-1}'+p_n'>V,\]
а
\[p_1+p_2+p_3+\dots+p_{n-1}<V,\]
то
\[V-(p_1+p_2+p_3+\dots+p_{n-1})<p_n'\]
то есть
\[x<p_n'.\]

Но 
\[p_n'=\text{площади}\, ABC\cdot \frac Hn\]
(если $ABC$ есть основание); поэтому
\[x<(\text{площади}\, ABC)\cdot \frac Hn.\]

При неограниченном возрастании числа $n$ величина $\frac Hn$, очевидно, может быть сделана как угодно малой (стремится к нулю).
Поэтому и произведение ($\text{площадь}\, ABC\cdot \frac Hn$), в котором множимое не изменяется, а множитель стремится к нулю, тоже стремится к нулю, и так как положительная величина $x$ меньше этого произведения, то она и подавно стремится к нулю.

То же самое рассуждение можно было бы повторить и о величине~$y$.

Мы доказали, таким образом, что при неограниченном увеличении числа призм объёмы вспомогательных ступенчатых тел приближаются к объёмам соответствующих пирамид как угодно близко.

III.
Заметив это, возьмём написанное выше равенство (2) и придадим ему такой вид: 
\[V - V_1 = x - y.\]
Докажем теперь, что это равенство возможно только тогда, когда $V\z=V_1$ и $x=y$.
Действительно, разность $V-V_1$, как всякая разность постоянных величин, должна равняться постоянной величине, разность же $x-y$, как всякая разность между переменными величинами, стремящимися к нулю, должна или равняться некоторой переменной величине (стремящейся к нулю), или равняться нулю.
Так как постоянная величина не может равняться переменной, то из двух возможностей надо оставить только одну: разность $x-y=0$;
но тогда $V=V_1$ и $x=y$.

Мы доказали, таким образом, что рассматриваемые пирамиды равновелики.

\medskip

\so{Замечание}
Необходимость столь сложного доказательства объясняется тем, что два равновеликих многогранника нельзя так легко преобразовывать один в другой, как это можно было делать с равновеликими многоугольниками на плоскости.
А именно, если даны два равновеликих многогранника, то в общем случае оказывается невозможным разбить один из них на такие части из которых можно было бы составить другой.
В частности, это невозможно для двух произвольных треугольных пирамид с равновеликими основаниями и равными высотами.
Эти утверждения были строго доказаны немецким математиком Максом Деном.\footnote{Доказательство приводится в книге «Третья проблема Гильберта» В. Г. Болтянского.}

Метод использованный нами в доказательстве называется \so{методом исчерпывания}; он является основным методом нахождения объёмов.
Этот метод лежит в основе понятия интеграла, которое изучается в курсе алгебры и начал анализа.

\paragraph{Принцип Кавальери.}\label{1938/s89} 
Следующее предложение высказано итальянским математиком Кавальери:
\emph{Если два тела} (ограниченные плоскостями или кривыми поверхностями — всё равно) \emph{могут быть помещены в такое положение, при котором всякая плоскость, параллельная какой-нибудь данной плоскости либо даёт в сечении с ними равновеликие фигуры либо вовсе не пересекает оба тела, то объёмы таких тел одинаковы.}

Это предложение может быть строго доказано тем же методом, которым мы воспользовались в §~\ref{1938/s90}.
Но это доказательство сложнее и потому мы ограничимся только формулировкой.

\begin{wrapfigure}{r}{60 mm}
\vskip-0mm
\centering
\includegraphics{mppics/s-ris-97}
\caption{}\label{1938/s-ris-97}
\bigskip
\includegraphics{mppics/s-ris-101}
\caption{}\label{1938/s-ris-101}
\bigskip
\includegraphics{mppics/s-ris-98}
\caption{}\label{1938/s-ris-98}
\end{wrapfigure}

Принцип Кавальери можно проверить на отдельных простых примерах.
Условиям принципа удовлетворяют, например, две прямые призмы (треугольные или многоугольные — всё равно) с равновеликими основаниями и равными высотами (рис.~\ref{1938/s-ris-97}).
Такие призмы, как мы знаем, равновелики.
Вместе с тем, если поставим такие призмы основаниями на какую-нибудь плоскость, то всякая плоскость, параллельная основаниям и пересекающая одну из призм, пересечёт и другую, причём в сечениях получатся равновеликие фигуры, так как фигуры эти равны основаниям, а основания равновелики.
Значит, принцип Кавальери подтверждается в этом частном случае.

Условиям принципа удовлетворяют также пирамиды с равновеликими основаниями и разными высотами.
Действительно, вообразим, что две такие пирамиды поставлены основаниями на какую-нибудь плоскость $P$ (рис.~\ref{1938/s-ris-101}), тогда всякая секущая плоскость $Q$, параллельная $P$, даёт в сечении с пирамидами равновеликие треугольники (§~\ref{1938/s77});
следовательно, пирамиды эти удовлетворяют условиям принципа Кавальери, и потому объёмы их должны быть одинаковы.
Заметим, что мы привели доказательство леммы в §~\ref{1938/s90}, используя принцип Кавальери;
однако это доказательство нельзя считать строгим, так как принцип Кавальери нами не был доказан.

В §~\ref{1938/s147} будет представлен ещё один пример применения принципа Кавальери.

Принцип Кавальери применим и к площадям фигур на плоскости, а именно:
\emph{если две фигуры могут быть помещены в такое положение, что всякая прямая, параллельная какой-нибудь данной прямой, либо даёт в сечении с ними равные отрезки либо не пересекает обе фигуры, то такие фигуры равновелики.}
Примером могут служить два параллелограмма или два треугольника с равными основаниями и равными высотами (рис.~\ref{1938/s-ris-98}).

}

\paragraph{}\label{1938/s91}
\so{Теорема}.
\textbf{\emph{Объём пирамиды равен произведению площади её основания на треть её высоты.}}

Сначала докажем эту теорему для пирамиды треугольной, а затем и многоугольной.

1) На основании треугольной пирамиды $SABC$ (рис.~\ref{1938/s-ris-102}) построим такую призму $ABCDES$, у которой высота равна высоте пирамиды, а одно боковое ребро совпадает с ребром $SB$.
Докажем, что объём пирамиды составляет третью часть объёма этой призмы.
Отделим от призмы данную пирамиду.
Тогда останется четырёхугольная пирамида $SADEC$ (которая для ясности изображена отдельно).
Проведём в ней секущую плоскость через вершину $S$ и диагональ основания $DC$.
Получившиеся от этого две треугольные пирамиды имеют общую вершину $S$ и равные основания $\triangle DEC$ и $\triangle DAC$, лежащие в одной плоскости;
значит, согласно доказанной выше лемме пирамиды эти равновелики.
Сравним одну из них, а именно $SDEC$, с данной пирамидой.
За основание пирамиды $SDEC$ можно взять $\triangle SDE$;
тогда вершина её будет в точке $C$ и высота равна высоте данной пирамиды.
Так как $\triangle SDE=\triangle ABC$, то согласно той же лемме пирамиды $SDEC$ и $SABC$ равновелики.

Призма $ABCDES$ нами разбита на три равновеликие пирамиды: $SABC$, $SDEC$ и $SDAC$.
(Такому разбиению, очевидно, можно подвергнуть всякую треугольную призму.) Таким образом, сумма объёмов трёх пирамид, равновеликих данной, составляет объём призмы;
следовательно,
\begin{align*}
\text{объём}\,SABC &= \tfrac13 \text{объёма}\, SDEABC =
\\
&=\frac{(\text{площади}\,ABC)\cdot H}3=
\\
&=(\text{площади}\,ABC)\cdot \frac{H}3,
\end{align*}
где $H$ есть высота пирамиды.

\begin{figure}[!ht]
\begin{minipage}{.58\textwidth}
\centering
\includegraphics{mppics/s-ris-102}
\end{minipage}\hfill
\begin{minipage}{.38\textwidth}
\centering
\includegraphics{mppics/s-ris-103}
\end{minipage}

\medskip

\begin{minipage}{.58\textwidth}
\centering
\caption{}\label{1938/s-ris-102}
\end{minipage}\hfill
\begin{minipage}{.38\textwidth}
\centering
\caption{}\label{1938/s-ris-103}
\end{minipage}
\vskip-4mm
\end{figure}

2) Через какую-нибудь вершину $E$ (рис.~\ref{1938/s-ris-103}) основания многоугольной пирамиды $SABCDE$ проведём диагонали $EB$ и $EC$.
Затем через ребро $SE$ и каждую из этих диагоналей проведём секущие плоскости.
Тогда многоугольная пирамида разобьётся на несколько треугольных, имеющих высоту, общую с данной пирамидой.
Обозначив площади оснований треугольных пирамид через $b_1$, $b_2$, $b_3$, высоту через $H$, будем иметь:
\begin{align*}
\text{объём}\,SABCDE &= \tfrac13b_1\cdot H + \tfrac13b_2\cdot H + \tfrac13b_3 \cdot H =
\\
&=(b_1+b_2+b_3)\cdot \frac H3=
\\
&=(\text{площади}\,ABCDE)\cdot \frac H3
\end{align*}
\so{Следствие}.
Если $V$, $B$ и $H$ означают числа, выражающие в соответствующих единицах объём, площадь основания и высоту какой угодно пирамиды, то
\[V = \tfrac13 B\cdot H.\]

\paragraph{}\label{1938/s92}
\so{Теорема}.
\textbf{\emph{Объём, усечённой пирамиды равен сумме объёмов трёх пирамид, имеющих высоту, одинаковую с высотой усечённой пирамиды, а основаниями: одна — нижнее основание данной пирамиды, другая — верхнее основание, а площадь основания третьей пирамиды равна среднему геометрическому площадей верхнего и нижнего оснований.}}

\begin{wrapfigure}{o}{38 mm}
\vskip-0mm
\centering
\includegraphics{mppics/s-ris-104}
\caption{}\label{1938/s-ris-104}
\vskip-0mm
\end{wrapfigure}

Пусть площади оснований усечённой пирамиды (рис.~\ref{1938/s-ris-104}) будут $B$ и $b$, высота $H$ и объём $V$ (усечённая пирамида может быть треугольная или многоугольная — всё равно).
Требуется доказать, что
\[V = \tfrac13B\cdot H + \tfrac13 b \cdot H + \tfrac13 H\cdot \sqrt{Bb} = \tfrac13 H(B + b + \sqrt{Bb}),\]
где $\sqrt{Bb}$ есть среднее геометрическое между $B$ и $b$.

Для доказательства на меньшем основании поместим малую пирамиду, дополняющую данную усечённую пирамиду до полной.
Тогда объём усечённой пирамиды $V$ мы можем рассматривать как разность двух объёмов полной пирамиды и верхней дополнительной.

Обозначив высоту дополнительной пирамиды буквой $x$, мы получим, что
\[V = \tfrac13B(H + x)- \tfrac13bx = \tfrac13(BH + Bx-bx) = \tfrac13 [BH + (B - b)x].\]

Для нахождения высоты $x$ воспользуемся теоремой §~\ref{1938/s74}, согласно которой мы можем написать уравнение:
\[\frac Bb=\frac{(H+x)^2}{x^2}.\]

Для упрощения этого уравнения извлечём из обеих частей его квадратный корень:
\[\frac {\sqrt{B}}{\sqrt{b}}=\frac{H+x}{x}.\]

Из этого уравнения (которое можно рассматривать как пропорцию) получим:
\[x\sqrt{B}=(H+x)\sqrt{b},\]
откуда
\[x(\sqrt{B}-\sqrt{b})=H\sqrt{b}\]
и, следовательно,
\[x=\frac{H\sqrt{b}}{\sqrt{B}-\sqrt{b}}.\]
Подставив это выражение в формулу, выведенную нами для объёма $V$, получим:
\[V=\frac13\left[BH+\frac{(B-b)H\sqrt{b}}{\sqrt{B}-\sqrt{b}}\right].\]
Так как $B - b= (\sqrt{B} + \sqrt{b})(\sqrt{B} - \sqrt{b})$, то по сокращении дроби на разность $\sqrt{B} - \sqrt{b}$ получим:
\begin{align*}
V&=\frac13\left[BH+(\sqrt{B}+\sqrt{b})H\sqrt{b}\right]=
\\
&=\frac13\left[BH+H\sqrt{Bb}+bH\right]=
\\
&=\frac13 H\left[B+\sqrt{Bb}+b\right],
\end{align*}
то есть получим ту формулу, которую требовалось доказать.

{\small

\medskip

\so{Замечание.}
Если основания усечённой пирамиды являются квадратами со сторонами $x$ и $y$, то согласно доказанной теореме, 
\[V=H(x^2+xy+y^2).\]
Эта формула для объёма квадратной усечённой пирамиды была известна в древнем Египте около 2 тысяч лет до нашей эры; в Московском математическом папирусе рассмотрен пример при $x=2$, $y=4$ и $H=6$.

}

{\small

\section{Правильные многогранники}

Многогранник называется правильным, если все его грани — равные правильные многоугольники и все многогранные углы равны (таков, например, куб).
Из этого определения следует, что в правильных многогранниках равны все плоские углы, все двугранные углы и все рёбра.

\paragraph{Перечисление правильных многогранников.}\label{1938/s97}
Примем во внимание, что в многогранном угле наименьшее число граней три, и что сумма плоских углов выпуклого многогранного угла меньше $360\degree$ (§~\ref{1938/s51}).

Каждый угол правильного треугольника равен $60\degree$.
Если повторим $60\degree$ слагаемым 3, 4 и 5 раз, то получим суммы, меньшие $360\degree$, а если повторим $60\degree$ слагаемым 6 раз или более, то получим в сумме $360\degree$ или более.
Поэтому из плоских углов, равных углам правильного треугольника, можно образовать выпуклые многогранные углы только трёх видов: трёхгранные, четырёхгранные и пятигранные.
Следовательно, если гранями правильного многогранника служат правильные треугольники, то в вершине многогранника могут сходиться или 3 ребра, или 4 ребра, или 5 рёбер.
Соответственно с этим имеется три вида правильных многогранников с треугольными гранями:

\begin{figure}[h]
\begin{minipage}{.32\textwidth}
\centering
\includegraphics{mppics/s-ris-107}
\end{minipage}\hfill
\begin{minipage}{.32\textwidth}
\centering
\includegraphics{mppics/s-ris-108}
\end{minipage}\hfill
\begin{minipage}{.32\textwidth}
\centering
\includegraphics{mppics/s-ris-109}
\end{minipage}

\medskip

\begin{minipage}{.32\textwidth}
\centering
\caption{}\label{1938/s-ris-107}
\end{minipage}\hfill
\begin{minipage}{.32\textwidth}
\centering
\caption{}\label{1938/s-ris-108}
\end{minipage}\hfill
\begin{minipage}{.32\textwidth}
\centering
\caption{}\label{1938/s-ris-109}
\end{minipage}
\vskip-4mm
\end{figure}

1) \rindex{тетраэдр}\textbf{Тетраэдр}, или правильный четырёхгранник, поверхность которого составлена из четырёх правильных треугольников (рис.~\ref{1938/s-ris-107}).
Он имеет 4 грани, 4 вершины и 6 рёбер.

2) \rindex{октаэдр}\textbf{Октаэдр}, или правильный восьмигранник, поверхность которого составлена из восьми правильных треугольников (рис.~\ref{1938/s-ris-108}).
Он имеет 8 граней, 6 вершин и 12 рёбер.

3) \rindex{икосаэдр}\textbf{Икосаэдр}, или правильный 20-гранник, образованный двадцатью правильными треугольниками (рис.~\ref{1938/s-ris-109}).
Он имеет 20 граней, 12 вершин и 30 рёбер.

Угол квадрата равен $90\degree$, а угол правильного пятиугольника равен $108\degree$, повторяя эти углы слагаемым 3 раза, получаем суммы, меньшие $360\degree$, а повторяя их 4 раза или более, получаем $360\degree$ или более.
Поэтому из плоских углов, равных углам квадрата или правильного пятиугольника, можно образовать только трёхгранные углы.

А поэтому, если гранями многогранника служат квадраты, то в каждой вершине могут сходиться лишь 3 ребра.
Имеется единственный правильный многогранник этого рода — это \rindex{куб}\textbf{куб}, или правильный шестигранник (рис.~\ref{1938/s-ris-110}).
Он имеет 6 граней, 8 вершин и 12 рёбер.

\begin{wrapfigure}{o}{35 mm}
\vskip-0mm
\centering
\includegraphics{mppics/s-ris-110}
\caption{}\label{1938/s-ris-110}
\bigskip
\includegraphics{mppics/s-ris-111}
\caption{}\label{1938/s-ris-111}
\vskip-0mm
\end{wrapfigure}

Если гранями правильного многогранника служат правильные пятиугольники, то в каждой вершине могут сходиться лишь 3 ребра.

Существует единственный правильный многогранник этого рода — правильный 12-гранник, или \rindex{додекаэдр}\textbf{додекаэдр}.
Он имеет 12 граней, 20 вершин и 30 рёбер (рис.~\ref{1938/s-ris-111}).

Угол правильного шестиугольника равен $120\degree$, поэтому из таких углов нельзя образовать даже трёхгранного угла.
Из углов правильных многоугольников, имеющих более 6 сторон, подавно нельзя образовать никакого выпуклого многогранного угла.

Отсюда следует, что гранями правильного многогранника могут служить лишь правильные треугольники, квадраты и правильные пятиугольники.

Таким образом, всего может существовать лишь пять видов правильных многогранников, указанных выше.

\paragraph{Построение правильных многогранников.}\label{1938/s98}
Изложенные выше рассуждения о возможных видах правильных многогранников доказывают, что может существовать не более пяти видов правильных многогранников.

Но из этих рассуждений ещё не вытекает, что все эти пять видов правильных многогранников действительно существуют, то есть что можно проведением плоскостей в пространстве осуществить построение каждого из этих пяти правильных многогранников.
Чтобы убедиться в существовании всех правильных многогранников, достаточно указать способ построения каждого из них.

Способ построения куба указать весьма легко.
Действительно, берём произвольную плоскость $P$ и в ней какой-либо квадрат;
через стороны этого квадрата проводим плоскости, перпендикулярные к плоскости~$P$.
Таких плоскостей будет четыре.
Далее проводим плоскость $Q$, параллельную $P$ и отстоящую от неё на расстоянии, равном стороне квадрата.
Шесть полученных плоскостей образуют грани куба;
двенадцать прямых — пересечения каждой пары пересекающихся плоскостей — являются рёбрами куба, а восемь точек пересечения каждой тройки пересекающихся плоскостей служат вершинами куба.

Умея построить куб, легко найти способ построения всех других правильных многогранников.

\textbf{Построение тетраэдра.}
Пусть дан куб (рис.~\ref{1938/s-ris-112}).
Возьмём какую-нибудь его вершину, например $A$.
В ней сходятся три грани куба, имеющие форму квадратов.
В каждом из этих квадратов берём вершину, противоположную точке $A$.
Пусть это будут вершины куба $B$, $C$ и $D$.

Точки $A$, $B$, $C$ и $D$ служат вершинами тетраэдра.
Действительно, каждый из отрезков $AB$, $BC$, $CD$, $AD$, $BD$ и $AC$, очевидно, служит диагональю одной из граней куба.
А потому все эти отрезки равны между собой.
Отсюда следует, что в треугольной пирамиде с вершиной $A$ и основанием $BCD$ все грани — правильные треугольники; следовательно, эта пирамида является тетраэдром.
Этот тетраэдр вписан в данный куб.

Полезно заметить, что оставшиеся четыре вершины куба служат вершинами другого тетраэдра, равного первому и также вписанного в данный куб.

\begin{figure}[h]
\begin{minipage}{.24\textwidth}
\centering
\includegraphics{mppics/s-ris-112}
\end{minipage}\hfill
\begin{minipage}{.24\textwidth}
\centering
\includegraphics{mppics/s-ris-113}
\end{minipage}\hfill
\begin{minipage}{.24\textwidth}
\centering
\includegraphics{mppics/s-ris-502}
\end{minipage}
\hfill
\begin{minipage}{.24\textwidth}
\centering
\includegraphics{mppics/s-ris-503}
\end{minipage}

\medskip

\begin{minipage}{.24\textwidth}
\centering
\caption{}\label{1938/s-ris-112}
\end{minipage}\hfill
\begin{minipage}{.24\textwidth}
\centering
\caption{}\label{1938/s-ris-113}
\end{minipage}\hfill
\begin{minipage}{.24\textwidth}
\centering
\caption{}\label{1938/s-ris-502}
\end{minipage}
\hfill
\begin{minipage}{.24\textwidth}
\centering
\caption{}\label{1938/s-ris-503}
\end{minipage}
\vskip-4mm
\end{figure}

\textbf{Построение октаэдра.}
Если в данном кубе построить центры всех его граней, то шесть полученных точек служат вершинами октаэдра (рис.~\ref{1938/s-ris-113}).

\textbf{Построение додекаэдра и икосаэдра.}
Если через каждое из 12 рёбер куба провести плоскость, не имеющую с поверхностью куба других общих точек, кроме точек того ребра, через которое она проведена, то полученные 12 плоскостей образуют грани некоторого 12-гранника.
Более подробное изучение формы этого многогранника показывает, что можно так подобрать наклон этих плоскостей к граням куба (рис.~\ref{1938/s-ris-502}), что полученный 12-гранник будет додекаэдром.
При этом куб окажется вписанным в додекаэдр; то есть 8 из 20 вершин додекаэдра являются вершинами куба.

Наконец, если мы умеем построить додекаэдр, то построение икосаэдра не представляет затруднений: центры граней додекаэдра служат вершинами икосаэдра (рис.~\ref{1938/s-ris-503}).

}

\section{Понятие о симметрии}

\begin{wrapfigure}{r}{35 mm}
\vskip-4mm
\centering
\includegraphics{mppics/s-ris-114}
\caption{}\label{1938/s-ris-114}
\vskip-0mm
\end{wrapfigure}

\paragraph{Центральная симметрия.}\label{1938/s99}
Две фигуры называются центрально симметричными относительно какой-либо точки $O$ пространства, если каждой точке $A$ одной фигуры соответствует в другой фигуре точка $A'$, расположенная на прямой $OA$ по другую сторону от точки $O$, на расстоянии, равном расстоянию точки $A$ от точки $O$ (рис.~\ref{1938/s-ris-114}).
Точка $O$ называется центром симметрии фигур.

Пример таких центрально симметричных фигур в пространстве мы уже встречали (§~\ref{1938/s53}), когда, продолжая за вершину ребра и грани многогранного угла, получали многогранный угол, центрально симметричный данному.
Соответственные отрезки и углы, входящие в состав двух центрально симметричных фигур, равны между собой.
Тем не менее фигуры в целом не могут быть названы равными: их нельзя совместить одну с другой вследствие того, что порядок расположения частей в одной фигуре иной, чем в другой, как это мы видели на примере центрально симметричных многогранных углов.

В отдельных случаях центрально симметричные фигуры могут совмещаться, но при этом будут совпадать несоответственные их части.

\begin{wrapfigure}{o}{40 mm}
\vskip-2mm
\centering
\includegraphics{mppics/s-ris-115}
\caption{}\label{1938/s-ris-115}
\vskip-0mm
\end{wrapfigure}

Например, возьмём прямой трёхгранный угол (рис.~\ref{1938/s-ris-115}) с вершиной в точке $O$ и рёбрами $OX$, $OY$, $OZ$.
Построим ему центрально симметричный угол $OX'Y'Z'$.
Угол $OXYZ$ можно совместить с $OX'Y'Z'$ так, чтобы ребро $OX$ совпало с $OY'$, а ребро $OY$ с $OX'$.
Если же совместить соответственные рёбра $OX$ с $OX'$ и $OY$ с $OY'$, то рёбра $OZ$ и $OZ'$ окажутся направленными в противоположные стороны.

Если центрально симметричные фигуры составляют в совокупности одно геометрическое тело, то говорят, что это геометрическое тело имеет центр симметрии.
Таким образом, если данное тело имеет центр симметрии, то всякой точке, принадлежащей этому телу, соответствует симметричная точка, тоже принадлежащая данному телу.
Из рассмотренных нами геометрических тел центр симметрии имеют, например: 1) параллелепипед, 2) призма, имеющая в основании правильный многоугольник с чётным числом сторон.

Тетраэдр не имеет центра симметрии.

\paragraph{Зеркальная симметрия.}\label{1938/s100}
Две пространственные фигуры называются зеркально симметричными относительно плоскости $P$, если каждой точке $A$ в одной фигуре соответствует в другой точка $A'$, причём отрезок $AA'$ перпендикулярен к плоскости $P$ и в точке пересечения с этой плоскостью делится пополам.

\medskip

\so{Теорема}.
\textbf{\emph{Всякие два соответственных отрезка в двух зеркально симметричных фигурах равны между собой.}}

\begin{wrapfigure}{r}{40 mm}
\vskip-0mm
\centering
\includegraphics{mppics/s-ris-116}
\caption{}\label{1938/s-ris-116}
\vskip-0mm
\end{wrapfigure}

Пусть даны две фигуры, зеркально симметричные относительно плоскости~$P$.
Выделим две какие-нибудь точки $A$ и $B$ первой фигуры, пусть $A'$ и $B'$ — соответствующие им точки второй фигуры (рис.~\ref{1938/s-ris-116}, на рисунке фигуры не изображены).
Пусть далее $C$ — точка пересечения отрезка $AA'$ с плоскостью $P$, $D$ — точка пересечения отрезка $BB'$ с той же плоскостью.
Соединив прямолинейным отрезком точки $C$ и $D$, получим два четырёхугольника $ABDC$ и $A'B'DC$.
Так как \[AC = A'C,\quad BD = B'D\] 
и 
\[\angle ACD = \angle A'CD=90\degree,\quad  \angle BDC = \angle B'DC=90\degree,\] то эти четырёхугольники равны (в чём легко убеждаемся наложением).
Следовательно, $AB=A'B'$.
Из этой теоремы непосредственно вытекает, что соответствующие плоские и двугранные углы двух зеркально симметричных фигур равны между собой.

Тем не менее совместить эти две фигуры одну с другой так, чтобы совместились их соответственные части, невозможно, так как порядок расположения частей в одной фигуре зеркальный тому, который имеет место в другой (это будет доказано ниже, §~\ref{1938/s102}).
Простейшим примером двух зеркально симметричных фигур являются: любой предмет и его отражение в плоском зеркале;
всякая фигура, симметричная со своим зеркальным отражением относительно плоскости зеркала.

Если какое-либо геометрическое тело можно разбить на две части, симметричные относительно некоторой плоскости, то эта плоскость называется плоскостью симметрии данного тела.

Геометрические тела, имеющие плоскость симметрии, чрезвычайно распространены в природе и в обыденной жизни.
Тело человека и животного имеет плоскость симметрии, разделяющую его на правую и левую части.

На этом примере особенно ясно видно, что зеркально симметричные фигуры нельзя совместить.
Так, кисти правой и левой рук симметричны, но совместить их нельзя, что можно видеть хотя бы из того, что одна и та же перчатка не может подходить и к правой и к левой руке.
Большое число предметов домашнего обихода имеет плоскость симметрии: стул, обеденный стол, книжный шкаф, диван и др.
Некоторые, как например обеденный стол, имеют даже не одну, а две плоскости симметрии (рис.~\ref{1938/s-ris-117}).

{

\begin{wrapfigure}{r}{43 mm}
\vskip-7mm
\centering
\includegraphics{mppics/s-ris-117}
\caption{}\label{1938/s-ris-117}
\vskip-0mm
\end{wrapfigure}

Если, рассматривая предмет, имеющий плоскость симметрии, занять по отношению к нему такое положение, чтобы плоскость симметрии нашего тела, или по крайней мере нашей головы, совпала с плоскостью симметрии самого предмета, то симметричная форма предмета становится особенно заметной.

}

\paragraph{Осевая симметрия.}\label{1938/s101}
Две фигуры называются симметричными относительно оси $l$ (ось — прямая линия), если каждой точке первой фигуры соответствует точка $A'$ второй фигуры, так что отрезок $AA'$ перпендикулярен к оси, пересекается с нею и в точке пересечения делится пополам.
Сама прямая $l$ называется осью симметрии второго порядка.

Из этого определения непосредственно следует, что если два геометрических тела, симметричных относительно какой-либо оси, пересечь плоскостью, перпендикулярной к этой оси, то в сечении получатся две плоские фигуры, центрально симметричные относительно точки пересечения плоскости с осью симметрии тел.

Отсюда далее легко вывести, что два тела, симметричных относительно оси, можно совместить одно с другим, вращая одно из них на $180\degree$ вокруг оси симметрии.
В самом деле, вообразим все возможные плоскости, перпендикулярные к оси симметрии.
Каждая такая плоскость, пересекающая оба тела, содержит две фигуры, центрально симметричные относительно точки встречи плоскости с осью симметрии тел.
Если заставить скользить секущую плоскость саму по себе, вращая её вокруг оси симметрии тела на $180\degree$, то первая фигура совпадает со второй.
Это справедливо для любой секущей плоскости.
Вращение же всех сечений тела на $180\degree$ равносильно повороту всего тела на $180\degree$ вокруг оси симметрии.
Отсюда и вытекает справедливость нашего утверждения.

Если после вращения пространственной фигуры вокруг некоторой прямой на $180\degree$ она совпадает сама с собой, то говорят, что фигура имеет эту прямую своею осью симметрии второго порядка.

Название «ось симметрии второго порядка» объясняется тем, что при полном обороте вокруг этой оси тело будет в процессе вращения дважды принимать положение, совпадающее с исходным (считая и исходное).
Примерами геометрических тел, имеющих ось симметрии второго порядка, могут служить:

1) правильная пирамида с чётным числом боковых граней;
осью её симметрии служит её высота;

2) прямоугольный параллелепипед;
он имеет три оси симметрии: прямые, соединяющие центры его противоположных граней;

3) правильная призма с чётным числом боковых граней.
Осью её симметрии служит каждая прямая, соединяющая центры любой пары противоположных граней (боковых граней и двух оснований призмы).
Если число боковых граней призмы равно $2k$, то число таких осей симметрии будет $k+1$.
Кроме того, осью симметрии для такой призмы служит каждая прямая, соединяющая середины её противоположных боковых рёбер.
Таких осей симметрии призма имеет $k$.

Таким образом, правильная $2k$-гранная призма имеет $2k+1$ осей симметрии.

\paragraph{Зависимость между различными видами симметрии в пространстве.}\label{1938/s102}
Между различными видами симметрии в пространстве — осевой, плоскостной и центральной — существует зависимость, выражаемая следующей теоремой.

\medskip

\so{Теорема}.
\textbf{\emph{Если фигура $F$ зеркально симметрична с фигурой $F'$ относительно плоскости $P$ и в то же время центрально симметрична с фигурой $F''$ относительно точки $O$, лежащей в плоскости $P$, то фигуры $F'$ и $F''$ симметричны относительно оси, проходящей через точку $O$ и перпендикулярной к плоскости~$P$.}}

\begin{wrapfigure}{o}{33 mm}
\vskip-0mm
\centering
\includegraphics{mppics/s-ris-118}
\caption{}\label{1938/s-ris-118}
\vskip-0mm
\end{wrapfigure}

Возьмём какую-нибудь точку $A$ фигуры $F$ (рис.~\ref{1938/s-ris-118}).
Ей соответствует точка $A'$ фигуры $F'$ и точка $A''$ фигуры $F''$ (сами фигуры $F$, $F'$ и $F''$ на чертеже не изображены).

Пусть $B$ — точка пересечения отрезка $AA'$ с плоскостью~$P$.
Проведём плоскость через точки $A$, $A'$ и $O$.
Эта плоскость будет перпендикулярна к плоскости $P$, так как проходит через прямую $AA'$, перпендикулярную к этой плоскости.
В плоскости $AA'O$ проведём прямую $OH$, перпендикулярную к $OB$.
Эта прямая $OH$ будет перпендикулярна и к плоскости~$P$.
Пусть далее $C$ — точка пересечения прямых $AA''$ и $OH$.

В треугольнике $AA'A''$ отрезок $BO$ соединяет середины сторон $AA'$ и $AA''$, следовательно, $BO\parallel A'A''$, но $BO\perp OH$, значит, $A'A''\perp OH$.
Далее, так как $O$ — середина стороны $AA''$ и $CO\parallel AA'$, то $A'C = A''C$.
Отсюда заключаем, что точки $A'$ и $A''$ симметричны относительно оси $OH$.
То же самое справедливо и для всех других точек фигуры.
Значит, наша теорема доказана.

Из этой теоремы непосредственно следует, что две зеркально симметричные фигуры не могут быть совмещены так, чтобы совместились их соответственные части.
В самом деле, фигура $F'$ совмещается с $F''$ путём вращения вокруг оси $OH$ на $180\degree$.
Но фигуры $F''$ и $F$ не могут быть совмещены как центрально симметричные, следовательно, фигуры $F$ и $F'$ также не могут быть совмещены.

\begin{wrapfigure}{r}{33 mm}
\vskip-8mm
\centering
\includegraphics{mppics/s-ris-119}
\caption{}\label{1938/s-ris-119}
\vskip-0mm
\end{wrapfigure}

\paragraph{Оси симметрии высших порядков.}\label{1938/s103}
Фигура, имеющая ось симметрии, совмещается сама с собой после поворота вокруг оси симметрии на угол в $180\degree$.
Но возможны случаи, когда фигура приходит к совмещению с исходным положением после поворота вокруг некоторой оси на угол, меньший $180\degree$.
Таким образом, если тело сделает полный оборот вокруг этой оси, то в процессе вращения оно несколько раз совместится со своим первоначальным положением.
Такая ось вращения называется осью симметрии высшего порядка, причём число положений тела, совпадающих с первоначальным, называется порядком оси симметрии.
Эта ось может и не совпадать с осью симметрии второго порядка.
Так, правильная треугольная пирамида не имеет оси симметрии второго порядка, но её высота служит для неё осью симметрии третьего порядка.
В самом деле, после поворота этой пирамиды вокруг высоты на угол в $120\degree$ она совмещается сама с собой (рис.~\ref{1938/s-ris-119}).
При вращении пирамиды вокруг высоты она может занимать три положения, совпадающие с исходным, считая и исходное.
Легко заметить, что всякая ось симметрии чётного порядка является в то же время осью симметрии второго порядка.

Примеры осей симметрии высших порядков:

1) Правильная $n$-угольная пирамида имеет ось симметрии $n$-го порядка.
Этой осью служит высота пирамиды.

2) Правильная $n$-угольная призма имеет ось симметрии $n$-го порядка.
Этой осью служит прямая, соединяющая центры оснований призмы.

\paragraph{Симметрии куба.}\label{1938/s104}
Как и для всякого параллелепипеда, точка пересечения диагоналей куба есть центр его симметрии.

Куб имеет девять плоскостей симметрии: шесть диагональных плоскостей и три плоскости, проходящие через середины каждой четвёрки его параллельных рёбер.

Куб имеет девять осей симметрии второго порядка — это шесть прямых, соединяющих середины его противоположных рёбер, и три прямые, соединяющие центры противоположных граней (рис.~\ref{1938/s-ris-120}).
Эти последние три прямые также являются осями симметрии четвёртого порядка.
Кроме того, куб имеет четыре оси симметрии третьего порядка, которые являются его диагоналями.
В самом деле, диагональ куба $AG$, очевидно, одинаково наклонена к рёбрам $AB$, $AD$ и $AE$, а эти рёбра одинаково наклонены одно к другому.
Если соединить точки $B$, $D$ и $E$, то получим правильную треугольную пирамиду $ADBE$, для которой диагональ куба AG служит высотой.
Когда при вращении вокруг высоты эта пирамида будет совмещаться сама с собой, весь куб будет совмещаться со своим исходным положением.
Других осей симметрии, как нетрудно убедиться, куб не имеет.

\begin{wrapfigure}{r}{40 mm}
\vskip-6mm
\centering
\includegraphics{mppics/s-ris-120}
\caption{}\label{1938/s-ris-120}
\vskip-0mm
\end{wrapfigure}

Посмотрим, сколькими различными способами куб может быть совмещён сам с собой.
Выберем пару соседних вершин куба, скажем $A$ и $B$ (рис.~\ref{1938/s-ris-120}).
При совмещении куба вершина $A$ может перейти в любую из 8 вершин, 
а вершина $B$ в любую из трёх её соседей.
Положение этих двух вершин полностью определяет положение остальных.
Значит всего существует $24=8\cdot 3$  способа совмещения куба с самим собой, включая исходное положение.

Заметим, что вращение вокруг обыкновенной оси симметрии даёт одно положение куба, отличное от исходного, при котором куб в целом совмещается сам с собой.
Вращение вокруг оси третьего порядка даёт два таких положения, и вращение вокруг оси четвёртого порядка — три таких положения.
Так как куб имеет шесть осей второго порядка (это обыкновенные оси симметрии), четыре оси третьего порядка и три оси четвёртого порядка, то имеются $6\cdot 1 + 4\cdot 2 + 3\cdot 3 = 23$ положения куба, отличные от исходного, при которых он совмещается сам с собой. 

Легко убедиться непосредственно, что все эти положения отличны одно от другого, а также и от исходного положения куба.
Вместе с исходным положением они составляют все 24 способа совмещения куба с самим собой.

{\small

\subsection*{Упражнения}

\begin{enumerate}[noitemsep]

\item
Ребро данного куба равно $a$.
Найти ребро другого куба, объём которого вдвое более объёма данного куба.

\medskip

\so{Замечание}.
Эта \so{задача об удвоении куба}, известная с древних времён, легко решается вычислением (а именно $x\z=\sqrt[3]{2a^3}\z=a\cdot\sqrt[3]{2}\z=a\cdot1{,}26992 \dots $), но построением с помощью циркуля и линейки она решена быть не может, так как формула для неизвестного содержит радикал третьей степени из числа, не являющегося кубом рационального числа.

\item
Вычислить площадь поверхности и объём прямой призмы, у которой основание — правильный треугольник, вписанный в круг радиуса $r \z= 2$~м, а высота равна стороне правильного шестиугольника, описанного около того же круга.

\item
Определить площадь поверхности и объём правильной восьмиугольной призмы, у которой высота $h$ = 6 м, а сторона основания $a$ = 8 см.

\item
Определить площадь боковой поверхности и объём правильной шестиугольной пирамиды, у которой высота равна 1 м, а апофема составляет с высотой угол в $30\degree$.

\item
Вычислить объём треугольной пирамиды, у которой каждое боковое ребро равно $l$, а стороны основания — $a$, $b$ и $c$.

\item
Дан трёхгранный угол $SABC$, у которого все три плоских угла прямые.
На его рёбрах отложены длины: $SA = a$;
$SB = b$ и $SC\z=c$.
Через точки $A$, $B$ и $C$ проведена плоскость.
Определить объём пирамиды $SABC$.

\item
Высота пирамиды равна $h$, а основание — правильный шестиугольник со стороной $a$.
На каком расстоянии $x$ от вершины пирамиды следует провести плоскость, параллельную основанию, чтобы объём образовавшейся усечённой пирамиды равнялся $V$?

\item
Определить объём тетраэдра с ребром $a$.

\item
Определить объём октаэдра с ребром $a$.

\item
Усечённая пирамида, объём которой $V = 1465 \text{см}^3$, имеет основаниями правильные шестиугольники со сторонами: $a = 23 \text{см}$ и $b = 17 \text{см}$.
Вычислить высоту этой пирамиды.

\item
Объём $V$ усечённой пирамиды равен $10{,}5\text{м}^3$, высота $h = \sqrt{3}\text{м}$ и сторона $a$ правильного шестиугольника, служащего нижним основанием, равна 2 м.
Вычислить сторону правильного шестиугольника, служащего верхним основанием.

\item
На каком расстоянии от вершины $S$ пирамиды $SABC$ надо провести плоскость, параллельную основанию, чтобы отношение объёмов частей, на которые рассекается этой плоскостью пирамида, равнялось $m$?

\item
Пирамида с высотой $h$ разделена плоскостями, параллельными основанию, на три части, причём объёмы этих частей находятся в отношении $m : n : p$.
Определить расстояние этих плоскостей до вершины пирамиды.

\item
Разделить усечённую пирамиду плоскостью, параллельной основаниям $B$ и $b$, на две части, чтобы объёмы находились в отношении $m:n$.

\item
Найти центр, оси и плоскости симметрии фигуры, состоящей из плоскости и пересекающей её прямой, не перпендикулярной к этой плоскости.

\so{Ответ}: центр симметрии — точка пересечения прямой с плоскостью;
плоскость симметрии — плоскость, перпендикулярная данной, проходящая через данную прямую;
осью симметрии служит прямая, лежащая в данной плоскости и перпендикулярная к данной прямой.

\item
Найти центр, оси и плоскости симметрии фигуры, состоящей из двух пересекающихся прямых.

\so{Ответ}: фигура имеет три плоскости симметрии и три оси симметрии (указать какие).

\item Доказать, что многогранник не может иметь двух различных центров симметрии.

\item Описать оси симметрии тетраэдра и указать их порядки.

\item Сколькими разными способами можно совместить тетраэдр с собой, включая его исходное положение?
\end{enumerate}

}

%% file: 3D/kruglye_tela.tex
\chapter{Круглые тела}

\section{Цилиндр и конус}

\begin{wrapfigure}{r}{43 mm}
\vskip-0mm
\centering
\includegraphics{mppics/s-ris-121}
\caption{}\label{1938/s-ris-121}
\vskip-0mm
\end{wrapfigure}

\paragraph{Поверхность вращения.}\label{1938/s105}
Поверхностью вращения называется поверхность, которая получается от вращения какой-нибудь линии ($MN$, рис.~\ref{1938/s-ris-121}), называемой образующей, вокруг неподвижной прямой ($AB$), называемой осью, при этом предполагается, что образующая ($MN$) при своём вращении неизменно связана с осью ($AB$).

Возьмём на образующей какую-нибудь точку $P$ и опустим из неё на ось перпендикуляр $PO$.
Очевидно, что при вращении не изменяются ни длина этого перпендикуляра, ни величина угла $AOP$, ни положение точки $O$.
Поэтому каждая точка образующей описывает окружность, плоскость которой перпендикулярна к оси $AB$ и центр которой лежит на пересечении этой плоскости с осью.
Отсюда следует:

\emph{Плоскость, перпендикулярная к оси, пересекаясь с поверхностью вращения, даёт в сечении окружность (или набор окружностей).} 
Такие окружности называются  \rindex{параллель}\textbf{параллелями} поверхности вращения.

Всякая секущая плоскость, проходящая через ось, называется \rindex{меридиональная плоскость}\textbf{меридиональной} плоскостью, а линия её пересечения с поверхностью вращения — \rindex{меридиан}\textbf{меридианом}.
Все меридианы равны между собой, потому что при вращении каждый из них проходит через то положение, в котором ранее был всякий другой меридиан.

\begin{wrapfigure}{r}{38 mm}
\vskip-0mm
\centering
\includegraphics{mppics/s-ris-122}
\caption{}\label{1938/s-ris-122}
\bigskip
\includegraphics{mppics/s-ris-123}
\caption{}\label{1938/s-ris-123}
\bigskip
\includegraphics{mppics/s-ris-124}
\caption{}\label{1938/s-ris-124}
\bigskip
\includegraphics{mppics/s-ris-125}
\caption{}\label{1938/s-ris-125}
\vskip-0mm
\end{wrapfigure}

\paragraph{Цилиндрическая поверхность.}\label{1938/s106}
Цилиндрической поверхностью называется поверхность, производимая движением прямой ($AB$, рис.~\ref{1938/s-ris-122}), перемещающейся в пространстве параллельно данной прямой и пересекающей при этом данную линию ($MN$).
Прямая $AB$ называется \rindex{образующая}\textbf{образующей}, а линия $MN$ — \rindex{направляющая}\textbf{направляющей}.

\paragraph{Цилиндр.}\label{1938/s107}
Цилиндром называется тело, ограниченное цилиндрической поверхностью и двумя параллельными плоскостями (рис.~\ref{1938/s-ris-123}).

Часть цилиндрической поверхности, заключённая между плоскостями, называется \rindex{боковая поверхность!цилиндра}\textbf{боковой поверхностью}, а части плоскостей, отсекаемые этой поверхностью, — \rindex{основание!цилиндра}\textbf{основаниями} цилиндра.
Расстояние между плоскостями оснований есть \rindex{высота!цилиндра}\textbf{высота} цилиндра.
Цилиндр называется \rindex{прямой цилиндр}\textbf{прямым} или \rindex{наклонный цилиндр}\textbf{наклонным}, смотря по тому, перпендикулярны или наклонны к основаниям его образующие.

Прямой цилиндр (рис.~\ref{1938/s-ris-124}) называется круговым, если его основания — круги.
Такой цилиндр можно рассматривать как тело, происходящее от вращения прямоугольника $OAA_1O_1$ вокруг стороны $OO_1$ как оси;
при этом сторона $AA_1$ описывает боковую поверхность, а стороны $OA$ и $O_1A_1$ — круги оснований.
Всякий отрезок $BC$, параллельный $OA$, описывает также круг, плоскость которого перпендикулярна к оси.
Отсюда следует:

\emph{Сечение прямого кругового цилиндра плоскостью, параллельной основаниям, есть круг.}

В элементарной геометрии рассматривается только прямой круговой цилиндр;
для краткости его называют просто цилиндром.
Иногда приходится рассматривать такие призмы, основания которых — многоугольники, вписанные в основания цилиндра или описанные около них, а высоты равны высоте цилиндра;
такие призмы называются \rindex{вписанная призма}\textbf{вписанными} в цилиндр или \rindex{описанная призма}\textbf{описанными} около него.

{

\begin{wrapfigure}{r}{38 mm}
\vskip-0mm
\centering
\includegraphics{mppics/s-ris-126}
\caption{}\label{1938/s-ris-126}
\bigskip
\includegraphics{mppics/s-ris-127}
\caption{}\label{1938/s-ris-127}
\bigskip
\includegraphics{mppics/s-ris-128}
\caption{}\label{1938/s-ris-128}
\vskip-0mm
\end{wrapfigure}

\paragraph{Коническая поверхность.}\label{1938/s108}
Конической поверхностью называется поверхность, производимая движением прямой ($AB$, рис.~\ref{1938/s-ris-125}), перемещающейся в пространстве так, что она при этом постоянно проходит через неподвижную точку ($S$) и пересекает данную линию ($MN$).
Прямая $AB$ называется \rindex{образующая}\textbf{образующей}, линия $MN$ — \rindex{направляющая}\textbf{направляющей}, а точка $S$ — \rindex{вершина!конической поверхности}\textbf{вершиной конической поверхности}.

\paragraph{Конус.}\label{1938/s109}
Конусом называется тело, ограниченное частью конической поверхности, расположенной по одну сторону от вершины, и плоскостью, пересекающей все образующие по ту же сторону от вершины (рис.~\ref{1938/s-ris-126}).
Часть конической поверхности, ограниченная этой плоскостью, называется \rindex{боковая поверхность!конуса}\textbf{боковой поверхностью}, а часть плоскости, отсекаемая боковой поверхностью, — основанием конуса.
Перпендикуляр, опущенный из вершины на плоскость основания, называется \rindex{высота!конуса}\textbf{высотой конуса}.

Конус называется \rindex{прямой круговой конус}\textbf{прямым круговым}, если его основание есть круг, а высота проходит через центр основания (рис.~\ref{1938/s-ris-127}).
Такой конус можно рассматривать как тело, происходящее от вращения прямоугольного треугольника $SOA$ вокруг катета $SO$ как оси.
При этом гипотенуза $SA$ описывает боковую поверхность, а катет $OA$ — основание конуса.
Всякий отрезок $BO$\, параллельный $OA$, описывает при вращении круг, плоскость которого перпендикулярна к оси.
Отсюда следует:

\emph{Сечение прямого кругового конуса плоскостью, параллельной основанию, есть круг.}

В элементарной геометрии рассматривается только прямой круговой конус, который для краткости называется просто конусом.

}

Иногда приходится рассматривать такие пирамиды, основаниями которых являются многоугольники, вписанные в основание конуса или описанные около него, а вершина совпадает с вершиной конуса.
Такие пирамиды называются \rindex{вписанная пирамида}\textbf{вписанными} в конус или \rindex{описанная пирамида}\textbf{описанными} около него.

\paragraph{Усечённый конус.}\label{1938/s110}
Так называется часть полного конуса, заключённая между основанием и секущей плоскостью, \so{параллельной основанию}.

Круги, по которым параллельные плоскости пересекают конус, называются \rindex{основание!усечённого конуса}\textbf{основаниями} усечённого конуса.

Усечённый конус (рис.~\ref{1938/s-ris-128}) можно рассматривать как тело, происходящее от вращения прямоугольной трапеции $OAA_1O_1$ вокруг стороны $OO_1$, перпендикулярной к основаниям трапеции.

\subsection*{Поверхность цилиндра и конуса}

\paragraph{Площадь кривой поверхности.}\label{1938/s111}
Боковые поверхности цилиндра и конуса принадлежат к поверхностям {}\textbf{кривым}, то есть к таким, никакая часть которых не может совместиться с плоскостью.
Поэтому мы должны особо оговорить, как измерить площадь \so{кривой} поверхности, используя \so{плоскую} единицу площади.

Напомним, что длину дуги мы определяли как предел длин вписанных в неё ломаных при условии, что стороны ломаных неограниченно уменьшаются.
Может показаться, что площадь кривой поверхности можно определить аналогично: рассмотреть вписанные в неё многогранные поверхности и перейти к пределу их площадей при условии, что размер их граней стремится к нулю.
Однако, как будет показано в §~\ref{fikhtengoltz/3-623b}, это определение непригодно даже для боковой поверхности прямого кругового цилиндра — даже в этом случае предел о котором идёт речь не существует.

На помощь приходит то обстоятельство, что цилиндры, конусы, усечённые конусы (а также шары, которые будут рассматриваться в~§~\ref{1938/s136}) являются \so{выпуклыми телами}; то есть с любой парой своих точек они содержат и отрезок между ними.
Доказано, что если рассматривать приближения выпуклых тел только \so{выпуклыми} вписанными многогранниками, то независимо от выбора последовательности вписанных многогранников, площади их поверхностей стремятся к одному и тому же пределу.
Этот предел и принимается за площадь поверхности исходного тела. 
Более того, если выпуклое тело приближается последовательностью выпуклых тел,
то площади поверхностей последовательности выпуклых тел стремятся к площади поверхности предельного тела.%
\footnote{Доказательство можно найти например в книге «Круг и шар» В. Бляшке.}
Можно также определить площадь частей поверхностей выпуклых тел так, чтобы удовлетворялись условия подобные описанным в~§~\ref{1938/243}.
В частности (1) площади двух равных поверхностей, равны между собой; (2) если поверхность разбита на несколько частей, то  площадь всей поверхности равна сумме площадей отдельных её частей.

Отсюда легко выводятся следующие два утверждения, которые можно принять за \so{определение} площади боковых поверхностей цилиндра и конуса:

1) \emph{За площадь боковой поверхности цилиндра принимают предел, к которому стремится площадь боковой поверхности вписанной в этот цилиндр правильной призмы, когда число сторон правильного многоугольника, вписанного в основание, неограниченно удваивается} (и, следовательно, площадь каждой боковой грани неограниченно убывает).

2) \emph{За площадь боковой поверхности конуса} (полного или усечённого) \emph{принимается предел, к которому стремится площадь боковой поверхности вписанной в этот конус правильной пирамиды} (полной или усечённой), \emph{когда число сторон правильного многоугольника, вписанного в основание, неограниченно удваивается} (и, следовательно, площадь каждой боковой грани неограниченно убывает).

{\small

\begin{wrapfigure}{r}{40 mm}
\centering
\includegraphics{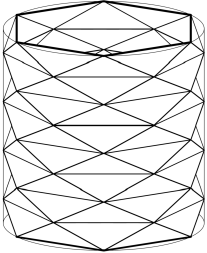}
\caption{}
\label{acy/schwarz}
\end{wrapfigure}

\paragraph{Сапог Шварца.}\label{fikhtengoltz/3-623b}
Мы построим приближение боковой поверхности прямого кругового цилиндра \so{невыпуклыми} вписанными многогранными поверхностями, площади которых не имеют предела.
Этот пример был приведён Карлом Шварцем; он демонстрирует несостоятельность определения площади поверхности тела как предела площадей поверхностей вписанных в него многогранников (без дополнительного предположения, что многогранники выпуклы).

Разделим данный цилиндр на $m$ равных цилиндров плоскостями, параллельными его основаниям.
Так на поверхности данного цилиндра получится $m+1$ окружностей, включая и окружности обоих оснований цилиндра.
Далее, каждую из этих окружностей разделим на $n$ равных дуг так, чтобы точки деления вышележащей окружности находились в точности над серединами дуг нижележащей окружности.
Рассмотрим все треугольники, образованные хордами всех этих дуг и отрезками, соединяющими концы хорд с теми точками деления выше- и нижележащих окружностей, которые расположены как раз над или под серединами соответствующих дуг.
Все эти $2mn$ равных треугольников образуют многогранную поверхность; 
она и называется \so{сапогом Шварца}.
Модель сапога Шварца при $m=8$ и $n=6$ показана на рис.~\ref{acy/schwarz}.

\begin{wrapfigure}{o}{40 mm}
\centering
\includegraphics{mppics/s-ris-501}
\caption{}\label{treug-sapog}
\end{wrapfigure}

Предположим, что основанием цилиндра является круг радиуса $R$.
Рассмотрим один из треугольников $ABC$, из которых состоит сапог Шварца, рис.~\ref{treug-sapog}.
По построению, основание $AC$ лежит на горизонтальной плоскости,
а проекция $B'$ вершины $B$ на эту плоскость делит дугу $AC$ на две равные части.
Значит вершины треугольника $AB'C$ это три последовательные вершины правильного $2n$-угольника вписанного в окружность радиуса~$R$.
Обозначим через $s_n$ площадь треугольника $AB'C$; очевидно она зависит только от $R$ и~$n$.
 
Оба треугольника $ABC$ и $AB'C$ равнобедренные и имеют общее основание $AC$;
при этом высота $BK$ больше высоты $B'K$ так как отрезок $BK$ длинней своей проекции $B'K$.
Отсюда получаем, что
\begin{align*}
\text{площадь}\, ABC&=\tfrac12 AC\cdot BK>
\\
&>\tfrac12 AC\cdot B'K=
\\
&=\text{площадь}\, AB'C=
\\
&=s_n.
\end{align*}
Значит общая площадь $S_{m,n}$ сапога Шварца превосходит $2mns_n$.

Заметим, что построенные многогранные поверхности приближают боковую поверхность цилиндра если $m$ и $n$ неограниченно возрастают.
Предположим, дополнительно, что $m$ возрастает существенно быстрее $n$,
а именно выполняется условие $m>\tfrac1{s_n}$.
Тогда
\[S_{m,n}> 2 m n s_n>2 \tfrac1{s_n} n s_n=2n.\]
Значит, при условии $m>\tfrac1{s_n}$, площадь $S_{m,n}$ неограниченно возрастает с ростом~$n$.
В частности $S_{m,n}$ не имеет предела при неограниченном возрастании $m$ и~$n$.%
\footnote{Следует отметить, что при дополнительном предположении $m=n$, площади $S_{m,n}$ стремятся к площади боковой поверхности цилиндра. Подробный анализ зависимости площади $S_{m,n}$ от $m$ и $n$ приведён в §~623, третьего тома  «Курса дифференциального и интегрального исчисления» Г.~М.~Фихтенгольца.}

}

\paragraph{}\label{1938/s112}
\so{Теорема}.
\textbf{\emph{Площадь боковой поверхности цилиндра равна произведению длины окружности основания на высоту.}}

\begin{wrapfigure}{o}{38 mm}
\vskip-0mm
\centering
\includegraphics{mppics/s-ris-129}
\caption{}\label{1938/s-ris-129}
\vskip4mm
\end{wrapfigure}

Впишем в цилиндр (рис.~\ref{1938/s-ris-129}) какую-нибудь правильную призму.
Обозначим буквами $p$ и $H$ числа, выражающие длины периметра основания и высоты этой призмы.
Тогда площадь её боковой поверхности выразится произведением $p\cdot H$ (§~\ref{1938/s79}).
Предположим теперь, что число сторон вписанного в основание многоугольника неограниченно возрастает.

Тогда периметр $p$ будет стремиться к пределу, принимаемому за длину $C$ окружности основания, а высота $H$ останется без изменения;
следовательно, площадь боковой поверхности призмы, равная произведению $p\cdot H$, будет стремиться к пределу $C\cdot H$.
Этот предел и принимается за площадь боковой поверхности цилиндра;
обозначив её буквой $S$, можем написать:
\[S = C\cdot H.\]

\paragraph{}\label{1938/s113}
\mbox{\so{Следствия}.}
1) Если $R$ обозначает радиус основания цилиндра, то $C= 2\pi R$, поэтому площадь боковой поверхности цилиндра выразится формулой:
\[S = 2\pi R \cdot H.\]

2) Чтобы получить площадь полной поверхности цилиндра, достаточно приложить к боковой поверхности сумму площадей двух оснований, поэтому, обозначая площадь полной поверхности через $T$, будем иметь:
\[T= 2\pi RH + \pi R^2 + \pi R^2 = 2\pi R(H + R).\]

\paragraph{}\label{1938/s114}
\so{Теорема}.
\textbf{\emph{Площадь боковой поверхности конуса равна произведению длины окружности основания на половину образующей.}}

\begin{wrapfigure}{o}{38 mm}
\vskip-0mm
\centering
\includegraphics{mppics/s-ris-130}
\caption{}\label{1938/s-ris-130}
\vskip-0mm
\end{wrapfigure}

Впишем в конус (рис.~\ref{1938/s-ris-130}) какую-нибудь правильную пирамиду и обозначим буквами $p$ и $l$ числа, выражающие длины периметра основания и апофемы этой пирамиды.
Тогда площадь её боковой поверхности выразится произведением $\tfrac12 p\cdot l$ (§~\ref{1938/s80}).

Предположим теперь, что число сторон вписанного в основание многоугольника неограниченно возрастает.
Тогда периметр $p$ будет стремиться к пределу, принимаемому за длину $C$ окружности основания, а апофема $l$ будет иметь пределом образующую конуса (так как из $\triangle SAK$ следует, что $SA\z-SK\z<AK$).
Значит, если образующую конуса обозначим буквой $L$, то площадь боковой поверхности вписанной пирамиды, равная $\tfrac12 p\cdot l$, будет стремиться к пределу $\tfrac12 C\cdot L$. 

Этот предел и принимается за площадь боковой поверхности конуса; обозначив её буквой $S$, можем написать:
\[S = \tfrac12 C\cdot L.\]

\paragraph{}\label{1938/s115}
\so{Следствия}. 1) Если $R$ обозначает радиус основания конуса, то $C= 2\pi R$, поэтому площадь боковой поверхности конуса выразится формулой:
\[S
= \tfrac12 \cdot 2\pi R \cdot L
=\pi RL.\]

2) Площадь полной поверхности конуса получим, если площадь боковой поверхности сложим с площадью основания;
поэтому, обозначая площадь полной поверхности через $T$:
\[T= \pi RL + \pi R^2 = \pi R(L + R).\]

\paragraph{}\label{1938/s116}
\so{Теорема}.
\textbf{\emph{Площадь боковой поверхности усечённого конуса равна произведению полусуммы длин окружностей оснований на образующую.}}

\begin{wrapfigure}{o}{44 mm}
\vskip-0mm
\centering
\includegraphics{mppics/s-ris-131}
\caption{}\label{1938/s-ris-131}
\vskip-0mm
\end{wrapfigure}

Впишем в усечённый конус (рис.~\ref{1938/s-ris-131}) какую-нибудь правильную усечённую пирамиду и обозначим буквами $p$, $p_1$ и $l$ числа, выражающие в одинаковых линейных единицах длины периметров нижнего и верхнего оснований и апофемы этой пирамиды.
Тогда площадь боковой поверхности вписанной пирамиды равна $\tfrac12 (p+p_1)l$ (§~\ref{1938/s81}).

При неограниченном возрастании числа боковых граней вписанной пирамиды периметры $p$ и $p_1$ стремятся к пределам, принимаемым за длины $C$ и $C_1$ окружностей оснований, а апофема $l$ имеет пределом образующую $L$ усечённого конуса.
Следовательно, площадь боковой поверхности вписанной пирамиды стремится при этом к пределу, равному $\tfrac12(C + C_1)L$.
Этот предел и принимается за площадь боковой поверхности усечённого конуса;
обозначив её буквой $S$, будем иметь:
\[S=\tfrac12 (C+C_1)L.\]

\paragraph{}\label{1938/s117}
\so{Следствия}. 1) Если $R$ и $R_1$ означают радиусы окружностей нижнего и верхнего оснований, то площадь боковой поверхности усечённого конуса будет:
\[S = \tfrac12(2\pi R + 2\pi R_1)L = \pi (R + R_1)L.\]

2) Если в трапеции $OO_1A_1A$ (рис.~\ref{1938/s-ris-131}), от вращения которой получается усечённый конус, проведём среднюю линию $BC$, то получим:
\[BC = \tfrac12(OA + O_1A_1) = \tfrac12(R + R_1),\]
откуда
\[R + R_1=2\cdot BC\]
Следовательно,
\[S = 2\pi\cdot BC\cdot L,\]
то есть площадь боковой поверхности усечённого конуса равна произведению длины окружности среднего сечения на образующую.

3) Площадь $T$ полной поверхности усечённого конуса выразится так:
\[T = \pi (R^2 + R_1^2 + RL + R_1L).\]

\paragraph{Развёртка цилиндра и конуса.}\label{1938/s118}
Впишем в цилиндр (рис.~\ref{1938/s-ris-132}) какую-нибудь правильную призму и затем вообразим, что боковая её поверхность разрезана вдоль бокового ребра.
Очевидно, что, вращая её грани вокруг рёбер, мы можем \so{развернуть} эту поверхность в плоскую фигуру без разрыва и без складок.
Тогда получится то, что называется \rindex{развёртка поверхности}\textbf{развёрткой} боковой поверхности призмы.
\begin{figure}[!ht]
\vskip-0mm
\centering
\includegraphics{mppics/s-ris-132}
\caption{}\label{1938/s-ris-132}
\vskip-0mm
\end{figure}
Она представляет собой прямоугольник $KLMN$, составленный из стольких отдельных прямоугольников, сколько в призме боковых граней.
Основание его $MN$ равно периметру основания призмы, а высота $KN$ есть высота призмы.

Вообразим теперь, что число боковых граней вписанной призмы неограниченно удваивается;
тогда её развёртка будет всё удлиняться, приближаясь к предельному прямоугольнику $KPQN$, у которого длина основания равна длине окружности основания цилиндра, а высота есть высота цилиндра.
Этот прямоугольник называется \rindex{развёртка поверхности}\textbf{развёрткой} боковой поверхности цилиндра.

Подобно этому вообразим, что в конус вписана какая-нибудь правильная пирамида (рис.~\ref{1938/s-ris-133}).
Мы можем разрезать её боковую поверхность по одному из рёбер и затем, поворачивая грани вокруг рёбер, получить её плоскую развёртку в виде многоугольного сектора $SKL$, составленного из стольких равнобедренных треугольников, сколько в пирамиде боковых граней.
\begin{figure}[!ht]
\vskip-0mm
\centering
\includegraphics{mppics/s-ris-133}
\caption{}\label{1938/s-ris-133}
\vskip-0mm
\end{figure}
Отрезки $SK$, $Sa$, $Sb,\dots$ равны боковому ребру пирамиды (или образующей конуса), а длина ломаной $Kab\dots L$ равна периметру основания пирамиды.
При неограниченном удвоении числа боковых граней вписанной пирамиды развёртка её увеличивается, приближаясь к предельному сектору $SKM$, у которого длина дуги $KM$ равна длине окружности основания, а радиус $SK$ равен образующей конуса.
Этот сектор называется \rindex{развёртка поверхности}\textbf{развёрткой} боковой поверхности конуса.

Подобно этому можно получить развёртку боковой поверхности усечённого конуса (рис.~\ref{1938/s-ris-133}) в виде части кругового кольца $KMNP$.
Легко видеть, что площадь боковой поверхности цилиндра или конуса равна площади соответствующей развёртки.

\subsection*{Объём цилиндра и конуса}

\paragraph{}\label{1914/470} \so{Лемма} 1. 
\textbf{\emph{Объём цилиндра есть общий предел объёмов правильных вписанных и описанных призм при неограниченном удвоении числа их боковых граней.}}

Впишем в цилиндр и опишем около него по правильной призме с правильным $n$-угольником как основание.
Обозначим объём, площадь основания и высоту соответственно; 
для цилиндра — $V$, $B$, $H$, 
для вписанной призмы — $v_n$, $b_n$, $H$ 
и для описанной призмы — $V_n$, $B_n$, $H$.
Тогда будем иметь (§~\ref{1938/s88}): 
\begin{align*}
v_n&=b_n\cdot H;
&
V_n&=B_n\cdot H.
\end{align*}
Откуда:
\[V_n-v_n=(B_n-b_n)H.\]

При неограниченном удвоении $n$ разность: $B_n-b_n$ стремится к нулю (§~\ref{1938/262}), а множитель $H$ есть число постоянное.
Поэтому правая часть последнего равенства, и следовательно, и его левая часть, стремятся к нулю.

Объём цилиндра, очевидно, больше объёма вписанной призмы, но меньше объёма описанной.
Поэтому каждая из разностей $V-v_n$ и $V_n-V$ меньше разности $V_n-v_n$. 
Но последняя, по доказанному, стремится к нулю.
Следовательно, и первые две разности стремятся к нулю.
Это, по определению предела, и означает, что $V$ есть предел обеих последовательностей $v_n$ и $V_n$.

\paragraph{}\label{1914/471} \so{Лемма} 2. 
\textbf{\emph{Объём полного конуса есть общий предел объёмов правильных вписанных и описанных пирамид при неограниченном удвоении числа их боковых граней.}}

Впишем в конус и опишем около него по какой-нибудь пирамиде с правильным $n$-угольником как основание.
Так же как в §~\ref{1914/470}, обозначим объём, площадь основания и высоту соответственно; 
для конуса — $V$, $B$, $H$, 
для вписанной пирамиды — $v_n$, $b_n$, $H$ 
и для описанной пирамиды — $V_n$, $B_n$, $H$.
Тогда будем иметь (§~\ref{1938/s91}): 
\begin{align*}
v_n&=\tfrac13 b_n\cdot H;
&
V_n&=\tfrac13 B_n\cdot H.\end{align*}
Откуда:
\[V_n-v_n=\tfrac13(B_n-b_n)H.\]

Отсюда, точно также как и в §~\ref{1914/470}, можно заключить, что $V$ есть предел обоих последовательностей $v_n$ и $V_n$.

{\small
\medskip

\so{Замечание}. В доказанных леммах вписанные и описанные призмы и пирамиды предполагаются правильными только ради простоты доказательства. Содержание этих лемм остается в полной силе и тогда, когда призмы и пирамиды будут неправильные, лишь бы боковые грани их неограниченно уменьшались.

}

\paragraph{}\label{1938/s120}
\so{Теоремы}.
1) \textbf{\emph{Объём цилиндра равен произведению площади основания на высоту.}}

2) \textbf{\emph{Объём конуса равен произведению площади основания на треть высоты.}}

Впишем в цилиндр какую-нибудь правильную призму, а в конус — какую-нибудь правильную пирамиду с правильным $n$-угольником как основание.
Обозначим площадь основания призмы или пирамиды буквой $b_n$, высоту их буквой $H$ и объём — $v_n$, получим:

\medskip

\columnratio{0.5}
\setlength{\columnseprule}{.2pt}
\begin{paracol}{2}
для призмы $v_n = b_n\cdot H.$
\switchcolumn
для пирамиды $v_n= \tfrac13 b_n\cdot H.$
\end{paracol}

\medskip

Вообразим теперь, что число $n$ неограниченно удваивается.
Тогда $b_n$ будет иметь пределом площадь $B$ основания цилиндра или конуса, а высота $H$ остаётся без изменения;
значит, произведения $b_nH$ и $\tfrac13b_nH$ будут стремиться к пределам $BH$ и $\tfrac13BH$.
Согласно доказанным леммам (§~\ref{1914/470}, \ref{1914/471}), объём $V$ цилиндра или конуса будет:

\medskip

\columnratio{0.5}
\setlength{\columnseprule}{.2pt}
\begin{paracol}{2}
для цилиндра $V = B\cdot  H$.
\switchcolumn
для конуса $V= \tfrac13 B\cdot  H$.
\end{paracol}

\medskip

\paragraph{}\label{1938/s121}
\so{Следствие}.
Если радиус основания цилиндра или конуса обозначим через $R$, то $B= \pi R^2$, поэтому 

\medskip

\columnratio{0.5}
\setlength{\columnseprule}{.2pt}
\begin{paracol}{2}
$\text{объём цилиндра} = \pi R^2H$.
\switchcolumn
$\text{объём конуса} = \tfrac13\pi R^2H.$
\end{paracol}

\medskip

\paragraph{}\label{1938/s122}
\so{Теорема}.
\textbf{\emph{Объём, усечённого конуса равен сумме объёмов трёх конусов, имеющих одинаковую высоту с усечённым конусом, а основаниями: один — нижнее основание этого конуса, другой — верхнее, третий — круг, площадь которого есть среднее геометрическое между площадями верхнего и нижнего оснований.}}

Теорему эту докажем совершенно так же, как раньше мы доказали теорему для объёма усечённой пирамиды (§~\ref{1938/s92}).

На верхнем основании усечённого конуса (рис.~\ref{1938/s-ris-134}) поместим такой малый конус (с высотой $h$), который дополняет данный усечённый конус до полного.
Тогда объём $V$ усечённого конуса можно рассматривать как разность объёмов полного конуса и дополнительного.
Поэтому
\[V =\tfrac13 \pi R^2(H + h)-\tfrac13 \pi r^2h=\tfrac13 \pi [R^2H+(R^2-r^2)h].\]

\begin{wrapfigure}[10]{o}{54 mm}
\vskip-0mm
\centering
\includegraphics{mppics/s-ris-134}
\caption{}\label{1938/s-ris-134}
\vskip-0mm
\end{wrapfigure}

Из подобия треугольников находим:
\[\frac Rr=\frac{H+h}h\]
откуда получаем:
\begin{align*}
Rh &= rH + rh;
\\
(R-r)h &= rH;
\\
h&= \frac{rH}{R-r}.
\end{align*}

Поскольку 
\[R^2-r^2=(R+r)(R-r),\] 
получаем
\begin{align*}
V &= \tfrac13 \pi [R^2H+(R^2-r^2)h]=
\\
&= \tfrac13\pi [R^2H + (R +r)rH]=
\\
&=
\tfrac13\pi H(R^2 + Rr + r^2)
= 
\\
&=\tfrac13\pi R^2 H+\tfrac13\pi Rr H+\tfrac13\pi r^2 H.
\end{align*}

Так как $\pi R^2$ выражает площадь нижнего основания, $\pi r^2$ — площадь верхнего основания и $\pi Rr \z= \sqrt{\pi R^2\cdot \pi r^2}$ есть среднее геометрическое между площадями верхнего и нижнего оснований, то полученная нами формула подтверждает теорему.

\subsection*{Подобные цилиндры и конусы}

\paragraph{}\label{1938/s123}
\so{Определение}.
Два цилиндра или конуса называются подобными, если они произошли от вращения подобных прямоугольников или прямоугольных треугольников вокруг сходственных сторон.

Пусть (рис.~\ref{1938/s-ris-135} и \ref{1938/s-ris-136}) $h$ и $h_1$ будут высоты двух подобных цилиндров или конусов, $r$ и $r_1$ — радиусы их оснований, $l$ и $l_1$ — образующие;
тогда согласно определению
\[\frac r{r_1}=\frac h{h_1}\quad\text{и}\quad\frac r{r_1}=\frac l{l_1}\]
откуда (по свойству равных пропорций) находим:
\[\frac {r+h}{r_1+h_1}=\frac r{r_1}
\quad\text{и}\quad
\frac {r+l}{r_1+l_1}=\frac r{r_1}\]

\begin{figure}[!ht]
\begin{minipage}{.48\textwidth}
\centering
\includegraphics{mppics/s-ris-135}
\end{minipage}\hfill
\begin{minipage}{.48\textwidth}
\centering
\includegraphics{mppics/s-ris-136}
\end{minipage}

\medskip

\begin{minipage}{.48\textwidth}
\centering
\caption{}\label{1938/s-ris-135}
\end{minipage}\hfill
\begin{minipage}{.48\textwidth}
\centering
\caption{}\label{1938/s-ris-136}
\end{minipage}
\vskip-4mm
\end{figure}

Заметив эти пропорции, докажем следующую теорему:

\paragraph{}\label{1938/s124}
\so{Теорема}.
\textbf{\emph{Площади боковых и полных поверхностей подобных цилиндров или конусов относятся, как квадраты радиусов или высот;
объёмы — как кубы радиусов или высот.}}

Пусть $S$, $T$ и $V$ будут площади соответственно боковой поверхности, полной поверхности и объём одного цилиндра или конуса;
$S_1$, $T_1$ и $V_1$ — те же величины для другого цилиндра или конуса, подобного первому.
Тогда будем иметь для цилиндров:
\begin{align*}
\frac S{S_1}&=
\frac{2\pi rh}{2\pi r_1h_1}=
\frac{rh}{r_1h_1}=
\frac{r}{r_1}\cdot \frac{h}{h_1}=
\frac{r^2}{r_1^2}=
\frac{h^2}{h_1^2};
\\
\frac T{T_1}&=
\frac{2\pi r(r+h)}{2\pi r_1(r_1+h_1)}=
\frac{r}{r_1}\cdot\frac{r+h}{r_1+h_1}=
\frac{r^2}{r_1^2}=
\frac{h^2}{h_1^2};
\\
\frac V{V_1}&=
\frac{\pi r^2h}{\pi r_1^2h_1}=
\frac{r^2}{r_1^2}\cdot \frac{h}{h_1}=
\frac{r^3}{r_1^3}=
\frac{h^3}{h_1^3};
\end{align*}
для конусов
\begin{align*}
\frac S{S_1}&=
\frac{\pi rl}{\pi r_1l_1}=
\frac{rl}{r_1l_1}=
\frac{r}{r_1}\cdot \frac{l}{l_1}=
\frac{r^2}{r_1^2}=
\frac{l^2}{l_1^2};
\\
\frac T{T_1}&=
\frac{\pi r(r+l)}{\pi r_1(r_1+l_1)}=
\frac{r}{r_1}\cdot\frac{r+l}{r_1+l_1}=
\frac{r^2}{r_1^2}=
\frac{l^2}{l_1^2};
\\
\frac V{V_1}&=
\frac{\frac13\pi r^2h}{\frac13\pi r_1^2h_1}=
\frac{r^2}{r_1^2}\cdot \frac{h}{h_1}=
\frac{r^3}{r_1^3}=
\frac{h^3}{h_1^3}.
\end{align*}

\section{Шар}

\subsection*{Сечение шара плоскостью}

\paragraph{}\label{1938/s125}
\so{Определение}.
Тело, происходящее от вращения полукруга вокруг диаметра, называется \rindex{шар}\textbf{шаром}, а поверхность, образуемая при этом полуокружностью, называется \rindex{сфера}\textbf{сферой} или шаровой поверхностью.
Можно также сказать, что сфера есть геометрическое место точек, одинаково удалённых от одной и той же точки (называемой \rindex{центр!сферы}\textbf{центром} сферы).

Отрезок, соединяющий центр с какой-нибудь точкой сферы, называется \rindex{радиус}\textbf{радиусом}, а отрезок, соединяющий две точки сферы и проходящий через её центр, называется \rindex{диаметр}\textbf{диаметром}.
Все радиусы одного шара равны между собой;
всякий диаметр равен двум радиусам.

Два шара одинакового радиуса равны, потому что при вложении они совмещаются.

\paragraph{}\label{1938/s126}
\so{Теорема}.
\textbf{\emph{Всякое сечение шара плоскостью есть круг.}}

\begin{wrapfigure}{o}{38 mm}
\vskip-0mm
\centering
\includegraphics{mppics/s-ris-137}
\caption{}\label{1938/s-ris-137}
\vskip-0mm
\end{wrapfigure}

1) Предположим сначала, что (рис.~\ref{1938/s-ris-137}) секущая плоскость $AB$ проходит через центр $O$ шара.
Все точки линии пересечения принадлежат сфере и поэтому одинаково удалены от точки $O$, лежащей в секущей плоскости;
следовательно, сечение есть круг с центром в точке $O$.

2) Положим теперь, что секущая плоскость $CD$ не проходит через центр.
Опустим на неё из центра перпендикуляр $OK$ и возьмём на линии пересечения какую-нибудь точку $M$.
Соединив её с $O$ и $K$, получим прямоугольный треугольник $MOK$, из которого находим:
\[MK = \sqrt{OM^2 - OK^2}. \eqno(1)\]

Так как длины отрезков $OM$ и $OK$ не изменяются при изменении положений точки $M$ на линии пересечения, то расстояние $MK$ есть величина постоянная для данного сечения;
значит, линия пересечения есть окружность, центр которой есть точка $K$.

\paragraph{}\label{1938/s127}
\so{Следствие}.
Пусть $R$ и $r$ будут длины радиуса шара и радиуса круга сечения, a $d$ — расстояние секущей плоскости от центра, тогда равенство (1) примет вид: 
\[r= \sqrt{R^2 - d^2}.\]

Из этой формулы выводим:

1) \emph{Наибольший радиус сечения получается при $d=0$, то есть когда секущая плоскость проходит через центр шара.}
В этом случае $r=R$.
Круг, получаемый в этом случае, называется \rindex{большой круг}\textbf{большим кругом}.

2) \emph{Наименьший радиус сечения получается при $d=R$.}
В этом случае $r=0$, то есть круг сечения вырождается в точку.

3) \emph{Сечения, равноотстоящие от центра шара, равны.}

4) \emph{Из двух сечений, неодинаково удалённых от центра шара, то, которое ближе к центру, имеет больший радиус.}

\paragraph{}\label{1938/s128}
\so{Теорема}.
\textbf{\emph{Всякая плоскость}} ($P$, рис.~\ref{1938/s-ris-138}), \textbf{\emph{проходящая через центр сферы, делит её на две зеркально симметричные и равные части.}}

\begin{wrapfigure}{o}{53 mm}
\vskip-0mm
\centering
\includegraphics{mppics/s-ris-138}
\caption{}\label{1938/s-ris-138}
\vskip-0mm
\end{wrapfigure}

Возьмём на сфере какую-нибудь точку $A$;
пусть $AB$ есть перпендикуляр, опущенный из точки $A$ на плоскость~$P$.
Продолжим $AB$ до пересечения со сферой в точке $C$.
Проведя $BO$, мы получим два равных прямоугольных треугольника $AOB$ и $BOC$ (общий катет $BO$, а гипотенузы равны, как радиусы сферы);
следовательно, $AB=BC$, таким образом, всякой точке $A$ сферы соответствует другая точка $C$ сферы, зеркально симметричная относительно плоскости $P$ с точкой $A$.
Значит, плоскость $P$ делит сферу на две зеркально симметричные части.

Эти части не только зеркально симметричны, но и равны, так как, разрезав сферу по плоскости $P$, мы можем совместить эти части повернув одну из них на угол $180\degree$ вокруг любой прямой плоскости $P$, проходящей через~$O$.

\begin{wrapfigure}[9]{r}{43 mm}
\vskip-8mm
\centering
\includegraphics{mppics/s-ris-139}
\caption{}\label{1938/s-ris-139}
\vskip-0mm
\end{wrapfigure}

\paragraph{}\label{1938/s129}
\mbox{\so{Теорема}.}
\textbf{\emph{Через две точки сферы, не лежащие на концах одного диаметра, можно провести окружность большого круга и только одну.}}

Пусть на сфере (рис.~\ref{1938/s-ris-139}), имеющей центр $O$, взяты какие-нибудь две точки, например $C$ и $N$, не лежащие на одной прямой с точкой $O$.
Тогда через точки $C$, $O$ и $N$ можно провести плоскость.
Эта плоскость, проходя через центр $O$, даст в пересечении со сферой окружность большого круга.

Другой окружности большого круга через те же две точки $C$ и $N$ провести нельзя.
Действительно, всякая окружность большого круга должна, по определению, лежать в плоскости, проходящей через центр сферы;
следовательно, если бы через $C$ и $N$ можно было провести ещё другую окружность большого круга, тогда выходило бы, что через три точки $C$, $N$ и $O$, не лежащие на одной прямой, можно провести две различные плоскости, что невозможно.

\paragraph{}\label{1938/s130}
\so{Теорема}.
\textbf{\emph{Окружности двух больших кругов при пересечении делятся пополам.}}

Центр $O$ (рис.~\ref{1938/s-ris-139}), находясь на плоскостях обоих больших кругов, лежит на прямой, по которой эти круги пересекаются;
значит, эта прямая есть диаметр того и другого круга, а диаметр делит окружность пополам.

\subsection*{Плоскость, касательная к сфере}

\paragraph{}\label{1938/s131}
\so{Определение}.
Плоскость, имеющая со сферой только одну общую точку, называется касательной плоскостью.
Возможность существования такой плоскости доказывается следующей теоремой.

\paragraph{}\label{1938/s132}
\so{Теорема}.
\textbf{\emph{Плоскость}} ($P$, рис.~\ref{1938/s-ris-140}), \textbf{\emph{перпендикулярная к радиусу}} ($AO$) \textbf{\emph{в конце его, лежащем на сфере, есть касательная плоскость.}}

\begin{wrapfigure}{r}{43 mm}
\vskip-2mm
\centering
\includegraphics{mppics/s-ris-140}
\caption{}\label{1938/s-ris-140}
\bigskip
\includegraphics{mppics/s-ris-504}
\caption{}\label{1938/s-ris-504}
\end{wrapfigure}

Возьмём на плоскости $P$ произвольную точку $B$ и проведём прямую $OB$.
Так как $OB$ — наклонная, а $OA$ — перпендикуляр к плоскости $P$, то $OB>OA$.
Поэтому точка $B$ лежит вне сферы;
следовательно, у плоскости $P$ есть только одна общая точка $A$ со сферой;
значит, эта плоскость касательная.

\paragraph{}\label{1938/s133}
\mbox{\so{Обратная теорема}.}
\textbf{\emph{Касательная плоскость}} ($P$, рис. \ref{1938/s-ris-140}) \textbf{\emph{перпендикулярна к радиусу}} ($OA$), \textbf{\emph{проведённому в точку касания.}}

Проведём плоскость $Q$ через центр $O$ и произвольную прямую $AB$ на плоскости~$P$.
Заметим, что $Q$ пересекает сферу по окружности большого круга и прямая $AB$ имеет с этим кругом ровно одну общую точку — точку $A$ (рис.~\ref{1938/s-ris-504}).
То есть прямая $AB$ касается окружности большого круга в плоскости $Q$ (§~\ref{1938/113}) и значит $AB\perp OA$.

Точно так же докажем, что $AC\z\perp OA$ для другой прямой $AC$ на плоскости $P$
и значит $P\perp OA$ (§~\ref{1938/s23}).

Прямая, имеющая ровно одну общую точку со сферой, называется \so{касательной} к сфере.
Легко видеть, что существует бесчисленное множество прямых, касающихся сферы в данной точке.
Действительно, всякая прямая ($AC$, рис.~\ref{1938/s-ris-140}), лежащая в плоскости, касательной к сфере, в данной точке ($A$) и проходящая через точку касания ($A$), есть касательная к сфере в этой точке.

\subsection*{Сфера и её части} 

\begin{wrapfigure}{r}{34 mm}
\vskip-6mm
\centering
\includegraphics{mppics/s-ris-141}
\caption{}\label{1938/s-ris-141}
\vskip-0mm
\end{wrapfigure}

\paragraph{}\label{1938/s134}
\mbox{\so{Определения}.}
1) Часть сферы (рис. \ref{1938/s-ris-141}), отсекаемая от неё какой-нибудь плоскостью ($AA_1$), называется \rindex{сегмент}\rindex{сферический!сегмент}\textbf{сферическим сегментом}.

Окружность $AA_1$ называется \rindex{основание!сферического сегмента}\textbf{основанием}, а отрезок $KM$ радиуса, перпендикулярного к плоскости сечения, — \rindex{высота!сферического сегмента}\textbf{высотой} сферического сегмента.

2) Часть сферы, заключённая между двумя параллельными секущими плоскостями ($AA_1$ и $BB_1$), называется \rindex{сферический!пояс}\textbf{сферическим поясом}.

Окружности сечения $AA_1$ и $BB_1$ называются \rindex{основание!сферического пояса}\textbf{основаниями}, а расстояние $KL$ между параллельными плоскостями — \rindex{высота!пояса}\textbf{высотой} пояса.

Сферический пояс и сегмент можно рассматривать как поверхности вращения, в то время как полуокружность $MABN$, вращаясь вокруг диаметра $MN$, описывает сферу, часть её $AB$ описывает пояс, а часть $MA$ — сегмент.

Для нахождения площади сферы и её частей мы докажем следующую лемму:

\paragraph{}\label{1938/s135}
\mbox{\so{Лемма}.}
\textbf{\emph{Площадь боковой поверхности каждого из трёх тел: конуса, усечённого конуса и цилиндра — равна произведению высоты тела на длину окружности, у которой радиус есть перпендикуляр, восстановленный к образующей из её середины до пересечения с осью.}}

1) Пусть конус образуется (рис.~\ref{1938/s-ris-142}) вращением треугольника $ABC$ вокруг катета $AC$.
Если $D$ есть середина образующей $AB$, то (§~\ref{1938/s115})
\[\text{боковая поверхность конуса} = 2\pi \cdot BC\cdot AD. \eqno(1)\]

Проведя $DE\perp AB$, получим два подобных треугольника $ABC$ и $AED$ (они прямоугольные и имеют общий угол $A$);
из их подобия выводим:
\[\frac{BC}{ED} = \frac{AC}{AD},\]
откуда
\[BC\cdot AD = ED\cdot AC,\]
и равенство (1) даёт:
\[\text{боковая поверхность конуса} = 2\pi \cdot ED\cdot AC,\]
что и требовалось доказать.

\begin{figure}
\begin{minipage}{.32\textwidth}
\centering
\includegraphics{mppics/s-ris-142}
\end{minipage}\hfill
\begin{minipage}{.32\textwidth}
\centering
\includegraphics{mppics/s-ris-143}
\end{minipage}\hfill
\begin{minipage}{.32\textwidth}
\centering
\includegraphics{mppics/s-ris-144}
\end{minipage}

\medskip

\begin{minipage}{.32\textwidth}
\centering
\caption{}\label{1938/s-ris-142}
\end{minipage}
\begin{minipage}{.32\textwidth}
\vfill
\centering
\caption{}\label{1938/s-ris-143}
\end{minipage}
\begin{minipage}{.32\textwidth}
\vfill
\centering
\caption{}\label{1938/s-ris-144}
\end{minipage}
\vskip-4mm
\end{figure} 

2) Пусть усечённый конус (рис.~\ref{1938/s-ris-143}) образуется вращением трапеции $ABCD$ вокруг стороны $AD$.

Проведя среднюю линию $EF$, будем иметь (§~\ref{1938/s117}):
\[\text{боковая поверхность усечённого конуса} = 2\pi \cdot EF\cdot BC. \eqno(2)\]

Проведём $EG\perp BC$ и $BH\perp DC$.
Тогда получим два подобных треугольника $EFG$ и $BHC$ (стороны одного перпендикулярны к сторонам другого);
из их подобия выводим:
\[\frac{EF}{BH} = \frac{EG}{BC},\]
откуда
\[EF\cdot BC = BH\cdot EG = AD \cdot EG.\]
Поэтому равенство (2) можно записать так:
\[\text{боковая поверхность усечённого конуса} = 2\pi \cdot EG\cdot AD,\]
что и требовалось доказать.

3) Теорема остаётся верной и в применении к цилиндру, так как окружность, о которой говорится в теореме, равна окружности основания цилиндра.

\paragraph{}\label{1938/s136} 
Используя рассуждения в §~\ref{1938/s111}, легко прийти к следующему утверждению, которое можно принять за \so{определение} площади сферического пояса:

\emph{За площадь поверхности сферического пояса, образуемого вращением} (рис.~\ref{1938/s-ris-144}) \emph{какой-нибудь дуги} ($BE$) \emph{полуокружности вокруг диаметра} ($AF$), \emph{принимают предел, к которому стремится площадь поверхности, образуемая вращением вокруг того же диаметра правильной вписанной ломаной линии} ($BCDE$), \emph{когда её стороны неограниченно уменьшаются} (и, следовательно, число сторон неограниченно увеличивается).

Это утверждение распространяется и на сферический сегмент, и на всю сферу;
в последнем случае ломаная линия вписывается в целую полуокружность.

\paragraph{}\label{1938/s137}
\mbox{\so{Теоремы}.}
1) \textbf{\emph{Площадь сферического сегмента равна произведению его высоты на длину окружности большого круга.}}

2) \textbf{\emph{Площадь сферического пояса равна произведению его высоты на длину окружности большого круга.}}

Впишем в дугу $AF$ (рис.~\ref{1938/s-ris-145}), образующую при вращении сферический сегмент, правильную ломаную линию $ACDEF$ с произвольным числом сторон.

Поверхность, получающаяся от вращения этой ломаной, состоит из частей, образуемых вращением сторон $AC$, $CD$, $DE$ и так далее.
Эти части представляют собой боковые поверхности или полного конуса (от вращения $AC$), или усечённого конуса (от вращения $CD$, $EF,\dots$), или цилиндра (от вращения $DE$), если $DE\parallel AB$.
Поэтому мы можем применить к ним лемму §~\ref{1938/s135}.
При этом заметим, что каждый из перпендикуляров, восстановленных из середин образующих до пересечения с осью, равен апофеме ломаной линии.
Обозначив эту апофему буквой $a$, получим:
\begin{align*}
\text{поверхность, образованная вращением}&\ AC = Ac \cdot 2\pi a;
\\
\text{—\textquotedbl—\qquad\qquad—\textquotedbl—\qquad\qquad—\textquotedbl—\qquad}&\ CD = cd \cdot 2\pi a;
\\
\text{—\textquotedbl—\qquad\qquad—\textquotedbl—\qquad\qquad—\textquotedbl—\qquad}&\ DE = de \cdot 2\pi a;
\\
\text{и так далее\qquad\qquad}&
\end{align*}

Сложив эти равенства почленно, получим:
\[\text{поверхность вращения}\ ACDEF = Af \cdot 2\pi a.\]

При неограниченном увеличении числа сторон вписанной ломаной апофема $a$ стремится к пределу, равному радиусу сферы $R$, а отрезок $Af$ остаётся без изменения;
следовательно, предел площади поверхности, образованной вращением $ACDEF$ равен $Af\cdot 2\pi R$.
Но предел площади поверхности, образованной вращением $ACDEF$, принимают за площадь сферического сегмента, а отрезок $Af$ есть высота $H$ сегмента;
поэтому
\[\text{площадь сферического сегмента} = H\cdot 2\pi R = 2\pi RH.\]

\begin{wrapfigure}{r}{34 mm}
\vskip-0mm
\centering
\includegraphics{mppics/s-ris-145}
\caption{}\label{1938/s-ris-145}
\vskip-0mm
\end{wrapfigure}

2) Предположим, что правильная ломаная линия вписана не в дугу $AF$, образующую сферический сегмент, а в какую-нибудь дугу $CF$, образующую сферический пояс (рис.~\ref{1938/s-ris-145}).
Это изменение, как легко видеть, нисколько не влияет на ход предыдущих рассуждений, поэтому и вывод остаётся тот же, то есть что
\[\text{площадь сферического пояса} = H \cdot 2\pi R = 2\pi RH,\]
где буквой $H$ обозначена высота $cf$ сферического пояса.

\paragraph{}\label{1938/s138}
\mbox{\so{Теорема}.}
\textbf{\emph{Площадь сферы равна произведению длины окружности большого круга на диаметр,}} или: \textbf{\emph{площадь сферы равна учетверённой площади большого круга.}}

Сфера, образуемая вращением полуокружности $ADB$ (рис. \ref{1938/s-ris-145}) состоит из поверхностей, образуемых вращением дуг $AD$ и $DB$.
Поэтому согласно предыдущей теореме можно написать:
\begin{align*}
\text{площадь сферы} &= 2\pi R\cdot Ad + 2\pi R\cdot dB =
\\
&=2\pi R(Ad + dB) =
\\
&= 2\pi R\cdot 2R = 
\\
&=4\pi R^2.
\end{align*}

\paragraph{}\label{1938/s139}
\so{Следствие}.
\emph{Площади сфер относятся, как квадраты их радиусов или диаметров}, потому что, обозначая через $R$ и $R_1$ радиусы, а через $S$ и $S_1$ площади двух сфер, будем иметь:
\begin{align*}
\frac{S}{S_1} &= \frac{4\pi R^2}{4\pi R^2_1} =\frac{R^2}{R^2_1} = 
\\&=\frac{4R^2}{4R_1^2} =\frac{(2R)^2}{(2R_1)^2}.
\end{align*}

\subsection*{Объём шара и его частей}

{

\begin{wrapfigure}{r}{38 mm}
\vskip-14mm
\centering
\includegraphics{mppics/s-ris-146}
\caption{}\label{1938/s-ris-146}
\vskip-0mm
\end{wrapfigure}

\paragraph{}\label{1938/s140}
\mbox{\so{Определение}.}
Тело, получаемое от вращения (рис.~\ref{1938/s-ris-146}) кругового сектора ($AOB$) вокруг одной из его сторон ($OA$) называется \rindex{сектор}\rindex{шаровой!сектор}\textbf{шаровым сектором}.
Это тело ограничено боковой поверхностью конуса и сферическим сегментом.

}

\paragraph{}\label{1938/s141}
Для нахождения объёма шарового сектора и целого шара мы предварительно докажем следующую лемму.

\medskip

\mbox{\so{Лемма}.}
\textbf{\emph{Если $\triangle ABC$}} (рис.~\ref{1938/s-ris-147}) \textbf{\emph{вращается вокруг оси $xy$, которая лежит в плоскости треугольника, проходит через его вершину $A$, но не пересекает стороны $BC$, то объём тела, получаемого при этом вращении, равен произведению площади поверхности, образуемой противоположной стороной $BC$, на одну треть высоты $h$, опущенной на эту сторону.}}

\begin{figure}[!ht]
\begin{minipage}{.32\textwidth}
\centering
\includegraphics{mppics/s-ris-147}
\end{minipage}
\hfill
\begin{minipage}{.32\textwidth}
\centering
\includegraphics{mppics/s-ris-148}
\end{minipage}
\hfill
\begin{minipage}{.32\textwidth}
\centering
\includegraphics{mppics/s-ris-149}
\end{minipage}

\medskip

\begin{minipage}{.32\textwidth}
\centering
\caption{}\label{1938/s-ris-147}
\end{minipage}
\hfill
\begin{minipage}{.32\textwidth}
\centering
\caption{}\label{1938/s-ris-148}
\end{minipage}
\hfill
\begin{minipage}{.32\textwidth}
\centering
\caption{}\label{1938/s-ris-149}
\end{minipage}
\vskip-4mm
\end{figure}

При доказательстве рассмотрим три случая:

1) Ось совпадает со стороной $AB$ (рис.~\ref{1938/s-ris-148}).
В этом случае искомый объём равен сумме объёмов двух конусов, получаемых вращением прямоугольных треугольников $BCD$ и $DCA$.
Первый объём равен $\tfrac13 \pi CD^2\cdot DB$, а второй $\tfrac13\pi CD^2\cdot DA$;
поэтому объём, образованный вращением $ABC$, равен 
\[\tfrac13\pi CD^2(DB + DA) = \tfrac13\pi CD\cdot CD\cdot BA.\]

Заметим, что $CD\cdot BA=BC\cdot h$, так как каждое из этих произведений выражает удвоенную площадь $\triangle ABC$;
поэтому
\[\text{объём}\, ABC = \tfrac13\pi CD\cdot BC\cdot h.\]

Но произведение $\pi CD\cdot BC$ равно площади боковой поверхности конуса $BDC$;
значит,
\[\text{объём}\, ABC = (\text{поверхность}\, BC)\cdot h.\]

2) Ось не совпадает с $AB$ и не параллельна $BC$ (рис.~\ref{1938/s-ris-149}).
В этом случае искомый объём равен разности объёмов тел, производимых вращением треугольников $AMC$ и $AMB$.
По доказанному в первом случае
\begin{align*}
\text{объём}\, AMC&= h\cdot (\text{поверхность}\, MC),
\\
\text{объём}\, AMB &= h\cdot (\text{поверхность}\, MB),
\intertext{следовательно,}
\text{объём}\, ABC &= h\cdot (\text{поверхность}\, MC-\text{поверхность}\, MB)=
\\&=h\cdot (\text{поверхность}\, BC).
\end{align*}

\begin{wrapfigure}{o}{33 mm}
\vskip0mm
\centering
\includegraphics{mppics/s-ris-150}
\caption{}\label{1938/s-ris-150}
\vskip-0mm
\end{wrapfigure}

3) Ось параллельна стороне $BC$ (рис.~\ref{1938/s-ris-150}).
Тогда искомый объём равен объёму цилиндра, производимому вращением прямоугольника $DEBC$ без суммы объёмов конусов, производимых вращением треугольников $AEB$ и $ACD$;
первый из них равен $\pi DC^2\cdot ED$;
второй — $\tfrac13\pi EB^2\cdot EA$ 
и третий — $\tfrac13\pi DC^2\cdot AD$.
Приняв во внимание, что $EB=DC$, получим:
\begin{align*}
\text{объём}\,ABC &= \pi DC^2[ED-\tfrac13(EA + AD)]=
\\
&=\pi DC^2[ED-\tfrac13 ED]=
\\
&= \tfrac23\pi DC^2\cdot ED.
\end{align*}
Произведение $2\pi DC\cdot ED$ выражает площади боковой поверхности цилиндра, образуемую стороной $BC$;
поэтому
\begin{align*}
\text{объём}\, ABC&= (\text{поверхность}\, BC)\cdot \tfrac13 DC
\\&=(\text{поверхность}\, BC)\cdot \tfrac13 h.
\end{align*}

\paragraph{}\label{1938/s143} 
\so{Теорема}.
\textbf{\emph{Объём шарового сектора равен произведению площади его сферического сегмента на треть радиуса.}}

Пусть шаровой сектор производится вращением вокруг стороны $OA$ (рис.~\ref{1938/s-ris-151}) сектора $AOD$.

\begin{wrapfigure}{o}{33 mm}
\vskip-0mm
\centering
\includegraphics{mppics/s-ris-151}
\caption{}\label{1938/s-ris-151}
\vskip-0mm
\end{wrapfigure}

Впишем в дугу $AD$ правильную ломаную линию $ABCD$ с произвольным числом сторон равным~$n$.
Затем, продолжив конечные радиусы $OA$ и $OD$, опишем около дуги $AD$ правильную ломаную $A_1B_1C_1D_1$ стороны которой параллельны сторонам вписанной ломаной.
Многоугольники $OABCD$ и $OA_1B_1C_1D_1$ произведут при вращении некоторые тела, объёмы которых обозначим: первого через $v_n$, а второго через $V_n$.

Докажем прежде всего, что при неограниченном удвоении числа $n$ разность $V_n-v_n$ стремится к нулю.

Объем $v_n$ есть сумма объёмов, получаемых вращением треугольников $OAB$, $OBC$, $OCD$ вокруг оси $OA$.
Объем $V_n$ есть сумма, объёмов, получаемых вращением вокруг той же оси треугольников $OA_1B_1$, $OB_1C_1$ $OC_1D_1$.
Применим к этим объёмам лемму в §~\ref{1938/s141}, причём заметим, что высоты первых треугольников равны апофеме $a_n$ вписанной ломаной, а высоты вторых треугольников равны радиусу $R$ шара.
Согласно этой лемме будем иметь:
\begin{align*}
v_n&=
(\text{поверхность}\, AB)\cdot\tfrac {a_n}3
+
(\text{поверхность}\, BC)\cdot\tfrac {a_n}3
+
\dots=
\\
&=
(\text{поверхность}\, ABCD)\cdot \tfrac {a_n}3.
\\
V_n&=
(\text{поверхность}\, AB)\cdot\tfrac R3
+
(\text{поверхность}\, BC)\cdot\tfrac R3
+
\dots=
\\
&=
(\text{поверхность}\, ABCD)\cdot \tfrac R3.
\end{align*}

Вообразим теперь, что число $n$ неограниченно удваивается.
При этом условии площади поверхностей $ABCD$ и $A_1B_1C_1D_1$ стремятся к общему пределу, именно к площади сферического сегмента $AD$ (§~\ref{1938/s136}), а апофема $a_n$ имеет пределом радиус $R$.
Следовательно, объёмы $v_n$ и $V_n$ стремятся при этом к общему пределу, именно к произведению 
\[(\text{поверхность сегмента}\,  AD)\cdot \tfrac R3.\]
То есть обе последовательности объёмов $v_n$ и $V_n$ приближаются к одной и той же постоянной величине как угодно близко; это возможно только тогда, когда разность $V_n-v_n$ стремится к $0$.

Обозначим буквою $V$ объём шарового сектора $OAD$.
Очевидно, что $V>v_n$ и $V<V_n$.
Значит, каждая из разностей $V_n-V$ и $V-v_n$ меньше разности $V_n-v_n$.
Ho эта разность, как мы видели, при неограниченном удвоении числа сторон ломаных стремится к $0$. Следовательно, разности $V_n-V$ и $V-v_n$ и подавно при этом стремятся к $0$.
Отсюда заключаем, что постоянная величина $V$ есть общий предел последовательностей объёмов $v_n$ и $V_n$.
Ho этот общий предел, как мы нашли, есть произведение $(\text{поверхность сегмента}\,  AD)\cdot \tfrac R3$.
Значит:
\[V=(\text{поверхность сегмента}\,  AD)\cdot \tfrac R3.\]

\begin{wrapfigure}{r}{48 mm}
\vskip-0mm
\centering
\includegraphics{mppics/s-ris-152}
\caption{}\label{1938/s-ris-152}
\vskip-0mm
\end{wrapfigure}

\paragraph{}\label{1938/s144}
\mbox{\so{Теорема}.}
\textbf{\emph{Объём шара равняется произведению площади его поверхности на треть радиуса.}}

Разбив полукруг $ABC$ (рис.~\ref{1938/s-ris-152}), производящий шар, на пару круговых секторов $AOB$ и $BOC$, мы заметим, что объём шара можно рассматривать как сумму объёмов шаровых секторов, производимых вращением этих круговых секторов.
Согласно предыдущей теореме
\begin{align*}
\text{объём}\, AOB &= (\text{поверхность}\, AB)\cdot \tfrac13 R,
\\
\text{объём}\, BOC &= (\text{поверхность}\, BC) \cdot \tfrac13 R.
\end{align*}
Значит
\begin{align*}
\text{объём шара} &= (\text{поверхность}\, AB+\text{поверхность}\, BC)\cdot \tfrac13 R=
\\
&=(\text{поверхность}\, ABC)\cdot \tfrac13 R.
\end{align*}

{\small
\medskip

\so{Замечание}.
Можно и непосредственно рассматривать объём шара как объём тела, образованного вращением вокруг диаметра кругового сектора, центральный угол которого равен $180\degree$.

В таком случае объём шара можно получить как частный случай объёма шарового сектора, чей сферический сегмент составляет всю сферу.
В силу предыдущей теоремы \emph{объём шара будет при этом равен площади его поверхности, умноженной на треть радиуса.}

}

\paragraph{}\label{1938/s145}
\so{Следствие} 1 .
Обозначим высоту сферического пояса или сегмента через $H$, радиус сферы — через $R$, а диаметр — через $D$;
тогда площадь пояса или сегмента выразится, как мы видели (§~\ref{1938/s137}), формулой $2\pi RH$, а площадь сферы (§~\ref{1938/s138}) — формулой $4\pi R^2$;
поэтому

\begin{align*}
\text{объём шарового сектора}&= 2\pi RH\cdot \tfrac13 R =
\\
&= \tfrac23 \pi R^2 H;
\\
\text{объём шара} &= 4\pi R^2\cdot \tfrac13 R = 
\\
&=\tfrac43\pi R^3=
\\
&= \tfrac43\pi\left(\frac D2\right)^3=
\\
&= \tfrac16\pi D^3.
\end{align*}
Отсюда видно, что объёмы шаров относятся, как кубы их радиусов или диаметров.

Объём шара может быть выведен (не вполне, впрочем, строго) следующим простым рассуждением.
Вообразим, что вся поверхность шара разбита на очень малые участки и что все точки контура каждого участка соединены радиусами с центром шара.
Тогда шар разделится на очень большое число маленьких тел, из которых каждое можно рассматривать как пирамиду с вершиной в центре шара.
Так как объём пирамиды равен произведению площади поверхности основания на третью часть высоты (которую можно принять равной радиусу шара), то объём шара, равный, очевидно, сумме объёмов всех пирамид, выразится так:
\[\text{объём шара} = S \cdot \tfrac13 R,\]
где $S$ — сумма площадей оснований всех пирамид.
Но эта сумма должна составить площадь сферы, и, значит,
\[\text{объём шара} = 4\pi R^2\cdot \tfrac13 R = \tfrac43 \pi R^3.\]
Таким образом, объём шара может быть найден через площадь его поверхности.
Обратно, площадь сферы может быть найдена с помощью формулы его объёма из равенства:
\[S\cdot \tfrac13 R = \tfrac43 \pi R^3,\quad\text{откуда}\quad S = 4\pi R^2.\]

\paragraph{}\label{1938/s146}
\so{Следствие} 2.
\emph{Площадь поверхности шара и объём шара соответственно составляют $\tfrac23$ площади полной поверхности и объёма цилиндра, описанного около шара.}

Действительно, у цилиндра, описанного около шара, радиус основания равен радиусу шара, а высота равна диаметру шара;
поэтому для такого цилиндра
\begin{align*}
\text{полная поверхность описанного цилиндра} &= 2\pi R \cdot 2R + 2\pi R^2 = 6\pi R^2,
\\
\text{объём описанного цилиндра} &= \pi R^2\cdot 2R = 2\pi R^3.
\end{align*}
Отсюда видно, что $\tfrac23$ площади полной поверхности этого цилиндра равны $4\pi R^2$, то есть равны площади поверхности шара, а $\tfrac23$ объёма цилиндра составляют $\tfrac43\pi R^3$, то есть объём шара.

Это предложение было доказано Архимедом (в III веке до начала нашей эры).
Архимед выразил желание, чтобы чертёж;
этой теоремы был изображён на его гробнице, что и было исполнено римским военачальником Марцеллом (Ф. Кэджори. История элементарной математики).

Предлагаем учащимся как полезное упражнение доказать, что площадь поверхности и объём шара составляют $\tfrac49$ соответственно площади полной поверхности и объёма описанного конуса, у которого образующая равна диаметру основания.
Соединяя это предложение с указанным в следствии 2, мы можем написать такое равенство, где $Q$ обозначает площадь поверхности или объём:
\[\frac{Q_{\text{шара}}}4
=\frac{Q_{\text{цилиндра}}}6
=\frac{Q_{\text{конуса}}}9.
\]

{\small

\paragraph{}\label{1938/s147}
\so{Замечание}.
Формулу для объёма шара можно весьма просто получить, основываясь на принципе Кавальери (§~\ref{1938/s89}), следующим образом.

\begin{figure}[!ht]
\vskip-0mm
\centering
\includegraphics{mppics/s-ris-153}
\caption{}\label{1938/s-ris-153}
\vskip-0mm
\end{figure}

Пусть на одной и той же плоскости $H$ (рис.~\ref{1938/s-ris-153}) помещены шар радиуса $R$ и цилиндр, радиус основания которого равен $R$, а высота $2R$ (значит, это такой цилиндр, который может быть описан около шара радиуса $R$).
Вообразим далее, что из цилиндра вырезаны и удалены два конуса, имеющие общую вершину на середине $a$ оси цилиндра, и основания — у одного верхнее основание цилиндра, у другого нижнее.
От цилиндра останется тогда некоторое тело, объём которого, как мы сейчас увидим, равен объёму нашего шара.

Проведём какую-нибудь плоскость, параллельную плоскости $H$ и которая пересекалась бы с обоими телами.
Пусть расстояние этой плоскости от центра шара будет $d$, а радиус круга, полученного в сечении плоскости с шаром, пусть будет $r$.
Заметим, что $R$, $d$ и $r$ являются сторонами прямоугольного треугольника.
По теореме Пифагора, имеем $r^2=R^2-d^2$.
Значит площадь этого круга окажется равной 
\[\pi r^2 = \pi (R^2 - d^2).\]

Та же секущая плоскость даст в сечении с телом, оставшимся от цилиндра, круговое кольцо (оно на чертеже покрыто штрихами), у которого радиус внешнего круга равен $R$, а внутреннего $d$ (прямоугольный треугольник, образованный этим радиусом и отрезком $am$, равнобедренный, так как каждый острый угол его равен $45\degree$).
Значит, площадь этого кольца равна 
\[\pi R^2 - \pi d^2 = \pi(R^2 - d^2).\]

Мы видим, таким образом, что секущая плоскость, параллельная плоскости $H$, даёт в сечении с шаром и телом, оставшимся от цилиндра, фигуры одинаковой площади, следовательно, согласно принципу Кавальери, объёмы этих тел равны.
Но объём тела, оставшегося от цилиндра, равен объёму цилиндра без удвоенного объёма конуса, то есть он равен
\[\pi R^2\cdot 2R - 2 \cdot \tfrac13\pi R^2\cdot R = 2\pi R^3 - \tfrac23\pi R^3 = \tfrac43\pi R^3,\]
значит, это и будет объём шара.

}

\begin{wrapfigure}{r}{38 mm}
\vskip-0mm
\centering
\includegraphics{mppics/s-ris-154}
\caption{}\label{1938/s-ris-154}
\vskip-0mm
\end{wrapfigure}

\paragraph{}\label{1938/s148}
\mbox{\so{Определения}.}
1) Часть шара ($ACC'$, рис.~\ref{1938/s-ris-154}), отсекаемая от него какой-нибудь плоскостью ($CC'$), называется \rindex{шаровой!сегмент}\textbf{шаровым сегментом}.
Круг сечения называется \rindex{основание!шарового сегмента}\textbf{основанием сегмента}, а отрезок $Am$ радиуса, перпендикулярного к основанию, — высотой сегмента.

2) Часть шара, заключённая между двумя параллельными секущими плоскостями ($CC'$ и $DD'$), называется \rindex{шаровой!слой}\textbf{шаровым слоем}.
Круги параллельных сечений называются основаниями слоя, а расстояние $mn$ между ними — его \rindex{высота!шарового слоя}\textbf{высотой}.

Оба эти тела можно рассматривать как происходящие от вращения вокруг диаметра $AB$ части круга $AmC$ или части $CmnD$.

\paragraph{}\label{1938/s149}
\mbox{\so{Теорема}.}
\textbf{\emph{Объём шарового сегмента равен объёму цилиндра, у которого радиус основания есть высота сегмента, а высота равна радиусу шара, уменьшенному на треть высоты сегмента,}} то есть
\[ V = \pi H^2(R-\tfrac13 H)\]
где $H$ есть высота сегмента, a $R$ — радиус шара.

\begin{wrapfigure}{o}{38 mm}
\vskip-0mm
\centering
\includegraphics{mppics/s-ris-155}
\caption{}\label{1938/s-ris-155}
\vskip-0mm
\end{wrapfigure}

Объём шарового сегмента, получаемого вращением вокруг диаметра $AD$ (рис.~\ref{1938/s-ris-155}) части круга $ACB$, найдётся, если из объёма шарового сектора, получаемого вращением кругового сектора $AOB$, вычтем объём конуса, получаемого вращением $\triangle COB$.
Первый из них равен $\tfrac23\pi R^2H$, а второй $\tfrac13\pi CB^2\cdot CO$.
Так как $CB$ есть средняя пропорциональная между $AC$ и $CD$, то $CB^2 = H(2R - H)$, поэтому
\begin{align*}
CB^2\cdot CO &= H(2R - H)(R - H)=
\\
&=2R^2H - RH^2 - 2RH^2 + H^3 =
\\
&= 2R^2H-3H^2R + H^3;
\end{align*}
следовательно, 
\begin{align*}
\text{объём}\, ABB_1
&= \text{объёму}\, OBAB_1- \text{объём}\, OBB_1=
\\
&=\tfrac23\pi R^2H -\tfrac13 \pi CB^2\cdot CO = 
\\
&=\tfrac23\pi R^2H - \tfrac23\pi R^2H + \pi RH^2 - \tfrac13\pi H^3 =
\\
&= \pi H^2(R-\tfrac13H).
\end{align*}

{\small

\subsection*{Упражнения}

\begin{enumerate}[noitemsep]

\item 
Объём цилиндра, у которого высота вдвое более диаметра основания, равен 1~м$^3$.
Вычислить его высоту.

\item
Вычислить площадь боковой поверхности и объём усечённого конуса, у которого радиусы оснований равны 27~см и 18~см, а образующая равна 21~см.

\item
На каком расстоянии от центра шара, радиус которого равен 2,425~м, следует провести секущую плоскость, чтобы отношение площади поверхности меньшего сегмента к площади боковой поверхности конуса, имеющего общее с сегментом основание, а вершину в центре шара, равнялось 7:4?

\item
Найти объём тела, происходящего от вращения правильного шестиугольника со стороной $a$ вокруг одной из его сторон.

\item
Вычислить радиус шара, описанного около куба, ребро которого равно 1 м.

\item
Вычислить объём тела, происходящего от вращения правильного треугольника со стороной $a$ вокруг оси, проходящей через его вершину и параллельной противоположной стороне.

\item
Дан равносторонний треугольник $ABC$ со стороной $a$;
на $AC$ строят квадрат $BCDE$, располагая его в противоположную сторону от треугольника.
Вычислить объём тела, происходящего от вращения пятиугольника $ABEDC$ вокруг стороны $AB$.

\item
Дан квадрат $ABCD$ со стороной $a$.
Через вершину $A$ проводят прямую $AM$, перпендикулярную к диагонали $AC$, и вращают квадрат вокруг $AM$.
Вычислить площадь поверхности, образуемую контуром квадрата, и объём, образуемый площадью квадрата.

\item
Дан правильный шестиугольник $ABCDEF$ со стороной $a$.
Через вершину $A$ проводят прямую $AM$, перпендикулярную к радиусу $OA$, и вращают шестиугольник вокруг $AM$.
Вычислить площадь поверхности, образуемую контуром, и объём, образуемый площадью правильного шестиугольника.

\item
Вычислить объём шара, который, будучи вложен в коническую воронку с радиусом основания $r = 5$ см и с образующей $l = 13$ см, касается основания воронки.

\item
Центр сферы радиуса $r$ располагается на другой сфере радиуса $R>2r$ и высекает из неё сегмент.
Доказать, что площадь сегмента равна площади круга радиуса $r$; в частности эта площадь не зависит от радиуса~$R$.

\item
В шаре радиуса $r$, вдоль его диаметра, просверлено цилиндрическое отверстие длины $l$.
Доказать, что объём оставшейся части равен объёму шара диаметра $l$; в частности этот объём не зависит от радиуса~$r$.

\item
Около круга радиуса $r$ описан равносторонний треугольник.
Найти отношение объёмов тел, которые производятся вращением круга и площади треугольника вокруг высоты треугольника.

\item
В цилиндрический сосуд, у которого диаметр основания равен 6 см, а высота 36 см, налита вода до половины высоты сосуда.
На сколько поднимается уровень воды в сосуде, если в него погрузить шар диаметром 5 см.

\item
Железный пустой шар, внешний радиус которого равен 0,154 м, плавает в воде, погружаясь в неё наполовину.
Вычислить толщину оболочки этого шара, зная, что плотность железа в 7,7 раза выше плотности воды.

\item
Диаметр Марса составляет половину земного.
Во сколько раз площадь поверхности и объём Марса меньше, чем соответственные величины для Земли?

\item
Диаметр Юпитера в 11 раз больше земного.
Во сколько раз Юпитер превышает Марс по площади поверхности и объёму? (Используйте условие предыдущей задачи.)
\end{enumerate}

}